\documentclass[12pt]{amsart}

\usepackage{times,amsfonts,amsmath,amstext,amsbsy,amssymb,
  amsopn,amsthm,upref,eucal, amscd} 
\usepackage[T1]{fontenc}

\newtheorem{theorem}{Theorem}[section]
\newtheorem{lemma}[theorem]{Lemma}
\newtheorem{corollary}[theorem]{Corollary}
\newtheorem{proposition}[theorem]{Proposition} 

\numberwithin{equation}{section}

\newtheorem*{tha1}{Theorem A1}
\newtheorem*{tha2}{Theorem A2}
\newtheorem*{thb1}{Theorem B1}
\newtheorem*{thb2}{Theorem B2}
\newtheorem*{thc}{Theorem C}

\theoremstyle{definition}
\newtheorem{definition}[theorem]{Definition}

\newcommand{\field}[1]{\mathbb{#1}}
\newcommand{\R}{\field{R}} 
\newcommand{\N}{\field{N}}
\newcommand{\C}{\field{C}} 
\newcommand{\Z}{\field{Z}}
 
\newcommand{\Cal}{\mathcal}

\newcommand{\<}{{\langle}} 
\renewcommand{\>}{{\rangle}}
\newcommand{\pref}[1]{(\ref{#1})}

\newcommand{\lyap}{l}

\begin{document}
\title[Sobolev regularity of solutions ]
{Sobolev regularity of solutions \\ of the cohomological equation}

\author{Giovanni Forni $^{\dag}$}

\address{Department of Mathematics, University of Toronto\\
40 St. George St., Toronto, Ontario,  CANADA M5S 2E4}

\email{forni@math.toronto.edu}

\address{Laboratoire de Math\'ematiques\\
        Universit\'e de Paris-Sud, B\^ atiment 425\\
       91405 Orsay Cedex, FRANCE}
        
\email{giovanni.forni@math.u-psud.fr}

\thanks{$\dag$ The author gratefully acknowledges  support from National Science Foundation  grant DMS-0244463.}

\keywords{Invariant distributions, cohomological equation, Teich\"muller flow, Kontsevich-Zorich
cocycle, distributional cocycles}

\subjclass{37A20, 37C10, 37C15}

\date{\today}

  


\maketitle
\tableofcontents

\section{Introduction} 

We prove the sharpest results available on the loss of regularity for solutions of the cohomological 
equation for translation flows. For \emph{any given translation surface} and for the directional flow in almost all directions the smallest loss Sobolev regularity available to the Fourier analysis methods developed in \cite{F97} is essentially $3+\epsilon$ (for any $\epsilon>0$).  We remark that this is the best result available for the flow of rational polygonal billiards in almost all directions. The motivation
for improving the estimate on the loss of regularity obtained in \cite{F97} was provided by a question
of Marmi, Moussa and Yoccoz \cite{MMY05}. We also remark that their results, on the related problem of
solutions of the cohomological equation for almost all interval exchange transformations (IET's), 
do not apply to rational billiards for the well-known reason that IET's induced by rational billiard flows form a zero measure set in the space of all IET's. 

 For \emph{almost all translation surfaces} in every stratum of the moduli space, we prove the 
 refined, optimal result that the loss of Sobolev regularity for the directional flow in almost all directions 
 is $1+\epsilon$ (for any $\epsilon>0$). In this case, in fact we prove that for any function of Sobolev order $s>1$, the solution and its derivatives up to order $k<s-1$ are $L^\infty$ functions on the surface. We also determine precisely the Sobolev orders of the distributional obstructions to the existence of solutions (first constructed in \cite{F97}) in terms of the Kontsevich-Zorich Lyapunov exponents \cite{Zo96}, \cite{Ko97}. As a consequence we are able to determine the exact codimension of coboundaries for every Sobolev regularity class of the transfer function. For instance, the codimension 
 of coboundaries with square-integrable transfer functions (in the space of functions of Sobolev order $s>1$) is exactly equal to the genus of the surface. For such coboundaries the transfer function is actually in $L^\infty$. 
 
These results implies quite immediately corresponding results for interval exchange transformations, which  improve on the loss of regularity established in \cite{MMY05}. We should point out that in that
 paper the authors are mostly concerned with \emph{Diophantine conditions} on interval exchange
 transformations for which the cohomological equation admits smooth solution, while we have not investigated this question at all. The reason is that the full measure sets of systems for which our
 results (as well as \cite{MMY05}) are determined by several conditions which always include 
 \emph{Oseledec regularity} (or rather a weaker  \emph{coherence} property) with respect to Kontsevich-Zorich renormalization cocycle. The Oseledec's theorem is invoked to ensure that the set of regular
 (coherent) IET's has full measure. Any substantial progress over \cite{MMY05} would have to succeed in characterizing explicitly a full measure set of regular (coherent) points without relying on the Oseledec's
 theorem. To the best of our knowledge this goal is still beyond reach.
 
 There are several motivations for this work. The study of cohomological equations is a relevant
 part of the theory of (smooth) dynamical systems directly connected to basic questions such
 as \emph{triviality of time-changes }for flows, \emph{asymptotic of ergodic averages} and the 
 \emph{smooth conjugacy problem} via linearization and Nash-Moser implicit function theorem.
 In the hyperbolic case (for dynamical systems with exponential divergence of nearby orbits) such there
are extensive, deep results on the cohomological equation pioneered in the work of Livsic \cite{L71}, later developed by several authors (see \cite{GK}, \cite{CEG}, \cite{LMM}). The completely different
case of (Diophantine) linear flows on the torus  is well-known, since the cohomological equation for such systems is closely related to the linearized equation in the classical KAM theory for Hamiltonian flows. This in an example of elliptic dynamics (no divergence of nearby orbits) which can be studied to a great extent by the classical theory of Fourier series. It is characterized by the `small divisors'  appearing in 
the Fourier coefficients of solutions, which lead to a loss of regularity. It is not difficult to see that the
optimal loss of Sobolev regularity for the full measure set of Roth flows is $1+\epsilon$ and that for
such flows any zero average function of Sobolev order $s>1$ is an $L^2$ coboundary. It can be proved by our methods (and by the Gottschalk-Hedlund theorem \cite{GH55}) that the transfer function is in fact continuous. We have not been able to locate this result in the literature, however it is well within reach of the methods of  \cite{He83}, Chap. VI, \S 3. However, only the measure zero case of  rotation numbers of \emph{constant type} seems to have been explicitly considered there.

For systems with intermediate behavior, that is, for elliptic systems with singularities or for parabolic systems (characterized by polynomial divergence of nearby orbits) much less is known. The author
discovered in \cite{F97} that the cohomological equation for generic translation flows (or equivalently
for generic  IET's) has finitely smooth solutions  for sufficiently smooth data under finitely many
distributional conditions. In other terms, on one hand the problem shares a typical feature of `small divisors' problems,  namely the finite loss of regularity of solutions with respect to the data; on the
other hand, a new phenomenon appears: the existence of infinitely many independent distributional
obstructions (of increasing order) which are not given by invariant measures. In \cite{F97} only
a rough estimate for the loss of derivatives is explicitly obtained ($\leq 9$). Our goal in this paper is
to improve such estimate as much as possible. In joint papers with L. Flaminio  
 the authors have investigated the existence of smooth solutions of the cohomological
equation for horocycle flows (on surfaces constant negative curvature) \cite{FF03}, for generic nilflows on quotients of the Heisenberg group \cite{FF06}  and generic nilflows on general nilmanifolds \cite{FF07}. In all cases the fundamental features of finite loss derivatives and of the existence of infinitely many independent distributional obstructions have been established (although the structure of the space of invariant distributions is significantly different for IET's, horocycle flows and nilflows). For horocycle flows and for Heisenberg nilflows it was possible to estimate that the loss of Sobolev regularity is $1+\epsilon$ (for any $\epsilon>0$) and to establish the conjectural relation that the Sobolev order of the distributional obstructions be related to the Lyapunov exponents of the distribution under the appropriate renormalization dynamics. In this paper we prove analogous results for generic translation flows. We should point out that for generic nilflows on general nilmanifolds the loss of regularity and the regularity of the distributional obstructions seem to depend on the depth and rank of the nilpotent group considered, although no lower bound was established in \cite{FF07}.

Let $q$ be a \emph{holomorphic orientable quadratic differential} on a Riemann surface $M$ of genus
$g \geq 1$.  The horizontal and vertical measured foliations (in the Thurston's sense) associated to a holomorphic quadratic differential  $q$ on $M$ are defined as $\Cal F_q =\{ \Im(q^{1/2}) =0 \}$ (the horizontal foliation) and $\Cal F_q =\{ \Re(q^{1/2}) =0 \}$. Such foliations are well-defined even in the case that there is no globally defined square root of the quadratic differential. The horizontal foliation is endowed with the transverse measure given by $\vert \Im(q^{1/2})\vert$, the vertical foliation is endowed with the transverse measure given by $\vert \Re(q^{1/2})\vert$. The quadratic differential is called
orientable if the horizontal and vertical foliations are both orientable. Orientability is equivalent to
the condition that the quadratic differential is globally the square of a holomorphic (abelian) differential.
The structure induced by an orientable holomorphic quadratic differential (or by a holomorphic
abelian differential) can also be described as follows. There is a flat metric $R_q$ associated with any
quadratic differential $q$ on $M$. Such a metric has conical singularities as the finite set $\Sigma_q
=\{p\in M \vert  q(p)=0\}$. If $q$ is orientable there exists a (positively oriented) parallel orthonormal frame $\{S_q, T_q\}$ of the tangent bundle $TM\vert M\setminus\Sigma_q$ such that $S_q$ is tangent
to the horizontal foliation $\Cal F_q$ and $T_q$ is tangent to the vertical foliation $\Cal F_{-q}$ everywhere on $M\setminus \Sigma_q$. In other terms, \emph{ the flat metric $R_q$ has trivial holonomy}. In another equivalent formulation, any orientable holomorphic quadratic differential determines a \emph{translation structure} on $M$, that is, an equivalence class of atlases with transition functions given by translation of the euclidean plane (see for instance the excellent survey \cite{MTHB}, \S 1.8. For a given orientable quadratic differential $q$ on a Riemann surface $M$, we will consider
the one-parameter family of vector fields on $M\setminus \Sigma_q$ defined as
\begin{equation}
\label{eq:vectorfields}
S_\theta := \cos \theta \, S_q  \, + \, \sin\theta\, T_q \,, \quad \theta \in S^1 .
\end{equation}
The vector field $S_\theta$ is a parallel normalized  vector field in the direction at angle $\theta\in S^1$
with the horizontal. We remark that it is not defined as the singular set $\Sigma_q$ of the flat metric. hence the flow it generates is defined (almost everywhere) on the complement of the union of all separatrices of the orbit foliation (a measure foliation). The singularities of the orbit foliation are all
saddle-like, but the saddles are degenerate if the order of zero of the quadratic differential at the
singularity is strictly greater than $2$. In fact, since the quadratic differential is supposed to be orientable
it has zeroes of even order and the orbit foliations of the vector fields \pref{eq:vectorfields} has $m$
stable and $m$ unstable separatrices at any zero of order $2m$.

\smallskip
\noindent  Our goal is to investigate the loss of (Sobolev) regularity of solutions of the \emph{cohomological equation} $S_\theta u =f$ for Lebesgue almost all $\theta \in S^1$. The author proved in \cite{F97} that  if the function $f$ is sufficiently regular, satisfies a finite number of independent distributional conditions (which include conditions on the jets at the singularities) then there exists a finitely smooth solution (unique up to additive constants). The loss of regularity was estimated in that paper to be no more than $9$ derivatives in the Sobolev sense. 

\noindent If $q$ is any orientable quadratic differential, the regularity of functions on the translation surface $(M,q)$ is expressed in terms of a family $\{H^s_q(M) \vert s\in \R\}$ of {\it weighted Sobolev spaces}. Such spaces were introduced in \cite{F97} for all $s\in \Z$ as follows. Let $\omega_q$ be the standard (degenerate) volume form on $M$ of the flat metric $R_q$. The space $H^0_q(M)$ is the space $L^2(M, \omega_q)$ of square-integrable functions.  For $k\in \N$, the space $H^k_q(M)$ is the subspace of functions $f \in H^0_q(M)$  such that the weak derivatives $S_q ^i T_q^j f \in H^0_q(M)$ [and $T_q ^i S_q^j f \in H^0_q(M)$] for all $i+j \leq k$ and the space $H^{-k}_q(M)$ is the dual Hilbert space $H^k_q(M)^\ast$. In  \S \ref{fracsobspaces} of this paper we introduce weighted \emph{Sobolev spaces with arbitrary (real) exponents }by methods of interpolation theory. Although the Sobolev 
norms we construct do not form an interpolation family in the sense of (holomorphic) interpolation theory, they do satisfy a standard interpolation inequality.  The weighted Sobolev spaces combine standard Sobolev smoothness conditions on $M\setminus\Sigma_q$ with restrictions on the jet of the functions at the singular set $\Sigma_q\subset M$.

\noindent As discovered in \cite{F97}, for functions $f\in C_0^\infty(M\setminus \Sigma_q)$ the space of all distributional obstructions to the existence of a solution $u\in C_0^\infty(M\setminus \Sigma_q)$ of the cohomological equation $S_\theta u=f$ coincides for almost all $\theta\in S^1$ with the infinite dimensional space of all \emph{$S_\theta$-invariant distributions}:
\begin{equation}
\Cal I_{q,\theta} (M\setminus\Sigma_q) := \{ \Cal D \in \Cal D'(M\setminus\Sigma_q) \, \vert \,   S_\theta \Cal D = 0 \text{ \rm in }  \Cal D'(M\setminus\Sigma_q)  \}\,.
\end{equation}
For data $f\in H^s_q(M)$ of finite Sobolev differentiability, a complete set of obstructions is given for almost all $\theta\in S^1$ by the finite dimensional subspace of invariant distributions
\begin{equation}
\Cal I^s_{q,\theta} (M) :=  \{ \Cal D \in H^{-s}_q(M) \, \vert \,   S_\theta \Cal D = 0  \text{ \rm in }
H^{-s}_q(M) \}\,.
\end{equation}
The goal of this paper is to prove \emph{optimal }estimates on the Sobolev regularity of solutions of the cohomological equation $S_\theta u=f$ and on the dimension of the spaces $\Cal I^s_{q,\theta} (M)$ 
of invariant distributions for all $s>0$. In \S \ref{ss:general}, Theorem \ref{thm:GCEsmooth}, we prove

\begin{tha1}  Let $q$ be any orientable holomorphic quadratic differential. Let 
$k\in \N$ be any integer such that $k\geq 3$ and let $s>k$ and $r<k-3$. For almost all $\theta\in S^1$ (with respect  to the Lebesgue measure), there exists a constant $C_{r,s}(\theta)>0$ such that the following holds.  If $f\in H^{s}_q(M)$ is such that  $\Cal D(f)=0$ for all $\Cal D \in \Cal I^s_{q,\theta}(M)$, the cohomological equation $S_{\theta}u=f$ has a  solution $u\in H^{r}_q(M)$ satisfying the following estimate:
\begin{equation}
\label{eq:ThA1}
\vert u\vert_{r}\leq C_{r,s}(\theta)\, \vert f \vert_{s}\,\,.
\end{equation}
\end{tha1}

\noindent The dimensions of the spaces of invariant distributions can be estimated as follows
(see \S \ref{ss:bc}, Corollary \ref{cor:basiccohom} and Theorem \ref{thm:bcstruct}):
\begin{tha2} Let $q$ be any orientable holomorphic quadratic differential. Let 
$k\in \N$ be any integer such that $k\geq 3$ and let $k<s \leq k+1$. For almost all $\theta\in S^1$ (with respect  to the Lebesgue measure), 
\begin{equation}
\label{eq:ThA2}
1 + 2(k-2)(g-1)  \leq  \text{\rm dim}\, \Cal I^s_{q,\theta}(M) \leq  1 + (2k-1)(g-1) \,.
\end{equation}
\end{tha2}
\noindent The proof of the above results is essentially based on the harmonic analysis methods
developed in \cite{F97}. We remark that no other methods are known for the case of an \emph{arbitrary }orientable quadratic differential. 

\smallskip
\noindent We prove much sharper results for \emph{almost all }orientable quadratic differentials.
The moduli space of orientable holomorphic quadratic differentials $q$ on some Riemann surface 
$M_q$ with a given pattern of zeroes, that is, with zeroes of (even) multiplicities $\kappa= (k_1, \dots, k_\sigma)$ at a finite set $\Sigma_q =\{ p_1, \dots, p_\sigma\}\subset M_q$ is a \emph{stratum} 
$\Cal M_{\kappa}$ of the moduli space $\Cal M_g$ of all holomorphic quadratic differential. 
Let $\Cal M^{(1)}_{\kappa} \subset \Cal M^{(1)}_g$ be the subsets of quadratic differential of total
area equal to $1$. It was proved by H. Masur \cite{Ma82} and W. Veech \cite{Ve86} that each 
stratum $\Cal M^{(1)}_{\kappa}$ carries an absolutely continuous probability measure $\mu^{(1)}_\kappa$, invariant under the action of the Teichm\"uller geodesic flow, which is ergodic 
when restricted to each connected component of $\Cal M^{(1)}_{\kappa}$ (the connected components 
of strata of orientable quadratic differentials were classified in \cite{KZ03}). In fact, there is natural action of the group $SL(2, \R)$ on the moduli space $\Cal M^{(1)}_g$ such that the Teichm\"uller geodesic flow corresponds to the action of the diagonal subgroup of $SL(2,\R)$ and the measure $\mu^{(1)}_\kappa$
is $SL(2, \R)$ invariant.

\noindent In \cite{Ko97} M. Kontsevich introduced a renormalization cocycle for translation flows, inspired to the Rauzy-Veech-Zorich cocycle for interval exchange transformations. The Kontsevich-Zorich cocycle is a dynamical system on an orbifold vector bundle over $\Cal M^{(1)}_\kappa$ with fiber the first cohomology $H^1(M_q,\R)$ of the Riemann surface carrying the orientable holomorphic quadratic differential $q\in \Cal M^{(1)}_{\kappa}$. The action of such a dynamical system  is (by definition of a cocycle) linear on the fibers and projects onto the Teichm\"uller geodesic flow on the base $\Cal M^{(1)}_\kappa$. Since the cocycle is symplectic, for any probability measure $\mu$ on a stratum $\Cal M^{(1)}_\kappa$, the \emph{Lyapunov spectrum} of the Kontsevich-Zorich cocycle takes the form:
\begin{equation}
\label{eq:introKZexp}\lambda^{\mu}_1\geq \dots \geq\lambda^{\mu}_g\geq 0\geq \lambda^\mu_{g+1}=-\lambda^{\mu}_g\geq \dots \geq \lambda^\mu_{2g}=-\lambda^{\mu}_1\,.
\end{equation}
In addition, it is not difficult to prove that $\lambda^{\mu}_1=1$. A probability measure $\mu$ on a stratum $\Cal M^{(1)}_\kappa$, invariant under the Teichm\"uller geodesic flow, will be called  $(a)$ \emph{$SO(2, \R)$-absolutely continuous} if it induces absolutely continuous measures on every orbit 
of the circle group $SO(2, \R) \subset SL(2, \R)$; $(b)$ \emph{KZ-hyperbolic} if all the Lyapunov exponents in \pref{eq:introKZexp} are non-zero. It is immediate that all $SL(2, \R)$-invariant measures are  $SO(2, \R)$-absolutely continuous. It was first proved in \cite{F02} that the measure $\mu^{(1)}_\kappa$ is KZ-hyperbolic. A different proof that also reaches the stronger conclusion that the exponents \pref{eq:introKZexp} are all distinct has been given more recently by A. Avila and M. Viana \cite{AV} who have thus completed the proof of the Zorich-Kontsevich conjectures \cite{Zo96}, \cite{Ko97} on the Lyapunov spectrum of the Kontsevich-Zorich cocycle (and its discrete counterparts).

\smallskip
\noindent Our sharpest results are proved for \emph{ almost all }quadratic differentials
with respect to any $SO(2, \R)$-absolutely continuous, KZ-hyperbolic, Teichm\"uller invariant, probability
measure on any stratum $\Cal M^{(1)}_\kappa$ of orientable quadratic differentials. The smoothness
informations on the solutions is stronger than just Sobolev $L^2$ regularity and it is naturally encoded
by the following spaces. For any $k\in \N$, let $B^k_q(M)$ be the space of all functions $u\in H^{k}_q(M)$ such that $S_q^i T_q^j u=T_q^i S_q^j u \in L^{\infty}(M)$ for all pairs of integers $(i,j)$ such that $0\leq i+j \leq k$. The space $B^k_q(M)$ is endowed with the norm defined as follows: for any $u \in B^k_q(M)$,
\begin{equation}
\label{eq:introuniformnorms1}
\vert u \vert_{k,\infty} \,:= \,  \left[\sum_{i+j \leq k}  \vert S_q^i T_q^j u\vert^2_{\infty}\right]^{1/2}
\,=\,\left[\sum_{i+j \leq k}  \vert T_q^iS_q^j u\vert^2_{\infty} \right]^{1/2} \,.
\end{equation}
For $s\in [k,k+1)$, let $B^s_q(M):= B^k_q(M) \cap H^s_q(M)$ endowed with the norm defined
as follows: for any  $u \in B^s_q(M)$, 
\begin{equation}
\label{eq:introuniformnorms2}
\vert u \vert_{s,\infty} :=  \left(\,\vert u \vert^2_{k,\infty}\,+\,  \vert u \vert^2_s \,\right)^{1/2}\,.
\end{equation}
In \S \ref{ss:gencase}, Theorem \ref{thm:CEsharpth}, we prove the following:

\begin{thb1}  Let $\mu$ be any $SO(2, \R)$-absolutely continuous, KZ-hyperbolic probability measure on any stratum $\Cal M^{(1)}_\kappa$ of orientable quadratic differentials. Let $s>1$ and let $r<s-1$.  For $\mu$-almost all $q\in \Cal M^{(1)}_\kappa$ and for almost all $\theta\in S^1$ (with respect  to the Lebesgue measure), there exists a constant $C_{r,s}(\theta)>0$ such that the following holds.  If $f\in H^{s}_q(M)$ is such that  $\Cal D(f)=0$ for all $\Cal D \in \Cal I^s_{q,\theta}(M)$, the cohomological equation $S_{\theta}u=f$ has a  solution $u\in B^{r}_q(M)$ satisfying the following estimate:
\begin{equation}
\label{eq:ThB1}
\vert u\vert_{r,\infty}\leq C_{r,s}(\theta)\, \vert f \vert_{s}\,\,.
\end{equation}
\end{thb1}

\noindent The regularity of invariant distributions can be precisely determined as follows
(see \S \ref{ss:Lyapexp}, Corollary \ref{cor:Lspectrum}, and \S \ref{ss:gencase}, Theorem
\ref{thm:CEsharpth}). For any $q\in \Cal M^{(1)}_\kappa$ and any distribution $\Cal D \in \Cal D'(M\setminus\Sigma_q)$, the \emph{weighted Sobolev order} is the number
\begin{equation}
\Cal O^H_q(\Cal D) =\inf \{  s \in \R^+ \vert   \Cal D \in H^{-s} _q(M) \}\,.
\end{equation}
Let $\Cal I_{q,\theta}(M):= \cup \{ \Cal I^s_{q,\theta}(M) \vert s\geq 0\}$ denote the space of all $S_\theta$-invariant distribution of finite Sobolev order and let $\hat{\Cal I}_{q,\theta}(M) \subset
\Cal I_{q,\theta}(M)$ be the subspace of invariant distributions vanishing on constant functions.
It follows immediately by the definitions that $\Cal I_{q,\theta}(M)= \C \oplus \hat{\Cal I}_{q,\theta}(M)$.

\begin{thb2}  Let $\mu$ be any $SO(2, \R)$-absolutely continuous, KZ-hyperbolic probability measure on any stratum $\Cal M^{(1)}_\kappa$ of orientable quadratic differentials.  For $\mu$-almost all $q\in \Cal M^{(1)}_\kappa$ and for almost all $\theta\in S^1$ (with respect  to the Lebesgue measure),  the
space $\hat {\Cal I}_{q,\theta}(M)$ has a basis $\{ \Cal D_{i,j}(\theta)\}$ such that 
\begin{equation}
O^H_q(\Cal D_{i,j}(\theta) )=   \lambda^\mu_i - (j+1)\,, \quad  i\in \{2,\dots, 2g-1\}, \,  j\in \N \cup \{0\} \,.
\end{equation}
In addition, the basis $\{ \Cal D_{i,j}(\theta)\}$ can be generated from the finite dimensional subsystem
$\{ \Cal D_2, \dots, \Cal D_{2g-1}\}$ by the following differential relations:
\begin{equation}
\Cal D_{i,j}(\theta) =  T_\theta^j \Cal D_{i,0} (\theta) \,, \quad  i\in \{2,\dots, 2g-1\}, \,  j\in \N \cup \{0\}\,.
\end{equation}
\end{thb2}

\noindent The above Theorems B1 and B2 are proved by methods based on renormalization, which were inspired by the work of Marmi-Moussa-Yoccoz \cite{MMY05}. However,
our approach differs from theirs since we explicitly assume that the Lyapunov exponents of the Kontsevich-Zorich cocycle are all non-zero.  The main idea of the argument, as in \cite{MMY05}, is
to prove uniform estimates for ergodic integrals of weakly differentiable functions, then apply
a version of Gottshalk-Hedlund theorem. The asymptotics of ergodic averages of functions in
$H^1(M)$ was studied by the author in \cite{F02}, where the Kontsevich-Zorich conjectures
on the deviation of ergodic averages for smooth functions (formulated in \cite{Ko97}) were proved.
The approach of \cite{F02} is based on the analysis of distributional cocycles over the Teichm\"uller
flow which extend the Kontsevich-Zorich cocycles.  The estimates proved in \cite{F02} are (barely) not strong enough to yield the required uniform boundedness of ergodic averages under the appropriate distributional conditions. In \S \ref{distcocycles} of this paper  we have recalled the definition of distributional cocycles, and we have strengthened the estimates proved in \cite{F02}  under the slightly stronger (and correct) assumption that the functions considered belong to $H_q^s(M)$ for some $s>1$. 

\noindent Another important techical issue that separates Theorems B1 and B2 from the less precise
Theorems A1 and A2 is related to interpolation theory in the presence of distributional obstructions.
In the general case, we have not been able to overcome the related difficulties, hence the lack of precision of Theorem A1 and A2 for intermediate Sobolev regularity. In the \emph{generic }case of Theorems B1 and B2 we have been able to prove a remarkable linear independence property of invariant distributions which makes interpolation possible in the construction of solutions of the cohomological equation. 

\smallskip
\noindent We introduce the following definition (see Definition \ref{def:regularsystem}). A finite system $\{\Cal D_1, \dots, \Cal D_J\} \subset H^{-\sigma}_q(M)$ of finite order distributions is called \emph{$\sigma$-regular} (with respect to the family $\{H^s_q(M)\}$ of weighted Sobolev spaces) if for any $\tau \in (0,1]$ there exists a dual system  $\{u_1(\tau), \dots, u_J(\tau)\} \subset H^\sigma_q(M)$ (that is, the identities $\Cal D_i( u_j(\tau)) =\delta_{ij}$ hold for all $i$, $j\in \{1, \dots, J\}$ and all $\tau \in (0,1]$) such that the following estimates hold. For all $0\leq r  \leq \sigma$  and all $\epsilon>0$, there exists a constant $C^\sigma_r(\epsilon)>0$ such that, for all $i$, $j\in \{1,\dots, J\}$,
\begin{equation}
\label{eq:introregularsystem}
\vert u_j(\tau) \vert_r \leq C^\sigma_r(\epsilon)\, \tau^{{\Cal O}^H(\Cal D_j)-r-\epsilon}\,.
\end{equation}
A finite system  $\{\Cal D_1, \dots, \Cal D_J\} \subset H^{-s}_q(M)$ of finite order distributions will be called \emph{regular} if it is $\sigma$-regular  for any $\sigma\geq s$. A finite dimensional subspace $\Cal I \subset H^{-s}_q(M)$ of finite order distributions will be called \emph{$\sigma$-regular [regular]} 
if it admits a $\sigma$-regular [regular] basis. 

\noindent We have proved that the spaces of distributional obstructions for the cohomological
equation are regular in the above sense (see Theorem \ref{thm:CEsharpfundth}). 

\begin{thc}
Let $\mu$ be any $SO(2, \R)$-absolutely continuous, KZ-hyperbolic probability measure on any stratum $\Cal M^{(1)}_\kappa$ of orientable quadratic differentials.  For $\mu$-almost all $q\in \Cal M^{(1)}_\kappa$, for almost all $\theta\in S^1$  and for all $s>0$, the space ${\Cal I}^s_{q,\theta}(M) \subset H^{-s}_q(M)$ of $S_\theta$-invariant distributions is regular. 
\end{thc}

\section{Fractional weighted Sobolev spaces}
\label{fracsobspaces}
 \noindent In \cite{F97} we have introduced a natural scale of weighted Sobolev spaces with integer exponent associated with any orientable holomorphic quadratic differential $q$ on a Riemann surface $M$ (of genus $g\ge 2$). In this section we extend  the definition of weighted Sobolev spaces  to arbitrary (real) exponents by methods of interpolation theory.

 \subsection{ Weighted Sobolev spaces}
 \label{WSS}
Let $\Sigma_{q}:=\{p_{1},\dots,p_{\sigma}\}\subset M$ be the set of zeros of the holomorphic quadratic differential~$q$, of even orders $(k_{1},\dots,k_{\sigma})$ respectively with~$k_{1} + \dots + k_{\sigma}=4g-4$.  Let  $R_{q}:=\vert q\vert^{1/2}$ be the flat metric with cone singularities at $\Sigma_{q}$ induced by the quadratic  differential $q$ on $M$.  With respect to a holomorphic local coordinate $z=x+iy$, the quadratic differential $q$ has the form $q=\phi(z)dz^{2}$, where $\phi$ is a locally defined holomorphic function, and, consequently, 
\begin{equation}
\label{eq:metric}
 R_q= |\phi(z)|^{1/2} (dx^2 +dy^2)^{1/2}\,,\quad  \omega_q=|\phi(z)|\,dx\wedge dy\,.
\end{equation}
\noindent The metric $R_{q}$ is flat, it is degenerate at the finite set $\Sigma_{q}$ of zeroes of $q$ 
and, if $q$ is orientable, it has trivial holonomy, hence $q$ induces a structure of {\it translation 
surface}  on $M$. 

 \smallskip
  \noindent The weighted $L^{2}$ space is the standard space $L^{2}_{q}(M):= L^{2}(M,\omega_{q})$  with respect to the area element $\omega_{q}$  of the metric $R_{q}$. Hence the weighted $L^{2}$ 
  norm $\vert \cdot\vert_{0}$ are induced by the hermitian product $\<\cdot, \cdot\>_{q}$ defined as follows: for all functions $u$,$v\in L^{2}_{q}(M)$,
 \begin{equation}
 \label{eq:0norm}
 \< u ,v\>_{q} :=  \int _{M} u\,\bar v \, \omega_{q}\,\,.
 \end{equation}
 Let $\Cal F_{q}$ be the {\it horizontal foliation},  $\Cal F_{-q}$ be the {\it vertical foliation} for
the holomorphic quadratic differential $q$ on $M$. The foliations $\Cal F_{q}$ and $\Cal F_{-q}$
are measured foliations (in the Thurston's sense):  $\Cal F_{q}$ is the foliation given (locally) by the equation $\Im (q^{1/2})=0$ endowed with the invariant transverse measure $\vert \Im (q^{1/2}) \vert$,  
$\Cal F_{-q}$ is the foliation given (locally) by the equation $\Re (q^{1/2})=0$ endowed with the 
invariant transverse measure $\vert \Re (q^{1/2}) \vert$.  If the quadratic differential $q$ is orientable,
since the metric $R_q$ is flat with trivial holonomy,  there exist  commuting vector fields $S_q$
and $T_q$ on $M\setminus \Sigma_{q}$ such that 
\begin{enumerate}
\item The frame $\{S_q,T_q\}$ is a parallel  orthonormal frame with respect to the metric $R_{q}$ for the restriction of the tangent bundle $TM$ to the complement $M\setminus \Sigma_{q}$  of the set of cone points;
\item the vector field $S_{q}$ is tangent to the horizontal foliation $\Cal F_{q}$, the vector field $T_{q}$
 is tangent to the vertical foliation $\Cal F_{-q}$ on $M\setminus \Sigma_{q}$ \cite{F97}. 
 \end{enumerate}
 In the following  we will often drop the dependence of the vector fields $S_{q}$, $T_{q}$ on the quadratic differential in order to simplify the notations. We have:
\begin{enumerate}
\item $\Cal L _{S} \omega_{q} = \Cal L_{T}\omega_{q} =0$ on $M\setminus \Sigma_{q}$ , that is, the 
area form $\omega_{q}$ is invariant with respect to the flows generated by $S$ and $T$;
\item $\imath_{S} \omega_{q}= \Re (q^{1/2})$ and $\imath_{T} \omega_{q}= \Im (q^{1/2})$, hence
the $1$-forms $\eta_{S} :=\imath_{S} \omega_{q}$,  $\eta_{T} :=-\imath_{T} \omega_{q}$ are smooth 
and closed on $M$ and $\omega_{q}= \eta_{T}\wedge \eta_{S}$.
\end{enumerate}
It follows from the area-preserving property $(1)$ that the vector field $S$, $T$ are anti-symmetric
as densely defined operators on $L^{2}_{q}(M)$, that is, for all functions $u$, $v \in C_0^{\infty} (M\setminus\Sigma_q)$,  (see \cite{F97}, $(2.5)$),
\begin{equation}
\label{eq:antisymm}
\< Su ,v\>_{q} = -\< u ,Sv\>_{q}\,\,, \quad \text{ respectively } \,\, \< Tu ,v\>_{q} =-\< u ,Tv\>_{q} \,\,.
\end{equation}
In fact, by Nelson's criterion~\cite{Ne59}, Lemma 3.10, the anti-symmetric operators $S$, $T$ are {\it essentially skew-adjoint} on the Hilbert space $L^{2}_{q}(M)$.

\smallskip
\noindent The {\it weighted Sobolev norms} $\vert \cdot\vert_{k}$, with integer exponent $k>0$, are the euclidean norms, introduced in \cite{F97}, induced by the hermitian product defined as follows: for all 
functions $u$, $v\in L^{2}_{q}(M)$,
\begin{equation}
 \label{eq:knorm}
 \< u,v \>_{k} :=   \frac{1}{2}\sum_{i+j\leq k}\<S^{i}T^{j}u, S^{i}T^{j}v\>_{q} + 
 \<T^{i}S^{j}u, T^{i}S^{j}v\>_{q}\,.
 \end{equation}
The  {\it weighted Sobolev norms }with integer exponent $-k<0$ are defined to be the dual norms.

\smallskip
\noindent The {\it weighted Sobolev space }$H^{k}_q(M)$, with integer exponent $k\in\Z$, is the 
Hilbert space obtained as the completion with respect to the norm $\vert \cdot  \vert _{k}$ of the 
maximal {\it common invariant domain}
 \begin{equation}
 \label{eq:cid}
 H^{\infty}_q(M):= \bigcap_{i,j\in \N}  D( \bar S^i \bar T^j) \cap D( \bar T^i \bar S^j)\,.
 \end{equation}
 of the closures $\bar S$, $\bar T$ of the essentially skew-adjoint operators $S$,  $T$ on $L^2_q(M)$.
 The weighted Sobolev space $H^{-k}_q(M)$ is isomorphic to  the dual space of the Hilbert space $H^{k}_q(M)$, for all $k\in \Z$. 
 
 \smallskip
 \noindent  Since the vector fields $S$, $T$ commute (infinitesimally) on $M\setminus\Sigma_{q}$, the following weak commutation identity holds on $M$. 
  \begin{lemma} 
   \label{lemma:commut}
  (\cite{F97}, Lemma 3.1)  
  For all functions $u$,$v  \in H^{1}_{q}(M)$,
 \begin{equation}
 \label{eq:commut}
\<Su,Tv\>_{q} =  \<Tu,Sv\>_{q}\,\,.
 \end{equation}
  \end{lemma}
  \smallskip
  \noindent By the anti-symmetry property~\pref{eq:antisymm} and the commutativity property 
  \pref{eq:commut},  the frame $\{S,T\}$ yields an essentially skew-adjoint action of the Lie 
  algebra $\R^{2}$ on the Hilbert space $L^{2}_{q}(M)$ with common domain $H^{1}_{q}(M)$.  
  If $\Sigma_{q}\not=\emptyset$, the (flat) Riemannian manifold $(M\setminus\Sigma_{q}, R_{q})$ 
  is not complete, hence its Laplacian $\Delta_q$ is not essentially self-adjoint on $C_{0}^{\infty}(M \setminus\Sigma_{q})$. By a theorem of Nelson~\cite{Ne59}, \S 9, this is equivalent to the non-integrability of the action of~$\R^{2}$ as a Lie algebra (to an action of $\R^{2}$ as a Lie group). 
 
 \medskip
 \noindent  Following \cite{F97}, the Fourier analysis on the flat surface $M_q$ will be based 
 on a canonical self-adjoint extension $\Delta_q^F$ of the Laplacian $\Delta_q$,  called  the 
 {\it Friedrichs extension},  which is uniquely determined by  the {\it Dirichlet hermitian form }$\Cal Q:H^{1}_{q}(M)\times H^{1}_{q}(M) \to \C$. We recall that, for all $u$, $v \in H^{1}_{q}(M)$,
 \begin{equation}
 \label{eq:Dirichlet}
 \Cal Q(u,v) := \<Su,Sv\>_{q} \,+\, \<Tu,Tv\>_{q} \,\,. 
 \end{equation}
   \begin{theorem} (\cite{F97}, Th. 2.3)
   \label{thm:Dirichlet}
  The hermitian form $\Cal Q$ on $L^{2}_{q}(M)$ has the following spectral properties:
  \begin{enumerate}
  \item $\Cal Q$ is positive semi-definite and the set $\text{\rm EV}(\Cal Q)$ of its eigenvalues is a  
  discrete subset of $[0,+\infty)$;
  \item Each eigenvalue has finite multiplicity, in particular $0\in \text{\rm EV}(\Cal Q)$ is simple and the kernel of $\Cal Q$ consists only of constant functions;
  \item The space $L^{2}_{q}(M)$ splits as the orthogonal sum of the eigenspaces. In addition,
  all eigenfunctions are $C^{\infty}$ (real analytic) on $M$.
   \end{enumerate} 
  \end {theorem} 
  \noindent The {\it Weyl asymptotics }holds for the eigenvalue  spectrum of the Dirichlet form . For any $\Lambda>0$, let $N_{q}(\Lambda):=\text{\rm card} \{ \lambda \in \text{\rm EV}(\Cal Q)\,/\, \lambda \leq \Lambda \}$, where each eigenvalue $\lambda \in \text{\rm EV}(\Cal Q)$ is counted according to its multiplicity. 
  \begin{theorem} (\cite{F97}, Th. 2.5) 
    \label{thm:Weyl}
  There exists a constant $C>0$ such that
   \begin{equation}
 \label{eq:Weyl}
\lim_{\Lambda\to + \infty} \,\frac {N_{q}(\Lambda)}{\Lambda} \,\, = \,\, \text{\rm vol}(M,R_{q})\,\,. 
 \end{equation}
 \end {theorem} 
 \medskip
 \noindent Let $\partial^{\pm}_{q}:=S_{q}\pm \imath\, T_{q}$ be the {\it Cauchy-Riemann 
 operators }induced by the holomorphic orientable quadratic differential $q$ on $M$, introduced in~\cite{F97}, \S 3.  Let $\Cal M^{\pm}_{q}\subset L^{2}_{q}(M)$ be the subspaces of meromorphic, respectively anti-meromorphic functions (with poles at $\Sigma_{q}$). By the Riemann-Roch 
 theorem, the subspaces $\Cal M^{\pm}_{q}$ have  the same complex dimension equal to 
 the genus $g\geq 1$ of the Riemann surface $M$. In addition,  $\Cal M^{+}_{q}\cap 
 \Cal M^{-}_{q}=\C$, hence 
  \begin{equation}
  \label{eq:H}
H_{q}:= \left(\Cal M^{+}_{q}\right)^{\perp} \oplus  \left(\Cal M^{-}_{q}\right)^{\perp} = 
 \{ u\in L^{2}_{q}(M)\,\vert \,  \int_{M} u\, \omega_{q} \,=\,0\,\}\,\,.
 \end{equation}
 Let $H^{1}_{q}: = H_{q} \cap H^{1}_{q}(M)$. By Theorem~\ref{thm:Dirichlet}, the restriction of the
 hermitian form to $H^{1}_{q}$ is positive definite, hence it induces a norm. By the Poincar\'e 
 inequality (see \cite{F97}, Lemma 2.2 or \cite{F02}, Lemma 6.9), the Hilbert space $(H^{1}_{q}, 
 \Cal Q)$ is isomorphic to the Hilbert space $(H^{1}_{q}, \<\cdot,\cdot\>_{1})$.
 
 \begin{proposition} (\cite{F97}, Prop. 3.2) 
 \label{prop:CR}
  The Cauchy-Riemann operators $\partial^{\pm}_{q}$ are closable operators on the common domain
 $C^{\infty}_{0}(M\setminus\Sigma_{q}) \subset  L^{2}_{q}(M)$ and their closures (denote by the same symbols) have the following properties:
 
\begin{enumerate}
\item the domains $D(\partial^{\pm}_{q}) = H^{1}_{q}(M)$ and the kernels $N(\partial^{\pm}_{q}) = \C$;
\item the ranges $R_q^{\pm} :=\text{\rm Ran}(\partial^{\pm}_{q}) =   \left(\Cal M^{\mp}_{q}\right)^{\perp} $
are closed in  $L^{2}_{q}(M)$;
\item the operators  $\partial^{\pm}_{q}: (H^{1}_{q},  \Cal Q) \to (R^{\pm}, \<\cdot,\cdot\>_{q})$ are isometric.
\end{enumerate}
\end{proposition} 
\medskip
\noindent Let $\Cal E= \{ e_{n} \,\vert \, n\in \N \}\subset H^{1}_{q}(M)\cap C^{\infty}(M)$ be an orthonormal basis of the Hilbert space $L^{2}_{q}(M)$ of eigenfunctions of the Dirichlet form~\pref{eq:Dirichlet}  and let $\lambda:\N\to \R^{+}\cup \{0\}$ be the corresponding sequence of eigenvalues:
\begin{equation}
\lambda_{n}:= \Cal Q(e_{n}, e_{n}) \,, \quad \text{\rm  for each }\,n \in \N\,.
\end{equation}
The {\it Friedrichs weighted Sobolev norm} $\Vert \cdot \Vert_{s}$ of order $s\in\R^{+}$ is the norm induced by the hermitian product defined as follows: 
for all $u$, $v\in L^{2}_{q}(M)$,
\begin{equation}
 \label{eq:extsnorm}
 (u,v) _{s} :=   \sum_{n\in \N} (1+\lambda_{n})^{s}\,
 \<u,e_{n}\>_{q}  \,  \<e_{n},v\>_{q} \,.
 \end{equation}
The inner products~\pref{eq:knorm} and~\pref{eq:extsnorm} induce {\it equivalent }Sobolev norms on the weighted Sobolev space $H^{k}_{q}(M)$, for all $k\in \Z^{+}$. In fact, the following result, a sharp 
version of Lemma 4.2 of \cite{F97},  holds:
 \begin{lemma} 
 \label{lemma:equivnorms}
 For each $k\in \Z^{+}$ there exists a constant~$C_{k}>1$ such that, for any orientable holomorphic  quadratic differential $q$ on $M$ and for all functions $u\in H^{k}_{q}(M)$,
 \begin{equation}
 \label{eq:equivnormstwo}
 C_{k}^{-1}\, \vert u \vert_{k}\,\, \leq \,\, \Vert u \Vert_{k} \,\, \leq \,\, C_{k} \, \vert u \vert_{k}\,.
 \end{equation}
\end{lemma}
\begin{proof}
By Proposition~\ref{prop:CR}, $(3)$, and Lemma~\ref{lemma:commut}, for all $u\in H^{k+1}_{q}(M)$,
$k\in \N$, the following identity holds (see  $(4.4)$ in~\cite{F97}):
\begin{multline} 
\label{eq:indstepone}
\vert u\vert_{k+1}^{2} =  \vert u\vert_{0}^{2}\,+\, \sum_{i+j\leq k} {\Cal Q} (S^{i}T^{j}u, S^{i}T^{j}u) 
= \\  = \vert u\vert_{0}^{2}\,+\, \sum_{i+j\leq k} \<S^{i}T^{j}\partial^{\pm}_{q}u, S^{i}T^{j}\partial^{\pm}_{q}u\>_{q}  =  \vert u\vert_{0}^{2} \, +\, \vert \partial^{\pm}_{q} u \vert_{k}^{2} \,\,.
\end{multline}
Hence in particular $u\in H^{k+1}_{q}(M)$, $k\in \N$, implies $\partial^{\pm}_{q} u\in H^{k}_{q}(M)$. 
If  $k\geq 1$,  a second application of the identity~\pref{eq:indstepone} yields
\begin{equation}
\label{eq:indsteptwo}
\vert u\vert_{k+1}^{2} = \vert u\vert_{1}^{2} \,+\, \vert \partial^{\mp}_{q} 
 \partial^{\pm}_{q} u \vert_{k-1}^{2}\,\,.
\end{equation}
 The statement then follows by induction on $k\in \N$. For $k=0$ it is immediate and for $k=1$,
 by the identity~\pref{eq:indstepone}, 
  \begin{equation}
  \label{eq:indhypone}
\vert u \vert_{1}^{2} = \vert u \vert_{0}^{2} \,+\, \Cal Q (u,u) = \sum_{n\in \N}(1+ \lambda_{n}) \, \vert \<u,e_{n}\>_{q}  \vert^{2}\,\,.
\end{equation}
For $k>1$, by induction hypothesis we can assume that the norms  $\vert \cdot \vert_{k-1}$
and $\Vert \cdot \Vert_{k-1}$ are equivalent, that is, there exists a constant~$C_{k-1} >1$ such 
that, for all $u\in H^{k-1}_{q}(M)$,
\begin{equation}
 \label{eq:indhyptwo}
 C_{k-1}^{-1}\, \vert u \vert_{k-1}\,\, \leq \,\, \Vert u \Vert_{k-1} \,\, \leq \,\, C _{k-1}\, \vert u \vert_{k-1}\,.
 \end{equation}
 Since by~\pref{eq:antisymm} for all $u$, $v \in H^{1}_{q}(M)$,  $\<\partial _{q}^{\pm}u,v\>_{q}=
 -\<u,\partial _{q}^{\mp}v\>_{q}$, the adjoint operator~$(\partial_{q}^{\pm})^{\ast} = -\partial_{q}^{\mp}$ 
 on~$H^{1}_{q}(M)$. By Proposition~\ref{prop:CR}, $(3)$, we have
\begin{multline}
\label{eq:CRest}
\Vert \partial^{\mp}_{q}  \partial^{\pm}_{q} u \Vert_{k-1}^{2} = \sum_{n\in \N}(1+ \lambda_{n})^{k-1} 
\, \vert \< \partial^{\pm}_{q}u,\partial^{\pm}_{q} e_{n}\>_{q}\vert^{2}=\\
  =\sum_{n\in \N}(1+ \lambda_{n})^{k-1} \, \vert \Cal Q (u,e_{n})\vert^{2}  = \sum_{n\in \N}(1+ \lambda_{n})^{k-1}  \lambda_{n}^{2} \, \vert \<u,e_{n}\>_{q}\vert^{2}\,\,.
\end{multline}
There exists a constant~$C_{k+1}>1$ such that, for all $\lambda\geq 0$,
\begin{multline}
\label{eq:constest}
C_{k+1}^{-2}\,(1+\lambda)^{k+1} \leq   1+\lambda +  C_{k-1}^{-2} \, \lambda^{2} (1+\lambda)^{k-1} \\
\leq 1+\lambda +  C^2_{k-1} \lambda^{2} \,(1+\lambda)^{k-1} \leq  C^2_{k+1}\,(1+\lambda)^{k+1} \,\,.
\end{multline}
By~\pref{eq:indsteptwo}, \pref{eq:indhypone},   \pref{eq:indhyptwo},  \pref{eq:CRest} and~\pref{eq:constest}, the estimate
\begin{equation}
 \label{eq:indhypthree}
 C_{k+1}^{-1}\, \vert u \vert_{k+1}\,\, \leq \,\, \Vert u \Vert_{k+1} \,\, \leq \,\, C _{k+1}\, \vert u \vert_{k+1}
 \end{equation}
follows, thereby completing the induction step.
\end{proof}

\subsection{Fractional Sobolev norms}
\label{FWSN}
Let $q$ be any orientable quadratic differential on $M$. For all $s\geq 0$, let 
\begin{equation}
\label{eq:Sobbar}
\bar H_{q}^{s} (M):= \{ u\in L^{2}_{q}(M)\,/\,  \sum_{n\in \N} (1+\lambda_{n})^{s}
\vert \<u,e_{n}\>_{q}\vert^{2} < +\infty\,\} ,
\end{equation}
endowed with the hermitian product given by \pref{eq:extsnorm} and, for any $s>0$, let 
$\bar H_{q}^{-s} (M)$ be the dual space of the Hilbert space $\bar H_{q}^{s} (M)$. The spaces 
$\bar H_{q}^{s} (M)$ will be called the {\it Friedrichs (fractional) weighted Sobolev spaces}.  

\smallskip
\noindent Let $H_{1}\subset H_{2}$ be Hilbert spaces such that $H_{1}$ embeds continuously into
$H_{2}$ with dense image. For all $\theta\in [0,1]$, let $[H_{1}, H_{2}]_{\theta}$ be the (holomorphic) interpolation space of  $H_{1}\subset H_{2}$ in the sense of Lions-Magenes  \cite{LiMa68}, Chap. 1, endowed with the canonical interpolation norm.  By the results of ~\cite{LiMa68}, Chap. 1, \S\S  2, 5, 6 
and 14, we have the following:
 
\begin{lemma} 
\label{lemma:extSobspaces} 
The Friedrichs weighted Sobolev spaces form an interpolation family  
$ \{ \bar H_{q}^{s} (M)\}_{s\in\R}$ of Hilbert spaces: for all $r$, $s\in \R$ with $r<s$,
\begin{equation}
\label{eq:IS}
 \bar H_{q}^{(1-\theta)r +\theta s} (M) \equiv  [ \bar H_{q}^{r} (M),  \bar H_{q}^{s} (M)]_{\theta} \,\,.
 \end{equation}
\end{lemma}

\medskip
\noindent The family $\{H^s_q(M)\}_{s\in \R}$ of {\it fractional weighted Sobolev spaces} will be 
defined as follows. Let $[s]\in \N$ denote the {\it integer part} and $\{s\} \in [0,1)$ the {\it fractional 
part }of any real number $s\geq 0$. 

\begin{definition}
\label{def:snorm}
\begin{enumerate}
\item The \emph{fractional weighted Sobolev norm} $\vert \cdot  \vert _{s}$ of order $s\geq 0$ is the euclidean norm induced by the hermitian product defined as follows: for all functions $u$, $v\in H^{\infty}_{q}(M)$,
\begin{equation}
 \label{eq:snorm}
 \< u,v \>_{s} :=   \frac{1}{2}\sum_{i+j\leq [s]} (S^{i}T^{j}u, S^{i}T^{j}v)_{\{s\}} + 
 (T^{i}S^{j}u, T^{i}S^{j}v)_{\{s\}}\,.
 \end{equation}
\item The  \emph{fractional weighted Sobolev norm }$\vert \cdot  \vert _{-s}$ of order $-s<0$ is defined as the dual norm of the weighted Sobolev norm $\vert \cdot  \vert _{s}$. 
\item The \emph{fractional weighted Sobolev space }$H^{s}_q(M)$ of order $s\in \R$ is defined as the completion with respect to the norm $\vert \cdot  \vert _{s}$ of the maximal common invariant domain $H^{\infty}_q(M)$.  
\end{enumerate}
\end{definition}
\noindent It can be proved that the weighted Sobolev space $H^{-s}_q(M)$ is isomorphic to the dual space of the Hilbert space $H^{s}_q(M)$, for all $s\in \R$. 
\begin{definition} 
\label{def:distsobord}
For any distribution $\Cal D$ on $M\setminus\Sigma_q$, the \emph{weighted Sobolev order }$
{\Cal O}^H_q(\Cal D)$ and the \emph{Friedrichs weighted Sobolev order }$\bar {\Cal O}^H_q(\Cal D)$ 
are the real numbers defined as follows:
\begin{equation}
\begin{aligned}
{\Cal O}^H_q(\Cal D):=& \inf\{ s\in\R \,\vert \, \Cal D \in H^{-s}_q(M)\}\,; \\
\bar {\Cal O}^H_q(\Cal D):=& \inf\{ s\in\R \,\vert \, \Cal D \in \bar H^{-s}_q(M)\}\,.
\end{aligned}
\end{equation}
\end{definition}

\smallskip
\noindent The definition of the fractional weighted Sobolev norms is  motivated by the following
basic result.

\begin{lemma} For all~$s\geq 0$, the restrictions of the Cauchy-Riemann operators 
$\partial^{\pm}_{q}:H^1_q(M) \to L^2_q(M)$ to the subspaces $H^{s+1}_q(M)\subset H^1_q(M)$ 
yield bounded operators 
$$
\partial^{\pm}_{s}:H_{q}^{s+1} (M) \to H_{q}^{s} (M)
$$
(which do not extend 
to operators  $\bar H_{q}^{s+1} (M) \to \bar H_{q}^{s} (M)$ unless $M$ is the torus). 
On the other hand, the Laplace operator 
\begin{equation}
\Delta_{q} =\partial^{+}_{q} \partial^{-}_{q} =\partial^{-}_{q} \partial^{+}_{q}: 
H^{2}_{q}(M) \to L^{2}_{q}(M)
\end{equation} 
yields a bounded operator $\bar \Delta_{s} :\bar H_{q}^{s+2}  (M) \to \bar H_{q}^{s} (M)$, defined
as the restriction of the Friedrichs extension $\Delta_q^F: \bar H_{q}^{2}  (M) \to L^2_q(M)$.
\end{lemma}
\begin{proof}  The restrictions  $\partial^\pm_s : H^{s+1}_q(M) \to H^s_q(M)$ of the Cauchy-Riemann
operators are well-defined and bounded   for all $s>0$ by definition of the Sobolev spaces 
$H^s_q(M)$.

\smallskip
\noindent The operators $\partial^{\pm}_s: H^{s+1}_q(M)\to H^s_q(M)$ do not extend to bounded operators $\bar H^{s+1}_q(M)\to \bar H^s_q(M)$ unless $M_q$ is a flat torus. In fact, every finite combination $f$ of eigenfunctions of the Dirichlet form belongs to $\bar H_{q}^{s} (M)$, for all $s\in \R$, but $\partial^{\pm}_{q}f \not\in H^{1}_{q}(M)$ in all cases because of the presence of obstructions in the Taylor expansion of eigenfunctions at the singular set $\Sigma_{q}\not=\emptyset$.  In fact, if $\partial^+_q e_n \in H^{1}_{q}(M)$ (or $\partial^-_q e_n \in H^{1}_{q}(M)$) for all $n\in \N$, then $\partial^+_q e_n \in H^{1}_{q}(M)$ and $\partial^-_q e_n \in H^{1}_{q}(M)$, for all  $n\in \N$, since the eigenfunctions $e_n$ can be chosen real. It follows that $e_n\in H^2_q(M)$, for all $n\in \N$,  hence $H_q^2(M)=\bar H^2_q(M)$ is the domain of the Friedrichs extension $\Delta^F_{q}$ of the Laplacian $\Delta _q$ of the metric $R_q$. Thus the Laplacian $\Delta _q$ is self-adjoint on the domain $H^2_q(M)$, hence the metric $R_q$ has no singularities and $M$ is the torus. In fact, the space $H^2_q(M)$ is the domain of the closure of the Laplacian $\Delta_q $ on the common invariant domain $H^{\infty}_q(M)$. If $H_q^2(M)=\bar H^2_q(M)$, then $\Delta_q$ is essentially self-adjoint on $H^{\infty}_q(M)$ and by \cite{Ne59}, Th. 5 or Cor. 9.1, the action of the commutative Lie algebra spanned by $\{S,T\}$ integrates to a Lie group action. Hence, the singularity set $\Sigma_q=\emptyset$ and $M$ is the torus.

\smallskip
\noindent Finally, the Friedrichs extension  $\Delta^F_{q}$, defined on $\bar H^2_q(M)$, has a 
bounded restriction $\bar \Delta_{s}:\bar H_{q}^{s+2} (M) \to \bar H_{q}^{s} (M)$,  for all $s\geq 0$. 
In fact, we have
\begin{equation}
\Delta_q^F u = \sum_{n\in \N}  \lambda_n  \, \<u,e_n\>_q \,e_n\,\,,  \quad \text { for all } 
u\in \bar H^2_q(M)\,,
\end{equation}
hence $\Delta_q^F u  \in \bar H_{q}^s(M)$ if  $u\in \bar H_{q}^{s+2} (M) $, by definition of the Friedrichs weighted Sobolev spaces $\bar H_{q}^s (M)$ (in terms of eigenfunction expansions for the Dirichlet form).
\end{proof}

\smallskip
\noindent 
\begin{lemma} 
\label{lemma:intineq} The fractional weighted Sobolev norms satisfy the following interpolation inequalities. For any $ 0\leq r < s $ there exists a constant $C_{r,s}>0$ such that, for any $\theta\in [0,1]$ and any function $u \in H^{s}_q(M)$,
\begin{equation}  
\label{eq:intineq}
\vert u \vert_{(1-\theta) r +\theta s} \,\,    \leq  \,\, C_{r,s}\, \vert u\vert_{r} ^{1-\theta} \, 
\vert u \vert_{s} ^{\theta}\,.
\end{equation}
\end{lemma}

\begin{proof}
The argument will be carried out in three steps: $(1)$ the open interval $(r,s)$ does not contain
integers; $(2)$ the open interval $(r,s)$ contains a single integer; $(3)$ the general case.

\smallskip
\noindent In case $(1)$ there exists $k\in \N$ such that $ k\leq r <s \leq k+1$. The interpolation 
inequality follows from the definition \pref{eq:snorm} of the euclidean product which induces the fractional Sobolev norms, from the interpolation inequality for Friedrichs weigthed Sobolev norms 
(which are by definition interpolation norms) and from the H\"older inequality. In fact,  since 
$0\leq r-k \leq s-k\leq 1$ and $\theta\in (0,1)$,  the fractional part
\begin{equation}
\{(1-\theta) r +\theta s\} = (1-\theta) (r-k)+\theta (s-k)\,,
\end{equation}
hence, by the interpolation inequality for Friedrichs norms (see for instance  \cite{LiMa68}, 
Chap. 1, \S 2.5), for all $i+j \leq k$ the following estimates hold:
\begin{equation}
\label{eq:intineqone}
\begin{aligned}
\vert S^{i}T^{j}u \vert_{\{(1-\theta) r +\theta s\}} &
\leq   \vert S^{i}T^{j}u \vert^{1-\theta}_{r-k} \, \vert S^{i}T^{j}u \vert^{\theta}_{s-k}\,; \\
\vert T^{i}S^{j}u \vert_{\{(1-\theta) r +\theta s\}} &
\leq   \vert T^{i}S^{j}u \vert^{1-\theta}_{r-k}  \, \vert T^{i}S^{j}u \vert^{\theta}_{s-k}\,.
\end{aligned}
\end{equation}
By the definition \pref{eq:snorm} of the Sobolev norms, the interpolation inequality \pref{eq:intineq} follows from \pref{eq:intineqone} by H\"older inequality.

\smallskip
\noindent In case $(2)$, there exists $k\in \N$ such that $k-1\leq r < k < s \leq k+1$. We claim that, 
for any $u\in H^s_q(M)$, 
\begin{equation}
\label{eq:intineqclaim}
\vert u \vert_k\,\,    \leq  \,\, C_{r,s} \, \vert u\vert_{r} ^{\frac{s-k}{s-r}} \, 
\vert u \vert_{s} ^{\frac{k-r}{s-r}}\,.
\end{equation}
Let us prove that step $(1)$ and the above claim \pref{eq:intineqclaim} imply step $(2)$.
Let $\sigma= (1-\theta) r +\theta s$. We will consider only the case when $\sigma \in (r,k)$ since
the case when $\sigma\in (k,s)$ is similar. By step $(1)$ we have the inequality
\begin{equation}
\label{eq:intineqtwo}
\vert u \vert_\sigma\,\,    \leq  \,\, C_r \,  \vert u\vert_{r} ^{\frac{k-\sigma}{k-r}} \, 
\vert u \vert_{k} ^{\frac{\sigma-r}{k-r}}\,.
\end{equation}
By the claim \pref{eq:intineqclaim} it then follows that
\begin{equation}
\vert u \vert_\sigma\,\,    \leq  \,\,  C^{(1)}_{r,s} \,  \vert u\vert_{r}^{\frac{k-\sigma}{k-r}+
 \frac{\sigma-r}{k-r} \frac{s-k}{s-r}   } \, \vert u \vert_{s} ^{\frac{\sigma-r}{k-r} \frac{k-r}{s-r}} \,.
\end{equation}
It is immediate to verify that
$$
\frac{k-\sigma}{k-r}+
 \frac{\sigma-r}{k-r} \frac{s-k}{s-r} = \frac{s-\sigma}{s-r}\,.
$$ 
Let us turn to the proof of the claim \pref{eq:intineqclaim}. Since $-1\leq r-k<0< s-k \leq 1$ and
the weighted Sobolev norm $\vert \cdot \vert_s$ coincides with the Friedrichs norm  $\Vert \cdot 
\Vert_s$ for any $s\in [-1, 1]$, by the interpolation inequality for the Friedrichs Sobolev norms, the following estimates hold:  for all $i$, $j\in \N$ with $i+j\leq k$ and for any function $u\in H^s_q(M)$,
\begin{equation}
\label{eq:intineqthree}
\begin{aligned}
\vert S^{i} T^j u\vert_0 \,\, &\leq   \,\,C^{(2)}_{r,s} \, \vert S^{i} T^j u\vert_{r-k}^{\frac{s-k}{s-r}} \,\,   
\vert S^{i} T^j u\vert_{s-k} ^{\frac{k-r}{s-r}}\,; \\
\vert T^{i} S^j u\vert_0 \,\, &\leq   \,\, C^{(2)}_{r,s}\,  \vert T^{i} S^j u\vert_{r-k} ^{\frac{s-k}{s-r}} \,\,   
\vert T^{i} S^j u\vert_{s-k} ^{\frac{k-r}{s-r}}\,.
\end{aligned}
\end{equation}
Since the operators $S$, $T: \bar H^1_q(M) \to L^2_q(M)$ are well-defined, bounded and
have bounded linear extensions $L^2_q(M) \to \bar H^{-1}_q(M)$, by the fundamental theorem of interpolation (see for instance \cite{LiMa68}, Chap. 1, \S 5.1), the operators $S$, $T: \bar H^s_q(M)
\to \bar H^{s-1}_q(M)$ are well-defined and bounded for any $s\in [0,1]$. It follows that there exists a constant $C_r>0$ such that 
\begin{equation}
\label{eq:intineqfour}
\begin{aligned}
\sum_{i+j\leq k} \vert S^{i} T^j u\vert^2_{r-k} & \leq (C'_r)^2 \,  \sum_{i+j\leq k-1} 
\vert S^{i} T^j u\vert^2_{r-(k-1)} \,; \\
\sum_{i+j\leq k} \vert T^{i} S^j u\vert^2_{r-k} & \leq (C'_r)^2 \,  \sum_{i+j\leq k-1} 
\vert T^{i} S^j u\vert^2_{r-(k-1)} 
\end{aligned}
\end{equation}
The claim \pref{eq:intineqclaim} then follows by H\"older inequality from  
\pref{eq:intineqthree} and \pref{eq:intineqfour}.

\smallskip
\noindent In general, let $k_1 <k_2$ be positive integers such that
$$
k_1-1 \leq r < k_1 <k_2 <s \leq k_2+1\,.
$$
By Lemma \ref{lemma:equivnorms}, since the Friedrichs norms are interpolation norms,
we have that there exists a constant $C_{k_1,k_2}>0$ such that, for all $k\in \N \cap [k_1, k_2]$,
\begin{equation}
\label{eq:intineqfive}
\vert u \vert_{k}\,\,    \leq  \,\, C_{k_1,k_2}\,  \vert u\vert_{k_1}^{ \frac{k_2-k}{k_2-k_1} } \,
\vert u \vert_{k_2} ^{\frac{k-k_1}{k_2-k_1}}\,.
\end{equation}
By step $(2)$ there exists a constant $C^{(3)}_{r,s}>0$ such that 
\begin{equation}
\label{eq:intineqsix}
\begin{aligned}
\vert u \vert_{k_1}\,\,    \leq  \,\, C^{(3)}_{r,s}\,   \vert u\vert_{r}^{ \frac{1}{k_1+1-r}} \,
\vert u \vert_{k_1+1} ^{ \frac{k_1-r}{k_1+1-r}}\,; \\
\vert u \vert_{k_2}\,\,    \leq  \,\,  C^{(3)}_{r,s}\,  \vert u\vert_{k_2-1}^{ \frac{s-k_2}{s-k_2+1}} \,
\vert u \vert_{s} ^{ \frac{1}{s-k_2+1}}\,.
\end{aligned}
\end{equation}
The estimates in \pref{eq:intineqsix} imply, by bootstrap-type estimates based on  \pref{eq:intineqfive} 
for $k=k_1+1$ and $k=k_2-1$, that there exists a constant $C^{(4)}_{r,s}>0$ 
 such that 
\begin{equation}
\label{eq:intineqseven}
\begin{aligned}
\vert u \vert_{k_1}\,\,    \leq  \,\, C^{(4)}_{r,s}\,   \vert u\vert_{r}^{ \frac{k_2-k_1}{k_2-r}} \,
\vert u \vert_{k_2} ^{ \frac{k_1-r}{k_2-r}}\,; \\
\vert u \vert_{k_2}\,\,    \leq  \,\,  C^{(4)}_{r,s}\,  \vert u\vert_{k_1}^{ \frac{s-k_2}{s-k_1}} \,
\vert u \vert_{s} ^{ \frac{k_2-k_1}{s-k_1}}\,.
\end{aligned}
\end{equation}
By \pref{eq:intineqseven} and again by bootstrap, there exists a constant $C^{(5)}_{r,s}>0$  
such that
\begin{equation}
\label{eq:intineqeight}
\begin{aligned}
\vert u \vert_{k_1}\,\,    \leq  \,\, C^{(5)}_{r,s}\,   \vert u\vert_{r}^{ \frac{s-k_1}{s-r}} \,
\vert u \vert_{s} ^{ \frac{k_1-r}{s-r}}\,; \\
\vert u \vert_{k_2}\,\,    \leq  \,\,  C^{(5)}_{r,s}\,  \vert u\vert_{r}^{ \frac{s-k_2}{s-r}} \,
\vert u \vert_{s} ^{ \frac{k_2-r}{s-r}}\,,
\end{aligned}
\end{equation}
Let $\sigma\in (r,s)$. We have proved the interpolation inequality for the
subcases $\sigma= k_1$ and $\sigma=k_2$. Let us prove that the general case can be
reduced to these subcases.  If $\sigma \in (r, k_1)$, by step $(1)$ there exists
$C''_r>0$ such that
\begin{equation}
\label{eq:intineqnine}
\vert u \vert_{\sigma}\,\,    \leq  \,\, C''_r  \vert u\vert_{r}^{ \frac{k_1-\sigma}{k_1-r}} \,
\vert u \vert_{k_1} ^{ \frac{\sigma-r}{k_1-r}}\,.
\end{equation}
The interpolation inequality in this case follows immediately from \pref{eq:intineqeight} and
\pref{eq:intineqnine}. If $\sigma\in (k_2,s)$, the argument is similar. If $\sigma\in (k_1,k_2)$, then 
by step $(1)$, there exists $C_{[\sigma]}>0$ such that 
\begin{equation}
\label{eq:intineqten}
\vert u \vert_{\sigma}\,\,    \leq  \,\, C_{[\sigma]}  \vert u\vert_{[\sigma]}^{ 1-\{\sigma\} } \,
\vert u \vert_{[\sigma]+1} ^{ \{\sigma\}}\,.
\end{equation}
The interpolation inequality then follows from \pref{eq:intineqfive}, \pref{eq:intineqeight} and \pref{eq:intineqten}.

\end{proof}

\smallskip
\noindent Let $H^{s}(M)$, $s\in \R$, denote a family of standard Sobolev spaces on the compact 
manifold $M$ (defined with respect to a Riemannian metric).  The comparison lemma below clarifies
to some extent the relations between the different scales of fractional Sobolev spaces. 

\begin{lemma}  
\label{lemma:comparison}
The following continuous embedding and isomorphisms of Banach spaces hold:
\begin{enumerate}
\item $ \,\, H^{s}(M) \,\, \subset   \,\, H_{q}^{s}(M) \,\, \equiv \,\,  \bar H_{q}^{s}(M) \,,  
\quad\text{for }0\leq s<1$;
\item $\,\,H^{s}(M) \,\, \equiv   \,\, H_{q}^{s}(M) \,\, \equiv \,\,  \bar H_{q}^{s}(M)\,,
\quad\text{for }s=1$;
\item $\,\,H_{q}^{s}(M) \,\, \subset  \,\, \bar H_{q}^{s}(M) \,\, \subset  \,\,  H^{s}(M)\,, 
\quad\text{for }s >1$.
\end{enumerate}
For $s \in [0,1]$, the space $H^{s}(M)$ is dense in $H_{q}^{s}(M)$ and, for $s >1$, the closure 
of $H_{q}^{s}(M)$ in $\bar H_{q}^{s}(M)$ or $H^s(M)$ has finite codimension.   
\end{lemma}
\begin{proof}
By definition $H^{0}(M)=L^{2}(M)$ and $H^{0}_{q} (M)=\bar H^{0}_{q} (M)=L^{2}_{q}(M)$. Since
the area form induced by any quadratic differential is smooth on $M$,  which is a compact surface, it follows that~$L^{2}(M)\subset L^{2}_{q}(M)$. The embedding $H^{1}_{q} (M) \subset \bar H^{1}_{q} (M)$
follows by Lemma~\ref{lemma:equivnorms} and the embedding $\bar H^{1}_{q} (M) \subset H^{1}_{q} (M)$ holds since the eigenfunctions of the Dirichlet form are in $H^{1}_{q} (M)$. The isomorphism
$H^{1}(M) \,\equiv   \,H_{q}^{1}(M)$ is proved in~\cite{F02}, \S 6.2. Hence $(2)$ is proved and $(1)$
follows by interpolation.

\smallskip
\noindent  Let $s>1$. If $[s]=2k$ is even, there exists a constant $A_k>0$ such that, for all functions 
$u\in H^{\infty}_q(M)$, we have
\begin{equation}
\Vert u \Vert^2_s = \Vert (I-\Delta^F_q)^k u\Vert^2_{\{s\}}  \leq A^2_k \sum_{i+j\leq 2k} 
\Vert S^{i}T^{j}u \Vert^2_{\{s\}} = C^2_k \,\vert u\vert_s\,.
\end{equation}
If $[s]=2k+1$ is odd, we argue as follows. The Cauchy-Riemann operators $\partial^\pm_q:H^1_q(M)
\to L^2_q(M)$ are bounded and extend by duality to bounded operators $\partial^\pm_0: L^2_q(M) \to H^{-1}_q(M)$. Hence, by the fundamental theorem of interpolation (see \cite{LiMa68}, Chap. 1, \S 5.1), for all $\sigma \in [0,1]$ the Cauchy-Riemann operators have bounded restrictions
\begin{equation}
\partial^{\pm} _\sigma: \bar H^\sigma_q(M) \to \bar H^{\sigma-1}_q(M)\,.
\end{equation}
It follows that there exists a constant $B_k>0$ such that, for all functions $u\in H^{\infty}_q(M)$, 
we have
\begin{equation}
\Vert u \Vert^2_s = \Vert (I-\Delta^F_q)^{k+1} u\Vert^2_{\{s\}-1}  \leq B^2_k \sum_{i+j\leq 2k+1} 
\Vert S^{i}T^{j}u \Vert^2_{\{s\}} = C^2_k \,\vert u\vert_s\,.
\end{equation}
Thus the embeddings $H^s_q(M) \subset \bar H^s_q(M)$ hold for all $s>1$.

\smallskip
\noindent It was proved in~\cite{F02}, \S 6.2, that $H_{q}^{k}(M) \subset  H^{k}(M)$, for all $k\in \Z^{+}$. 
We prove below the stronger statement that $\bar H_{q}^{s}(M) \subset  H^{s}(M)$, for all $s\in \R^{+}$. Let $R$ be a  smooth Riemannian metric on $M$ conformally equivalent to the degenerate metric $R_q$ and let $H_R^s(M)$, $s\geq 0$, denote the Sobolev spaces of the Riemannian manifold 
$(M,R)$ which are defined as the domains of the powers of  the essentially self-adjoint Laplacian $\Delta_R$ of the metric.  Since $M$ is compact, the Sobolev spaces $H_R^s(M)\equiv H^s(M)$ are independent, as topological vector spaces, of the choice of the Riemannian metric $R$, for all $s\in\R$. We claim that  $\bar H_{q}^{2k}(M) \subset  H_R^{2k}(M)$,  for all $k\in \Z^+$. In fact, there exists a smooth non-negative real-valued function $W$ on $M$ (vanishing only at $\Sigma_q$) such that  
$W  \Delta^F_q \subset \Delta_R$. Let $W$ be the unique function such that the area forms of the metrics are related by the identity $\omega_q= W \omega_R$. If  $u\in \bar H_{q}^{2} (M)$, then $ \Delta^F_q u \in L^2_q(M)$, so that 
\begin{equation}
\Delta_R u = W \Delta^F_q u \in L^2(M,\omega_R)\,.
\end{equation}
Let us assume that $\bar H_{q}^{2k-2}(M) \subset  H_R^{2k-2}(M)$ and let $u\in \bar H_{q}^{2k}(M)$.
We have
\begin{equation}
\Delta_R ^k u = \Delta_R ^{k-1} W \Delta^F_q u = [ \Delta_R ^{k-1}, W]  \Delta^F_q u + W 
 \Delta_R ^{k-1}  \Delta^F_q u\,.
 \end{equation}
 Since the commutator $[ \Delta_R ^{k-1}, W]$ and $\Delta_R ^{k-1}$ are differential operators of 
 order $2k-2$ on $M$ and $\Delta^F_q u \in H_R^{2k-2}(M)$ by the induction hypothesis, the function
$\Delta_R ^k u \in L^2(M,\omega_R)$. The claim is therefore proved. It follows by interpolation that  $\bar H_{q}^{s}(M)  \subset  H_R^{s}(M)$,  for all $s\geq 1$. Thus $(3)$ is proved. 

\smallskip
\noindent For $s\in [0,1]$, the space $C_0^{\infty}(M\setminus\Sigma_q)\subset H^s_q(M)$ is dense
in $H^s_q(M)$. For $s>1$, the subset $C^{\infty}(M)\cap \bar H^{s}_{q}(M)$ is dense in $\bar H^{s}_{q}(M)$, since the eigenfunctions of the Dirichlet form (hence all finite linear combinations) 
belong to $C^{\infty}(M)$ and the space $C^{\infty}(M)$ is dense in $H^s(M)$. Finally, the subspace~$C^{\infty}(M) \cap H^{k}_{q}(M) \subset C^{\infty}(M)$ can be described, for any $k\in \N$, 
as the kernel of a finite number of distributions of finite order supported on the finite set $\Sigma_{q}$ 
(see~\cite{F02}, $(7.9)$), hence for any $k>s$ the closure of $H_{q}^{k}(M)\subset H^s_q(M)$ in~$\bar H_{q}^{s}(M)$ or in $H^{s}(M)$ has finite codimension. 
\end{proof}

\subsection{Local analysis}
\label{LA}
For each $p\in M$ and all $k\in \Z^+$, let $H^k_q(p)$, $\bar H^k_q(p)$, and $H^k(p)$ 
the spaces of germs of fuctions at $p$ which belong to $H^k_q(M)$, $\bar H^k_q(M)$ 
and $H^k(M)$, endowed with the respective direct limit topologies. More precisely, a germ of 
function $f$ at $p$ belongs to the space $H^k_q(p)$, $\bar H^k_q(p)$ or $H^k(p)$ iff it can be
realized by a function $F$ on $M$ which belongs to the space $H^k_q(M)$, $\bar H^k_q(M)$ or $H^k(M)$ respectively and the open sets in $H^k_q(p)$, $\bar H^k_q(p)$ or $H^k(p)$  are defined 
as the images of open sets in $H^k_q(M)$, $\bar H^k_q(M)$ or $H^k(M)$ under the natural
maps $H^k_q(M)\to H^k_q(p)$, $\bar H^k_q(M) \to \bar H^k_q(p)$ or $H^k(M)\to H^k(p)$.

\noindent By Lemma \ref{lemma:comparison} we have the inclusions 
\begin{equation}
\label{eq:germsinclusions}
\begin{aligned}
H^0(p) &\subset H^0_q(p) \subset \bar H^0_q(p)\,; \\
H^1_q(p)&=\bar H^1_q(p)=H^1(p)\,;\\
H^k_q(p)&\subset \bar H^k_q(p) \subset H^k(p)\,.
\end{aligned}
\end{equation}
If $p\not\in \Sigma_q$, since there is an open neighbourhood $D_p$ of $p$ in $M$ isomorphic 
to a flat disk and the operator $\Delta^F_q$ is elliptic of order $2$ on $D_p$ (isomorphic to
the flat Laplacian), all the inclusions in \pref{eq:germsinclusions} are identities. We will describe
precisely the inclusions $H^k_q(p)\subset \bar H^k_q(p) \subset H^k(p)$ for $k>1$.

\smallskip
\noindent Let $p\in \Sigma_q$ be a zero of (even) order $2m$ of the (orientable) quadratic differential $q$ on $M$. There exists a unique canonical holomorphic coordinate $z:D_p\to \C$, defined 
on a neighbourhood $D_p$ of $p\in M$, such that $z(p)=0$ and $q(z)=z^{2m} dz^2$.  With 
respect to the canonical coordinate the Cauchy-Riemann operators $\partial^\pm_q$ can 
be written in the following form:
\begin{equation}
\label{eq:CRlocal}
\partial^+_q = \frac{2}{\bar z^{m}}\, \frac{\partial}{\partial \bar z} \qquad \text{ and } \qquad
\partial^-_q =\frac{2}{z^{m}}\, \frac{\partial}{\partial z}\,\,.
\end{equation}
Let $C^{\infty}(p)$ be the space of germs at $p\in M$ of smooth complex-valued functions on $M$
and for any $u\in C^{\infty}(p)$ let  
\begin{equation}
u(z,\bar z) = \sum_{i,j\in \N}  a_{ij}(u,p) z^{i} \bar z^{j}\,.
\end{equation}
be its (formal) Taylor series at $p$ (with respect to the canonical coordinate).

\begin{lemma} 
\label{lemma:localbarH}
Let $p\in \Sigma_q$ be a zero of (even) order $2m$ of the (orientable) quadratic differential $q$ on
$M$. For any $k\in \N$, a germ 
$$
u \in C^{\infty}(p) \cap \bar H^k_q(p) \Leftrightarrow a_{ij}(u,p) =0\,, 
\quad\text{\rm for all } i+j \leq (k-1)(m+1)\,,
$$ 
except all pairs $(i,j)$ for which one of the following conditions holds:
\begin{equation}
\begin{aligned}
&(1)\,\, i \in \N\cdot (m+1)\,, \,\,j \in \N\cdot (m+1)\,;\\
&(2)\,\, i \in \N\cdot (m+1)\, ,\,\, j\not \in  \N\cdot (m+1) \,\,\,\,\text {\rm and }\,\,i<j\,;\\
&(3)\,\, i\not \in  \N\cdot (m+1)\, , \,\,j\in \N\cdot (m+1) \,\,\,\,\text {\rm and }\,\,i>j\,.
\end{aligned}
\end{equation}
\end{lemma}
\begin{proof}
Let $u\in C^{\infty}(p)$. For any $n\in \N$, there is a local Taylor expansion 
$$
u(z,\bar z) = \sum_{i+j\leq n} a_{ij}(u,p) z^{i} \bar z^{j}  \, +\, R^u_n(z,\bar z)
$$
where the remainder $R^{u}_n$ is a smooth function vanishing at order $n$ at $p$. A straightforward
calculation (based on formulas \pref{eq:CRlocal} yields that any smooth function $R$ vanishing 
at $p$ at order $n$ belongs to the space $H^k_q(p)\subset \bar H^k_q(p)$ if $n>(k-1)m$. It follows 
that $u\in \bar H^k (p)$ iff its Taylor polynomial of any order $n>(k-1)m$ does.  The argument
can therefore be reduced to the case of polynomials.

\smallskip
\noindent  It follows from formulas \pref{eq:CRlocal} that, for all $\ell\in \N$ and all $(i,j)\in \N\times \N$, there exists a complex constant $c^{m,\ell}_{ij}$ such that 
\begin{equation}
\label{eq:lapllocfor}
\Delta^{\ell}_q (z^{i} \bar z^{j}) = c^{m,\ell}_{ij}\,  z^{i-\ell(m+1)} \, \bar z^{j-\ell(m+1)} \,.
\end{equation}
The area form of the quadratic differential  $q$ can be written as 
\begin{equation}
\label{eq:locarea}
\omega_q= \vert z \vert ^{2m} dx\wedge dy
\end{equation}
with respect to the canonical coordinate $z:=x+\imath y$. Hence, straightforward computations 
in polar coordinates yield that, if $c^{m,\ell}_{ij}\not =0$,
\begin{equation}
\label{eq:germLaplacian}
\begin{aligned}
\Delta^{\ell}_q (z^{i} \bar z^{j})  \in H^0_q(p)  &\Leftrightarrow    i+j- 2\ell (m+1) > -(m+1)\,;\\
 \Delta^{\ell}_q (z^{i} \bar z^{j})  \in H^1_q(p)  &\Leftrightarrow     i+j- 2\ell (m+1) > 0\,.
\end{aligned}
\end{equation}
If $c^{m,\ell}_{ij}=0$, then either $i\in \N\cdot (m+1)$ and $i<\ell(m+1)$ or $j\in \N\cdot (m+1)$ and $j<\ell(m+1)$. It follows that, if $i$, $j\not\in  \N\cdot (m+1)$, then
\begin{equation}
\label{eq:ijkcondone}
z^{i} \bar z^{j} \in \bar H^k_q(p)  \Leftrightarrow   i+j- (k-1)(m+1) > 0\,.
\end{equation}
If $i \in  \N\cdot (m+1)$, $i=h(m+1)$, and $j\not\in  \N\cdot (m+1)$, then conditions 
\pref{eq:germLaplacian} apply for all $\ell \leq h$, hence \pref{eq:ijkcondone} holds 
 if $k\leq 2h$, while if $k>2h$,
\begin{equation}
\label{eq:ijkcondtwo}
z^{i} \bar z^{j} \in \bar H^k_q(p)  \Leftrightarrow  \Delta^{h}_q (z^{i} \bar z^{j})  \in H^1_q(p)
\Leftrightarrow j> i\,.
\end{equation}
Similarly, if $j\in  \N\cdot (m+1)$, $j=h(m+1)$, and $i\not\in  \N\cdot (m+1)$, then 
 \pref{eq:ijkcondone} holds  if $k\leq 2h$, while if $k>2h$,
\begin{equation}
\label{eq:ijkcondthree}
z^{i} \bar z^{j} \in \bar H^k_q(p)  \Leftrightarrow  \Delta^{h}_q (z^{i} \bar z^{j})  \in H^1_q(p)
\Leftrightarrow i> j\,.
\end{equation}
It follows immediately from \pref{eq:ijkcondone}, \pref{eq:ijkcondtwo} and \pref{eq:ijkcondthree} 
that the conditions listed in the statement of the lemma are sufficient. The necessity follows from
the following argument. For any $r_1<r_2$, let $D(r_1,r_2)\subset M$ be the annulus (centered 
at $p$) defined by the inequalities $r_1<\vert z\vert <r_2$.  The system of Laurent monomials 
$\{z^{i}\bar z^{j}\vert i,j \in \Z\}$ is orthogonal in $\bar H^k(D_{r_1,r_2})$. In fact, a computation 
in polar coordinates shows that, for all $(i,j)\not= (i',j')\in \Z\times \Z$,
\begin{equation}
\int_{D(r_1,r_2)}  \Delta_q^\ell (z^{i}\bar z^{j})   \Delta_q^{\ell'} (\bar z^{i'}z^{j'})\,\,\omega_q \,=\,0\,,
\end{equation}
hence 
\begin{equation}
\Vert \sum_{i+j\leq n} a_{ij}  z^{i} \bar z^{j}  \Vert^2_k =  \sum_{i+j\leq n} \vert a_{ij} \vert^2 
\, \Vert  z^{i} \bar z^{j}\Vert^2_k\,.
\end{equation}
It follows that  only Laurent monomials $z^{i}\bar z^{j}\in \bar H^k(p)$ can appear in the Taylor expansion of a function $f\in C^{\infty}(p)\cap \bar H^k(p)$.
\end{proof}

\begin{lemma} 
\label{lemma:localH}
Let $p\in \Sigma_q$ be a zero of order $2m$ of the quadratic differential $q$ on $M$. For any 
$k\in \N$, a germ 
$$
u \in C^{\infty}(p) \cap H^k_q(p) \Leftrightarrow a_{ij}(u,p) =0\,, 
\quad\text{\rm for all } i+j \leq (k-1)(m+1)\,,
$$ 
except all pairs $(i,j) \in \N\cdot (m+1) \times \N\cdot (m+1)$.
\end{lemma}

\begin{proof} The proof is similar to that of Lemma \ref{lemma:localbarH} above. In fact, 
formulas \pref{eq:lapllocfor} are replaced by the following formulas. For all $\alpha\in \N$
and all $i\in \N$, there exists a complex constant $c^{m,\alpha}_{i}$ such that 
\begin{equation}
\label{eq:CRplusminus}
(\partial_q^+)^\alpha (\partial_q^-)^\beta z^{i} \bar z^{j}  =     c^{m,\alpha}_{i} 
\bar c^{m,\beta}_{j}\,  z^{i-\alpha(m+1)} \, \bar z^{j-\beta(m+1)} \,.
\end{equation}
As in the proof of Lemma \ref{lemma:localbarH}, it follows by a straightforward computation in polar coordinates that, if $c^{m,\alpha}_{i} \bar c^{m,\beta}_{j}\not =0$, 
\begin{equation}
\label{eq:germCR}
(\partial^+)^\alpha (\partial^-)^\beta z^{i} \bar z^{j}  \in H^0_q(p) \Leftrightarrow
i+j -(\alpha+\beta)(m+1) > -(m+1)\,.
\end{equation}
Since there exists $\alpha\in \N$ such that $c^{m,\alpha}_{i} =0$ iff $i\in \N\cdot (m+1)$,
either  $i$, $j \in \N\cdot (m+1)$, in which case $z^{i} \bar z^{j} \in H^k_q(p)$, or there 
exists $(\alpha, \beta)$ such that $\alpha+\beta =k$ and $c^{m,\alpha}_{i} 
\bar c^{m,\beta}_{j}\not=0$, in which case $z^{i} \bar z^{j}\in H^k_q(p)$ iff $i+j>(k-1)(m+1)$.
\end{proof}


\noindent  Let $p\in \Sigma_q$ be a zero of order $2m_p$ of the orientable quadratic differential 
$q$. Let $\Cal T_p \subset \N\times\N$ be the set of $(i,j)$ such that
\begin{equation}
\label{eq:CalT}
\begin{aligned}
&i\in \N\cdot (m_p+1)\,,\,\, j\not \in \N\cdot (m_p+1) \quad \text{and} \quad i< j \quad \text{or} \\
&i\not \in \N\cdot (m_p+1)\,,\,\, j \in \N\cdot (m_p+1)  \quad \text{and} \quad i>j , 
\end{aligned}
\end{equation}
For any $(i,j)\in \Cal T_p$, let $\delta^{ij}_p$ be the linear functional (distribution) on $C^{\infty}(M)$ defined as follows. Let $u(z,\bar z) = \sum a_{ij}(u,p) z^{i} \bar z^{j} $ denote the Taylor expansion of $u\in C^{\infty}(M)$ at $p\in \Sigma_q$ with respect to the canonical coordinate $z:D_p\to \C$ for the differential $q$ at $p$. Let 
\begin{equation}
\delta^{ij}_p(u) :=  a_{ij}(u,p)\,.
\end{equation}
\noindent It is clear from the definition that $\delta^{ij}_p= \overline{\delta^{ji}_p}$ for all $(i,j)\in \Cal T_p$.  A calculation shows that, for any $h\in\N\setminus \N\cdot (m_p+1)$ we have the following representation in terms of the Cauchy principal value: for any $u\in C^{\infty}(p)$, 
\begin{equation}
\label{eq:deltazeros}
\begin{aligned}
\delta^{h0}_p(u) &=-\frac{1}{4\pi h}\, \text{\rm PV} \int_M   \frac{\Delta_q u}{ z^h} \, \omega_q\, ;\\
\quad\delta^{0h}_p(u) &= -\frac{1}{4\pi h} \, \text{\rm PV} \int_M  
 \frac{\Delta_q u}{ \bar z^h} \, \omega_q\,  .
\end{aligned}
\end{equation}
(The above formulas can be derived one from the other by conjugation).  In addition, for any $\ell\in \N$, by formulas~\pref{eq:lapllocfor} there exist complex constants $c^{m_p,\ell}_{0,h}\not=0$ and $c^{m_p,\ell}_{h,0}\not=0$  such that the following identities hold in the sense of distributions:
\begin{equation}
\label{eq:Laplondeltas}
\begin{aligned}
c^{m_p,\ell}_{0,h}\,\delta_p^{\ell(m_p+1),\ell(m_p+1)+h} &=  \Delta_q^\ell \left(\delta_p^{0,h}\right)\,;  \\  c^{m_p,\ell}_{h,0}\,\delta_p^{\ell(m_p+1)+h,\ell(m_p+1)} &=  \Delta_q^\ell  \left( \delta_p^{h,0} \right)\,.
\end{aligned}
\end{equation}

\smallskip
\noindent Let $\Cal T_p^k \subset \Cal T_p$ be the subset of $(i,j)$ such that $i+j\leq (k-1)(m_p+1)$.
\begin{lemma}
\label{lemma:continuity} For each $(i,j)\in \Cal T_p^k$, the functional $\delta^{ij}_p$ has a unique (non-trivial) continuous extension to the space $\bar H^k_q(p)$ and the following holds:
\begin{equation}
\label{eq:Hk}
H^k_q(p) =  \{ u\in \bar H^k_q(p)\,\vert \, \delta^{ij}_p(u)=0 \quad \text{\rm for all } (i,j)\in \Cal T_p^k\,\}\,.
\end{equation}
\end{lemma}
\begin{proof} The functions $z^{-h}$ and $\bar z^{-h} \in L^2_q(D) $ for all $1\leq h \leq m_p$ and the operator $\Delta_q^F: \bar H^2_q(p) \to L^2_q(p)$ is bounded. Hence the linear functionals $\delta^{0h}_p$ and $\delta^{h0}_p$ are continuous on $\bar H^2_q(p)$ for all $1\leq h \leq m_p$.  
 Similarly, the distributions  $\text{\rm PV} (z^{-h})$ and $\text{\rm PV}(\bar z^{-h}) \in 
H^{-1}(D_p)$ for all $m_p < h <2(m_p+1)$ and the operator $\Delta_q^F: \bar H^3_q(p) \to H^1(p)$ is bounded. Hence the linear functionals $\delta^{0h}_p$ and $\delta^{h0}_p$ are continuous on $\bar H^3_q(p)$ for all $m_p< h < 2(m_p+1)$.  Since the space $H^k_q(p)$ is equal to the closure in $\bar H^k_q(p)$ of the subspace $C^{\infty}(p)\cap H^k_q(p)$, the statement for $k=2$, $k=3$ follows from Lemmas~\ref{lemma:localbarH} and~\ref{lemma:localH}.

\smallskip
\noindent We complete the argument by induction on $k\in \N$. The Friedrichs extension defines
bouned operators $\Delta^F_q : \bar H^{k+1}_q(p) \to  \bar H^{k-1}_q(p)$ and its dual 
$\Delta^F_q : \bar H^{-k+1}_q(p) \to  \bar H^{-k-1}_q(p)$. By the induction hypothesis, since
all functionals in $\Cal T_p^{k-1}$ extend (uniquely) to bounded functionals in $ \bar H^{-k+1}_q(p)$,
it follows that all functionals in $\Delta^F_q (\Cal T_p^{k-1})$ extend (uniquely) to bounded functionals 
in $ \bar H^{-k-1}_q(p)$ and the following holds. For any $u\in \bar H^{k+1}_q(p)$,
\begin{equation}
 \Delta^F_q u \in H^{k-1}_q(p) \, \Leftrightarrow  \, \delta^{ij}_p(u) =0  \quad \text{\rm for all } \delta^{ij}_p
 \in  \Delta^F_q (\Cal T_p^{k-1})\,.
\end{equation}
Let $E^{k+1}_q \subset  \bar H^{k+1}_q(p)$ the closed finite-codimensional subspace defined as
\begin{equation}
E^{k+1}_q:= \{u\in \bar H^{k+1}_q(p)\,\vert\, \delta^{ij}_p(u) =0  \quad \text{\rm for all } \delta^{ij}_p
 \in  \Delta^F_q (\Cal T_p^{k-1})\,\}\,.
\end{equation}
By formulas~\pref{eq:Laplondeltas}, any distribution $\delta_p\in \Cal T_p^{k+1} \setminus \Delta^F_q (\Cal T_p^{k-1})$ is of the form $\delta_p= \delta^{0h}_p$ or $\delta_p= \delta^{h0}_p$ with $1\leq h
\leq k(m_p+1)$. By formulas~\pref{eq:deltazeros}, such distributions have a (unique) continuous extension to the subspace $E^{k+1}_q \subset  \bar H^{k+1}_q(p)$. In fact, the distributions $\text{\rm PV}(z^{-h})$ and  $\text{PV} (\bar z^{-h}) \in H^{-k+1}_q(D_p)$, for all $1\leq h < k(m_p+1)$, and $\Delta^F_q (f)\in H^{k-1}_q(D)$ for all $f\in E^{k+1}_q$. Hence all distributions in the set 
$ \Cal T_p^{k+1}$ have a continuous extension to the space $\bar H^{k+1}_q(p)$ and the characterization~\pref{eq:Hk} of the subspace $H^{k+1}_q(p) \subset \bar H^{k+1}_q(p)$ follows 
by Lemmas~\ref{lemma:localbarH} and~\ref{lemma:localH}.

\end{proof}

\noindent Let $\Cal D^k_q$ be the set of all continuous extensions to the space $\bar H^k_q(M)$
of the functionals $\delta^{ij}_p$ for all $p\in \Sigma_q$ and all $(i,j)\in \Cal T^k_p$. Let 
\begin{equation}
\label{eq:deltasing}
\Cal D_q := \bigcup_{k\in \N}  \Cal D^k_q \,\,.
\end{equation}

\begin{theorem}  
\label{thm:Hk}
The (closed) kernel of the system $\Cal D^k_q$ on $\bar H^k_q(M)$ coincides 
with the subspace $H^k_q(M)$, that is,
$$
H^k_q(M)=\{u\in \bar H^k_q(M)\,\vert\, \delta(u)=0\,,\,\,  \text{ \rm for all } \,\delta\in \Cal D^k_q\}\,.
$$
\end{theorem}

\subsection{Smoothing operators}
\label{SO}
We will establish below finer results on the Sobolev regularity of the distributions in $\Cal D_q$. 
The key step will be to construct smoothing operators for the scale of Sobolev spaces 
$\{ H_q^k(p)\,\vert \, k\in \N\}$.

\smallskip
\noindent For any $p\in \Sigma_q$, let $z: D_p \to \C$ be a canonical coordinate for  the (orientable) quadratic differential $q$ such that $p\in D_p$ and $z(p)=0$.  For any $(i,j)\in \N\times \N$, let 
$Z_p^{ij} \in C^{\infty}(M)$ be a function such that 
\begin{equation}
\label{eq:Zp}
Z_p^{ij} (z) \equiv z^{i} \bar z^{j}  \quad \text{ \rm on } \, D_p \,.
\end{equation}

\begin{lemma}
\label{lemma:SO} Let $p \in \Sigma_q$ be a zero of order $2m$ of the quadratic differential $q$ on $M$. There exists a one-parameter family $\{K_p(\tau)\,\vert\, \tau\in (0,1]\}$ of bounded operators $K_p(\tau) : C^{\infty}(M) \to H^{\infty}_q(M)$ such that the following estimates hold. For each $k\in \N$, there exists a constant $C_k>0$ such that, for any $(i,j)\in \N\times \N$ and for all $\tau\in (0,1]$:
\begin{equation}
\label{eq:SOest}
\begin{aligned}
&\vert  K_p(\tau) (Z_p^{ij} ) -Z_p^{ij}  \vert _k  \leq  C_k \,  \tau^{1+\frac{i+j}{m+1} -k } \,,
\quad &\text{for }\, k <  1+\frac{i+j}{m+1}  \,;   \\
&\vert K_p(\tau) (Z_p^{ij}  ) \vert _k  \leq C_k \,  \vert \log \tau \vert^{1/2}  \,,
\quad &\text{for }\, k =1+\frac{i+j}{m+1} \,;\\
&\vert K_p(\tau) (Z_p^{ij} ) \vert _k \leq C_k \,  \tau^{- [k-(1+\frac{i+j}{m+1})] }  \,,
\quad &\text{for }\, k > 1+\frac{i+j}{m+1} \,.
\end{aligned}
\end{equation}
\end{lemma}
\begin{proof} 
Let $z:D_p \to \C$ be a canonical coordinate at $p\in \Sigma_q$ defined on an open neighbourhood
$D_p \subset M$ such that $D_p\cap \Sigma_q=\{p\}$. There exists $r_1>0$ such that $D(r_1) \subset\subset  z(D_p)$, where $D(r_1)$ is the euclidean disk centered at the origin of radius $r_1>0$ .  Let $0<r_2 <r_1$ and let $\phi: \C \to \R$ be any non-negative smooth function identically zero on the closure of $D(r_2)$ and identically equal to $1$ outside $D(r_1)$. Let $D'_p \subset \subset D_p$
be any relatively compact neighbourhood of $p$ in $D_p$ such that $D(r_1) \subset \subset z(D'_p)$. 
For any $\tau\in (0,1]$, let us define, for all $F\in C^{\infty}(M)$,
\begin{equation}
\label{eq:Kp}
K_p(\tau) \left(F \right)(x) := 
\begin{cases}
 \phi ( \tau^{-\frac{1}{m+1}} z(x))\,  F \left(z(x) \right)\, ,\quad &\text{ for }\,x\in D_p\,; \\
F (x)            \, ,\quad &\text{ for }\,x\not \in  D'_p \,.      
\end{cases}                                   
\end{equation}
\noindent Clearly, the definition \pref{eq:Kp} is well-posed for all $\tau\in (0,1]$ and the functions $K_p(\tau) \left (F  \right)\in H^{\infty}_q(M)$, since the rescaled functions $\phi_\tau: \C \to \R$ defined as
$$
\phi_\tau (z) := \phi ( \tau^{-\frac{1}{m+1}} z)\,, \quad \text{ for }\,z\in \C\,,
$$ 
are smooth and have support away from the origin. Since the functions $\phi_\tau$
are bounded uniformly with respect to the parameter $\tau\in (0,1]$, for any function $F\in C^{\infty}(M)
\subset L^2_q(M)$, by the dominated convergence theorem,
\begin{equation}
\label{eq:Ltwoconv}
\vert K_p(\tau)\left(F\right)  - F \vert_0 \to 0 \,, \quad \text{ \rm as } \tau\to 0^+\,   \,.
\end{equation}

\smallskip
\noindent Let us denote for convenience, for any  $(a,b)\in \N \times \N$,
$$
\phi^{ab}_\tau (z):= [(\partial_q^+)^{a} (\partial_q^-)^{b} \phi] (\tau^{-\frac{1}{m+1}} z) \,, 
\quad \text{for }\, z\in \C\,.
$$
The functions $\phi^{ab}_\tau $ are smooth and bounded on $\C$, uniformly with 
respect to $\tau>0$. For any $(a,b)$ such that $(a,b)\not=(0,0)$, the function $\phi^{ab}_\tau$ 
has compact support contained in the euclidean annulus centered at the origin of 
inner radius $r_1\,\tau^{\frac{1}{m+1}}$ and outer radius $r_2 \tau^{\frac{1}{m+1}}$ .  
For $a=b=0$, $\phi^{ab}_\tau = \phi_\tau$ has support outside the euclidean disk of
radius $r_1\,\tau^{\frac{1}{m+1}}$ and is identically equal to $1$ outside the disk of radius
$r_2 \,\tau^{\frac{1}{m+1}}$. 

\smallskip
\noindent For any $(i,j)\in \N\times \N$, let $K^{ij}_\tau:= K_p(\tau) \left ( Z_p^{ij} \right)\in H^{\infty}_q(M)$. 
A calculation shows that  the following formulas hold. For all $m$, $i$ and $\alpha \in \N$ and all 
$a\in \N$ such that $0\leq a \leq \alpha$, let 
\begin{equation}
C^{m,i}_{\alpha,a} = \tbinom{\alpha}{ a}\, \prod_{\ell=1}^{\alpha-a} [ i- \ell(m+1)] \,.
\end{equation}
For any $(\alpha,\beta)\in \N\times \N$, the derivative $(\partial_q^+)^{\alpha} 
(\partial_q^-)^{\beta} \left(  K^{ij}_\tau \right)(z)$ is given on $D_p\setminus\{p\}$ by 
the sum
\begin{equation}
\label{eq:CRder}
\sum_{a=0}^{\alpha} \sum_ {b=0}^{\beta}  C^{m,i}_{\alpha,a} \, C^{m,j}_{\beta,b}
\,\phi^{ab}_\tau (z) \tau^{-(a+b)}\, z^{i- (\alpha-a)(m+1)} \bar z ^{j - (\beta-b)(m+1)}\,.
\end{equation}
If $i+j >(\alpha+\beta-1)(m+1)$, since the functions 
$$z^{i- (\alpha-a)(m+1)} \bar z ^{j - (\beta-b)(m+1)} \in L^2_q(D_p)\, ,$$ 
for all $0\leq a \leq \alpha$ and $0\leq b \leq \beta$, and $\phi^{ab}_\tau$ is bounded, by 
change of variables we obtain that for $(a,b)\not=(0,0)$ there exists a constant 
$C_{a,b}>0$ such that
\begin{equation}
\label{eq:SOone}
\vert \phi^{ab}_\tau (z) \tau^{-(a+b)}\, z^{i- (\alpha-a)(m+1)} \bar z ^{j - (\beta-b)(m+1)}\vert_0
\leq C_{a,b} \, \tau^{1+\frac{i+j}{m+1} -(\alpha +\beta)}\,,
\end{equation}
and that similarly, for $a=b=0$, there exists a constant $C_0>0$
\begin{equation}
\label{eq:SOtwo}
\vert (\phi_\tau (z) -1)\, z^{i- \alpha(m+1)} \bar z ^{j - \beta(m+1)}\vert_0
\leq C_0 \, \tau^{1+\frac{i+j}{m+1} -(\alpha +\beta)}\,.
\end{equation}
If $i+j < [(\alpha-a)+(\beta-b)-1] (m+1) $,  there exists a constant $C'_{a,b}>0$ such that
\begin{equation}
\label{eq:SOthree}
\vert \phi^{ab}_\tau (z) \tau^{-(a+b)}\, z^{i- (\alpha-a)(m+1)} \bar z ^{j - (\beta-b)(m+1)}\vert_0
\leq C'_{a,b} \, \tau^{1+\frac{i+j}{m+1} -(\alpha +\beta)}\,,
\end{equation}
since the function $\phi^{ab}_\tau $ is bounded and supported outside the euclidean disk 
of radius $r_1\,\tau^{\frac{1}{m+1}}$ centered at the origin. 

\smallskip
\noindent If $i+j =[(\alpha-a)+(\beta-b)-1] (m+1)$, a similar calculation yields
\begin{equation}
\label{eq:SOfour}
\begin{aligned}
\vert \phi^{ab}_\tau (z) \tau^{-(a+b)}\, z^{i- (\alpha-a)(m+1)} \bar z ^{j - (\beta-b)(m+1)}\vert_0& \\
\leq C'_{a,b} \, &\tau^{1+\frac{i+j}{m+1} -(\alpha +\beta)} \, \vert \log\tau  \vert^{\frac{1}{2}}\,.
\end{aligned}
\end{equation}
\noindent By formula \pref{eq:CRder} for the Cauchy-Riemann iterated derivatives, the required estimates \pref{eq:SOest} follow immediately from estimates \pref{eq:SOone},  \pref{eq:SOtwo}, 
\pref{eq:SOthree} and \pref{eq:SOfour}.
\end{proof} 

\noindent We derive below estimates for the local smoothing family  $\{K_p(\tau)\}$ constructed in Lemma \ref{lemma:SO} with respect to the fractional weighted Sobolev norms. Let $p\in \Sigma_q$ be any zero of order $2m_p$ of the quadratic differential $q$ on $M$. For each pair $(i,j)\in \N\times\N$, let $e^{ij}_p:\R^+ \to [0,1]$ denote the function defined as follows:
\begin{equation}
\label{eq:eij}
e^{ij}_p(s):= 
\begin{cases}
 \{s\} \,,   \quad  &\text{if }\, [s] = 1 +  \frac{i+j}{m_p+1} \,; \\
 1-\{s\} \,,  \quad  &\text{if }\, [s] = \frac{i+j}{m_p+1} \,;\\
 0    \,,   \quad  &\text{otherwise}\,.
\end{cases}
\end{equation}

\begin{theorem}
\label{thm:SO} The family $\{K_p(\tau)\,\vert\, \tau\in (0,1]\}$ of local smoothing operators $K_p(\tau) : C^{\infty}(M) \to H^{\infty}_q(M)$, defined in \pref{eq:Kp}, has the following properties. For each $s\in \R^+$, there exists a constant $C_s>0$ such that, for any pair $(i,j)\in \N\times \N$ and all $\tau\in (0,1]$:
\begin{equation*}
\begin{aligned}
&\vert  K_p(\tau) (Z_p^{ij} ) -Z_p^{ij}  \vert _s  \leq  C_s \tau^{1+\frac{i+j}{m_p+1} -s } \,
\vert \log \tau\vert^{\frac{e_p^{ij}(s)}{2}},
\,&\text{for }\, s <  1+\frac{i+j}{m_p+1} \, ;   \\
&\vert K_p(\tau) (Z_p^{ij}  ) \vert _s  \leq C_s  \vert \log\tau\vert  
\,\vert \log \tau\vert^{\frac{e_p^{ij}(s)}{2}},
\,&\text{for }\, s =1+\frac{i+j}{m_p+1}\, ;\\
&\vert K_p(\tau) (Z_p^{ij} ) \vert _s\leq C_s  \tau^{- [s-(1+\frac{i+j}{m_p+1})] }\,  \vert \log \tau\vert^{\frac{e_p^{ij}(s)}{2}} ,
\,&\text{for }\, s > 1+\frac{i+j}{m_p+1}\, .
\end{aligned}
\end{equation*}
\end{theorem}
\begin{proof}
\noindent For any $k\in \N$, the function $Z_p^{ij}\in  H^k_q(M)$, if $k<1 +(i+j)/(m+1)$, since $z^{i}\bar z^{j} \in H^k_q(p)$. By Lemma \ref{lemma:SO}, there exists a constant $C_k>0$ such that, for all 
$\tau\in (0,1]$,
\begin{equation}
\begin{aligned}
\vert & K^{ij}_p(\tau)- K^{ij}_p(\tau/2)\vert_{k} \leq  \vert ( K^{ij}_p(\tau)-Z_p^{ij})-
( K^{ij}_p(\tau/2)-Z_p^{ij})\vert_{k} \,; \\
& \leq \vert  K^{ij}_p(\tau)- Z_p^{ij} \vert_{k} + \vert  K^{ij}_p(\tau/2) - Z_p^{ij} \vert_{k} 
 \leq  C_{k}  \,\tau^{1+\frac{i+j}{m+1} -k}\,. 
\end{aligned}
\end{equation}
If $k\geq 1 +(i+j)/(m+1)$,
\begin{equation}
\vert  K^{ij}_p(\tau)- K^{ij}_p(\tau/2)\vert_{k}  \leq  \vert  K^{ij}_p(\tau)\vert_k \,+\, 
\vert K^{ij}_p(\tau/2)\vert_{k} 
\end{equation}
hence, by Lemma \ref{lemma:SO},
\begin{equation}
\begin{aligned}
\vert  K^{ij}_p(\tau)- K^{ij}_p(\tau/2)\vert_{k}  &\leq  2C_k \,\vert \log\tau \vert^{1/2} , \quad &\text{\rm if }  \,\,k=1 +\frac{i+j}{m+1}\,;\\
\vert  K^{ij}_p(\tau)- K^{ij}_p(\tau/2)\vert_{k}  &\leq  2C_k \,\tau^{1+\frac{i+j}{m+1} -k} , 
\quad &\text{\rm if } \,\, k>1 +\frac{i+j}{m+1} \,.
\end{aligned}
 \end{equation}
 As a consequence, by the interpolation inequality (Lemma \ref{lemma:intineq}), for every $s\in \R^+$ there exists a constant $C_s>0$ such that, for all $\tau\in (0,1]$,
 \begin{equation*}
\begin{aligned}
&\vert  K^{ij}_p(\tau)- K^{ij}_p(\tau/2)\vert_{s}  \leq  C_s\, \tau^{1+\frac{i+j}{m+1} -s}\,, 
\qquad \text{\rm if }  \,\,[s], [s]+1 \not = 1 +\frac{i+j}{m+1}\,;\\
&\vert  K^{ij}_p(\tau)- K^{ij}_p(\tau/2)\vert_{s}  \leq  C_s\, \tau^{1+\frac{i+j}{m+1} -s} 
\vert \log\tau \vert^{\frac{1-\{s\}}{2}} \,, 
\quad \text{\rm if } \,\, [s]=1 +\frac{i+j}{m+1} \,; \\
&\vert  K^{ij}_p(\tau)- K^{ij}_p(\tau/2)\vert_{s}  \leq  C_s\, \tau^{1+\frac{i+j}{m+1} -s} \vert \log\tau\vert^{\frac{\{s\}}{2}}\, , 
\qquad \text{\rm if } \,\, [s]=\frac{i+j}{m+1}\,.
\end{aligned}
 \end{equation*}
 If $s<1 + (i+j)/(m+1)$ and $[s] \not=(i+j)/(m+1)$, for all $n\in \N$,
\begin{equation}
\label{eq:2n}
\vert  K^{ij}_p(\tau/2^n)- K^{ij}_p(\tau/2^{n+1})\vert_{s} \leq C_s \,\tau^{1 + \frac{i+j}{m+1} -s}
2^{-n(1 + \frac{i+j}{m+1} -s)} \,,
\end{equation}
hence, for any fixed $\tau\in (0,1]$, the sequence $\{ K^{ij}_p(\tau/2^n)\}_{n\in \N}$ is Cauchy, 
and therefore convergent, in the Hilbert space $H^s_q(M)$. By Lemma \ref{lemma:SO} 
$\{ K^{ij}_p(\tau/2^n)\}_{n\in \N}$ converges to $Z_p^{ij}$ in $H^{[s]}_q(M)$. Since 
 $H^{[s]}_q(M)\subset H^s_q(M)$, by uniqueness of the limit $Z_p^{ij} \in H^s_q(M)$
and $\{  K^{ij}_p(\tau/2^n)\}_{n\in \N}$ converges to $Z_p^{ij}$ in $H^s_q(M)$.
The estimate \pref{eq:2n} also implies that
\begin{equation}
\vert  K^{ij}_p(\tau )- Z^{ij}_p \vert_{s} \leq  \sum_{n\in \N}
 \vert  K^{ij}_p(\tau/2^n)- K^{ij}_p (\tau/2^{n+1})\vert_{s} \leq 
 C'_s \,\tau^{1 + \frac{i+j}{m+1} -s}\,.
\end{equation}
If $s<1 + (i+j)/(m+1)$ and $[s] =(i+j)/(m+1)$, by a similar argument we again get that 
$Z_p^{ij} \in H^s_q(M)$ and 
\begin{equation}
\vert  K^{ij}_p(\tau )- Z^{ij}_p \vert_{s}  \leq C'_s \,\tau^{1 + \frac{i+j}{m+1} -s} \, 
\vert \log \tau \vert ^{\frac{\{s\}}{2}}\,.
\end{equation}

\smallskip
\noindent  If $s\geq 1 + (i+j)/(m+1)$  and $[s]$, $[s] +1 \not = 1+(i+j)/(m+1)$ we argue as follows.
For each $\tau \leq 1/2$, we have
 \begin{equation}
 \vert  K^{ij}_p(2\tau)- K^{ij}_p(\tau)\vert_{s}  \leq  C_s\,(2 \tau)^{1+\frac{i+j}{m+1} -s} \,,
  \end{equation}
  hence if $\tau\leq  2^{-n}$, for all $0\leq k < n$,
  \begin{equation}
 \vert  K^{ij}_p(2^{k+1} \tau)- K^{ij}_p(2^k\tau)\vert_{s}  \leq  C_s\,2 ^{(k+1)(1+\frac{i+j}{m+1} -s)}
\tau^{1+\frac{i+j}{m+1} -s} \,.
  \end{equation}
  It follows that, there exists a constant $C'_s>0$ such that
  \begin{equation}
   \vert  K^{ij}_p(2^n \tau)- K^{ij}_p(\tau)\vert_{s} \leq  C'_s\,2 ^{n(1+\frac{i+j}{m+1} -s)}
\tau^{1+\frac{i+j}{m+1} -s} \,.
 \end{equation}
 For every $\tau\in (0,1]$, let $n(\tau)$ be the maximum $n\in \N$ such that $2^n \tau \leq 1$. By
 this definition it follows that $1/2 < 2^{n(\tau)} \tau \leq 1$. Since
 $$
 \sup_{1/2 \leq \tau \leq 1}  \vert K^{ij}_p(\tau) \vert_s  \leq 
 \sup_{1/2 \leq \tau \leq 1}  \vert K^{ij}_p(\tau) \vert_{[s]+1} \,\,< \,\, +\infty\,\,,
 $$  
 it follows that, there exists a constant $C''_s >0$ such that 
  \begin{equation}
  \begin{aligned}
   \vert K^{ij}_p(\tau)\vert_{s} &\leq  C''_s\, \tau^{1+\frac{i+j}{m+1} -s}  \,, 
   \quad &\text{ \rm if }  s> 1+\frac{i+j}{m+1} \,;\\
   \vert K^{ij}_p(\tau)\vert_{s} &\leq  C''_s\,  \vert \log\tau\vert   \,, 
   \quad &\text{ \rm if }  s= 1+\frac{i+j}{m+1} \,.
   \end{aligned}
 \end{equation}
 By a similar argument, for $s>1+(i+j)/(m+1)$ we have
  \begin{equation}
  \begin{aligned}
 \vert K^{ij}_p(\tau)\vert_{s} &\leq  C''_s\, \tau^{1+\frac{i+j}{m+1} -s} \vert \log\tau \vert ^{\frac{1-\{s\}}{2}}\,,
  \quad \text{\rm if }  \,\, [s]= 1+\frac{i+j}{m+1}\,; \\  
  \vert K^{ij}_p(\tau)\vert_{s} &\leq  C''_s\, \tau^{1+\frac{i+j}{m+1} -s} \vert \log\tau \vert ^{\frac{\{s\}}{2}}\,,
  \qquad \text{\rm if }  \,\, [s]= \frac{i+j}{m+1}  \,,
   \end{aligned}
 \end{equation}
while for $s=1+(i+j)/(m+1)$ we have
\begin{equation}
  \begin{aligned}
 \vert K^{ij}_p(\tau)\vert_{s} &\leq  C''_s\,  \vert \log\tau \vert\, \vert \log\tau \vert ^{\frac{1-\{s\}}{2}}\,,
  \quad \text{\rm if }  \,\, [s]= 1+\frac{i+j}{m+1}\,; \\  
  \vert K^{ij}_p(\tau)\vert_{s} &\leq  C''_s\,  \vert \log\tau \vert\, \vert \log\tau \vert ^{\frac{\{s\}}{2}} \,,
  \qquad \text{\rm if }  \,\, [s]= \frac{i+j}{m+1}  \,.
   \end{aligned}
 \end{equation}
\end{proof}

\noindent Theorem \ref{thm:SO} implies in particular the following smoothness results.
\begin{corollary} 
\label{cor:ijreg}
Let $z:D_p \to C$ be a canonical coordinate for an orientable  quadratic differential $q$ at a zero 
$p\in \Sigma_q$ of order $2m$.  For each $(i,j)\in \N\times \N$, the function 
\begin{equation}
z^{i} \bar z^{j} \in H^s_q(p)\,\,,  \quad \text{\rm for all }\, s< 1+ \frac{i+j}{m_p+1} \,\,.
\end{equation}
\end{corollary}

\begin{corollary} 
\label{cor:deltareg}
Let $p\in \Sigma_q$ be a zero of order $2m_p$ and let $(i,j) \in {\Cal T}_p$. 
The distribution $\delta_p^{ij}$ has the following regularity properties:
\begin{equation}
\begin{aligned}
&\delta_p^{ij} \in \bar H^{-s}_q(p) \quad \text{ \rm for }  s> 1 + \frac{i+j}{m_p+1}\,, \\
&\delta_p^{ij} \not\in H^{-s}_q(p) \quad \text{ \rm for } s <1 + \frac{i+j}{m_p+1} \,.
\end{aligned}
\end{equation}
\end{corollary}
\begin{proof} By the formulas \pref{eq:deltazeros} and by Lemma \ref{lemma:continuity}, for any 
$h \in \N\setminus \N\cdot (m_p+1)$ and for any $\ell\in \N$, there exist constants $C^{m_p,\ell}_{h0} \not=0$ and $C^{m_p,\ell}_{0h} \not=0$ such that the following  identities hold in the dual Hilbert space 
$\bar H_q^{-k}(p)$ for any integer $k \geq 1 + h/(m_p+1)$:
\begin{equation}
\label{eq:hlapl}
\begin{aligned}
C^{m_p,\ell}_{h0}  \,\delta_p^{h0}& = \Delta^{\ell+1}_q \left ( z^{\ell(m_p+1)-h}
 \bar z^{\ell(m_p+1)}\right) \,,\\
C^{m_p,\ell}_{0h}  \,\delta_p^{0h}&  = \Delta^{\ell+1}_q \left ( z^{\ell(m_p+1)} 
\bar z^{\ell(m_p+1)-h}\right)\,.
\end{aligned}
\end{equation}
By Corollary \ref{cor:ijreg}, if $\ell(m_p+1)-h>0$, for all $s < 2\ell +1 - h/(m_p+1)$,
$$
z^{\ell(m_p+1)-h} \bar z^{\ell(m_p+1)}\,\,  \text{ and }\,\, z^{\ell(m_p+1)} \bar z^{\ell(m_p+1)-h} 
\in H^s_q(p) \subset \bar H^s_q(p)\,.
$$
Hence $\delta_p^{h0}$, $\delta_p^{0h} \in \bar H^{-s}_q(p)$ for all $s> 1+h/(m_p+1)$.
By formulas \pref{eq:Laplondeltas} it then follows that $\delta_p^{ij} \in \bar H^{-s}_q(p)$ for all 
$s> 1+(i+j)/(m_p+1)$ as claimed.

\smallskip
\noindent  Let $(i,j)\in \Cal T_p$ and $s<1 +(i+j)/(m_p+1)$. By Corollary \ref{cor:ijreg}, the 
function $z^{i} \bar z^{j}\in H^s_q(p)$. Since by definition $H^{\infty}_q(p)$ is dense in 
$H^s_q(p)$ for any $s>0$, the functional $\delta_p^{ij} \equiv 0$ on $H^{\infty}_q(p)$ and $\delta_p^{ij}(z^{i} \bar z^{j}) =1$, it follows that  $\delta_p^{ij}$ does not extend to a bounded 
functional on $H^s_q(p)$.
\end{proof}

\smallskip
\noindent Let $p\in \Sigma_q$ be a zero of (even) order $2m_p$ the quadratic differential $q$ 
on $M$. For every $s\in \R^+$, let $\Cal T ^s_p \subset \Cal T_p$ be the subset defined as
\begin{equation}
\label{eq:Tsp}
\Cal T ^s_p:=\{ (i,j)\in \Cal T_p \,\vert \, i+j < (s-1) (m_p+1)\}\,.
\end{equation}
Let $\Cal D^s_q \subset \bar H^{-s}_q(M)$ be the set of distributions defined as follows:
\begin{equation}
\label{eq:Dsq}
\Cal D^s_q := \{ \delta_p^{ij} \,\vert\,  p\in \Sigma_q\, \text{ \rm and }\, (i,j)\in \Cal T^s_p\}\,.
\end{equation}
\begin{corollary}
\label{cor:Hsclosure}
The closure of the subspace $H^s_q(M)$ in $\bar H^s_q(M)$ is a subset of the (closed) kernel 
of the system $\Cal D^s_q$ on $\bar H^s_q(M)$, that is,
\begin{equation}
\overline{H^s_q(M)} \subset \{u\in \bar H^s_q(M)\,\vert\, \delta(u)=0\,,\,\,  \text{ \rm for all } \,\delta\in
 \Cal D^s_q\}\,
\end{equation}
The reverse inclusion holds if the following sufficient condition is satisfied:
\begin{equation}
\label{eq:scond}
s \not\in \{ 1+ (i+j)/(m_p+1)\,\vert \, p\in \Sigma_q\, \text{ \rm and }\, (i,j)\in \Cal T_p\}\,.
\end{equation}
\end{corollary}
\begin{proof} Since $H^{\infty}_q(M)$ is dense in $H^s_q(M) \subset \bar H^s_q(M)$,
$$
H^{\infty}_q(M) \subset \{u\in \bar H^s_q(M)\,\vert\, \delta(u)=0\,,\,\,  \text{ \rm for all } \,\delta\in
 \Cal D^s_q\}
 $$
 and $ \Cal D^s_q \subset  \bar H^{-s}_q(M)$, it follows that
 $$
 \overline{H^s_q(M)} \subset \{u\in \bar H^s_q(M)\,\vert\, \delta(u)=0\,,\,\,  \text{ \rm for all } \,\delta\in
 \Cal D^s_q\}\,.
$$
Conversely, if condition \pref{eq:scond} is satisfied, by Corollary \ref{cor:ijreg} the subspace 
$$
\{u\in C^{\infty}(M)\,\vert\, \delta(u)=0\,,\,\,  \text{ \rm for all } \,\delta\in
 \Cal D^s_q\} \subset H^s_q(M) \,.
 $$
 Since $C^{\infty}(M) \cap  \bar H^s_q(M)$ is dense in  $\bar H^s_q(M)$, the result follows.
\end{proof}

\noindent The regularity result proved in Corollary \ref{cor:ijreg} extends to a certain subset of all pairs 
$(i,j)\in \Z  \times \Z$ if the functions $z^{i} \bar z^{j}$ are interpreted as distributions in the sense of the Cauchy principal value:
\begin{equation}
\label{eq:PVij}
\text{ \rm PV}\left(z^{i} \bar z^{j}\right)(v) := \text{ \rm PV} \int_ M  z^{i} \bar z^{j}\,v\,\omega_q \,, \quad 
 \text{ \rm for all } v\in C^{\infty}(p)\,.
 \end{equation}
The most general regularity result for the distributions \pref{eq:PVij} is based on the
following generalization of Corollary \ref{cor:ijreg} to include logarithmic factors.

\begin{lemma}
Let $z:D_p \to \C$ be a canonical coordinate for an orientable  holomorphic quadratic differential 
$q$ at a zero $p\in \Sigma_q$ of order $2m$. For each $(i,j,h)\in \N\times \N\times\N$, the function 
$$
 z^{i} \bar z^{j} \log^h \vert z\vert \in H^s_q(p)  \,, \quad \text{ \rm for all } s< 1+ \frac{i+j}{m_p+1} \,.
$$
\end{lemma}
\begin{proof} Simple calculations show that $\log \vert z\vert  \in L^2_q(M)$ and that by formulas
\pref{eq:CRlocal} the following identities hold on $D_p\setminus\{p\}$:
\begin{equation}
\label{eq:CRlog}
\partial^+ \log \vert z\vert = \frac{ 1}{\bar z^{m+1}} \quad \text{ \rm and } \quad 
\partial^- \log \vert z\vert = \frac{ 1}{z^{m+1}} \,.
\end{equation}
It follows that, for each $(i,j,h)\in \N\times \N\times\N$ and each $(\alpha,\beta)\in \N\times\N$, there 
exists a finite sequence of non-zero constants $C_1, \dots, C_h$, which depend on $(i,j,h,\alpha,\beta, m)$, such that the following identity holds on $D_p\setminus\{p\}$:
\begin{equation}
\label{eq:CRlogbis}
(\partial^+)^\alpha (\partial^-)^\beta \left(z^{i} \bar z^{j} 
\log^h \vert z\vert \right)= z^{i-\alpha(m+1)} \bar z^{j-\beta(m+1)}  \sum_{\ell=0}^h C_\ell \, \log^\ell \vert z\vert\,.
\end{equation}
For all $(i,j,h)\in \N\times \N \times\N$, let $L^{ijh}_p\in C^{\infty}(M\setminus\{p\})$ be any function such that $L^{ijh}_p(z)=z^{i} \bar z^{j} \log^h \vert z\vert$ for all $z\in D_p$. By \pref{eq:CRlogbis}, the function
$$
L^{ijh}_p \in H^k_q(M)  \,, \quad \text{ \rm if } k\in \N\,\, \text{ and }\,\,k< 1+ \frac{i+j}{m_p+1} 
$$
Let $\{K_p(\tau)\,\vert \, \tau\in (0,1]\}$ be the family of local smoothing operators defined by formulas \pref{eq:Kp}.  By computations similar to those carried out in the proof of Lemma \ref{lemma:SO}, based on formulas \pref{eq:CRlogbis}, it is possible to prove that for each $k\in \N$, there exists a constant $C_k>0$ such that for all $\tau\in (0,1]$:
\begin{equation*}
\begin{aligned}
&\vert  K_p(\tau) (L_p^{ijh} ) -L_p^{ijh}  \vert _k  \leq  C_k \, \tau^{1+\frac{i+j}{m+1} -k } \, 
\vert\log\tau \vert^h ,
&\text{for }k <  1+\frac{i+j}{m+1}  ;   \\
&\vert K_p(\tau) (L_p^{ijh}  ) \vert _k  \leq C_k \,  \vert \log\tau \vert ^{1/2} \vert\log\tau \vert^h  ,
&\text{for } k =1+\frac{i+j}{m+1} ;\\
&\vert K_p(\tau) (L_p^{ijh} ) \vert _k \leq C_k \, \tau^{- [k-(1+\frac{i+j}{m+1})] } \, \vert\log\tau \vert^h    ,
&\text{for } k > 1+\frac{i+j}{m+1} .
\end{aligned}
\end{equation*}
Reasoning as in the proof of Theorem \ref{thm:SO}, we can derive similar estimates for fractional Sobolev norms. For each $(i,j)\in \N\times \N$, let $e_p^{ij}:\R^+ \to [0,1]$ be the function defined in formula \pref{eq:eij}. By the interpolation Lemma \ref{lemma:intineq}, for any $s<1+(i+j)/(m+1)$ there exists a constant $C_s>0$ such that
$$
\vert  K_p(\tau) (L_p^{ijh} ) - K_p(\tau/2) (L_p^{ijh})  \vert _s  \leq  C_s \, \tau^{1+\frac{i+j}{m+1} -s } \, \vert\log\tau \vert^{h+\frac{e_p^{ij}(s)}{2}} \,.
$$
It follows that the sequence $\{ K_p(\tau/2^n) (L_p^{ijh} )\}_{n\in\N}$ is Cauchy and therefore converges in $H^s_q(M)$. By uniqueness of the limit 
$$
L^{ijh}_p \in H^s_q(M)  \,, \quad \text{ \rm for all } s< 1+ \frac{i+j}{m_p+1} \,.
$$
In addition, the following estimates hold. For each $s\in \R^+$ there exists a constant $C'_s>0$ such 
that  for all $\tau\in (0,1]$:
\begin{equation*}
\begin{aligned}
&\vert K_p(\tau) (L_p^{ijh}  ) \vert _s  \leq C'_s\,  \vert \log\tau  \vert   \,\vert \log\tau\vert^{h+\frac{e_p^{ij}(s)}{2}} \,,
\,\,&\text{for } s =1+\frac{i+j}{m_p+1}\,;\\
&\vert K_p(\tau) (L_p^{ijh} ) \vert _s\leq C'_s \, \tau^{- [s-(1+\frac{i+j}{m_p+1})] }
 \vert \log\tau\vert^{h+\frac{e_p^{ij}(s)}{2}} \,,
\,\,&\text{for } s > 1+\frac{i+j}{m_p+1} \,.
\end{aligned}
\end{equation*}
\end{proof}
\begin{theorem}
 \label{thm:ijreg}
Let $z:D_p \to \C\,$ be a canonical coordinate for an orientable  holomorphic quadratic differential 
$q$ at a zero $p\in \Sigma_q$ of order $2m_p$. For each $(i,j)\in \Z\times \Z$ such that $(1)$\, $i-j \not
\in \Z \cdot (m_p+1)$ or $(2)$\, $i> -(m_p+1)$ or $(3)$ \,$j>-(m_p+1)$, the distribution
\begin{equation}
\text{ \rm PV} \left(z^{i} \bar z^{j}\, \log^h \vert z\vert\right) \in H^s_q(p)\,\,,  \quad \text{\rm for all }\, 
s< 1+ \frac{i+j}{m_p+1} \,\,.
\end{equation}
\end{theorem}
\begin{proof}
For all $(i,j)\in \Z\times \Z$ such that $i-j \not \in \Z \cdot (m_p+1)$ and for any $h\in \Z$ the following formulas hold for all functions $v\in H^{\infty}_q(p)$:
\begin{equation}
\label{eq:PVformulas}
\begin{aligned}
(a) \text{ \rm PV} \int_M   \partial^+ (z^{i} \bar z^j \log^h \vert z\vert) v\,\omega_q 
 &= - \text{\rm PV} \int_M   z^{i} \bar z^j \log^h \vert z\vert \,\partial^+v\,\omega_q\,;\\
(b) \text{ \rm PV} \int_M   \partial^- (z^{i} \bar z^j \log^h \vert z\vert) v\,\omega_q &=
- \text{\rm PV} \int_M   z^{i} \bar z^j \log^h \vert z\vert \,\partial^-v\,\omega_q \,.
\end{aligned}
\end{equation}
Formulas \pref{eq:PVformulas} also hold in case $(a)$ if $i>-(m_p+1)$, $j\in \Z$, 
and in case $(b)$ if $j>-(m_p+1)$, $i\in \Z$, for all germs  $v\in C^{\infty}(p)$. 

\smallskip
\noindent By taking into account the formulas
 \pref{eq:CRlocal} for the Cauchy-Riemann operators with respect  to a canonical coordinate, 
 it follows from formulas \pref{eq:PVformulas} by induction on $h\in \N$ that if 
 $$
 \text{ \rm PV} (z^{i} \bar z^{j} \log^h\vert z\vert)  \in H^s_q(p)\,, \quad \text{ \rm for all }\,h\in \N\,,
 $$ 
 then, if $i-j\not \in \Z \cdot (m_p+1)$ or $i>-(m_p+1)$ and $j\in \Z$,
 $$
 \text{ \rm PV} (z^{i} \bar z^{j-(m_p+1)} \log^h\vert z\vert) \in 
 H^{s-1}_q(p)\,, \quad \text{ \rm for all }\,h\in \N\,,
 $$
 and, if  $i-j\not \in  \Z \cdot (m_p+1)$ or $j>-(m_p+1)$ and $i\in \Z$, 
 $$
 \text{ \rm PV} (z^{i-(m_p+1)} \bar z^{j}  \log^h\vert z\vert) \in H^{s-1}_q(p)\,,
 \quad \text{ \rm for all }\,h\in \N\,.
 $$
Thus, the statement of the theorem can be derived from Corollary \ref{cor:ijreg}, 
by an induction argument based on formulas \pref{eq:PVformulas}.
\end{proof}

\begin{corollary}
 \label{cor:ijregbis}
Let $z:D_p \to \C$ be a canonical coordinate for an orientable  holomorphic quadratic differential 
$q$ at a zero $p\in \Sigma_q$ of order $2m_p$. If $\,(i,j)\not\in  \Z \cdot (m_p+1) \times \Z \cdot  
(m_p+1)$,  the distribution
\begin{equation}
\text{ \rm PV} \left(z^{i} \bar z^{j}\, \log^h \vert z\vert \right) \not  \in H^s_q(p)\,\,,  \quad \text{\rm for }\, 
s> 1+ \frac{i+j}{m_p+1} \,,
\end{equation}
and, if $\,i-j\in \Z\cdot (m_p+1)\,$ and both $\,i\leq -(m+1)\,$ and $\,j\leq -(m+1)$, 
\begin{equation}
\text{ \rm PV} \left(z^{i} \bar z^{j}\, \log^h \vert z\vert \right) \not  \in H^{-\infty}_q(p)\,\,.
\end{equation}
\end{corollary}
\begin{proof}
We argue by contradiction. Assume there exists $(i,j)\in \Z  \times \Z$ such that $(i,j)\not\in  \Z \cdot (m_p+1) \times \Z \cdot  (m_p+1)$ and $\text{ \rm PV} \left(z^{i} \bar z^{j}\right) \in H^r_q(p)$ for some
$r > 1+ (i+j)/(m_p+1)$. By taking  Cauchy-Riemann derivatives if necessary, we can
assume that $i \leq 0$ and $j\leq 0$. By Theorem \ref{thm:ijreg}, the distribution 
$$
\text{ \rm PV} \left(z^{-i-(m+1)} \bar z^{-j-(m+1)}\right) \in H^s_q(p)\,,
\quad \text{\rm for all }\, s< -1- \frac{i+j}{m_p+1} \,\,.
$$
It follows that, for any positive smooth function $\phi \in C^{\infty}_0(D_p)$ identically equal to $1$ 
on a disk $D'_p \subset \subset D_p$ centered at $p\in \Sigma_q$, the principal value
$$
\text{ \rm PV}\int_M  \phi(z) \,z^{i} \bar z^{j}\, z^{-i-(m+1)} \bar z^{-j-(m+1)}\, \omega_q
$$
is finite, which can be proved to be false by a simple computation in (geodesic) polar
coordinates. This contradiction proves the first part of the statement.

\smallskip
\noindent If $i-j\in \Z\cdot (m_p+1)$ and both $i\leq -(m+1)$ and $j\leq -(m+1)$,  we argue as follows.
It is not restrictive to assume that $i\geq j$, hence the function ${\bar z}^{i-j} \in H_q^{\infty}(M)$. 
However, by a computation in polar coordinates, since $i\leq -(m_p+1)$, 
$$
\text {\rm PV} \int_M  \phi(z) \,z^{i} \bar z^{j} \log \vert z \vert \,{\bar z}^{i-j} \, \, \omega_q = +\infty\,.
$$
It follows that $\text{ \rm PV} \left(z^{i} \bar z^{j}\, \log^h \vert z\vert \right) \not  \in H^{-\infty}_q(p)$,
hence the second part of the statement is also proved. 
\end{proof}

\noindent We conclude this section with a fundamental smoothing theorem for the $1$-parameter family of weighted Sobolev spaces.

\begin{theorem}  
\label{thm:smoothingop}
For each $k\in \N$, there exists a family $\{ \Cal S^k(\tau) \,\vert \, \tau\in (0,1]\}$ of bounded operators $\Cal S^k(\tau): L^2_q(M) \to H^k_q(M)$ such that the following estimates hold. For any $s$, 
$r\in [0, k]$ and for any $\epsilon>0$, there exists a constant $C^k_{r,s}(\epsilon)>0$ such that, for all 
$u \in H^s_q(M)$ and for all $\tau\in (0,1]$:
\begin{equation}
\begin{aligned}
\vert \Cal S^k(\tau) (u)-u\vert_r   &\leq  C^k_{r,s}(\epsilon) \Vert u \Vert _s  \, \tau^{s-r-\epsilon} \,,
\quad &\text{ \rm if } \, s>r  \,; \\
\vert \Cal S^k(\tau) (u)\vert_r    &\leq  C^k_{r,s}(\epsilon)   \Vert u \Vert _s  \, \tau^{s-r-\epsilon}  \,,
  \quad &\text{\rm  if } \, s\leq r  \,.
\end{aligned}
\end{equation}
\end{theorem} 
\begin{proof} For each $p\in \Sigma_q$, let $z: D_p \to\C$ be a canonical coordinate defined
on a disk $D_p$ (centered at $p$) such that $D_p\cap \Sigma_q=\{p\}$. For each $(i,j)\in \Cal T_p$, let 
$Z^{ij}_p\in C^{\infty}(M)$ be a (fixed) smooth extension, as in \pref{eq:Zp}, of the locally defined
function $z^{i}\bar z^{j} \in C^{\infty}(p)$. Let $P^k$ be the linear operator defined as follows:
\begin{equation}
P^k(f) := f   - \sum_{p\in \Sigma_q} \sum_{(i,j)\in \Cal T_p^k}   \delta_p^{ij} (f)\, Z^{ij}_p\,\,,
\quad \text{ for all }\, f\in \bar H^k_q(M)\,.
\end{equation}
The operator $P_k: \bar H^k_q(M) \to H^k_q(M) $  is well-defined and bounded. It
is well-defined  by Lemma \ref{lemma:continuity} and Theorem \ref{thm:Hk}. It is bounded
since,  for all $p\in \Sigma_q$ and all $(i,j)\in \Cal T_p$,  the functions $Z^{ij}_p \in \bar 
H_q^\infty(M)$ by Lemma \ref{lemma:localbarH} and, for all $(i,j)\in \Cal T^\ell_p$,  the distributions $\delta_p^{ij}\in \bar H^{-k}_q(M)$ by Corollary \ref{cor:deltareg}. In fact, the condition 
$(i,j)\in \Cal T_p^\ell$ implies $i+j < (k-1)(m_p+1)$. 

\smallskip
\noindent For each $p\in \Sigma_q$, let $\{ K_p(\tau) \,\vert \, \tau\in (0,1]\}$ be the family of local smoothing operators constructed in Lemma \ref{lemma:SO}.  Let $ \{S^k_\tau \vert \tau\in (0,1]\}$ 
be the one-parameter family of bounded linear operators $S^k_\tau: \bar H^k_q(M) \to H^k_q(M)$ defined as follows. For all $f\in \bar H^k_q(M)$, we let 
\begin{equation}
\label{eq:SkNt}
S^k_\tau (f) :=  P^k (f)   \, +\, 
\sum_{p\in \Sigma_q} \sum_{(i,j)\in \Cal T_p^k}   \delta_p^{ij} (f)\,K_p(\tau) \left( Z^{ij}_p\right)\,.
\end{equation}
By definition the following identity holds for all $f\in  \bar H^k_q(M)$:
$$
S^k_\tau (f) -f = \sum_{p\in \Sigma_q} \sum_{(i,j)\in \Cal T_p^k}   \delta_p^{ij} (f)\,
\left[ K_p(\tau) \left( Z^{ij}_p\right) - Z^{ij}_p\right]\,.
$$
Since for all $p\in \Sigma_q$ the condition $(i,j)\in\Cal T_p$ implies $i+j\not\in \N\cdot (m_p+1)$,
by Lemma \ref{lemma:SO} the following estimate holds. For each $\ell\in \N$ such that $\ell \leq k$, 
there exists a constant $C^k_\ell>0$ such that, for any $f\in \bar H^k_q(M)$, 
\begin{equation}
\label{eq:Sk1}
\Vert S^k_\tau (f) -f\Vert_\ell \leq  C^k_\ell \, \sum_{p\in\Sigma_q}\sum_{(i,j)\in \Cal T_p^k} 
\tau^{1+\frac{i+j}{m_p+1}-\ell}\, \vert \delta_p^{ij} (f) \vert \,.
\end{equation}
In fact, for each $p\in \Sigma_q$ and each $(i,j)\in \Cal T_p$, since $ Z^{ij}_p\in H^\ell_q(M)$, which implies $K_p(\tau) \left( Z^{ij}_p\right) - Z^{ij}_p\in  H^\ell_q(M)$, if $\ell < 1 + (i+j)/(m_p+1)$, by Lemma \ref{lemma:SO}, there exist constants $C'_\ell>0$, $C''_\ell>0$ such that
$$
\Vert K_p(\tau) \left( Z^{ij}_p\right) - Z^{ij}_p \Vert_\ell  \leq C'_\ell \vert K_p(\tau) \left( Z^{ij}_p\right) - Z^{ij}_p \vert_\ell  \leq  C''_\ell  \,\tau^{1+\frac{i+j}{m_p+1}-\ell}\,,
$$
while, if $\ell > 1 + (i+j)/(m_p+1)$, since $Z^{ij}_p  \in \bar H^{\infty}_q(M)$, 
$$
\Vert K_p(\tau) \left( Z^{ij}_p\right) - Z^{ij}_p \Vert_\ell  \leq
\Vert K_p(\tau) \left( Z^{ij}_p\right)\Vert_\ell + \Vert  Z^{ij}_p \Vert_\ell   \leq 
C''_\ell \,\tau^{1+\frac{i+j}{m_p+1}-\ell}\,.
$$

\smallskip
\noindent   
The scale of Friedrichs Sobolev spaces admits a standard family of smoothing operators
$\{T_\tau \,\vert \, \tau>0\}$ such that the operator  $T_\tau : L^2_q(M) \to \bar H^{\infty}_q(M)$ is defined, 
for each $\tau>0$, by the following truncation of Fourier series. 
Let $\{e_n \,\vert \, n\in \N\}$ be an orthonormal basis of eigenfunctions of the Friedrichs Laplacian $\Delta^F_q$ and let $\lambda:\N\to \R^{+}\cup \{0\}$ be the corresponding sequence of eigenvalues. Then
$$
T_\tau (u) :=   \sum_{\tau^2 \lambda_n\leq 1}   \<u,e_n\>_q   \, e_n \,,     
\quad \text {\rm  if } \, u=  \sum_{n\in \N} \<u,e_n\>_q \,   e_n  \,.
$$
If $u \in \bar H^s_q(M)$, then $T_\tau(u) \to u$ in $ \bar H^{s}_q(M) $ (as $\tau \to 0^+$) and the following estimates hold. For all $r\in \R^+$, there exists a constant $C_{r,s}(q)>0$ such that 
\begin{equation}
\label{eq:trunc}
\begin{aligned}
\Vert T_\tau (u)-u\Vert_r  & \leq  C_{r,s}(q)\, \Vert u \Vert_s\,\tau^{s-r}   \,, \quad &\text{ \rm if }  \,\, s\geq r\,; \\
\Vert T_\tau (u)\Vert_r &\leq  C_{r,s}(q)\,\Vert u \Vert_s\,\tau^{-(r-s)}   \,, \quad &\text{ \rm if } \,\,  r\geq s\,. 
\end{aligned}
\end{equation}
If $p\in \Sigma_q$ and $(i,j)\in \Cal T_p$, the distribution $\delta_p^{ij}\in \bar H^{-s_{ij}}_q(M)$ for any $s_{ij}>0$ such that $i+j < (s_{ij}-1)(m_p+1)$. Hence there exists a constant $C_p^{ij}(q)>0$ such that, by estimates \pref{eq:trunc},   for all $u\in \bar H^s_q(M)$,
\begin{equation}
\label{eq:deltasboundone}
\vert \delta_p^{ij} \left(T_\tau(u)\right)\vert  \leq C_p^{ij}(q) \,\Vert u \Vert _s \, \max \{1,\tau^{s-s_{ij}}\}\,.
\end{equation}
If $u \in H^s_q(M)$, then $\delta^p_{ij} (u)=0$, for all $p\in \Sigma_q$ and all $(i,j)\in \Cal T^s_p$. Hence, if $s_{ij} \leq s$ and  $i+j < (s_{ij}-1)(m_p+1)$, 
\begin{equation}
\label{eq:deltasboundtwo}
\vert \delta_p^{ij} \left(T_\tau(u)\right)\vert  = \vert \delta_p^{ij} \left(T_\tau(u)-u\right)\vert \leq
C^{ij}_p (q) \, \Vert u \Vert _s \,\tau^{s-s_{ij}} \,.
\end{equation}
\noindent The following estimates hold. Let $s\in \R^+$ and $\ell \in \N$. For any $\epsilon>0$, there exists a constant $C_{\ell,s}(\epsilon)>0$ such that,  for all $\tau\in (0,1]$ and all $u\in H^s_q(M)$,
\begin{equation}
\label{eq:Sktestone}
\Vert S^k_\tau \circ T_\tau (u) - T_\tau(u) \Vert_\ell \leq C_{\ell,s}(\epsilon)\, \Vert u \Vert_s  
\,\tau^{s-\ell -\epsilon} \,.
\end{equation}
In fact, if $p\in \Sigma_q$ and $(i,j)\in   \Cal T_p^k$, 
\begin{equation}
 \vert \delta_p^{ij} \left (T_\tau(u)\right) \vert  \,\leq \,  C_p^{ij}(q) \,\Vert u \Vert _s \, \,\tau^{s-s_{ij}}\,;
\end{equation}
for any $s_{ij}> 1+ (i+j)/(m_p+1)\geq s$, if $(i,j)\in   \Cal T_p^k\setminus \Cal T^s_p$, and
for any $ 1+ (i+j)/(m_p+1) < s_{ij} \leq s$, if $(i,j)\in   \Cal T_p^s$. The claim \pref{eq:Sktestone} then follows from \pref{eq:Sk1}. By estimates \pref{eq:trunc} and \pref{eq:Sktestone}, for any 
$\epsilon>0$, there exists a constant $C'_{\ell,s}(\epsilon)>0$ such that, for all $\tau\in (0,1]$ and for 
all $u\in H^s_q(M)$,
\begin{equation}
\label{eq:Sktesttwo}
\Vert S^k_\tau \circ T_\tau (u) - u \Vert_\ell \leq C'_{\ell,s}(\epsilon)\, \Vert u \Vert_s  \,\tau^{s-\ell-\epsilon}\,\,.
\end{equation}
Let $\{ \Cal S^k(\tau) \,\vert \, \tau\in (0,1]\}$ be the family of operators $\Cal S^k(\tau): L^2_q(M) \to  H^k_q(M)$ defined as follows: for each $\tau\in (0,1]$,
\begin{equation}
\label{eq:CalSkt}
 \Cal S^k(\tau) := S^k_\tau \circ T_\tau  \,.
\end{equation}
By estimate \pref{eq:Sktesttwo}, for any $\epsilon>0$, there exists a constant $C''_{\ell,s}(\epsilon)>0$ such that, for all $\tau\in (0,1]$ and for all $u\in H^s_q(M)$,
\begin{equation}
\label{eq:Sktestthree}
\Vert \Cal S^k(\tau) (u) -  \Cal S^k(\tau/2) (u)\Vert_\ell \leq C''_{\ell,s}(\epsilon)\, \Vert u \Vert_s  
\,\tau^{s-\ell-\epsilon}\,\,.
\end{equation}
Since $\Cal S^k(\tau) (u)  \in H^k_q(M)$ for all $\tau\in (0,1]$, by the interpolation inequality proved
in Lemma \ref{lemma:intineq} it follows that, for any $r\in [0,k]$ and for any $\epsilon>0$ there exists 
$C_{r,s}(\epsilon)>0$ such that, for all $\tau\in (0,1]$ and for all $u\in H^s_q(M)$,
\begin{equation}
\label{eq:Sktestfour}
\vert \Cal S^k(\tau) (u) -  \Cal S^k(\tau/2) (u)\vert_r \leq C_{r,s}(\epsilon)\, \Vert u \Vert_s  \,\tau^{s-r-\epsilon}\,\,.
\end{equation}
It follows that, for every $n\in \N$ and for every $\tau\in (0,1]$,
\begin{equation}
\label{eq:Sktestfive}
\vert \Cal S^k(\tau/2^n) (u) -  \Cal S^k(\tau/2^{n+1}) (u)\vert_r \leq C_{r,s}(\epsilon)\, \Vert u \Vert_s  
2^{-n(s-r-\epsilon)}\,\tau^{s-r-\epsilon}\,.
\end{equation}
If $s>r$ and $0<\epsilon<s-r$, the sequence $\{ S^k(\tau/2^n) (u)\}_{n\in \N}$ is Cauchy and therefore convergent  to the function $u\in H^s_q(M) \subset H^r_q(M)$ in the Hilbert space $H^r_q(M)$. It follows that, for all $\tau\in (0,1]$ and for all $u\in H^s_q(M)$,
\begin{equation}
\label{eq:Sktestsix}
\vert \Cal S^k(\tau) (u) -  u\vert_r \leq C_{r,s}(\epsilon)\, \Vert u \Vert_s \,\tau^{s-r-\epsilon}\,.
\end{equation}
If $s\leq r$ and $s-r\leq 0<\epsilon$, we argue as follows. By estimate \pref{eq:Sktestfour}, for each 
$\tau \leq 1/2$ we have
\begin{equation}
\label{eq:Sktestseven}
\vert \Cal S^k(2\tau) (u) - \Cal S^k(\tau)(u) \vert_r \leq C_{r,s}(\epsilon)\, \Vert u \Vert_s \,
 (2\tau)^{s-r-\epsilon}\,,
\end{equation}
hence if $\tau\leq  2^{-n}$, for all $0\leq k < n$,
\begin{equation}
\label{eq:Sktesteight}
\vert \Cal S^k(2^{k+1}\tau) (u) - \Cal S^k(2^k \tau)(u) \vert_r \leq C_{r,s}(\epsilon)\, \Vert u \Vert_s\, 
2 ^{(k+1)(s-r-\epsilon)} \,\tau^{s-r-\epsilon}\,,
\end{equation}
It follows that, there exists a constant $C'_{r,s}(\epsilon)>0$ such that
\begin{equation}
\label{eq:Sktestnine}
\vert \Cal S^k(2^n \tau) (u) - \Cal S^k(\tau)(u) \vert_r \leq C'_{r,s}(\epsilon)\, \Vert u \Vert_s\, 
2 ^{n (s-r-\epsilon)} \,\tau^{s-r-\epsilon}\,,
\end{equation}
 For every $\tau\in (0,1]$, let $n(\tau)$ be the maximum $n\in \N$ such that $2^n \tau \leq 1$. By
 this definition it follows that $1/2 < 2^{n(\tau)} \tau \leq 1$, hence by estimate \pref{eq:Sktesttwo} there
 exists a constant $C_k>0$ such that, for all $u\in H^s_q(M)$, 
 $$
  \vert \Cal S^k(2^{n(\tau)} \tau)(u) \vert_r \leq 
 \sup_{1/2 \leq \tau  \leq 1}  \Vert \Cal S^k(\tau)(u) \Vert_k\,\,< \,\, C_k \Vert u\Vert_s\,.
 $$  
It follows that there exists $C''_{r,s}(\epsilon)>0$ such that, for all $\tau\in (0,1]$ and for all $u\in H^s_q(M)$,
\begin{equation}
\label{eq:Sktestten}
\vert \Cal S^k(\tau)(u) \vert_r \leq C''_{r,s}(\epsilon)\, \Vert u \Vert_s \,\tau^{s-r-\epsilon}\,,
\end{equation}
 \end{proof}
 \noindent By Lemma \ref{lemma:comparison} and Theorem \ref{thm:smoothingop}, we have the
 following comparison estimate for the (Friedrichs) weigthed Sobolev norms :
 
 \begin{corollary} For any $0<r<s$ there exists constants $C_r>0$ and $C_{r,s} >0$ such
  that, for all $u\in H^s_q(M)$, the following inequalities hold:
 $$
C_r^{-1}  \Vert u \Vert_r  \leq  \vert u\vert_r \leq C_{r,s} \, \Vert u \Vert_s\,.
 $$
 \end{corollary}
 
 \noindent Finally, we derive a crucial interpolation estimate for the dual weighted Sobolev norms:
 \begin{corollary} 
 \label{cor:dualintineq} Let $0\leq s_1 < s_2$. For any $s_1 \leq r<s \leq s_2$ there exists a constant $C_{r,s}>0$ such that  for any distribution $u\in H^{-s_1}_q(M)$ the following interpolation inequality holds:
 \begin{equation}
 \label{eq:dualintineq}
 \vert u \vert _{-s } \,\, \leq \, \, C_{r,s}  \,\, \vert u \vert _{-s_1} ^{ \frac{s_2-r}{s_2-s_1}} \,\,  
 \vert u \vert _{-s_2}^{\frac{r-s_1}{s_2-s_1}} \,.
 \end{equation} 
  \end{corollary}
  \begin{proof} Let $k\in \N$ be any integer larger than $s_2 >s_1$ and let $\Cal S^k(\tau) :L^2_q(M)
  \to H^k_q(M)$ be the family of smoothing operators constructed above. By Theorem \ref{thm:smoothingop}, since $0\leq r-s_1 < s-s_1$ and any $r-s_2 < s-s_2\leq 0$, there exists a constant $C^k_{r,s}>0$ such  that the following holds: for any $u\in H^{-s_1}_q(M)\setminus\{0\}$, 
  any  $v \in H^s_q(M)$ and for all $\tau\in (0,1]$, 
 \begin{equation}
 \begin{aligned}
  \vert \<u, v\>\vert   &\leq \vert u \vert_{-s_1}  \vert v -\Cal S^k(\tau)(v) \vert_{s_1}  \, + \, 
  \vert u\vert_{-s_2} \vert \Cal S^k(\tau) (v) \vert_{s_2} \\
  &\leq C^k_{r,s} \, \{ \tau^{ r-s_1} \vert u \vert_{-s_1} +
  \tau^{r-s_2}\vert u \vert_{-s_2}  \}\,  \vert v \vert_s\,.  
   \end{aligned}
  \end{equation}
  The interpolation inequality \pref{eq:dualintineq} then follows by taking 
  $$
  \tau = \left( \frac{ \vert u \vert_{-s_2}}{ \vert u \vert_{-s_1} }\right) ^{\frac{1}{s_2-s_1}} \in (0,1] \,.
  $$
  \end{proof}

\section{The cohomological equation}
\label{distsol}

\subsection{Distributional solutions}
\label{DS}
\noindent In this section we give a streamlined version of the main argument of  \cite{F97} 
(Theorem 4.1) with the goal of establishing the sharpest bound on the loss of Sobolev regularity 
within the reach of the methods of \cite{F97}. We were initially motivated by a question of Marmi, 
Moussa and Yoccoz who found for \emph{almost all }orientable quadratic differentials a loss of regularity
of  $1+ {\rm BV}$  (they find bounded solutions for absolutely continuous data with first derivative of bounded variation under finitely many independent compatibility conditions and corresponding 
results for higher smoothness) \cite{MMY03},  \cite{MMY05} . The results of this section, as those
of \cite{F97}, hold for \emph{all} orientable quadratic differentials.

\medskip
\noindent There is a natural action of the circle group $S^1\equiv SO(2,\R)$ on the space $Q(M)$
of holomorphic quadratic differentials on a Riemann surface $M$:
$$
r_\theta(q):=  e^{i\theta} q \,, \quad \text{ for all } \,\, (r_\theta,q) \in SO(2,\R)\times Q(M)\,.
$$
Let $q_\theta$ denote the quadratic differential $r_\theta(q)$ and let $\{S_\theta, T_\theta\}$ denote the frame (introduced in \S \ref{WSS}) associated to the quadratic differential $q_\theta$ for any $\theta \in S^1$. We have the following formulas:
\begin{equation}
\begin{aligned}
S_{\theta} &= \cos \left(\frac{\theta}{2}\right) \,S_q  +  \sin  \left(\frac{\theta}{2}\right)\, T_q =  \frac{  e^{- i\frac{\theta}{2}} }{2} \,\partial_q^+
\,+ \, \frac{ e^{i \frac{\theta}{2}}}{2} \,\partial_q^- \,; \\
T_{\theta} &= -\sin \left(\frac{\theta}{2}\right)\, S_q  +  \cos \left(\frac{\theta}{2}\right) \,T_q =  \frac{ e^{- i \frac{\theta}{2}}}{2i} \,\partial_q^+
\,- \, \frac{ e^{i \frac{\theta}{2}}}{2i} \,\partial_q^- \,;
\end{aligned}
\end{equation}

\medskip
\begin{definition} 
\label{def:distsol}
Let $q$ be an orientable  quadratic differential. A distribution $u\in \bar H^{-r}_q(M)$ will be called a \emph {(distributional) solution} of the cohomological equation $S_q u =f$ for a given function $f\in  \bar H^{-s}_q(M)$  if 
$$
\<u, S_q v\> = - \< f, v\>\, ,\quad \text{ \rm for all } \,\, v\in H^{r+1}_q(M) \cap \bar H^{s}_q(M)\,.
$$
\end{definition}

\smallskip
\noindent Let ${\bar {\Cal H}}_q^s(M) \subset {\bar H}^s_q(M)$, ${\Cal H}^s_q(M) \subset {H}^s_q(M)$ (for any $s\in \R$) be the subspaces orthogonal to constant functions, that is 
\begin{equation}
\begin{aligned}
{\bar {\Cal H}}_q^s(M) &:= \{f\in {\bar H}^s_q(M)\,\vert \, \<f,1\>_s=0\}\,, \\
{\Cal H}_q^s(M) &:= \{f\in H^s_q(M)\,\vert \, (f,1)_s=0\}\,.
\end{aligned}
\end{equation}
The spaces ${\bar {\Cal H}}_q^s(M) \subset {\bar H}^s_q(M)$ and ${\Cal H}^s_q(M) \subset H^s_q(M)$ coincide with the subspaces of functions of zero average for $s \geq 0$, and with the subspaces of distributions vanishing on constant functions for $s<0$. 

\begin{theorem}
\label{thm:distsol} Let $r>2$ and $p\in (0,1)$ be such that and $rp>2$. There exists a bounded 
linear operator
$$
\Cal U: {\bar {\Cal H}}_q^{-1}(M) \to L^p\left(S^1,   {\bar H}^{-r}_q(M)\right)
$$
such that the following holds. For any $f\in {\bar {\Cal H}}_q^{-1}(M)$ there exists a full measure subset 
${\Cal F}_r(f)\subset S^1$ such that  $u:=\Cal U(f)(\theta)\in {\bar H}^{-r}_q(M)$ is a distributional solution of the cohomological equation $S_{\theta} u= f$ for all $\theta \in  {\Cal F}_r(f)$.
\end{theorem}
\begin{proof} We claim that for any $r>2$, any $p\in (0,1)$ such that $pr>2$ and any $f\in 
{\bar {\Cal H}}_q^{-1}(M)$, there exists a measurable function $A_q:=A_q(p,r,f) \in L^p(S^1,\R^+)$ such that the following estimates hold. Let $\theta\in S^1$ be such that $A_q(\theta)<+\infty$.  For all $v\in H_q^{r+1}(M)$ 
we have
\begin{equation}
\label{eq:aprioribound}
\vert \<f,v\> \vert \leq A_q(\theta) \,\Vert S_{\theta}v \Vert_r\,\,.  
\end{equation}
In addition, the following bound for the $L^p$ norm of the function $A_q$ holds. There exists a constant $B_q(p)>0$ such that
\begin{equation}
\label{eq:Lpbound}
\vert A_q\vert_p \leq B_q(p) \, \Vert f \Vert_{-1} \,\,.
\end{equation}
Assuming the claim, we prove the statement of the theorem.  In fact, by the estimate 
\pref{eq:aprioribound} the linear map given by
\begin{equation}
\label{eq:functional}
S_{\theta}v\to -\<f,v\>\,\,, \quad \hbox{ for all } v\in H_q^{r+1}(M)\,\,,
\end{equation}
is well defined and extends by continuity to the closure of the range ${\bar R}_r(\theta)$ of the linear operator $S_{\theta}$ in ${\bar H}^r_q(M)$. Let $\Cal U(f)(\theta)$ be the extension uniquely defined 
by the condition that $\Cal U(f)(\theta)$ vanishes on the orthogonal complement of ${\bar R}_r(\theta)$
in ${\bar H}^r_q(M)$. By construction, for almost all $\theta\in S^1$ the linear functional $u:=\Cal U(f)(\theta)\in {\bar H}^{-r}_q(M)$ yields a distributional solution of the cohomological equation 
$S_{\theta}u=f$ whose norm satisfies the bound 
$$
\Vert \Cal U(f)(\theta) \Vert_{-r} \leq A_q(\theta)\,.
$$
By \pref{eq:Lpbound} the $L^p$ norm of the measurable function $\Cal U(f):S^1 \to {\bar H}^{-r}_q(M)$ satisfies the required estimate
$$
 \vert \Cal U(f) \vert_{p}:= \left(\int_{S^1} \Vert \Cal U(f)(\theta) \Vert_{-r}^p \,d\theta \right)^{1/p}
 \,\, \leq \,\,  B_q(p) \, \Vert f \Vert_{-1} \,\,.
 $$
 We turn now to the proof of the above claim. Let $R_q^{\pm} =  \left(\Cal M^{\mp}_{q}\right)^{\perp}$ 
 be the (closed) ranges of the Cauchy-Riemann operators  $\partial^{\pm}_{q}: H^{1}_{q}(M)\to L^2_q(M)$ (see Proposition \ref{prop:CR}). Following \cite{F97}, we introduce the linear operator 
 $U_q: R^-_q\to R^+_q$ defined as 
\begin{equation}
\label{eq:partiso}
U_q:=\partial^-_q (\partial^+_q)^{-1}\,\,.
\end{equation}
By Proposition \ref{prop:CR}, $(3)$, the operator $U_q$ is a partial isometry on $L^2_q(M)$, hence
by the standard theory of partial isometries on Hilbert spaces, it has a family of unitary extensions 
$U_J:L^2_q(M)\to L^2_q(M)$ parametrized by isometries $J: \Cal M^{+}_{q}\to \Cal M^{-}_{q}$ (see formulas $(3.10)$-$(3.12)$ in \cite{F97}). By definition the following identities hold on $H^1_q(M)$ (see formulas $(3.13)$ in \cite{F97}) :
\begin{equation}
\label{eq:CEidentity}
S_{\theta}=\frac{e^{i\frac{\theta}{2}}}{2}\bigl(U_J+e^{-i\theta}\bigr)\partial^+_q
=\frac{e^{-i\frac{\theta}{2}}}{2}\bigl(U_J^{-1}+e^{i\theta}\bigr)\partial^-_q\,\,.
 \end{equation}
The proof of estimate \pref{eq:aprioribound} is going to be based on properties of the resolvent 
of the operator $U_J$. In fact, the proof of \pref{eq:aprioribound} is based on the results, 
summarized in \cite{F97}, Corollary 3.4, concerning the non-tangential boundary behaviour 
of the resolvent of a unitary operator on a Hilbert space, applied to the operators $U_J$,  
$U_J^{-1}$ on $L^2_q(M)$. The Fourier analysis of \cite{F97}, \S 2, also plays a relevant role 
through Lemma 4.2 in \cite{F97} and the Weyl's asymptotic formula (Theorem \ref{thm:Weyl}). 

\smallskip
\noindent Following \cite{F97}, Prop. 4.6A, or \cite{F02}, Lemma 7.3, we prove that there exists a constant $C_q>0$ such that the following holds. For any distribution $f\in H^{-1}_q(M)$ there exist  (weak) solutions $F^{\pm} \in L^2_q(M)$ of the equations  $\partial_q^{\pm} F^{\pm}=f$ such that 
\begin{equation}
\label{eq:CRsolbound}
\vert F^{\pm}\vert_0 \leq C_q \, \Vert f \Vert_{-1}\,.
\end{equation}
In fact, the maps given by
\begin{equation}
\label{eq:CRsol}
 \partial^{\pm}_q v  \to  - \<f,v\> \, \,, \quad \hbox{ for all }v\in H^1_q(M)\,,
\end{equation}
are bounded linear functionals on the (closed) ranges $R^{\pm}_q\subset L^2_q(M)$ (of the Cauchy-Riemann operator $\partial^{\pm}_q:H^1_q(M)\to L^2_q(M)$. In fact, the functionals are well-defined since $f$ vanishes on constant functions, that is, on the kernel of the Cauchy-Riemann operators, and it is bounded since by Poincar\'e inequality (see \cite{F97}, Lemma 2.2 or \cite{F02}, Lemma 6.9) there exists a constant $C_q>0$ such that, for any $v\in {\Cal H}^1_q(M) \subset H^1_q(M)$,
\begin{equation}
\label{eq:solboundone}
\vert \<f,v\> \vert \,\leq\, \Vert f\Vert_{-1} \Vert v \Vert_1 \leq C_q \, \Vert f\Vert_{-1} \, 
\vert \partial^{\pm}_q v \vert_0\,\,.
\end{equation}
Let $\Phi^{\pm}$ be the unique linear extension of the linear map \pref{eq:CRsol} to $L^2_q(M)$ 
which vanishes on the orthogonal complement of $R^{\pm}_q$ in $L^2_q(M)$. By \pref{eq:solboundone}, the functionals $\Phi^{\pm}$  are bounded on $L^2_q(M)$ with norm
$$
\Vert \Phi^{\pm} \Vert  \leq   C_q \, \Vert f\Vert_{-1}\,.
$$
By the Riesz representation theorem, there exist two (unique) functions $F^{\pm}\in L^2_q(M)$ such
that
$$
\<v, F^{\pm}\>_q  \,=\,   \Phi^{\pm} (v)\,, \quad \text{ \rm for all } \, v\in L^2_q(M)\,.
$$
The functions $F^{\pm}$ are by construction (weak) solutions of the equations $\partial_q^{\pm} F^{\pm}=f$ satisfying the required bound \pref{eq:CRsolbound}. 

\smallskip
\noindent The identities \pref{eq:CEidentity}  immediately imply that
\begin{equation}
\label{eq:keyidentity}
\begin{aligned}
\<{\partial}^{\pm}_qv,F^{\pm}\>_q&=2 e^{\mp i\frac{\theta}{2}}\, \<{\Cal R}^{\pm}_J(z)S_{\theta}v,F^{\pm}\>_q \\ 
&-(z+e^{\mp i\theta})\<{\Cal R}^{\pm}_J(z){\partial}^{\pm}_qv,F^{\pm}\>_q\,,
\end{aligned}
\end{equation}
where ${\Cal R}^{+}_J(z)$ and ${\Cal R}^{-}_J(z)$ denote the resolvents of the unitary operators 
$U_J$ and $U_J^{-1}$ respectively, which yield holomorphic families of bounded operators on
the unit disk $D\subset \C$.  

\smallskip
\noindent Let $r>2$ and let $p\in (0,1)$ be such that $pr>2$. Let $\Cal E=\{e_k\}_{k\in {\N}}$ be the orthonormal Fourier basis of the Hilbert space $L^2_q(M)$ described in \S \ref{fracsobspaces}. By Corollary 3.4 in \cite{F97} all holomorphic 
functions
\begin{equation}
\label{eq:matrixelements}
{\Cal R}^{\pm}_k(z):=\<{\Cal R}^{\pm}_J(z) e_k,F^{\pm}\>_q\,\,,\quad k \in {\N}\,,
\end{equation}
belong to the Hardy space $H^p(D)$, for any $0<p<1$. The corresponding non-tangential maximal functions $N_k^{\pm}$ (over cones of arbitrary fixed aperture $0<\alpha<1$) belong to the space $L^p(S^1,d\theta)$ and for all $0<p<1$  there exists a constant $A_{\alpha,p}>0$ such that 
the following inequalities hold:
\begin{equation}
\label{eq:Nkbound}
\vert N_k^{\pm} \vert_p\leq A_{\alpha,p} \vert e_k \vert_0\, \vert F^{\pm} \vert_0=
A_{\alpha,p} \, \vert F^{\pm} \vert_0 \leq  A_{\alpha,p}\, C_q \, \Vert f\Vert _{-1}  \,\,. 
\end{equation}
Let $\{\lambda_k\}_{k\in {\N}}$ be the sequence of the eigenvalues of the Dirichlet form $\Cal Q$ introduced in \S \ref{fracsobspaces}. Let $w\in  \bar H^r_q(M)$. We  have
\begin{equation}
\label{eq:Rsum}
\<{\Cal R}^{\pm}_J(z)w,F^{\pm}\>_q = \sum_{k=0}^{\infty}\<w,e_k\>_q\,\,{\Cal R}^{\pm}_k(z) \,\,,
\end{equation}
hence, by the Cauchy-Schwarz inequality,
\begin{equation}
\label{eq:Rbound}
\vert\<{\Cal R}^{\pm}_J(z)w,F^{\pm}\>_q\vert \leq \Bigl( \sum_{k=0}^{\infty}
\frac{\vert {\Cal R}^{\pm}_k(z)\vert^2}{ (1+\lambda_k)^r}\,  \Bigr) ^{1/2}
 \Vert w\Vert_r \,\,,
\end{equation}
Let $N^{\pm}(\theta)$ be the functions defined as 
\begin{equation}
\label{eq:Ndef}
N^{\pm}(\theta):= \Bigl( \sum_{k=0}^{\infty} \frac{\vert N_k^{\pm}(\theta)\vert^2}{ (1+\lambda_k)^r}
\Bigr)^{1/2}\,\,.
\end{equation}
Let $N^{\pm}(w)$ be the non-tangential maximal function for the holomorphic function $\<{\Cal R}^{\pm}_J(z)w,F^{\pm}\>_q$.  By formulas \pref{eq:Rbound} and \pref{eq:Ndef},  for all $\theta\in S^1$ and all functions $w \in {\bar H}^r_q(M)$, we have
\begin{equation}
\label{eq:Nwbound}
N^{\pm}(w)(\theta) \leq  N^{\pm}(\theta) \,  \Vert w\Vert_r  \,\,.
\end{equation}
The functions $N^{\pm} \in L^p(S^1,d\theta)$ for any $0<p<1$. In fact, by formula
\pref{eq:Nkbound} and (following a suggestion of Stephen Semmes) by the `triangular  inequality'
for $L^p$ spaces with  $0<p<1$, we have
\begin{equation}
\label{eq:NLpbound}
\vert N^{\pm}  \vert^p_p  \leq (A_{\alpha,p} \, C_q )^p \bigl(\sum_{k=0}^{\infty}
\frac{1}{ (1+\lambda_k)^{pr/2}}\bigr) \, \Vert f\Vert _{-1} ^p \,<\, +\infty \,\,.
\end{equation}
The series in formula \pref{eq:NLpbound} is convergent by the Weyl asymptotics
(Theorem \ref{thm:Weyl}) since $pr/2>1$.

\smallskip
\noindent By taking the non-tangential limit as $z\to -e^{\mp i\theta} $ in the identity 
\pref{eq:keyidentity}, formula \pref{eq:Nwbound} implies that, for all $\theta\in S^1$
such that $N^{\pm}(\pi\mp \theta)<+\infty$, 
$$
\vert \<{\partial}^{\pm}v,F^{\pm}\>_q \vert \leq N^{\pm}(\pi\mp \theta) \,
 \Vert S_{\theta}v \Vert_r  \,\,,
 $$
hence the required estimates \pref{eq:aprioribound} and \pref{eq:Lpbound} are proved
wih the choice of the function $A_q(\theta):= N^+(\pi-\theta)$ or  $A_q(\theta):= 
N^-(\pi+\theta)$ for all $\theta\in S^1$. Since the claim is proved the result follows.
\end{proof}

\begin{theorem} 
\label{thm:CEdistribution} 
Let $r>2$. For almost all $\theta\in S^1$ (with respect  to the Lebesgue measure), there exists a 
constant $C_r(\theta)>0$ such that, for all $f\in {\bar H}^{r-1}_q(M)$ such that $\int_M f\,\omega_q=0$, 
the cohomological equation $S_{\theta}u=f$ has a distributional solution $u\in {\bar H}^{-r}_q(M)$ 
satisfying the following estimate:
$$\Vert u\Vert_{-r}\leq C_r(\theta)\, \Vert f \Vert_{r-1}\,\,.$$
\end{theorem}
\begin{proof} Let $\Cal E=\{e_k\}_{k\in {\N}}$ be the orthonormal Fourier basis of the Hilbert space $L^2_q(M)$ described in \S \ref{fracsobspaces}. Let $r>2$ and $p\in (0,1)$ be such that $pr>2$.
By Theorem \ref{thm:distsol}, for any $k\in \N\setminus\{0\}$ there exists a function with distributional values $u_k := \Cal U(e_k)\in L^p\left(S^1, {\bar H}^{-r}_q(M)\right)$ such that the following holds. There 
exists a constant $C_q:=C_q(p,r) >0$ such that
\begin{equation}
\label{eq:kpnorm}
\left( \int_{S^1}  \Vert u_k (\theta) \Vert_{-r}^p \, d\theta \right)^{1/p} \,\leq\,  C_q\, \Vert e_k \Vert_{-1}
\,= \,C_q \,(1+\lambda_k)^{-1/2}  \,\,.
\end{equation}
In addition, for any $k\in \N\setminus\{0\}$, there exists a full measure set $\Cal F_k\subset S^1$ such that, for all $\theta\in {\Cal F}_k$, the distribution $u:=u_k(\theta) \in {\bar H}^{-r}_q(M)$ is a (distributional) solution of the cohomological equation $S_{\theta}\,u  = e_k$. 

\smallskip
\noindent Any function $f\in {\bar H}^{r-1}_q(M)$ such that $\int_M f\,\omega_q=0$ has a Fourier decomposition in $L^2_q(M)$:
$$
f = \sum_{k\in\N\setminus\{0\}} \<f,e_k\>_q \, e_k\,\,.
$$
A (formal) solution of the cohomological equation $S_{\theta} u =f$ is therefore given by the
series
\begin{equation}
\label{eq:solseries}
u_{\theta}:=  \sum_{k\in\N\setminus\{0\}} \<f,e_k\>_q\, u_k(\theta)\,\,.
\end{equation}
By the triangular inequality in ${\bar H}^{-r}_q(M)$ and by H\"older inequality, we have
$$
\Vert u_{\theta}\Vert_{-r} \leq  \left(\sum_{k\in\N\setminus\{0\}}   
\frac{ \Vert u_k(\theta)\Vert_{-r}^2}{(1+\lambda_k)^{r-1}}\right)^{1/2} \, \Vert f\Vert_{r-1}\,\,,
$$
hence by the `triangular inequality'  for $L^p$ spaces (with $0<p<1$) and by the estimate
\pref{eq:kpnorm},
\begin{equation}
\label{eq:Lpsol}
\int_{S^1}  \Vert u_{\theta} \Vert_{-r}^p \, d\theta  \,\,\leq \,\, C_q^p \, \left( \sum_{k\in\N\setminus\{0\}}  \frac{1}{(1+\lambda_k)^{pr/2}} \right)  \, \Vert f\Vert_{r-1}^p\,\,.
\end{equation}
Since $pr/2>1$ the series in \pref{eq:Lpsol} is convergent, hence by Chebyshev inequality 
for the space $L^p(S^1)$, there exists a full measure set $\Cal B \subset S^1$ such that, for all 
$\theta\in \Cal B$, formula \pref{eq:solseries} yields a well-defined distribution $u_{\theta}\in
 {\bar H}^{-r}_q(M)$ and there exists a constant $C_q(\theta)>0$ such that
\begin{equation}
\label{eq:solboundtwo}
\Vert u_{\theta} \Vert_{-r}\,  \leq  \, C_q(\theta) \,  \Vert f \Vert_{r-1}\,\,.
\end{equation}
The set ${\Cal F}= \cap {\Cal F}_k \cap \Cal B$ has full measure and for all $\theta\in \Cal F$,
for all $k\in \N\setminus\{0\}$, the distribution $u_k(\theta)\in  {\bar H}^{-r}_q(M)$ is a solution of the equation $S_{\theta} u=e_k$. It follows that $u_{\theta}\in {\bar H}^{-r}_q(M)$ is a solution of the cohomological equation $S_{\theta}u=f$ which satisfies the required bound \pref{eq:solboundtwo}.
\end{proof}

\noindent We finally derive a result on distributional solutions of the cohomological equation for 
distributional data of arbitrary regularity:
\begin{corollary}
\label{cor:CEdistribution} 
For any  $s\in \R$ there exists $r>0$ such that the following holds. For almost all $\theta\in S^1$ (with respect  to the Lebesgue measure), there exists a  constant $C_{r,s}(\theta)>0$ such that, for all $F\in {\bar H}^{s}_q(M)$ orthogonal to constant functions, the cohomological equation $S_{\theta}U=F $ has a distributional solution $U\in {\bar H}^{-r}_q(M)$ satisfying the following estimate:
$$\Vert U\Vert_{-r}\leq C_{r,s}(\theta)\, \Vert F \Vert_{s}\,\,.$$
\end{corollary}   
\begin{proof} Since $F\in \bar H^s_q(M)$ is orthogonal to constant functions, for every $k\in \N$ there exists $f_k\in \bar H^{s+2k}_q(M)$, orthogonal to constant functions, such that $(I-\Delta_q^F)^k f_k = F$. In fact, the family of Friedrichs Sobolev spaces $\{ \bar H^s_q(M) \vert s\in \R\}$ is defined in terms of the Friedrichs extension $\Delta_q^F$ of the Laplace operator $\Delta_q$ of the flat metric determined by the quadratic differential. Let $n \in \N$ be the minimum integer such that $\sigma:=s+2n +1> 2$. By Theorem \ref{thm:CEdistribution}, for almost all $\theta\in S^1$ (with respect  to the Lebesgue measure) the cohomological equation $S_{\theta}u=f_n$ has a distributional solution $u_\theta \in
 {\bar H}^{-\sigma}_q(M)$ satisfying the following estimate:
\begin{equation}
\label{eq:unest}
\Vert u_\theta \Vert_{-\sigma}\leq C_{\sigma}(\theta)\, \Vert f_n \Vert_{\sigma-1}\,.
\end{equation}
Let $U_\theta := (I-\Delta^F_q)^n u_\theta \in  {\bar H}^{-\sigma-2n}_q(M)$. It follows immediately from the estimate \pref{eq:unest} and from the definitions that
\begin{equation*}
\Vert U_\theta \Vert_{-\sigma-2n} = \Vert u_n \Vert_{-\sigma}\leq C_{\sigma}(\theta)\, \Vert f_n \Vert_{\sigma-1}
=  C_{\sigma}(\theta) \, \Vert F \Vert_{s}\,.
\end{equation*}
Finally $U_\theta\in {\bar H}^{-\sigma-2n}_q(M)$ is a distributional solution of the cohomological equation $S_\theta U=F$, for almost all $\theta\in S^1$. In fact, for any $v \in H^{\sigma+2n+1}_q(M)$, the function $S_\theta v \in H^{\sigma+2n}_q(M)$, hence $\,(I-\Delta_q)^n  S_\theta v = S_\theta (I-\Delta_q)^n v\,$ and, since the distribution $u_\theta\in H^{-\sigma}_q(M)$ is a solution of the cohomological equation $S_\theta u =f_n$, for almost all $\theta\in S^1$, and $(I-\Delta^F_q)^n v\in H^{\sigma+1}_q(M)$,
 \begin{equation*}
 \begin{aligned}
\<U_\theta, S_\theta v\> &= \<(I-\Delta^F_q)^n u_\theta, S_\theta v\> 
=  \<u_\theta, S_\theta (I-\Delta_q)^n v\>  \\
 &= - \<f_n, (I-\Delta^F_q)^n v\> 
 = - \<(I-\Delta^F_q)^n f_n, v\> 
 =  -\<F, v\>\,,
 \end{aligned}
\end{equation*}
as required by the definition of distributional solution of the cohomological equation (Definition \ref{def:distsol}).
\end{proof}

\subsection{Invariant distributions and basic currents}
\label{invdistandbasic}

\noindent Invariant distributions yield obstructions to the existence of smooth solutions of the cohomological equation. We derive below from Theorem \ref{thm:distsol} a sharp version of the main results of  \cite{F02}, \S 6, about the Sobolev regularity of invariant distributions. We then recall the structure theorem proved in that paper on the space of invariant distributions (see \cite{F02}, Th. 7.7).

\smallskip
\noindent Invariant distributions for the horizontal [respectively vertical] vector field of an orientable quadratic differential $q$ are closely related to {\it basic currents} (of dimension and degree equal to $1$) for the horizontal [vertical] foliation $\Cal F_q$ [ $\Cal F_{-q}$].  The notion of a basic current for a measured foliation on a Riemann surface has been studied in detail in \cite{F02}, \S 6, in the context
of weighted Sobolev spaces with integer exponent. We outline below some of the basic constructions and results on basic currents and invariant distributions which carry over without modifications to
the more general context of fractional weighted Sobolev spaces. Finally, we derive from Theorem 
\ref{thm:distsol} a result on the Sobolev regularity of basic currents (or invariant distributions) which improves upon a similar result proved in \cite{F02} (see Theorem 7.1 (i)).

\smallskip 
\noindent Let $\Sigma\subset M$ be a finite subset. The space $\Cal D(M\setminus\Sigma)$ will
denote the standard space of de Rham currents on the open manifold $M\setminus \Sigma$, that
is the dual of the Fr\'echet space $\Omega_c(M\setminus\Sigma)$ of differential forms with 
compact support in $M\setminus\Sigma$.  A homogenous current of dimension $d\in \N$ 
(and degree $2-d$) on $M\setminus \Sigma$ is a continuous linear functional on the subspace $\Omega^d_c(M\setminus \Sigma)$ of diferential forms of degree $d$. The subspace of homogeneous currents of dimension $d$ on $M\setminus\Sigma$ will be denote by $\Cal D^d(M\setminus\Sigma)$.

\smallskip
\noindent Let $q$ be an orientable quadratic differential on a Riemann surface $M$. Let $\Sigma_q$ 
be the (finite) set of its zeroes. In \cite{F02} we have introduced the following space $\Omega_q(M)$ of smooth test forms on $M$. 

\begin{definition}
For any $p\in M$ of (even) order $k=2m\in \N$ ($m=0$ if $p\not\in \Sigma_q$), let  $z:{\Cal U}_p\to \C$ be a canonical complex coordinate on a neighbourhood ${\Cal U}_p$ of $p\in M$, that is a complex coordinate  such that $z(p)=0 $ and $q=z^k dz^2$ on ${\Cal U}_p$. Let $\pi_p:{\Cal U}_p \to \C$ be the (local) covering map defined by 
 $$
 \pi_p(z):=\frac{z^{m+1}}{m+1}\,, \quad z\in \C\,.
 $$ 
 The space $\Omega_q(M)$ is defined as the space of smooth forms $\alpha$ on $M$ such that
 the following holds: for all $p\in M$, there exists a smooth form $\lambda_p$ on a neighbourhood of $0\in \C$ such that $\alpha=\pi_p^{\ast}(\lambda_p)$ on ${\Cal U}_p'\subset {\Cal U}_p$.  The space $\Omega_q(M)$ is the direct sum of the subspaces  $\Omega^d_q(M)$ of homogeneous
 forms of degree $d\in \{0,1,2\}$. The spaces $\Omega^d_q(M)$, for any $d\in \{0,1,2\}$,  and 
 $\Omega_q(M)$ can be endowed with a natural Fr\'echet topology modeled on the smooth topology in every coordinate neighbourhood. 
 \end{definition}

\begin{lemma}
For any orientable quadratic differential $q\in Q(M)$, the space of functions $\Omega^0_q(M)$ 
is dense in the space $H^\infty_q(M)$ endowed with the inverse limit Fr\'echet topology 
induced by the family of weighted Sobolev norms.
\end{lemma}
\begin{proof}
By definition, the MacLaurin series of any $f\in \Omega_q(M)$ with respect to a canonical complex coordinate $z$ for $q$ at every $p\in \Sigma_q$ (of order $2m_p$) has the following form:
\begin{equation}
\label{eq:Texp}
f(z)= \sum_{h,k\in \N}  f_{hk} z^{h(m_p+1)} {\bar z}^{k(m_p+1)} \,.
\end{equation}
By Lemmas \ref{lemma:localbarH} and \ref{lemma:localH}, $f \in H^{\infty}_q(M)$. Thus
$ \Omega_q(M) \subset H^{\infty}_q(M)$. 

\smallskip
\noindent Let $F\in H^{\infty}_q(M)$. By Lemma \ref{lemma:comparison} the function
 $F\in C^{\infty}(M)$ and by Lemmas  \ref{lemma:localbarH} and \ref{lemma:localH} its  
 MacLaurin series has the form \pref{eq:Texp} at every $p\in \Sigma_q$ (of order $2m_p$).
By Borel's thereorem and by a partition of unity argument, there exists a function $f \in
 \Omega_q(M)$ such that $F-f \in C^{\infty}(M)$ vanishes at infinite order at $\Sigma_q$. 
 Let ${\Cal U}_\tau$ be the open neighbourhood of $\Sigma_q$  which is the union of a finite number 
 of {\it disjoint }geodesic disks $D_\tau(p)$ of radius $\tau\in (0,\tau_0)$, each centered at a point 
 $p\in \Sigma_q$. Let $\phi_\tau :M \to [0,1]$ be a smooth function such that  $(a)$\,$\phi_\tau\in C_0^{\infty}(M\setminus\Sigma_q)$, $(b)$\, $\phi_\tau\equiv 1$ on $M \setminus {\Cal U}_\tau$ and 
 $(c)$\, for each $(i,j)\in \N\times\N$ there exists a constant $C_{ij}>0$ such that, for all 
 $\tau\in (0,\tau_0)$, 
 $$
 \max_{x\in M} \vert S^{i} T^{j}\phi_\tau (x)\vert   \,\, \leq   \,\, \frac{C_{ij}}{\tau^{i+j} }\,.
 $$
 If can be proved that, since $F-f$ vanishes at infinite order at $\Sigma_q$,
 $$
 f+ \phi_\tau (F-f) \, \to F   \quad \text{ \rm in } \,\,H^{\infty}_q(M)\,, \quad \text{ as }\,\, \tau\to 0^+\,,
 $$
 which implies, since by construction $f+ \phi_\tau (F-f)  \in \Omega_q^0(M)$, that $F$ belongs
 to the closure of $\Omega_q^0(M)$ in $H^{\infty}_q(M)$.
\end{proof}

\begin{definition}  The space $\Cal S_q(M) \subset \Cal D(M\setminus\Sigma_q)$ of  $q$\emph{-tempered currents }(introduced in \cite{F02}, \S 6.1) is the dual space of the Fr\'echet space $\Omega_q(M)$. A {\it homogeneous }$q$-tempered current of dimension $d$ (and degree $2-d$) is a continuous functional on the subspace $\Omega^d_q(M) \subset\Omega_q(M)$ of homogeneous forms of degree $d\in \{0,1,2\}$. The space of homogeneous currents of dimension $d$ (and degree $2-d$) will be denoted by $\Cal S^d_q(M)$.
\end{definition}

\noindent For any quadratic differential $q$ on $M$, there is a natural operator $^\ast$, which maps
 the space $\Cal D^0(M\setminus\Sigma_q)$ of currents of dimension $0$ and degree $2$ on the
 non-compact manifold $M \setminus \Sigma_q$ (which is naturally identified with the space of distributions on $M \setminus \Sigma_q$) bijectively onto the space $\Cal D^2(M\setminus\Sigma_q)$ 
 of currents of dimension  $2$ and degree $0$ on $M \setminus \Sigma_q$. The operator 
 $$
 ^\ast: \Cal D^0(M\setminus\Sigma_q)\to \Cal D^2(M\setminus\Sigma_q)
 $$ 
 is defined as follows. Let $\omega_q$  be the smooth area form associated with the (orientable) quadratic differential $q$ on $M$. It is a standard fact in the theory of currents that any distribution $U$ on the $2$-dimensional surface $M\setminus\Sigma_q$ can be written as $U=U^\ast \omega_q$ for a unique current $U^\ast$ of dimension $2$ and degree $0$.  Since $\omega_q\in \Omega^2_q(M)$, the
map $^\ast$ extends to a bijective map 
$$
^\ast: \Cal S^0_q(M) \to \Cal S^2_q(M)\,.
$$

\begin{definition} A distribution $\Cal D \in \Cal D^2(M\setminus\Sigma_q)$ is \emph{horizontally
[vertically] quasi-invariant }if $S\Cal D =0$ [$T\Cal D=0$] in  $\Cal D^2(M\setminus\Sigma_q)$.
A distribution $\Cal D \in \Cal S_q^2(M)$ is \emph{horizontally [vertically] invariant }if $S\Cal D =0$ [$T\Cal D=0$] in  $\Cal S_q^2(M)$. The space of horizontally [vertically] quasi-invariant distributions will be denoted by $\Cal I_q(M\setminus\Sigma_q)$ [$\Cal I_q(M\setminus\Sigma_q)$]
and the subspace  of horizontally [vertically] invariant distributions will be denoted by $\Cal I_q (M)$ [$\Cal I_{-q}(M)$]. 
\end{definition}

\begin{definition} For any $s\in \R^+$, let
\begin{equation}
\begin{aligned}
\Cal I^s_{\pm q}(M\setminus\Sigma_q) &:= \Cal I_ {\pm q}(M\setminus\Sigma_q)  \cap H^s_q(M) \,;\\
\Cal I^s_{\pm q}(M) &:= \Cal I_{\pm q}(M)  \cap H^s_q(M) \,.
\end{aligned}
\end{equation}
The subspaces $\Cal I^s_{\pm q}(M) \subset  H^{-s}_q(M)$ of \emph{ horizontally [vertically] invariant}  distributions can also be defined as follows: 
\begin{equation}
\begin{aligned}
\Cal I^s_q(M)&:= \{ \Cal D\in H^s_q(M) \,\vert \, S\Cal D = 0 \,\, \quad   \text{ in } \,\,  H^{-s-1}_q(M)\}\,;\\
[\Cal I^s_{-q}(M)&:= \{ \Cal D\in H^s_q(M) \,\vert \, T\Cal D = 0 \,\, \quad   \text{ in } \,\,  H^{-s-1}_q(M)\}]\,.
\end{aligned}
\end{equation}
The subspaces of horizontally [vertically] invariant distributions which can be extended to bounded functionals on Friedrichs weighted Sobolev spaces will be denoted by
\begin{equation}
\begin{aligned}
\bar {\Cal I}^s_{\pm q}(M\setminus\Sigma_q)&:= \Cal I _{\pm q}(M\setminus\Sigma_q) \cap   
\bar H^{-s}_q(M)\,; \\
\bar {\Cal I}^s_{\pm q}(M)&:= \Cal I_{\pm q}(M) \cap   \bar H^{-s}_q(M)\,.
\end{aligned}
\end{equation}
\end{definition}

\smallskip
\noindent Let ${\Cal V}_q(M)$ be the space of vector fields $X$ on $M\setminus\Sigma_q$ such that the contraction $\imath_X \alpha$ and the Lie derivative ${\Cal L}_X\alpha\in \Omega_q(M)$ for all $\alpha\in \Omega_q(M)$. 

\begin{definition}  A current $C\in \Cal D^1(M\setminus\Sigma_q)$ is \emph{horizontally [vertically] quasi-basic}, that is basic for ${\Cal F}_q$ [${\Cal F}_{-q}$] in the standard sense on $M\setminus\Sigma_q$, if the identities 
\begin{equation}
\label{eq:defbasic}
\imath_X C \,=\, \Cal L_X C \, =\, 0 
\end{equation}
hold in ${\Cal D}(M\setminus\Sigma_q)$ for all smooth vector fields $X$ tangent to ${\Cal F}_q$ 
[${\Cal F}_{-q}$] with compact support on $M\setminus \Sigma_q$. A $q$-tempered current $C\in {\Cal S}^1_q(M)$ is \emph{horizontally [vertically] basic }if the identities \pref{eq:defbasic} holds in $\Cal S_q(M)$ for all vector fields $X\in {\Cal V}_q(M)$, tangent to ${\Cal F}_q$ [${\Cal F}_{-q}$]  on $M\setminus\Sigma_q$. The vector spaces of horizontally [vertically] quasi-basic (real) currents will be denoted by ${\Cal B}_{q}(M \setminus\Sigma_q)$ [${\Cal B}_{-q}(M\setminus\Sigma_q)$] and the subspace of horizontally [vertically] basic (real) currents will be denoted by ${\Cal B}_{q}(M)$ 
[${\Cal B}_{-q}(M)$]. 
\end{definition}

\begin{definition}
For any $s\in \R$, the Friedrichs weighted Sobolev space of $1$-currents $\bar {W}^s_q(M)
\subset \Cal S_q(M)$ and the weighted Sobolev space of $1$-currents $W^s_q(M)\subset 
\Cal S_q(M)$ are defined  as follows: 
 \begin{equation}
\begin{aligned}
\bar {W}^s_q(M)&:= \{ \alpha \in \Cal S_q(M) \,\vert \, (\imath_S\alpha, \imath_T\alpha) \in 
\bar H_q^s(M) \times \bar H_q^s(M) \}\,;\\
W^s_q(M)&:= \{ \alpha \in \Cal S_q(M) \,\vert \, (\imath_S\alpha, \imath_T\alpha) \in 
H_q^s(M) \times H_q^s(M) \}\,.
\end{aligned}
\end{equation}
\end{definition}

\begin{definition} 
\label{def:currsobord}
For any $1$-current $C\in \Cal D(M\setminus\Sigma_q)$, the \emph{weighted Sobolev order }${\Cal O}_q^W(C)$ and the \emph{Friedrichs weighted Sobolev order }$\bar {\Cal O}_q^W(C)$ are the real 
numbers defined as follows:
\begin{equation}
\begin{aligned}
{\Cal O}_q^W(C):=& \inf\{ s\in\R \,\vert \, \Cal D \in W^{-s}_q(M)\}\,; \\
\bar {\Cal O}_q^W(C):=& \inf\{ s\in\R \,\vert \, C  \in \bar W^{-s}_q(M)\}\,.
\end{aligned}
\end{equation}
\end{definition}

\begin{definition} For any $s\in \R$, let 
\begin{equation}
\label{eq:Hsbcurr}
\begin{aligned}
\Cal B^s_{\pm q}(M\setminus\Sigma_q) &:= \Cal B_{\pm q}(M\setminus\Sigma_q)  
\cap {W}^s_q(M) \,;\\
\Cal B^s_{\pm q}(M) &:= \Cal B_{\pm q}(M)  \cap {W}^s_q(M) \,.
\end{aligned}
\end{equation}
The subspaces $\Cal B^s_{\pm q}(M) \subset  W^{-s}_q(M)$ of \emph{horizontally [vertically] basic currents}  can also be defined as follows: 
\begin{equation}
\begin{aligned}
\imath_S C =0 \,\,  \text{ in } \,\,  H^{-s}_q(M) \quad &\text{\rm and } \quad \Cal L_S C=0  \,\, 
 \text{ in } \,\,  W^{-s-1}_q(M)\,; \\
 [\imath_T C =0 \,\,  \text{ in } \,\,  H^{-s}_q(M) \quad  &\text{\rm and } \quad \Cal L_T C=0  \,\, 
 \text{ in } \,\,  W^{-s-1}_q(M)]\,.
\end{aligned}
\end{equation} 
The subspaces of basic currents which can be extended to bounded functionals on Friedrichs weighted Sobolev spaces will be denoted by
\begin{equation}
\label{eq:FHsbcurr}
\begin{aligned}
\bar {\Cal B}^s_{\pm q}(M\setminus\Sigma_q)&:= \Cal B_{\pm q}(M\setminus\Sigma_q) \cap  
 \bar {W}^{-s}_q(M)\,;\\
\bar {\Cal B}^s_{\pm q}(M)&:= \Cal B_{\pm q}(M) \cap   \bar {W}^{-s}_q(M)\,.
\end{aligned}
\end{equation}
\end{definition}

\noindent According to Lemma 6.5 of \cite{F02}, the notions of invariant distributions and basic currents are related (see also Lemma 6.6 in \cite{F02}):
\begin{lemma}
\label{lemma:DtoC}
A current $C\in {\Cal B}^s_{q}(M\setminus \Sigma_q)$ [$C\in {\Cal B}^s_{ q}(M)$] if and only 
if the distribution $C\wedge\Re(q^{1/2}) \in {\Cal I}^s_q(M\setminus\Sigma_q)$ [$C\wedge\Re(q^{1/2})
\in {\Cal I}^s_q(M)$]. A current $C\in {\Cal B}^s_{-q}(M\setminus \Sigma_q)$ [$C\in {\Cal B}^s_{-q}(M)$] if and only if the distribution $C\wedge\Im(q^{1/2})\in {\Cal I}^s_q(M\setminus\Sigma_q)$ [$C\wedge\Im(q^{1/2})\in {\Cal I}^s_q(M)$].  In addition, 
the map
\begin{equation}
\label{eq:DtoC}
\begin{aligned}
&\Cal D_q: C  \to -C\wedge\Re(q^{1/2})\,\,; \\
[&\Cal D_{-q}: C  \to C\wedge \Im(q^{1/2})] \,\,;  
\end{aligned}
\end{equation}
is a bijection from the space ${\Cal B}^s_q(M\setminus\Sigma_q)$ [${\Cal B}^s_{-q}(M\setminus \Sigma_q)$] onto the space ${\Cal I}^s_q(M \setminus\Sigma_q)$ [${\Cal I}^s_{-q}(M \setminus\Sigma_q)$], which maps the subspace ${\Cal B}^s_q(M)$ [ ${\Cal B}^s_{-q}(M)$]  onto the subspace ${\Cal  I}^s_q(M)$ [${\Cal  I}^s_{-q}(M)$]. The map \pref{eq:DtoC} also maps the space $\bar {\Cal B}^s_{q}(M)$ [$\bar {\Cal B}^s_{-q}(M)$]  onto $\bar {\Cal I}^s_{q}(M)$ [$\bar {\Cal I}^s_{-q}(M)$].
\end{lemma}

\subsection{Basic cohomology}
\label{ss:bc}
\noindent Let $\Cal Z (M\setminus \Sigma)\subset \Cal D^1(M\setminus\Sigma)$ denote the 
subspace of all (real) closed currents, that is, the space of all (real) de Rham currents $C\in \Cal D^1(M\setminus\Sigma)$ such that the exterior derivative $dC=0$ in $\Cal D(M\setminus\Sigma)$. 
Let $\Cal Z_q (M)\subset \Cal S^1_q(M)$ be the subspace of all (real) closed $q$-tempered currents, that is, the space of all $q$-tempered (real) currents $C$ such that $dC=0$ in $\Cal S_q(M)$.  It was proved in \cite{F02}, Lemma 6.2,  that the natural de Rham cohomology map
\begin{equation}
\label{eq:cohomapone}
j_q: \Cal Z (M\setminus \Sigma_q) \to H^1(M\setminus\Sigma_q, \R)
\end{equation}
has the property that the subspace of closed $q$-tempered currents is mapped onto the absolute
real cohomology of the surface, that is,
\begin{equation}
\label{eq:cohomaptwo}
j_q: \Cal Z_q (M) \to H^1(M,\R) \subset H^1(M\setminus\Sigma_q, \R)\,.
\end{equation}
It was also proved in \cite{F02}, Lemma 6.2', that quasi-basic and basic currents are closed, in the 
sense that the following inclusions hold:
\begin{equation}
\begin{aligned}
\Cal B_{\pm q} (M\setminus\Sigma_q) &\subset  \Cal Z (M\setminus \Sigma_q)\,, \\
\Cal B_{\pm q} (M) &\subset  \Cal Z_q (M)\,.
\end{aligned}
\end{equation}
 The images of the restrictions of the natural cohomology map to the various spaces of basic currents 
 are called the\emph{ horizontal [vertical] basic cohomologies }, namely the spaces
  \begin{equation}
  \begin{aligned}
 H^1_{\pm q}(M\setminus\Sigma_q,\R)&:=  j_q\left( \Cal B_{\pm q}(M\setminus\Sigma_q) \right) 
 \,\subset \, H^1(M\setminus\Sigma_q,\R)\,;Ê\\
  H^{1,s}_{\pm q}(M\setminus\Sigma_q,\R) &:=  j_q\left( \Cal B^s_{\pm q}(M\setminus\Sigma_q) \right) \,\subset \, H^1_{\pm q}(M\setminus\Sigma_q,\R)\,; \\
  H^1_{\pm q}(M,\R)&:=  j_q\left( \Cal B_{\pm q}(M) \right) 
 \,\subset \, H^1(M,\R)\,; \\
H^{1,s}_{\pm q}(M,\R) &:=  j_q\left( \Cal B^s_{\pm q}(M) \right) \,\subset \, H^1_{\pm q}(M,\R)\,.
\end{aligned}
\end{equation}
Following \cite{F02}, Theorem 7.1, we give below a description of the horizontal [vertical]  basic cohomologies for the orientable quadratic differential $q_\theta$, for any orientable holomorphic quadratic differential $q$ on $M$ and for almost all $\theta\in S^1$. The result we obtain below is stronger than Theorem 7.1 of \cite{F02} since it requires weaker Sobolev regularity assumptions.

\smallskip
\noindent  (Absolute) real cohomology classes on $M$ can be represented in terms of meromorphic 
(or anti-meromorphic) functions in $L^2_q(M)$ (see \cite{F02}, \S 2). In fact, by the Hodge theory on Riemann surfaces \cite{FK92}, III.2, all real cohomology classes can be represented as the real (or imaginary) part of a holomorphic (or anti-holomorphic) differential on $M$. In turn, any orientable
holomorphic quadratic differential induces an isomorphism between the space $\text{\rm Hol}^+(M)$ 
[$\text{\rm Hol}^-(M)$] of holomorphic [anti-holomorphic] differentials and the space of square-integrable meromorphic [anti-meromorphic] functions. Let ${\Cal M}^+_q$ [${\Cal M}^{-}_q$] be the space of meromorphic [anti-meromorphic] functions on $M$ which belong to the Hilbert space $L^2_q(M)$ (see Proposition \ref{prop:CR}). Such spaces can be characterized as the spaces of all meromorphic  [anti-meromorphic] functions with poles at $\Sigma_{q}=\{q=0\}$ of orders bounded in terms of the  multiplicity of the points $p\in\Sigma_{q}$ as zeroes of the quadratic differential $q$. In fact, if $p\in \Sigma_q$ is
a zero of $q$ of order $2m$, that $p$ is a pole of order at most $m$ for any $m^{\pm} \in {\Cal M}^\pm_q$.

\smallskip
\noindent Let $q^{1/2}$ be a holomorphic square root of $q$ on $M$. Holomorphic [anti-holomorphic]
differentials $h^+$ [$h^-$]  on $M$ can be written in terms of meromorphic [anti-meromorphic] functions
in $L^2_q(M)$ as follows:
\begin{equation}
\label{eq:holmer}
\begin{aligned}
h^+ &:=  m^+  q^{1/2} \,, \quad m^+ \in \Cal M^+_q \,; \\
h^- &:=  m^-  \bar q^{1/2} \,, \quad m^- \in \Cal M^-_q \,.
\end{aligned}
\end{equation}
\noindent The following \emph{representations of real cohomology classes} therefore hold:
\begin{equation}
\label{eq:reprs}
\begin{aligned}
c\in H^1(M,\R) &\Longleftrightarrow  c=[ \Re( m^+  q^{1/2})]  \,, \quad m^+ \in \Cal M^+_q \,; \\
c\in H^1(M,\R) &\Longleftrightarrow  c=[ \Re( m^- \bar q^{1/2})]  \,, \quad m^- \in \Cal M^-_q\,.
\end{aligned}
\end{equation}
The maps $c^\pm_q:{\Cal M}^\pm_q \to H^1(M,{\R})$ given by the representations \pref{eq:reprs} are bijective and it is in fact {\it isometric }if the spaces ${\Cal M}^\pm_q$ are endowed with the euclidean structure induced by $L^2_q(M)$ and $H^1(M,{\R})$ with the {\it Hodge product }relative to the complex structure of the Riemann surface $M$. In fact, the Hodge norm  $\Vert c \Vert^2_H$ of a cohomology class $c\in H^1(M,\R)$ is defined as follows: 
\begin{equation}
\label{eq:Hodgenorm}
\Vert c \Vert^2_H := \frac{\imath}{2}\,  \int_M  h^\pm \wedge  \overline {h^\pm}  \quad \text{ \rm if }\,\,
c=[\Re(h^\pm)]\,, \,\, h^\pm\in \text{\rm Hol}^\pm(M)\,.
\end{equation}
\noindent We remark that the Hodge norm is defined in terms of the complex structure of the Riemann
surface $M$ (carrying a holomorphic quadratic differential $q\in Q(M)$) but does not depend
on the quadratic differential.  If $q\in Q(M)$ is any orientable quadratic differential on $M$, by the representation \pref{eq:reprs}, we can also write:
\begin{equation}
\label{eq:Hodgenormbis}
\begin{aligned}
\Vert c^+_q(m^+) \Vert^2_H &:=  \int_M  \vert m^+ \vert ^2 \,\omega_q \, ,  \quad \text{ \rm for all }\,\,
m^+\in \Cal M^+_q\,; \\
\Vert c^-_q(m^-) \Vert^2_H &:=  \int_M  \vert m^- \vert ^2 \,\omega_q \, ,  \quad \text{ \rm for all }\,\,
 m^-\in \Cal M^-_q\,.
\end{aligned}
\end{equation}

\smallskip
\noindent The representation \pref{eq:holmer}-\pref{eq:reprs} can be extended to the punctured cohomology $H^1(M\setminus\Sigma_q,\R)$ as follows. For any finite set $\Sigma\subset M$, let $\text{\rm Hol}^+(M\setminus\Sigma)$ [$\text{\rm Hol}^-(M\setminus\Sigma)$] be the space of meromorphic [anti-meromorphic] differentials with at most simple poles at $\Sigma$. By Riemann surface theory, any 
real cohomology class $c\in H^1(M\setminus\Sigma,\R)$ can be represented as the real (or imaginary) part of a differential $h^+ \in \text{\rm Hol}^+(M\setminus\Sigma)$ or $h^- \in \text{\rm Hol}^-(M\setminus\Sigma)$. Let $\Sigma\subset M$ be a finite set and let $\Cal M^{+}(\Sigma)$ [$\Cal M^{+}(\Sigma)$] be the space of all meromorphic [anti-meromorphic] functions which are holomorphic [anti-holomorphic] on $M\setminus\Sigma$. The spaces $\Cal M^{\pm}(\Sigma)$ can be identified with a subspace of the distributional space $\Cal D^2(M\setminus\Sigma)$. In fact, if $q$ is any orientable holomorphic quadratic differential on $M$, the spaces $\Cal M^{\pm}(\Sigma_q)$ identify with subspaces of the space $\Cal S^2_q(M)$ of $q$-tempered distributions. The distribution determined by a function $\phi \in \Cal M^+(\Sigma_q)$ or $\Cal M^-(\Sigma_q)$ is defined by integration (in the standard way) as a linear functional on $C_0^{\infty}(M \setminus\Sigma_q)$, which can be extended to the space $\Omega^0_q(M)$ as follows:
$$
\phi(v) := \text{\rm PV }\int_M  \phi \,v\,\, \omega_q  \,,\quad \text{ for all }  v\in \Omega^0_q(M)\,.
$$
The Sobolev regularity of a distributions $\phi\in \Cal M^{\pm}(\Sigma_q)$ depend on the order
of its poles. In fact, by Theorem \ref{thm:ijreg} we have the following:

\begin{lemma} 
\label{lemma:meromreg}
Let $\phi \in \Cal M^{+}(\Sigma_q)$ [ $\phi\in \Cal M^{+}(\Sigma_q)$ ] be a meromorphic
[anti-meromorphic] function with poles at $\Sigma_q$. For any $s\in \R$, the associated distribution $\phi\in H^{-s}_q(M)$ if at every $p\in\Sigma_q$ of order $2m_p$ the function  $\phi$ has a pole of 
order $<(m_p+1)(s+1)$.
\end{lemma}

\noindent We introduce the following notation: for all $s>0$,
\begin{equation}
 \Cal M^{\pm}_s(\Sigma_q):= \Cal M^{\pm}(\Sigma_q)\cap H^{-s}_q(M) \,.
\end{equation}

\smallskip
\noindent There exist natural maps $\phi_q^\pm : \text{\rm Hol}^\pm(M\setminus\Sigma_q) \to \Cal M^\pm
(\Sigma_q)$ defined as follows: for all $h^\pm  \in \text{\rm Hol}^\pm(M\setminus\Sigma_q)$, 
$$
\phi_q^+( h^+) = h^+/q^{1/2}\quad [ \phi_q^-( h^-) = h^-/{\bar q}^{1/2} ]\,.$$
By Lemma \ref{lemma:meromreg} the range of the maps $\phi_q^\pm$ is contained in the weighted
Sobolev space $H^{-s}_q(M)$ for all $s>0$, hence there are well-defined maps
\begin{equation}
\label{eq:punctcohomaps}
\phi_q^\pm : \text{\rm Hol}^\pm(M\setminus\Sigma_q) \to \Cal M^\pm
_s(\Sigma_q)\, \quad \text{ \rm for all } s>0\,.
\end{equation}
The maps \pref{eq:punctcohomaps} are clearly injective and by Corollary 
\ref{cor:ijregbis} there exists $s_q >0$ such that, for any $s\in (0, s_q)$, they are
also surjective. Let
\begin{equation}
\Cal M_{*}^{\pm}(\Sigma_q)= \bigcap_{s>0} \,\Cal M_{s}^{\pm}(\Sigma_q) =
\Cal M_{s}^{\pm}(\Sigma_q) \,, \quad \text{ for any } s\in (0,s_q)\,.
\end{equation}
\noindent The representation \pref{eq:reprs} of the absolute real cohomology generalizes to the punctured real cohomology as follows. 
\begin{equation}
\label{eq:punctreprs}
\begin{aligned}
c\in H^1(M\setminus\Sigma_q,\R) &\Longleftrightarrow  c=[ \Re( m^+  q^{1/2})]  \,, \quad 
m^+ \in \Cal M^+_*(\Sigma_q) \,; \\
c\in H^1(M\setminus\Sigma_q,\R) &\Longleftrightarrow  c=[ \Re( m^- \bar q^{1/2})]  \,, \quad 
m^- \in \Cal M^-_*(\Sigma_q)\,.
\end{aligned}
\end{equation}
\noindent The following lemma, proved in \cite{F02}, Lemma 7.6, for weighted Sobolev spaces
with integer exponent, holds:
\begin{lemma} 
\label{lemma:dUhC}
Let $s\in \R^+$. Let $C\in W_q^{-s}(M)$ be any real current of dimension (and degree) equal to
$1$, closed in the space $\Cal D(M\setminus\Sigma_q)$ of  currents on $M\setminus\Sigma_q$. 
There exists a distribution $U\in H^{-s+1}_q(M)$ and a meromorphic differential $h^+ \in \text{\rm Hol}^+(M\setminus\Sigma_q)$ such that 
\begin{equation}
\label{eq:dUhC}
dU^\ast \,=\, \Re ( h^+) \,-\, C \quad \text {\rm in }\, W^{-s}_q(M)\,.
\end{equation}
If $C$ is closed in the space $\Cal S_q(M)$ of $q$-tempered currents there exists a distribution 
$U\in H^{-s+1}_q(M)$ and a holomorphic differential $h^+ \in \text{\rm Hol}^+(M)$ such
that the identity \pref{eq:dUhC} holds. 
\end{lemma}
\noindent The argument given in \cite{F02}, Lemma 7.6, in the case of integer order $s\in \N$ extends the general case of order $s\in \R^+$. In fact, it follows from the distributional identity \pref{eq:dUhC} in $\Cal S_q(M)$ that the current $U^\ast\in H^{-s+1}_q(M)$ if and only if the current $C\in W^{-s}_q(M)$, for any $s\in \R^+$. Hence, Lemma \ref{lemma:dUhC} follows immediately from \cite{F02}, Lemma 7.6.

\smallskip
\noindent The construction of basic currents (or, equivalently, of invariant distributions) is based on
the following method.

\begin{lemma} 
\label{lemma:bcconstr}
Let $q$ be an orientable holomorphic quadratic differential on a Riemann surface $M$. Let $m^+\in
 \Cal M^+_s(\Sigma_q)$ be a meromorphic function with poles at $\Sigma_q\subset M$. A distribution $U\in H^{-s+1}_q(M)$ is a (distributional) solution in $\Cal D(M\setminus\Sigma_q)$ of the cohomological equation 
\begin{equation}
\label{eq:bccohomeq}
SU = \Re(m^+) \,\,  [TU = - \Im(m^+) ]  \quad \text{\rm in }\,\,
 \Cal D(M\setminus\Sigma_q)\,,
\end{equation}
if and only if the current $C\in W^{-s}_q(M)$ uniquely determined by the identity
\begin{equation}
\label{eq:dUmC}
dU^\ast = \Re ( m^+q^{1/2} ) \, + \,  C 
\end{equation}
is horizontally [vertically] quasi-basic. If $\,\Re ( m^+q^{1/2} ) \in \text{\rm Hol}^+(M)$, the distribution $U\in H^{-s+1}_q(M)$ is a solution of  the cohomological equation \pref{eq:bccohomeq} in the space $\Cal S_q(M)$ of $q$-tempered currents  if and only if the current $C\in W^{-s}_q(M)$ uniquely determined by formula \pref{eq:dUmC} is horizontally [vertically] basic.
\end{lemma}
\begin{proof}
If formula \pref{eq:dUmC} holds in $\Cal D(M\setminus\Sigma_q)$, then $C$ is closed in $\Cal D(M\setminus\Sigma_q)$. If the differential $\Re ( m^+q^{1/2} )$ is holomorphic and 
formula \pref{eq:dUmC} holds in $\Cal S_q(M)$, then $C$ is closed in $\Cal S_q(M)$.  
The standard formula for the Lie derivative of a current,
\begin{equation}
\Cal L_{X}C = \imath_{X}dC \,+ \, d \imath_{X}C\,=\, 0\,,
\end{equation}
holds in $\Cal D(M\setminus\Sigma_q)$ for any vector field $X$ with compact support contained 
in $M\setminus \Sigma_q$ and it holds in $\Cal S_q(M)$ for any vector field $X\in \Cal V_q(M)$
It follows that a current $C\in \Cal D(M\setminus\Sigma_q)$ is horizontally [vertically] quasi-basic if and only if it is closed and $\imath_SC=0$ [$\imath_TC=0$] in $\Cal D(M\setminus\Sigma_q)$ and it is
horizontally [vertically] basic if and only if $C\in \Cal S_q(M)$ is closed and $\imath_S C=0$ [$\imath_TC=0$] in $\Cal S_q(M)$. The distribution $U\in H^{-s}_q(M)$ in formula \pref{eq:dUmC} is a solution of the cohomological equation \pref{eq:bccohomeq} in $\Cal D(M\setminus\Sigma_q)$ or $\Cal S_q(M)$ if and only if $\imath_S C=0$ [$\imath_TC=0$] in $\Cal D(M\setminus\Sigma_q)$ or $\Cal S_q(M)$ respectively.  As a consequence, the lemma is proved.
\end{proof}

\smallskip
\noindent Let $q$ be an orientable holomorphic quadratic differential on a Riemann surface $M$
(of genus $g\geq 1$).  Let $\Pi_{\pm q}(M\setminus\Sigma_q,\R) \subset H^1(M\setminus\Sigma_q,\R)$ be the codimension $1$ subspaces
defined as follows:
\begin{equation}
\label{eq:Pipunct}
\begin{aligned}
 \Pi^1_{q}(M\setminus\Sigma_q,\R) &:= \{ c\in H^1(M\setminus\Sigma_q,\R)\,\vert\,  
 c\wedge [ \Im (q^{1/2})] =0 \}\,; \\
 \Pi^1_{-q}(M\setminus\Sigma_q,\R) &:= \{ c\in H^1(M\setminus\Sigma_q,\R)\,\vert\,  
 c\wedge [ \Re (q^{1/2})] =0 \}\,.
\end{aligned}
\end{equation}
Since the absolute cohomology can be regarded as a subspace of the punctured
cohomology it also is possible to define the subspaces
\begin{equation}
\label{eq:Piabs}
\begin{aligned}
 \Pi^1_{q}(M,\R) &:=  \Pi^1_{q}(M\setminus\Sigma_q,\R) \cap H^1(M,\R)   \,; \\
 \Pi^1_{-q}(M,\R) &:= \Pi^1_{-q}(M\setminus\Sigma_q,\R) \cap H^1(M,\R) \,.
\end{aligned}
\end{equation}

\begin{theorem}
\label{thm:bc}
For any $s>3$ there exists a full measure set $\Cal F_s \subset S^1$ such that the following holds. For
any $\theta\in \Cal F_s$, the following inclusions hold
\begin{equation}
\label{eq:cohomincl}
\begin{aligned}
\Pi^1_{\pm q_\theta}(M\setminus\Sigma_q,\R)  \,\, &\subset \,\,   H^{1,s}_{\pm q_\theta} (M\setminus\Sigma_q,\R) \,, \\
 \Pi^1_{\pm q_\theta}(M,\R)    \,\, &\subset \,\, H^{1,s}_{\pm q_\theta}(M,\R)\,.
\end{aligned}
\end{equation}
\end{theorem}
\begin{proof}
Let $m^+\in \Cal M_{-1}^+(\Sigma_q)$ be any meromorphic function such that the induced 
distribution $\text{\rm PV}(m^+) \in H^{-1}_q(M)$. A computation shows that, for all $\theta\in S^1$, 
\begin{equation}
\label{eq:zeroav}
\text{\rm PV } \int_M \Re( m^+)\,\omega_q =0 \Longleftrightarrow  [\Re(m^+q_\theta^{1/2})] \in \Pi^1_{q_\theta}(M\setminus\Sigma_q,\R)\,.
\end{equation}
Under the zero-average condition \pref{eq:zeroav}, by Theorem \ref{thm:distsol} for any $r>2$
there exists a full measure set $\Cal F_r (m^+) \subset S^1$ such that, for all $\theta\in S^1$, the cohomological equation $S_\theta U= \Re(m^+)$ has a distributional solution  $U_\theta(m^+) \in \bar H^{-r}_q(M)$. Let $U^\ast_\theta(m^+)$ be the current of dimension $2$ corresponding to the distribution $U_\theta(m^+)$. Let $C_\theta(m^+) \in W_q^{-r-1}(M)$ be the $1$-dimensional current determined by the identity
\begin{equation}
\label{eq:dUast}
dU^\ast_\theta(m^+) \,=\, \Re(m^+q_\theta^{1/2}) \,-\, C_\theta(m^+)\,.
\end{equation}
By Lemma \ref{lemma:bcconstr}, we have thus proved that, for all meromorphic functions
$m^+\in \Cal M_{-1}^+(\Sigma_q)$ and for all $\theta \in \Cal F_r(m^+)$, there exists a quasi-basic current $C_\theta(m^+)\in\Cal B^{r+1}_{q_\theta}(M\setminus \Sigma_q)$ such that 
$$
[C_\theta(m^+)] = [\Re(m^+q_\theta^{1/2})] \in \Pi^1_{q_\theta}(M\setminus\Sigma_q,\R)\,;
$$ 
in addition, whenever $m^+\in \Cal M^+_q$, the current $C_\theta(m^+)\in \Cal B^{r+1}_{q_\theta}(M)$
is basic and has a cohomology class
$$
[C_\theta(m^+)] = [\Re(m^+q_\theta^{1/2})] \in \Pi^1_{q_\theta}(M,\R)\,.
$$
\smallskip
\noindent  Let $\sigma\in \N$ be the cardinality of the set $\Sigma_q\subset M$ and let 
$$
\{m^+_1, \dots , m^+_{2g-1}, \dots, m^+_{2g+\sigma-1} \}
$$ 
be a basis (over $\R$) of the real subspace of $\Cal M^+_*(\Sigma_q)$ defined by the zero-average condition \pref{eq:zeroav}. For any $s>3$, let
$$
\Cal F^+_s := \bigcap_{i=1}^{2g+\sigma-1} \Cal F_{s-1}(m^+_i)\,.
$$
Clearly he set $\Cal F^+_s$ has full Lebesgue measure. We claim that for all $\theta\in \Cal F^+_s$ 
the following inclusions hold:
\begin{equation}
\label{eq:cohominclplus}
\begin{aligned}
\Pi^1_{ q_\theta}(M\setminus\Sigma_q,\R) \,\,& \subset \,\, 
H^{1,s}_{q_\theta} (M\setminus\Sigma_q,\R)  \,, \\
\Pi^1_{q_\theta}(M,\R) \,\, &\subset\,\, H^{1,s}_{q_\theta}(M\setminus\R)  \,.
\end{aligned}
\end{equation}
The claim is proved as follows. For any $c\in H^1(M\setminus\Sigma_q,\R)$ there 
exists a unique meromorphic function $m^+\in \Cal M^+_*(\Sigma_q)$ such that $c=[\Re(m^+q_\theta^{1/2})]$. The function $m^+\in \Cal M^+_q$ for any $c\in H^1(M,\R)$. If $c\in \Pi^1_{q_\theta}(M,\R)$, the distribution $\Re(m^+)$ vanishes on constant functions
as in \pref{eq:zeroav}, hence for all $\theta \in \Cal F^+_s$, there exists a solution $U\in H^{-s+1}_q(M)$ 
of the cohomological equation $S_\theta U = \Re(m^+)$. The current $C \in \Cal B^s_{q_\theta}(M\setminus\Sigma_q)$ such that $[C] =c\in  \Pi^1_{q_\theta}(M\setminus\Sigma_q,\R)$ 
is then given by the identity \pref{eq:dUast}. By the above discussion the current $C\in \Cal B^s_{q_\theta}(M)$  for all $c\in \Pi^1_{q_\theta}(M,\R)$.

\smallskip
\noindent By a similar  argument  it is possible to construct a full measure set $\Cal F^-_s$ such that,
for all $\theta\in \Cal F^-_s$, the following inclusions hold:
\begin{equation}
\label{eq:cohominclminus}
\begin{aligned}
\Pi^1_{- q_\theta}(M\setminus\Sigma_q,\R) \,\,& \subset \,\, H^{1,s}_{-q_\theta} 
(M\setminus\Sigma_q,\R)  \,, \\
\Pi^1_{-q_\theta}(M,\R) \,\, &\subset\,\, H^{1,s}_{-q_\theta}(M,\R)  \,.
\end{aligned}
\end{equation}
Thus the set $\Cal F_s:= \Cal F^+_s \cap \Cal F^-_s$ has the required properties
since it has full measure and the inclusions \pref{eq:cohomincl} hold.
\end{proof}

\noindent By Lemma \ref{lemma:DtoC} and Theorem \ref{thm:bc} the following holds:

\begin{corollary} 
For any $s>3$ there exists a full measure set $\Cal F_s \subset S^1$ such that, for all
$\theta\in \Cal F_s $, the spaces $\Cal I_{\pm q_\theta}^s(M) \subset H^{-s}_q(M)$ of horizontally
or vertically quasi-invariant distributions have dimension at least $2g+\sigma-1$ and the spaces 
$\Cal I_{\pm q_\theta}^s(M)  \subset H^{-s}_q(M)$ of horizontally or vertically invariant distributions have dimension at least $2g-1$.
\end{corollary}

\begin{corollary} 
\label{cor:basiccohom}
For any $s>3$ and for almost all $\theta\in S^1$, 
\begin{equation}
\label{eq:bcone}
H^{1,s}_{\pm q_\theta} (M,\R)=\Pi^1_{\pm q_\theta}(M,\R) \,.
\end{equation}
For any $s>4$ and for almost all $\theta\in S^1$,
\begin{equation}
\label{eq:bctwo}
H^{1,s}_{\pm q_\theta} (M\setminus\Sigma_q,\R)=H^1(M\setminus\Sigma_q,\R) \,.
\end{equation}
\end{corollary}
\begin{proof} The inclusions $H^{1,s}_{\pm q} (M,\R)\subset \Pi^1_{\pm q}(M,\R)$ hold
for any orientable quadratic differential $q$ on $M$ and for any $s>0$. In fact,
\begin{equation}
\label{eq:wedgecon}
\begin{aligned}
&[C\wedge  \Im (q^{1/2})](1) = \,\,\,\,\imath_SC(\omega_q) =0 \,, \quad \text{ if } \,C\in \Cal B_q(M)\,;\\
&[C\wedge  \Re (q^{1/2})](1) = -\imath_TC(\omega_q) =0 \,, \quad \text{ if } \,
C\in \Cal B_{-q}(M)\,.
\end{aligned}
\end{equation}
Thus, the identity \pref{eq:bcone} follows immediately from Theorem \ref{thm:bc}.

\smallskip
\noindent By Theorem \ref{thm:bc}, in order to prove the identity \pref{eq:bctwo} it is enough to
prove that for almost all $\theta\in S^1$ the cohomology class $[\Re(q_\theta)] \in \Cal 
B_{q_\theta}(M\setminus \Sigma_q)$ and the cohomology class $[\Im(q_\theta)] \in \Cal 
B_{-q_\theta}(M\setminus \Sigma_q)$.  By Lemma \ref{lemma:bcconstr}  the argument is therefore
reduced to the construction, for any $s>3$ and for almost all $\theta\in S^1$, of a solution $U\in 
H^{-s}_q(M)$ of the cohomological equation $S_\theta U=1$ [$T_\theta U=1$] in $\Cal D(M\setminus \Sigma_q)$. Such a construction can be carried out as follows. Let $\delta_p$ be the Dirac mass at any point $p\in \Sigma_q$. The distribution $F=1-\delta_p \in H^{-s}(M) \subset H^{-s}_q(M)$ for any $s>1$.
We claim that for any $s>3$ and for almost all $\theta\in S^1$ there exists a distributional solution 
$U\in H^{-s}_q(M)$ of the cohomological equation $S_\theta U = F$ [$T_\theta U=F$] in 
$H^{-s-1}_q(M)$. It follows that $U$ is a solution of the cohomological equation $S_\theta U=1$ [$T\theta U=1$] in $\Cal D(M\setminus \Sigma_q)$. The above claim is proved as follows. By \cite{F97}, Prop. 4.6, or \cite{F02}, Lemma 7.3, since $F\in H^{-2}_q(M)$ vanishes on constant functions, there exists a distribution $f \in H_q^{-1}(M)$ such that $\partial_q^+ f= F$ (as well as a distribution $f' \in 
H^{-1}_q(M)$ such that $\partial_q^-f'=F$).  By Theorem \ref{thm:distsol}, for almost all $\theta\in S^1$ and for all $s>2$, there exists a solution $u \in H^{-s}_q(M)$ of the equation $S_\theta u = f$ [$ T_\theta u = f$ ]. The distribution $U:= \partial_q ^+u \in H^{-s}_q(M)$ for any $s>3$ and solves the cohomological equation $S_\theta U=F$ [$T_\theta U=F$] in $H^{-s-1}_q(M)$.
\end{proof}

\noindent The structure of the space of basic currents with vanishing cohomology class with respect to the filtration induced by weighted Sobolev spaces with integer exponent was described 
in \cite{F02}, \S 7. We extend below the results of \cite{F02} to fractional weighted Sobolev spaces.

\smallskip
\noindent Let $\delta_{\pm q}: \Cal B_{\pm q}(M\setminus\Sigma_q) \to \Cal B_{\pm q}(M\setminus
\Sigma_q)$ be the linear maps defined as follows (see \cite{F02}, formula $(7.18')$): 
\begin{equation}
\label{eq:deltas}
\begin{aligned}
\delta_{q}(C) &:=  -\,d \left(C\wedge \Re(q^{1/2}) \right)^\star \,, \quad &\text{\rm for }\,\,C\in  
\,\Cal B_{q}(M\setminus\Sigma_q)   \,; \\
\delta_{-q}(C) &:=  \,\,\,\,\,d \left(C\wedge \Im(q^{1/2}) \right)^\star  \,,  \quad  &\text{\rm for }\,\,C\in  
\,\Cal B_{-q}(M\setminus\Sigma_q) \,.
\end{aligned}
\end{equation}
It can be proved by Lemma \ref{lemma:DtoC} and by the definition of the weighted Sobolev spaces $H^s_q(M)$ and $W^s_q(M)$ that the above formulas \pref{eq:deltas} define, for all $s\in \R^+$,
bounded linear maps 
\begin{equation}
\label{eq:deltasbis}
\begin{aligned}
\delta^s_{\pm q}: {\Cal B}_{\pm q}^{s-1}(M\setminus\Sigma_q)\,\, &\to\,\, {\Cal B}_{\pm q}^{s}(M\setminus\Sigma_q) \,; \\
\delta^s_{\pm q}: {\Cal B}_{\pm q}^{s-1}(M) \,\,&\to \,\,{\Cal B}_{\pm q}^{s}(M)\,.
\end{aligned}
\end{equation}
We remark that a similar statement is false in general for the Friedrichs Sobolev spaces of currents. 
The following result extends Theorem 7.7 of \cite{F02} to fractional weighted Sobolev spaces.  

 \begin{theorem} \label{thm:bcstruct}
For all $s\in \R^+$ there exist exact sequences
\begin{equation}
\label{eq:exactsequences}
\begin{aligned}
0\to {\R} \to {\Cal B}_{\pm q}^{s-1}(M\setminus\Sigma_q) &^{\underrightarrow{\,\,\,\,\delta^s_{\pm q}\,\,\,\,}} 
{\Cal B}_{\pm q}^s(M\setminus\Sigma_q) ^{\underrightarrow{\,\,\,\,j_{q}\,\,\,\,}}  H^1(M\setminus\Sigma_q,{\R})\,\,;\\
0\to {\R} \to {\Cal B}_{\pm q}^{s-1}(M) &^{\underrightarrow{\,\,\,\,\delta^s_{\pm q}\,\,\,\,}} 
{\Cal B}_{\pm q}^s(M) ^{\underrightarrow{\,\,\,\,j_{q}\,\,\,\,}}  H^1(M,{\R})\,\,.
\end{aligned}
\end{equation} 
  \end{theorem}
  \begin{proof} 
  The maps $i_{\pm q}: {\R} \to {\Cal B}_{\pm q}^{s-1}(M) \subset {\Cal B}_{\pm q}^{s-1}(M \setminus\Sigma_q)$ defined  as
  $$
  i_q (\tau) = \tau \eta_S   \quad \text { \rm and } \quad  i_{-q} (\tau) = \tau \eta_T \,, \qquad \text{ \rm for all }\,\tau\in \R\,,
  $$
 are clearly injective and the kernels $\text{\rm Ker}(\delta^s_{\pm q}) =i_{\pm q}(\R)$, 
 for all $s\in \R^+$. In fact, if a current $C\in \text{\rm Ker}(\delta^s_{+ q})$, then the distribution $\left(C\wedge \Re(q^{1/2})\right)^\star \in \R$, hence $C\wedge \Re
 (q^{1/2}) \in \R\cdot  \omega_q$ and $\,C \in \R\cdot \Re(q^{1/2}) + \R\cdot 
 \Im(q^{1/2})$. It follows that $C\in \R\cdot \Im(q^{1/2})$ since  $C\in \Cal B_q(M\setminus\Sigma_q)$. Similarly, if the current $C\in \text{\rm Ker}(\delta^s_{-q})$, it
 follows that $C\wedge \Im(q^{1/2}) \in \R\cdot \omega_q$ and $C\in \R\cdot  
 \Re(q^{1/2})$ since  $C\in  \Cal B_{-q}(M\setminus\Sigma_q)$. This proves the inclusions $\text{\rm Ker}(\delta^s_{\pm q})  \subset i_{\pm q}(\R)$. The opposite inclusions are immediate. 
 
 \smallskip
 \noindent By Lemma \ref{lemma:dUhC} a current $C \in {\Cal B}_{\pm q}^s(M\setminus\Sigma_q)$ has zero cohomology class, that is, it is in the kernel of the cohomology map, if and only if there exists a distribution $U_C \in H^{-s+1}_q(M)$ such that $dU^\ast_C= C$ in $\Cal D(M\setminus\Sigma_q)$ and this identity holds in $\Cal S_q(M)$, hence in $W^{-s}_q(M)$, if $C \in {\Cal B}_{\pm q}^s(M)$. It is immediate to verify that $U_C \in \Cal I_{\pm q}^{s-1}(M\setminus \Sigma_q)$ if and only  if 
 $C \in {\Cal B}_{\pm q}^s(M\setminus\Sigma_q)$ and that $U_C \in \Cal I_{\pm q}^{s-1}(M)$ if and  only if $C \in {\Cal B}_{\pm q}^s(M)$. By Lemma~\ref{lemma:DtoC} we have thus proved that the map $j_q:\Cal B^{s}_{\pm q}(M\setminus\Sigma_q)\to H^1(M\setminus\Sigma_q,\R)$ has kernel equal to the range $\delta^s_{\pm q}\left( {\Cal B}_{\pm q}^{s-1}(M) \right)$ and that  the map $j_q:\Cal B^{s}_{\pm q}(M\setminus)\to  H^1(M,\R)$ has kernel equal to the range $\delta^s_{\pm q}\left( {\Cal B}_{\pm q}^{s-1}(M) \right)$. 
\end{proof}

\noindent Finer results on invariant distributions and on smooth solutions of the cohomological
equation for directional flows can be obtained by combining the results of this section with the
renormalization method based on the Teichm\"uller flow and related cocycles, such as the 
Kontsevich-Zorich cocycle. Our goal is to improve upon the results of Marmi, Moussa and Yoccoz 
\cite{MMY05} who have studied the cohomological equation for interval exchange transformations  
solely by methods based on the renormalization dynamics (the Rauzy-Veech-Zorich induction).

\section{Cocycles over the Teichm\"uller flow}

\subsection{ The Kontsevich-Zorich cocycle. }
\label{KZcocycle}

\noindent The Kontsevich-Zorich cocycle is a multiplicative cocycle over the Teichm\"uller
geodesic flow on the moduli space of (orientable) holomorphic quadratic differentials on
compact Riemann surfaces. This cocycle appears in the study of the dynamics of interval 
exchange transformations and of (translation) flows on surfaces  for which it represents a renormalization dynamics and of the Teichm\"uller flow itself. In fact, the study tangent cocycle 
of the Teichm\"uller flow can be reduced to that of the Kontsevich-Zorich cocycle.

\smallskip
\noindent Let $T_g$ and  $Q_g$ be respectively the {\it Teichm\"uller spaces} of complex (conformal) structures  and of holomorphic quadratic differentials on a surface of genus $g\geq 1$. We recall that
 the spaces $T_g$ and $Q_g$ can be described as follows. Let $\text{Diff}_0^+(M)$ is the group of orientation preserving diffeomorphisms of the surface $M$ which are isotopic to the identity (equivalently, it is the connected component of the identity in the Lie group of all orientation preserving diffeomorphisms of $M$).  By definition
\begin{equation}
\begin{aligned}
T_g & := \{\text{ complex (conformal) structures }\}/ \text{Diff}_0^+(M) \,\,,\\ 
Q_g & :=\{\text{ holomorphic quadratic differentials }\}/ \text{Diff}_0^+(M)\,\,.
\end{aligned}
\end{equation}
A theorem of L. Ahlfors, L. Bers and S. Wolpert states that $T_g$ has a complex structure holomorphically equivalent to that of a Stein (strongly pseudo-convex) domain in ${\C}^{3g-3}$ \cite{Ber74}, \S 6, or  \cite{Na88}, Chap. 3, 4 and Appendix \S 6. The space $Q_g$ of holomorphic quadratic differentials is a complex vector bundle over $T_g$ which can be identified to the cotangent bundle of $T_g$. Let $\Gamma_g := \hbox{Diff}^+(M)/\hbox{Diff}_0^+(M)$ be the {\it mapping class group }and let $R_g$, ${\Cal M}_g$ be respectively the {\it moduli spaces }of complex (conformal) structures and of holomorphic quadratic differentials on a surface of genus  $g\geq 1$. The spaces $R_g$ and ${\Cal M}_g$ can be described as the quotient spaces: 
\begin{equation}
R_g := T_g/\Gamma_g\,\,,\quad 
 {\Cal M}_g := Q_g/\Gamma_g\,\,,
\end{equation}
In case $g=1$, the Teichm\"uller space $T_{1}$ of elliptic curves (complex structures on $T^2$) is isomorphic to the upper half plane ${\C}^+$ and the Teichm\"uller space $Q_1$ of holomorphic quadratic differentials on elliptic curves is a complex line bundle over $T_1$  (see \cite{Na88}, Ex. 2.1.8). The mapping class group can be identified with the lattice $SL(2,\Z)$ which acts on the
upper half plane ${\C}^+$ in the standard way. The moduli space $R_1:={\C}^+/ SL(2,\Z)$ is a non-compact finite volume surface with constant negative curvature, called the {\it modular surface}. The moduli space ${\Cal M}_1$ can be identified to the cotangent bundle of the modular surface.

\smallskip
\noindent The {\it Teichm\"uller (geodesic) flow}  on $\Cal M_g$ can be defined as the geodesic flow for a natural metric on $R_g$ called the {\it Teichm\"uller metric}. Such a metric measures the amount of {\it quasi-conformal distorsion }between two different (equivalent classes of) complex structures in $R_g$. In the higher genus case, the Teichm\"uller metric is not Riemannian, but only {\it Finsler }(that is, the norm on each tangent space does not come from an euclidean product) and, as H. Masur proved, does not have negative curvature in any reasonable sense \cite{Ber74}, \S 3 (E). If $g=1$, the Teichm\"uller metric coincides with the Poincar\'e metric on the modular surface $R_1$ \cite{Na88}, 2.6.5, in particular it is Riemannian with constant negative curvature. 
 
\smallskip
\noindent There is a natural action of the Lie group $GL(2,{\R})$ on $Q_g$ (see also \cite{HSHB},
 \S 1.4 or \cite{MaHB}, \S 3) which is defined as follows. The map 
 $$
 q \to (\Cal F_q, \Cal F_{-q}) \,, \quad q\in Q_g\,,
 $$
 is a bijection between the space $Q_g$ and the space of all pairs of transverse measured foliations.
 We recall that transversality for a pair $(\Cal F^, \Cal F^\perp)$ of measured foliations is taken in the sense that ${\Cal F}$ and ${\Cal F}^{\perp}$ have a common set $\Sigma$ of (saddle) singularities, 
 have the same index at each singularity and are transverse in the standard sense on $M\setminus
\Sigma$.  The set $\Sigma$ of common singularities coincides with the set $\Sigma_{q}$ of zeroes 
of the holomorphic quadratic differential $q\equiv ({\Cal F},{\Cal F}^{\perp})$. Any pair of transverse measured foliations is determined locally by a pair $(\eta,\eta^\perp)$ of (locally defined) transverse real-valued closed $1$-forms. The group $GL(2,{\R})$ acts naturally by left multiplication on the space 
of (locally defined) pairs of transverse real-valued closed $1$-forms, hence it acts on the space of all pairs of  transverse measured foliations an on the space of $Q_g$ of holomorphic quadratic differentials. Such an action is equivariant with respect to the action of the mapping class group, hence it passes to the quotient ${\Cal M}_g$. 
\smallskip
\noindent The Teichm\"uller flow $\{G_t\}_{t\in \R}$ is given by the action of the diagonal subgroup $\{\hbox{diag}\,(e^{-t},e^t)\}_{t\in \R}$ on $Q_g$ (on ${\Cal M}_g$). In other terms, if we identify holomorphic quadratic differentials with pairs of transverse measured foliations as explained above, 
we have:
\begin{equation}
\label{eq:Teichflow}
G_t ({\Cal F}_q,{\Cal F}_{-q}) := (e^{-t} {\Cal F}_q, e^t {\Cal F}_{-q})\,\,. 
\end{equation}
\noindent In geometric terms, the action of the Teichm\"uller flow on quadratic differentials 
induces a one-parameter family of deformations of the conformal structure which consist in 
contracting along vertical leaves (with respect to the horizontal length) and expanding along  
horizontal leaves (with respect to the vertical length)  by reciprocal (exponential) factors. The 
reader can compare the definition in terms of the $SL(2,{\R})$-action with the analogous 
description of the geodesic flow on a surface of constant negative curvature (such as the modular surface). In fact, in case $g=1$ the above definition reduces to the standard Lie group presentation of the geodesic flow on the modular surface: the unit sub-bundle ${\Cal M}_1^{(1)} \subset {\Cal M}_1$ 
of all holomorphic quadratic differentials of unit total area on elliptic curves can be identified with 
the homogeneous space $SL(2,{\R})/SL(2,\Z)$ and the geodesic flow on the modular surface is then identified with the action of the diagonal subgroup of $SL(2,{\R})$ on $SL(2,{\R})/SL(2,\Z)$.

\medskip
\noindent We list below, following \cite{Ve90}, \cite{Ko97}, the main structures carried by the moduli space ${\Cal M}_g$ of quadratic differentials:

\begin{enumerate}
\item the moduli space ${\Cal M}_g$ is a stratified analytic orbifold; each stratum ${\Cal M}_{\kappa} \subset {\Cal M}_g$ (determined by the multiplicities $\kappa=(k_1,\dots,k_{\sigma})$ 
of the zeroes $\{p_1,\dots,p_{\sigma}\}$ of quadratic differentials) is $GL(2,{\R})$-invariant, hence in particular $G_t$-invariant;

\item The total area function $A:{\Cal M}_g\to {\R}^+$,
\begin{equation*}
A(q):=\int_M |q|\,\,,
\end{equation*}
is $SL(2,{\R})$-invariant; hence the {\it unit bundle }${\Cal M}^{(1)}_g:=A^{-1}(\{1\})$ and its strata
${\Cal M}^{(1)}_{\kappa}:= {\Cal M}_{\kappa}\cap {\Cal M}^{(1)}_g$ are $SL(2,{\R})$-invariant and, in particular, $G_t$-invariant.

\smallskip
\noindent Let ${\Cal M}_{\kappa}$ be a stratum of {\it orientable }quadratic differentials, that is, 
quadratic differentials which are squares of holomorphic $1$-forms. In this case, the natural
numbers $(k_1,\dots,k_{\sigma})$ are all even.

\item  The stratum of squares ${\Cal M}_{\kappa}$ has a locally affine structure modeled 
on the affine space $H^1(M,\Sigma_{\kappa};{\C})$, with $\Sigma_{\kappa}:=\{p_1,\dots,p_{\sigma}\}$.
Local charts are given by the period map $q\to [q^{1/2}]\in H^1(M,\Sigma_{\kappa};{\C})$. 

\item The Lebesgue measure on the euclidean space $H^1(M,\Sigma_{\kappa};{\C})$, 
normalized so that the quotient torus 
\begin{equation*}
H^1(M,\Sigma_{\kappa};{\C})\,/\,H^1
(M,\Sigma_{\kappa};{\Z}\oplus i {\Z})
\end{equation*}
has volume $1$, induces an absolutely continuous $SL(2,{\R})$-invariant measure 
$\mu_{\kappa}$ on ${\Cal M}_{\kappa}$. The conditional measure $\mu^{(1)}_{\kappa}$ induced 
on the stratum ${\Cal M}^{(1)}_{\kappa}$  is $SL(2,{\R})$-invariant, hence $G_t$-invariant.
\end{enumerate}
\noindent All the above structures (the stratification, the area function, the locally affine structure on
the strata of squares)  lift to corresponding structures at the level of the Teichm\"uller space of
quadratic differentials, equivariant under the action of the mapping class group.

\noindent It was discovered by W. Veech \cite{Ve90} that ${\Cal M}^{(1)}_{\kappa}$ has in general several connected components. The connected components for the strata of abelian differentials
(or equivalently of orientable quadratic differentials) have been classified completely by M. Kontsevich and A. Zorich \cite{KZ03}. A similar classification for the case of strata of non-orientable holomorphic quadratic differentials has been recently obtained by E. Lanneau in his thesis \cite{La05}.  Taking this phenomenon into account, the following result holds:

\begin{theorem} 
\label{thm:erg}
\cite{Ma82}, \cite{Ve86}  The total volume of the measure $\mu_{\kappa}^{(1)}$ on 
${\Cal M}^{(1)}_{\kappa}$ is finite and the Teichm\"uller geodesic flow $\{G_t\}_{t\in \R}$ is ergodic on 
each connected component of ${\Cal M}^{(1)}_{\kappa}$. 
\end{theorem}

\smallskip
\noindent  We will describe below several results about the Lyapunov structure of various coycles 
over the Teichm\"uller geodesic flow (including the tangent cocycle). We refer the reader to the
recent and excellent survey \cite{BPHB} by L. Barreira and Ya. Pesin (and references therein), which covers all the relevant results on the {\it theory of Lyapunov exponents}, including Oseledec's  
multiplicative ergodic theorem and the Oseledec-Pesin reduction theorem.

\smallskip
\noindent The Lyapunov spectrum of the Teichm\"uller flow, with respect to any ergodic invariant probability measure $\mu$ on the moduli space, has symmetries. In fact, there exists 
non-negative numbers $\lambda^\mu_1=1 \geq \lambda^\mu_2 \geq \dots \geq \lambda^\mu_g$ 
such that the Lyapunov spectrum  of the Teichm\"uller flow has the following form (see \S 5 in \cite{Zo96}, \S 7 in \cite{Ko97} or \S 2.3 in \cite{Zo99}):
\begin{equation}
\label{eq:Tflowexp}
\begin{aligned}
 2 & \geq (1+\lambda^{\mu}_2)\geq \cdots\geq (1+\lambda^{\mu}_g) \geq 
\overbrace{1=\cdots= 1}^{\#(\Sigma_{\kappa})-1}  \geq  (1-\lambda^{\mu}_g) \geq \\
&\geq \cdots\geq (1-\lambda^{\mu}_2)\geq 0=0 \geq -(1-\lambda^{\mu}_2)\geq \cdots 
\geq -(1-\lambda^{\mu}_g) \geq \\
&\geq\underbrace{-1=\cdots= -1}_{\#(\Sigma_{\kappa})-1}\geq -(1+\lambda^{\mu}_g)
\geq \cdots \geq -(1+\lambda^{\mu}_2)\geq  -2 \,\,. 
\end{aligned}
\end{equation}
\noindent In \cite{Ve86} Veech proved that the Teichm\"uller flow is {\it non-uniformly hyperbolic}, in the sense that all of its Lyapunov exponents, except one corresponding to the flow direction, are non-zero. By formulas \pref{eq:Tflowexp} Veech's theorem can be formulated as follows:

\begin{theorem} \label{thm:l2ub} \cite{Ve86} The inequality
\begin{equation}
\label{eq:l2ub}
\lambda^{\mu}_2 <\lambda^{\mu}_1=1 \,\,. 
\end{equation}
holds if $\mu$ is the absolutely continuous $SL(2,\R)$-invariant ergodic probability measure on any connected component  of a stratum ${\Cal M}^{(1)}_{\kappa}\subset {\Cal M}^{(1)}_g$ of orientable quadratic differentials. 
\end{theorem}

\noindent  M. Kontsevich and A. Zorich have interpreted the non-negative numbers 
\begin{equation}
\label{eq:KZnonnegexps}
\lambda^{\mu}_1=1\geq \lambda^{\mu}_2\geq \dots \geq \lambda^{\mu}_g 
\end{equation}
as the non-negative Lyapunov  exponents of a symplectic cocycle over the Teichm\"uller flow, which we
now describe.

\smallskip
\noindent The {\it Kontsevich-Zorich cocycle }$\{\Phi_t\}_{t\in\R}$ is a cocycle over the Teichm\"uller flow 
$\{G_t\}_{t\in \R}$ on the moduli space ${\Cal M}_g$, defined as the projection of the trivial cocycle
\begin{equation}
\label{eq:KZcocycle}
G_t \times \hbox{id}: Q_g\times H^1(M,{\R})\to Q_g\times H^1(M,{\R}) 
\end{equation}
onto the orbifold vector bundle ${\Cal H}^1_g(M,{\R})$ over ${\Cal M}_g$ defined as 
\begin{equation}
\label{eq:CB}
{\Cal H}^1_g(M,{\R}):=\bigl(Q_g\times H^1(M,{\R})\bigr) /\Gamma_g\,.
\end{equation}
The mapping class group $\Gamma_g$ acts naturally on the cohomology by pull-back. 

\smallskip
\noindent  The Kontsevich-Zorich cocycle was introduced in \cite{Ko97} as a continuous-time 
version of the Zorich cocycle. The Zorich cocycle was introduced earlier by A. Zorich \cite{Zo96}, \cite{Zo97} in order to explain polynomial deviations in  the homological asymptotic behavior of 
 typical leaves of orientable measured foliations on compact surfaces, a phenomenon he had discovered in numerical experiments \cite{Zo94}.  We recall that the real homology $H_{1}(M,\R)$ 
 and the real cohomology $H^{1}(M,\R)$ of an orientable closed surface $M$ are (symplectically) isomorphic by Poincar\'e duality.

\smallskip 
\noindent Let ${\Cal M}_{\kappa}\subset {\Cal M}_g$ be a stratum of orientable holomorphic quadratic
differentials and $Q_{\kappa}\subset Q_g$ the pull-back of the stratum ${\Cal M}_{\kappa}$
to the Teichm\"uller space of quadratic differentials. The {\it relative }Kontsevich-Zorich coycle $\{\hat{\Phi}_t\}_{t\in\R}$ is defined as the projection of the trivial cocycle
\begin{equation}
\label{eq:KZrelcocycle}
G_t \times \hbox{id}: \bigcup_{q\in Q_{\kappa}} \{q\} \times H^1(M,
\Sigma_q;{\R})\to \bigcup_{q\in Q_{\kappa}} \{q\} \times H^1(M,
\Sigma_q;{\R})
\end{equation}
on the orbifold vector bundle ${\Cal H}^1_{\kappa}(M,\Sigma_{\kappa};{\R})$ over ${\Cal M}_{\kappa}$
defined as
\begin{equation}
\label{eq:RCB}
{\Cal H}^1_{\kappa}(M,\Sigma_{\kappa};{\R}):=\Bigl(\bigcup_{q\in Q_{\kappa}} 
\{q\} \times H^1(M,\Sigma_q;{\R})\Bigr)/\Gamma_g\,.
\end{equation}
By a similar construction it is possible to define a {\it punctured }Kontsevich-Zorich coycle $\{\Psi_t\}_{t\in\R}$ defined as the projection of the trivial cocycle 
\begin{equation}
\label{eq:KZpunctcocycle}
G_t \times \hbox{id}: \bigcup_{q\in Q_{\kappa}} \{q\} \times H^1(M\setminus
\Sigma_q;{\R})\to \bigcup_{q\in Q_{\kappa}} \{q\} \times H^1(M\setminus
\Sigma_q;{\R})
\end{equation}
onto the orbifold vector bundle ${\Cal H}^1_{\kappa}(M\setminus \Sigma_{\kappa};{\R})$ over 
${\Cal M}_{\kappa}$ defined as
\begin{equation}
\label{eq:PCB}
{\Cal H}^1_{\kappa}(M\setminus \Sigma_{\kappa};{\R}):=\Bigl(\bigcup_{q\in Q_{\kappa}} 
\{q\} \times H^1(M\setminus\Sigma_q;{\R})\Bigr)/\Gamma_g\,.
\end{equation}
However, by Poincar\'e-Lefschetz duality, there exists a natural isomorphism
\begin{equation}
\label{eq:PLD}
H^1(M, \Sigma_\kappa;\R) \equiv  H^1(M \setminus\Sigma_\kappa;\R)^\ast\,,
\end{equation}
hence the punctured Kontsevich-Zorich cocycle $\{\Psi_t\}_{t\in\R}$ is isomorphic to the dual $\{\hat{\Phi}^\ast_t\}_{t\in\R}$ of the relative Kontsevich-Zorich cocycle.

\smallskip
\noindent Let ${\Cal H}^1_{\kappa}(M,{\R})$ be the restriction to the stratum $\Cal M^{(1)}_\kappa$
of the bundle ${\Cal H}^1_g(M,{\R})$ defined in \pref{eq:CB}. Let $\Cal K^1_{\kappa} (M,\Sigma_{\kappa};{\R})$ be the bundle defined as the kernel of the natural bundle map
${\Cal H}^1_{\kappa}(M,\Sigma_{\kappa};{\R}) \,\, \to \,\, {\Cal H}^1_{\kappa}(M,{\R})$.

\begin{lemma} 
\label{lemma:cocycleiso}
The sub-bundle $\Cal K ^1_{\kappa}(M,\Sigma_{\kappa};{\R})$ is invariant under the relative Kontsevich-Zorich cocycle $\{\hat {\Phi}_t\}_{t\in\R}$. The Lyapunov spectrum of the restriction $\{\hat \Phi_t \vert \Cal K ^1_{\kappa}(M,\Sigma_{\kappa};{\R})\}$ consists of the single exponent $0$ with multiplicity $ \#(\Sigma_{\kappa}) -1$. In fact, there exists a Lyapunov norm on the bundle $\Cal K ^1_{\kappa}(M,\Sigma_{\kappa};{\R})$ for which the cocycle is isometric.
\end{lemma}
\begin{proof} By de Rham theorem for the relative cohomology, the relative cohomology complex $H^*(M,\Sigma_q; \R)$ is isomorphic to the cohomology of the complex of {\it relative }differential forms, which are defined as all differential forms vanishing at $\Sigma_q$. The map $H^*(M,\Sigma_q; \R) \to H^*(M, \R)$ is naturally defined on the de Rham cohomology since every closed, exact relative form is also a closed, respectively exact form on $M$. It follows that any class $c\in H^1(M,\Sigma_q; \R)$ which belongs to the kernel $K^1(M,\Sigma_q; \R)$ of the map $H^1(M,\Sigma_q; \R) \to H^1(M, \R)$ can be represented by a differential $1$-form exact on $M$. Hence for any $c\in K^1(M,\Sigma_q; \R)$ there exists a smooth function $f_c:M\to \R$, uniquely determined up to the addition of any function vanishing at $\Sigma_q$, such that $c=[df_c]$ as a relative de Rham cohomology class. 

\smallskip
\noindent As a consequence of the above discussion, the following formula yields a well-defined euclidean norm $\Vert \cdot \Vert_K$ on $K^1(M,\Sigma_q; \R)$:
\begin{equation}
\label{eq:Knorm}
\Vert c \Vert^2_K := \sum_{p_1, p_2\in \Sigma_q  }  \vert f_c(p_1) - f_c(p_2)\vert^2
\,, \quad c =[df_c] \in  K^1(M,\Sigma_q; \R)\,.
\end{equation}
The norm \pref{eq:Knorm} induces an euclidean norm on the bundle $\Cal K^1_{\kappa} (M,\Sigma_{\kappa};{\R})$ which is invariant under the relative Kontsevich-Zorich cocycle.
\end{proof}

\smallskip
\noindent Since the vector bundle ${\Cal H}^1_g(M,\R)$ has a symplectic structure, given
by the intersection form on its fibers, which are isomorphic to the cohomology $H^1(M,{\R})$, the Lyapunov spectrum of the cocycle $\{\Phi_t\}_{t\in \R}$ with respect to any $G_{t}$-invariant ergodic probability measure $\mu$ on ${\Cal M}^{(1)}_g$ is {\it symmetric}:
\begin{equation}
\label{eq:KZexps}
\lambda^{\mu}_1\geq \dots \geq\lambda^{\mu}_g\geq 0\geq \lambda^\mu_{g+1}=-\lambda^{\mu}_g\geq \dots \geq \lambda^\mu_{2g}=-\lambda^{\mu}_1\,.
\end{equation}
\noindent The non-negative Kontsevich-Zorich exponents  coincide with the numbers \pref{eq:KZnonnegexps} which appear in the Lyapunov spectrum \pref{eq:Tflowexp} of the Teichm\"uller flow. This relation is explained for instance in \cite{FHB}.

\smallskip
\noindent Zorich conjectured (see \cite{Zo96}) that the exponents \pref{eq:KZexps} are all distinct 
and different from zero when $\mu$ is the canonical absolutely continuous ergodic
invariant probability measure on any connected component of a stratum $\Cal M^{(1)}_\kappa$ of orientable quadratic differentials. In other terms he conjectured that the canonical measures are
$KZ$-hyperbolic and and $KZ$-simple, according to the following:

\begin{definition}
\label{def:nuhsimple} A $G_t$-invariant ergodic probability measure $\mu$ on a stratum of orientable quadratic differentials will be called \emph{KZ-hyperbolic }if the Kontsevich-Zorich cocycle $(\{\Phi_t\}, \mu)$ is \emph{non-uniformly hyperbolic}, in the sense that its Lyapunov exponents satisfy the inequalities 
\begin{equation}
\label{eq:nuh}
\lambda^{\mu}_1=1\geq \lambda^{\mu}_2\geq \dots \geq\lambda^{\mu}_g > 0 \,.
\end{equation}
A KZ-hyperbolic measure $\mu$ on a stratum of orientable quadratic differentials  will be called \emph{KZ-simple} if the Kontsevich-Zorich cocycle  $(\{\Phi_t\}, \mu)$ is \emph{simple}, in the sense that  all inequalities \pref{eq:nuh} are strict.
\end{definition}

\noindent In \cite{F02} we have proved the following result:

\begin{theorem} \label{thm:nuh}(\cite{F02}, Th. 8.5)  The absolutely continuous, $SL(2,\R)$-invariant, ergodic probability measure on any connected component  of a stratum ${\Cal M}^{(1)}_{\kappa}
\subset {\Cal M}^{(1)}_g$ of orientable quadratic differentials is KZ-hyperbolic.
\end{theorem}
\noindent In \cite{FHB}, \S 7, we have given an example of a $SL(2,\R)$-invariant
measure $\mu$ on $\Cal M_3$, supported on the closed $SL(2,\R)$-orbit of a particular genus $3$ branched cover of the $2$-torus, such that $\lambda^{\mu}_2=\lambda^{\mu}_3=0$.

\medskip
\noindent A proof of  the simplicity of the Zorich and Kontsevich-Zorich cocycles, which yields in particular a new independent proof of Theorem \ref{thm:nuh}, has been recently obtained by A. Avila and M. Viana \cite{AV} by methods completely different from ours. 

\begin{theorem} 
\label{thm:simplicity} \cite{AV} The absolutely continuous, $SL(2,\R)$-invariant, ergodic probability measure on any connected component  of a stratum ${\Cal M}^{(1)}_{\kappa}
\subset {\Cal M}^{(1)}_g$ of orientable quadratic differentials is KZ-simple.
\end{theorem}
\noindent The results of this paper do not depend in any way on the simplicity of the Kontsevich-Zorich cocycle, while the non-uniform hyperbolicity is crucial to the sharp estimates on the regularity
of solutions of the cohomological equation proven in \S \ref{ss:gencase}. However, our results can
be refined by taking into account that the Kontsevich-Zorich  exponents \pref{eq:nuh} are all distinct.

\smallskip
\noindent Our proof of Theorem \ref{thm:nuh} in \cite{F02} yields in particular a new independent 
proof of a strong version of Veech's Theorem \ref{thm:l2ub}. In fact, we have proved that the strict inequality \pref{eq:l2ub} holds for an arbitrary probability $G_t$-invariant measure on any stratum of orientable quadratic differentials. By combining our methods with a recent result of J. Athreya \cite{At06} on large deviations of the Teichm\"uller flow, it is possible to prove a similar strict upper bound for the second exponent for Lebesgue almost all quadratic differentials in {\it any }orbit of the circle group $SO(2,\R)$ on any stratum $\Cal M_\kappa$ of orientable quadratic differentials.

\smallskip
\noindent We recall below the variational formulas for the evolution of Hodge norm of absolute cohomology classes under the Kontsevich-Zorich cocycle. Following \S 2 in \cite{F02}, such formulas 
can be written in terms of a natural $\R$-linear extension $\hat U_q$ of the partial isometry $U_q$,  defined in formula \pref{eq:partiso}, which plays a crucial role in the construction of distributional solutions of the cohomological equation. Let $\hat U_q : L^2_q(M) \to L^2_q(M)$ be the $\R$-linear isometry defined as follows in terms of the partial isometry $U_q$ and of the orthogonal projections
$\pi^{\pm}_q: L^2_q(M) \to \Cal M^{\pm}_q$: 
\begin{equation}
\label{eq:Rliniso}
\hat U_q := U_q \circ (I -\pi_q^-) \, -\, \overline{ \pi^-_q} \,.
\end{equation}

\smallskip
\noindent  Let $\{q_t\}_{t\in \R}$ be the orbit  of an orientable quadratic differential $q\in \Cal Q^{(1)}_\kappa$ under the Teichm\"uller flow $\{G_t\}_{t\in \R}$.  We remark that by the definition 
of the Teichm\"uller flow $\{G_t\}_{t\in \R}$, the area form $\omega_t$ of the metric induced by the quadratic differential $q_t$ is {\it constant} equal to $\omega_q$ for all $t\in \R$. Hence the Hilbert space $L^2_q(M)$ is invariant under the action of the Teichm\"uller flow (in fact, it is invariant under the
full $SL(2,\R)$ action). For each $t\in \R$, let $M_t$ the Riemann surface carrying $q_t\in {\Cal Q}^{(1)}_{\kappa}$ and let ${\Cal M}^{\pm}_t \subset L^2_q(M)$ be the space of meromorphic, respectively anti-meromorphic, functions on the Riemann surface $M_t$. Such spaces are respectively the kernels of the adjoints of the Cauchy-Riemann operators $\partial^{\mp}_t$, associated to the holomorphic quadratic differential $q_t$. The dimension of ${\Cal M}^{\pm}_t$ is constant equal to the genus $g\geq 1$ of $M$ (see \S \ref{WSS}) and it can be proved that $\{{\Cal M}^{\pm}_t\}_{t\in {\R}}$ are smooth families of $g$-dimensional subspaces of the fixed Hilbert space $L^2_q(M)$. 

\smallskip
\noindent By \pref{eq:reprs} there exists a one-parameter family $\{m^+_t\}_{t\in \R} \subset 
\Cal M^+_t$ such that
\begin{equation}
\label{eq:ct}
c_t=c^+_{q_t} (m^+_t):=[\Re(m_t^+q_t^{1/2})] \in H^1(M_t,{\R})\,.
\end{equation}

\begin{lemma} (\cite{F02}, Lemma 2.1) 
\label{lemma:varfor}
The ordinary differential equation 
\begin{equation}
\label{eq:KZODE}
u'=\hat U_{q_t}(u) 
\end{equation} 
is well defined in $L^2_q(M)$ and satifies the following properties:
\begin{enumerate}
\item Solutions of the Cauchy problem for \pref{eq:KZODE} exist for all 
times and are uniquely determined by the initial condition;
\item  If $u_t\in L^2_q(M)$ is any solution of \pref{eq:KZODE} such that the initial
condition $u_0\in {\Cal M}^+_q$, then $u_t \in {\Cal M}^+_t$ for all 
$t\in {\R}$.
\item Let $m^+_t\in {\Cal M}^+_t$ be the unique solution of \pref{eq:KZODE} 
with initial condition $m^+_0=m^+\in {\Cal M}^+_q$. For all $t\in {\R}$, we have
\begin{equation}
\label{eq:varfor}
\Phi_t \bigl( c^+_q(m^+)\bigr) \,=\, c^+_{q_t}(m^+_t)\,\,. 
\end{equation}
\end{enumerate}
\end{lemma}
\noindent  It follows immediately from Proposition \ref{prop:CR} that, for every $u\in L^2_q(M)$, there exist functions $v^{\pm} \in H^1_q(M)$ such that
\begin{equation}
\label{eq:L2split}
u\,= \, \partial ^{+}_{q} v^{+} \, + \, \pi^{-}_{q}(u) \,=\,  \partial ^{-} _{q}v^{-} \, + \, \pi^{+}_{q}(u)\,.
\end{equation}
The O. D. E.  in formula \pref{eq:KZODE}  can be written explicitly, in terms of the orthogonal decompositions \pref{eq:L2split}, as follows. Let $\pi^{\pm}_t :  L^2_q(M) \to  {\Cal M}^{\pm}_t$ denote the orthogonal projections in the (fixed) Hilbert space $L^2_q(M)$. By definition, the projections $\pi^{\pm}_t$ coincide with the projections $\pi^{\pm}_q$ for $q=q_{t}$, for any $t\in \R$. A function
$u\in C^1\left(\R,L^2_q(M)\right)$ satisfies equation \pref{eq:KZODE} iff
\begin{equation}
\label{eq:varforbis}
\left \{ \begin{array}{cl}
u_t = \partial^+_t v_t \,+\, \pi^-_t (u_t)\,\,; \\
     { \frac{d}{dt}}  u_t =\partial^-_t v_t \,-\,{\overline{ \pi^-_t (u_t)}}\,\,. 
 \end{array} \right.  
\end{equation}

\noindent  An immediate consequence of Lemma \ref{lemma:varfor} is the following result on the variation of the Hodge norm of cohomology classes under the action of the Kontsevich-Zorich cocycle. Let $B_{q}:L^{2}_{q}(M)\times  L^{2}_{q}(M) \to \C$ be the complex bilinear form given by
\begin{equation}
\label{eq:Bform}
B_{q}(u,v) := \int_{M} u \,v \, \omega_{q}\, ,\quad \hbox{ for all } u,v \in L^{2}_{q}(M)\,.
\end{equation}
\begin{corollary}  
\label{corollary:normder} (\cite{F02}, Lemma 2.1') The variation of the Hodge
norm $\Vert c_{t} \Vert_H$, which coincides with the $L^{2}_{q}$-norm $\vert m^{+}_{t}\vert_{0}$
under the identification \pref{eq:ct}, is given by the following formulas:
\begin{equation}
\label{eq:normder}
\begin{aligned}
&(a)\,\,{ \frac{d}{dt}}|m^+_t|_0^2=-2\,\Re [B_q(m^+_t)]
                                 = -2\,\Re\left[\int_M (m^+_t)^2 \omega_q\right]\,; \\ 
&(b)\,\,{ \frac{d^2}{dt^2}}|m^+_t|_0^2 = 4\left\{|\pi^-_t(m^+_t)|_0^2 -\Re
             \left[\int_M ( \partial^+_t v_t )\,(\partial^-_t v_t) \,\omega_q\,\right] \right\}\,.
\end{aligned}
\end{equation}
\end{corollary}
\noindent The second order variational formula \pref{eq:normder}, $(b)$, is crucial in our proof
of {\it lower bounds }for the Kontsevich-Zorich exponents (see Theorem \ref{thm:nuh}).
The first order variational formula  \pref{eq:normder}, $(a)$, implies quite immediately
an effective upper bound for the second exponent, which yields in particular 
the {\it average spectral gap} result proved in \cite{F02}, Corollary 2.2 (a generalization of Veech's 
Theorem \ref{thm:l2ub} to arbitrary $G_t$-invariant ergodic probability measures):

\begin{theorem}  (\cite{F02}, Corollary 2.2)  
\label{thm:spectralgapone}
The inequality
\begin{equation}
\label{eq:l2mub}
\lambda^{\mu}_2 <\lambda^{\mu}_1=1 \,\,. 
\end{equation}
holds for any  $G_t$-invariant ergodic probability measure $\mu$ on any connected component  of any stratum ${\Cal M}^{(1)}_{\kappa}$ of orientable quadratic differentials. 
\end{theorem}

\noindent The above spectral gap result implies the following unique ergodicity theorem for measured foliations:
\begin{corollary} For any stratum $\Cal M^{(1)}_\kappa$ of orientable quadratic differentials, the set  
of quadratic differentials $q\in \Cal M^{(1)}_\kappa$ with minimal but not uniquely ergodic horizontal [vertical] foliation has zero measure with respect to any $G_t$-invariant probability measure $\mu$ on $\Cal M^{(1)}_\kappa$.
\end{corollary}

\noindent In the remainder of this section we prove a {\it pointwise spectral gap }result which holds for almost all quadratic differentials in {\it any }orbit of the circle group $SO(2,\R)$ on any stratum $\Cal M^{(1)}_\kappa$ of orientable quadratic differentials. The argument is based on formula  \pref{eq:normder}, $(a)$ and, as mentioned above, on a result of J. Athreya \cite{At06} on large deviations 
of the Teichm\"uller flow. 

\smallskip
\noindent The upper second Lyapunov exponent of the Kontsevich-Zorich cocycle at any (orientable) quadratic differentials $q\in \Cal M^{(1)}_\kappa$ is defined as follows. Let $I_q(M,\R) \subset H^1(M,\R)$ be the subspace of real dimension $2$ defined as 
\begin{equation}
\label{eq:Iq}
I_q(M,\R):= \R \cdot \Re(q^{1/2}) + \R \cdot \Im(q^{1/2}) 
\end{equation}
and let $I_q^\perp (M,\R)$ be the symplectic orthogonal of $I_q(M,\R)$ in $H^1(M,\R)$,
with respect to the symplectic structure induced by the intersection form:
\begin{equation}
\label{eq:Iqperp}
I_q^\perp (M,\R):=\{ c\in H^1(M,\R)   \vert    c \wedge [q^{1/2}] =0\}\,.
\end{equation}
\noindent The complementary sub-bundles $I_\kappa(M,\R)$ and $I^\perp_\kappa(M,\R) \subset \Cal H^1(M,\R)$, with fibers at any $q\in \Cal M^{(1)}_\kappa$ respectively equal to $I_q(M,\R)$ and $I^\perp_q(M,\R)$, are invariant under the Kontsevich-Zorich cocycle. In fact, it is immediate to verify that the sub-bundle $I_\kappa(M,\R)$ is invariant under the Kontsevich-Zorich cocycle and that the Lyapunov spectrum of the restriction of the Kontsevich-Zorich cocycle to $I_\kappa(M,\R)$ equals $\{1,-1\}$ (both exponents with multiplicity $1$). Since the Kontsevich-Zorich cocycle is symplectic,  the symplectic orthogonal bundle $I^\perp_\kappa(M,\R)$ is also invariant. In addition, it is not difficult to verify that $1$ is the top (upper) exponent for the cocycle on the full cohomology bundle 
$\Cal H^1_\kappa(M, \R)$.

\smallskip
\noindent The {\it second upper (forward) exponent }of the Kontsevich-Zorich cocycle is the top 
upper (forward) Lyapunov exponent at any quadratic differential $q\in\Cal M^{(1)}_\kappa$ of the restriction of the cocycle to the sub-bundle $I^\perp_\kappa(M,\R)$:

\begin{equation}
\label{eq:l2plus}
\lambda^+_2(q) :=  \limsup_{t\to +\infty}  \frac{1}{t} \log  \Vert
 \Phi_t \vert I_q^\perp(M,R) \Vert_H \,.
\end{equation}

\begin{theorem}
\label{thm:spectralgaptwo}
 For any stratum $\Cal M^{(1)}_\kappa$ of orientable 
quadratic differentials, there exists a measurable function $L_\kappa: \Cal M^{(1)}_\kappa \to [0,1)$ 
such that for any (orientable) quadratic differential $q\in \Cal M^{(1)}_\kappa$, 
\begin{equation}
\label{eq:specgap}
 \lambda^+_2(q_\theta)   \, \leq  \, L_\kappa (q)  <1\,,  \quad \text {\rm for almost all } \, \,\theta\in S^1\,.
\end{equation}
\end{theorem}

\begin{proof} The argument follows closely the proof of Corollary 2.2 in \cite{F02}.
Under the isomorphism \pref{eq:reprs}, the vector space $I^\perp_q(M,\R)$ is represented by meromorphic functions with {\it zero average }(orthogonal to constant functions). 
It can be seen that the subspace of zero average meromorphic functions is invariant under the flow of equation \pref{eq:varfor} or, equivalently \pref{eq:varforbis}. By formula $(a)$ in \pref{eq:normder},
\begin{equation}
\label{eq: derlog}
\frac{d}{dt} \log |m^+_t|_0^2=-2\, { \frac{\Re \,B_q(m^+_t) }{\vert m^+_t\vert_0^2}}\,.
\end{equation}
Following \cite{F02}, we introduce a continuous function $\Lambda_\kappa^+:  \Cal M^{(1)}_\kappa \to 
\R^+$ defined as follows: for any $q\in  \Cal M^{(1)}_\kappa$, 
\begin{equation}
\label{eq:Lambdaplus}
\Lambda_\kappa^+(q):= \max \{ { \frac{\vert B_q(m^+)\vert}{\vert m^+\vert_0^2}}\,|\, m^+\in {\Cal M}^+_q\setminus\{0\}\,,\,\,\int_M m^+\,\omega_q=0\,\} \,. 
\end{equation}
Since by the Schwarz inequality,
\begin{equation}
\label{eq:Schwarz}
\vert B_q(m^+_t)\vert=\vert(m^+_t, {\overline{m^+_t}})_q \vert \leq \vert m^+_t\vert_0^2\,\,, 
\end{equation}
the image of the function $\Lambda_\kappa^+$ is contained in the interval $[0,1]$. We claim
that $\Lambda_\kappa^+(q)<1$ for all $q\in \Cal M^{(1)}_\kappa$.  In fact, $\Lambda^+(q)=1$ if 
and only if there exists a {\it non-zero }meromorphic function with zero average $m^+\in {\Cal M}^+_q$ such that $\vert(m^+,{\overline {m^+}})_q\vert=\vert m^+\vert_0^2$. A well-known property of the Schwarz inequality then implies that there exists $u \in {\C}$ such that $m^+=u \,{\overline{m^+}}$. However, it cannot be so, since $m^+$ would be meromorphic and anti-meromorphic, hence constant, and by the zero average condition it would be zero. 

\smallskip
\noindent It follows from formula \pref{eq: derlog} that, for any $q\in \Cal M^{(1)}_\kappa$,
\begin{equation}
\label{eq:logbound}
 \frac{1}{t} \log  \Vert  \Phi_t \vert I_q^\perp(M,\R) \Vert_H \leq   \frac{1}{t}  \int_0^t  \Lambda_\kappa^+\left(G_s(q)\right)\, ds \,.
\end{equation}
Let $q_0\in \Cal M^{(1)}_\kappa$. By the large deviation result of J. Athreya (see \cite{At06}, 
Corollary 2.4) the following holds. For any $\lambda<1$ there exists a compact
set $K \subset  \Cal M^{(1)}_\kappa$ such that, for almost all $q\in SO(2,\R)\cdot q_0$,
\begin{equation}
\label{eq:largedev}
\limsup_{t\to \infty} \frac{1}{t} \, \vert\{ 0 \leq s \leq t\, \vert \, G_s(q) \not\in K \} \vert
\,\,  \leq \,\, \lambda \,.
\end{equation}
Let $\Lambda_K:= \max \{ \Lambda_\kappa(q) \,\vert \, q\in K\}$ and, for any $(t,q)\in \R^+\times
\Cal M^{(1)}_\kappa$, let
$$
\Cal E_K(t,q):= \vert\{ 0 \leq s \leq t\, \vert \, G_s(q) \not\in K \} \vert \,.
$$ 
Since the function $\Lambda^+_\kappa$ is continuous and $\Lambda^+_\kappa(q)<1$ for
all $q\in \Cal M^{(1)}_\kappa$, its maximum on any compact set is $<1$, in particular 
$\Lambda_K <1$. The following immediate inequality holds:
\begin{equation}
\label{eq:logintsplit}
\int_0^t  \Lambda_\kappa^+\left(G_s(q)\right)\, ds  \leq  
(1 -  \Lambda_K) \,  {\Cal E}_K(t,q) + t \, \Lambda_K\,.
\end{equation}
It follows from \pref{eq:logbound}, \pref{eq:largedev} and \pref{eq:logintsplit} that, for
almost all $q\in SO(2,\R)\cdot q_0$, 
\begin{equation}
\limsup _{t\to +\infty} \frac{1}{t} \log  \Vert  \Phi_t \vert I_{q}^\perp(M,\R) \Vert_H
\,\leq (1 -  \Lambda_K) \lambda \,+\, \Lambda_K \,<\,1\,.
\end{equation}
The function $L_\kappa :\Cal M^{(1)}_\kappa \to [0,1)$ can be defined for every $q_0\in 
\Cal M^{(1)}_\kappa$ as the essential supremum of the second upper Lyapunov exponent
over the orbit $SO(2,\R)\cdot q_0$ of the circle group. Such a function is measurable by definition 
and it is everywhere $<1$ by the above argument.  
\end{proof}
\noindent For any Oseledec regular point $q\in Q^{(1)} _{\kappa}$ of the Kontsevich-Zorich 
cocycle, let $E^{+}_q(M,\R) [E^{-}_q(M,\R)]  \subset H^1(M,\R)$ be the unstable [stable] subspace of the
Kontsevich-Zorich cocycle. Homology cycles which are Poincar\'e duals to cohomology classes 
in $E^\pm_q(M,\R)$ are called {\it Zorich cycles }for the foliation $\Cal F_{\pm q}$.  It follows from Theorem \ref{thm:nuh} and from the symplectic property of the cocycle that $E^\pm_q$ are transverse Lagrangian subspaces (with respect to the intersection form), as conjectured by Zorich in  \cite{Zo96}. In \cite{F02}, Theorem 8.3 (see also \cite{FHB}, Theorem 8.2) we have proved the following representation theorem:
\begin{theorem}
\label{thm:Zcyclesrepr}
For Lebesgue almost all $q\in {\Cal M}^{(1)}_{\kappa}$, 
we have
\begin{equation}
\label{eq:basicrepr}
E^+_q(M,\R) = H^{1,1}_q(M,{\R}) \,, \quad 
E^-_q(M,\R) = H^{1,1}_{-q}(M,{\R})\,.
\end{equation}
(The Poincar\'e duals of Zorich cycles for a generic orientable measured foliation $\Cal F$ are represented by basic currents for $\Cal F$ of Sobolev order $\leq 1$).
\end{theorem}

\noindent We prove below a conjectural relation between Lyapunov exponents of a cohomology class under the Kontsevich-Zorich cocycle and the Sobolev regularity of the basic current representing the cohomology class. This result answers in the affirmative a question posed by the author in \cite{F02} (Question 9.9) and, independently by M. Kontsevich (personal communication). In order to formulate
our results in the greatest possible generality,  we introduce the following class of measures on the moduli space.

\begin{definition} 
\label{def:circleabs}
A measure $\mu$ on the moduli space $\Cal M_g$ of holomorphic quadratic differentials will be called $SO(2, \R)$-\emph{absolutely continous}  if induces absolutely continous
conditional measures on $\pi_\ast(\mu)$-almost every fiber of the fibration $\pi: \Cal M_g \to \Cal M_g/SO(2,\R)$ (with respect to the Haar/Lebesgue measure class on each fiber). 
\end{definition} 
\noindent It is immediate that any $SO(2,\R)$-invariant, hence {\it a fortiori} any $SL(2,\R)$-invariant measure is $SO(2, \R)$-absolutely continous.

\begin{lemma} 
\label{lemma:reglowbound}
Let $\mu$ be any  $SO(2, \R)$-absolutely continous, $G_t$-invariant ergodic probability measure on a stratum $\Cal M^{(1)}_\kappa$ of orientable quadratic differentials. For $\mu$-almost all $q\in \Cal M^{(1)}_\kappa$, the unique basic current $C\in \Cal B^1_q(M)$ [$C\in \Cal B^1_{-q}(M)$]  which represents a cohomology class $c\in E^+_q(M,\R)$ [$c\in E^-_q(M,\R)]$ of Lyapunov exponent $\lambda(c)>0$  [$\lambda(c)<0$] 
under the Kontsevich-Zorich cocycle, has the following Sobolev regularity: 
$$
C \in W^{-s}_q(M)\, , \quad \text{ \rm for all }\,\, s> 1-\vert \lambda(c)\vert\,\,.
$$
\end{lemma}

\begin{proof}
The argument follows the proof of Lemma 8.2 in \cite{F02}. The interpolation inequality \pref{eq:intineq} for fractional weighted Sobolev spaces allows us to estimate the Sobolev regularity of the basic
current constructed there. 

\smallskip
\noindent Let $q\in \Cal M^{(1)} _{\kappa}$ be any Oseledec regular point of the Kontsevich-Zorich cocycle and let $c\in H^1(M,\R)$ be a cohomology class.  Let $\{q_t\}_{t\in \R}\subset Q^{(1)} _{\kappa}$
denotes the lift to the Teichm\"uller space of the orbit $\{G_t(q)\}_{t\in \R}$ of $q$ under the Teichm\"uller flow. Let ${\Cal M}^+_t \subset L^2_q(M)$ be the space of meromorphic functions on the Riemann surface $M_t$ carrying the quadratic differential $q_t \in Q^{(1)} _{\kappa}$. According to the representation formula  \pref{eq:reprs},  there exists a (smooth) family $\{m^+_t\}_{t\in \R} \subset {\Cal M}^+_t$  such that, for each $t\in \R$, 
\begin{equation}
\label{eq:oneparrep}
\Phi_t(c) = \Re [m^+_t q_t^{1/2}]\,\,.
\end{equation}
 \noindent By the variational formulas \pref{eq:varfor} and \pref{eq:varforbis}, there exists a (smooth) 
 family $\{v\}_{t\in\R} \subset H^1(M)$ of {\it zero-average functions} such that 
\begin{equation}
\label{eq:varfortris}
\left \{ \begin{array}{cl}
m^+_t=\partial_t^+ v_t  + \pi^-_t (m^+_t) \,,\\
 {\frac{d}{dt}} m^+_t = \partial_t^- v_t  -\overline{ \pi^-_t (m^+_t)}\,.
 \end{array} \right.
\end{equation}
\noindent Since $\Phi_t(c)\equiv c \in H^1(M,\R)$ (by the definition \pref{eq:KZcocycle} of the 
Kontsevich-Zorich cocycle $\{\Phi_t\}$ {\it over the Teichm\"uller space}), for each $t\in \R$ there exists a unique zero average function $U_t\in L^2_q(M)$ such that
\begin{equation}
\label{eq:dUt}
dU_t = \Re[m^+_t q_t^{1/2}]\,-\,\Re[m^+_0q^{1/2}]\,.
\end{equation}
It follows that the family $\{U_t\}_{t\in\R}$ is smooth and satisfies, by the variational formulas \pref{eq:varfortris}, the following Cauchy problem in $L^2_q(M)$:
 \begin{equation}
\label{eq:Udot}
\left \{ \begin{array}{cl}
{\frac{d}{dt}}U_t &= 2\,\Re(v_t)\,,\\
                    U_0 &=  0\,,
 \end{array} \right.
 \end{equation}
We claim that, if $c\in E^+_q(M,\R)$ [$c\in E^-q(M,\R)$] has Lyapunov exponent 
$\lambda(c)>0$ [$\lambda(c)<0$], the set $\{U_t\, \vert \,t \leq 0\}$ [ $\{U_t\, \vert \,t \geq 0\}$]
is a bounded subset of the Hilbert space $H^s_q(M)$ for any $s<\vert \lambda(c)\vert$.

\smallskip
\noindent By Oseledec's theorem, for any  $0<\lambda< \vert\lambda(c)\vert$, there exist a measurable function $K_\lambda>0$ on ${\Cal M}^{(1)}_{\kappa}$ and  such that
\begin{equation}
\label{eq:expbound}
\Vert \Phi_t(c) \Vert_H =\vert m^+_t \vert_0\leq K_\lambda(q) \vert m^+_0\vert_0\,\exp (-\lambda |t|)\,, \quad t\leq 0 \,\, [t\geq 0]\,.
\end{equation}

\noindent For any (orientable) quadratic differential $q\in Q^{(1)} _{\kappa}$, let $\Vert q \Vert$ 
denote the length of the shortest geodesic segment with endpoints in the set $\Sigma_q$ (of 
zeroes of $q$) with respect to the induced metric. By the Poincar\'e inequality proved in \cite{F02}, Lemma 6.9, there exists a constant $K_{g,\sigma}>0$ (depending on the genus $g\geq 2$ of the Riemann surface $M_q$ and on the cardinality $\sigma:= \#(\Sigma_q)$ of any $q\in \Cal Q^{(1)}_\kappa$) such that
\begin{equation}
\vert v-\int_M v\,\omega_q \vert _0 \leq    \frac{ K_{g,\sigma}}{ \Vert q \Vert} \,\, \Cal Q(v,v)\,\,, 
\quad\text{ for all } v\in H^1_q(M)\,.
\end{equation}
By the commutativity property  \pref{eq:commut} of the horizontal and vertical vector fields, the Dirichlet form $\Cal Q$ of the quadratic differential can be written as
$$
{\Cal Q}(v,v)=|\partial^{\pm}_q v|_0^2\,\,, \quad\text{ for all } v\in H^1_q(M)
$$ 
(where $\partial^{\pm}_q$ are the Cauchy-Riemann operators introduced in \S \ref{fracsobspaces}).

\smallskip
\noindent By the Poincar\'e inequality and by the orthogonality of the decompositions 
in \pref{eq:varfortris} with respect to the invariant euclidean structure on $L^2_q(M)$, we have
\begin{equation}
\label{eq:vt}
\vert v_t\vert_0\leq K_{g,\sigma}\Vert q_t \Vert^{-1} \vert \partial^+_t v_t\vert _0 \leq 
K_{g,\sigma} \Vert q_t \Vert ^{-1}  \vert m^+_t\vert_0 \,.
\end{equation}
\noindent It follows by  \pref{eq:Udot}, \pref{eq:expbound} and \pref{eq:vt} that there exists a measurable function $K'_{\lambda}>0$ on ${\Cal M}^{(1)}_{\kappa}$ such that if $c\in E^+_q(M,\R)$ [$c\in E^-_q(M,\R)$] the following inequality holds for all $t\leq 0$  [ $t\geq 0$]:
\begin{equation}
\label{eq:Udotbound}
\vert \frac{d}{dt} U_t \vert _0\leq 2\,\vert v_t \vert _0\leq K'_\lambda(q)\, \vert m^+_0 \vert_0\, 
\Vert q_t \Vert ^{-1}\, e^{- \lambda |t|} \,.
\end{equation}
By \pref{eq:varfortris} and \pref{eq:Udot} we also have the following straightforward estimate for
the norm of the funcion $\frac{d}{dt} U_t$ in the fixed Hilbert space $H^1_q(M)$:
\begin{equation}
\label{eq:H1est}
\vert \frac{d}{dt} U_t \vert _1 \leq e^{\vert t \vert} \, \vert m^+_t\vert_0 \leq K_\lambda(q) \vert m^+_0\vert_0\, e^{(1-\lambda) \vert t \vert }\,.
\end{equation}
It follows from the inequalities \pref{eq:Udotbound} and \pref{eq:H1est}, by the interpolation 
inequality proved in Lemma \ref{lemma:intineq}, that for any $s\in (0,1)$ there exists a measurable function $K^s_{\lambda}>0$ on the stratum $Q^{(1)}_\kappa$ such that, if $c\in E^+_q(M,\R)$ 
[$c\in E^-_q(M,\R)$] the following inequality holds for all $t\leq 0$  [ $t\geq 0$]:
\begin{equation}
\label{eq:Hsest}
\vert \frac{d}{dt} U_t \vert _s \leq e^{\vert t \vert} \, \vert m^+_t\vert_0 \leq K^s_\lambda(q) \vert m^+_0\vert_0\, e^{(s-\lambda) \vert t \vert }\,.
\end{equation}
Since $U_0=0$, by Minkowski's integral inequality we finally obtain the estimate:
\begin{equation}
\label{eq:Ubound}
\vert U_t \vert _s \leq K^s_\lambda (q) \vert m^+_0 \vert _0\int_0^{|t|} 
e^{(s-\lambda) \tau } \Vert q_\tau \Vert ^{-1}\, d\tau\,.
 \end{equation}
By the  {\it logarithmic law }for the Teichm\"uller geodesic flow on the moduli space, proved by 
H. Masur in \cite{Ma93},  the following estimate holds for almost all quadratic differentials $q\in 
{\Cal M}^{(1)}_{\kappa}$ (see \cite{Ma93}, Prop. 1.2):
\begin{equation}
\label{eq:Loglaw}
 \limsup_{\tau\to \pm \infty} {\frac{-\log\Vert q_\tau \Vert}{\,\log |\tau|}}\,\leq \, \frac{1}{2}\,. 
\end{equation}
As a consequence, for $\mu$-almost all $q\in {\Cal M}^{(1)}_{\kappa}$ and for any $s<\lambda$, 
the integral in formula \pref{eq:Ubound} is uniformly bounded for $t\leq 0$ [$t \geq 0$]. Since for any $s<\vert \lambda(c)\vert$, there exists $\lambda\in (s,\vert \lambda(c)\vert)$, it follows that 
the family of functions $\{U_t\, \vert \,t \leq 0\}$ [$\{U_t\, \vert \,t \geq 0\}$] is uniformy bounded in the 
Sobolev space $H^s_q(M)$ for any $s<\vert \lambda(c)\vert$, as claimed.

\smallskip
\noindent For any $s<\vert\lambda(c)\vert $, let $U^+_s \in H^s_q(M)$ [$U^-_s \in H^s_q(M)$] 
be any weak limit of the family $\{U_t\}$ as $t\to -\infty$  [as $t\to -\infty$], which exists since all
bounded subsets of the separable Hilbert space $H^s_q(M)$ are sequentially weakly compact. 
Since the functions $U_t$ have zero average for all $t\in \R$ and the subspace of zero average
functions is closed in $H^s_q(M)$ for all $s>0$, the weak limit $U_s^+$ [$U_s^-$] has zero average.
By contraction of the identity \pref{eq:dUt} with the horizontal vector field $S$ [with the vertical 
vector field $T$] we have:
\begin{equation} 
\label{eq:SU}
\begin{aligned}
SU_t &= -\Re(m^+_0)\,+\, e^t\,\Re(m^+_t)\,, \quad t\leq 0\,,\\
[\,TU_t &= \Im(m^+_0)\,+\, e^{-t}\,\Re(m^+_t)\,, \quad t\geq 0\,,]
\end{aligned} 
\end{equation}
and by taking the limit as $t\to -\infty$ [as $t\to +\infty$],
\begin{equation} 
\label{eq:SUbis}
\begin{aligned}
SU^+_s  & =-\Re(m^+_0)\,,\\
[\,TU^-_s  & =\Im(m^+_0)\,.]
\end{aligned} 
\end{equation}
Since for almost all quadratic differential $q\in Q^{(1)}_\kappa$ the horizontal foliation
[the vertical foliation] is ergodic, the solution $U^+\in L^2_q(M)$ [$U^- \in L^2_q(M)$]
of the cohomological equation \pref{eq:SUbis} is unique (if it exists). Hence there exists
a unique zero-average function $U^+\in L^2_q(M)$ [$U^- \in L^2_q(M)$] , which solves the cohomological equation \pref{eq:SUbis}, such that $U^+_s=U^+$ [$U^-_s=U^-$] for all 
$s<\vert \lambda(c)\vert$. As a consequence, $U^+ \in H^s_q(M)$ [$U^- \in H^s_q(M)$] 
for all $s<\vert \lambda(c)\vert$.  The current $C^+\in W^{s-1}_q(M)$ [$C^- \in W^{s-1}_q(M)$] uniquely determined by the identity
\begin{equation}
\label{eq:Cpm}
\begin{aligned}
dU^+  &= C^+  \,-\,  \Re[m^+_0q^{1/2}] \,,\\
[\,dU^-&= C^-  \,-\,  \Re[m^+_0q^{1/2}] \,,]
\end{aligned}
\end{equation}
is basic for the horizontal [vertical] foliation by Lemma \ref{lemma:bcconstr} and represents the cohomology class $c\in E^+_q(M,\R)$ [$c\in E^-_q(M,\R)$] of Lyapunov exponent $\lambda(c)>0$ [$\lambda(c)<0$].
\end{proof}

\subsection{Distributional cocycles}
\label{distcocycles}

\noindent  Let $\Cal Q_\kappa(M)$ be the space of all orientable quadratic differentials, holomorphic with respect to some complex structure on a closed surface $M$, with zeros of multiplicities $\kappa=(k_1, \dots, k_\sigma)$. For any $s\in \R$, there is a natural action of the group $\text{\rm Diff}^+(M)$ of orientation preserving diffeomorphisms on the trivial bundles
\begin{equation}
\label{eq:distbundleone}
\bigcup_{q\in  \Cal Q_\kappa(M)}  \{q\} \times H_q^{s}(M)\,\,  \subset \,\, 
 \bigcup_{q\in  \Cal Q_\kappa(M)}  \{q\} \times \bar H_q^{s}(M)\ \,.
\end{equation}
In fact, any diffeomorphism $f\in \text{\rm Diff}^+(M)$ defines by pull-back an isomorphism
$f^* :  \bar H_q^{s}(M) \to  \bar H_{f^*(q)}^{s}(M)$ which maps the subspace $H_q^{s}(M) 
\subset  \bar H_q^{s}(M)$ onto $ H_{f^*(q)}^{s}(M) \subset \bar H_{f^*(q)}^{s}(M)$. The 
quotient bundles 
\begin{equation}
\label{eq:distbundletwo}
\begin{aligned}
\Bigl(\bigcup_{q\in  \Cal Q_\kappa(M)}  \{q\} \,&\times H_q^{s}(M)\Bigr) /
 \text{\rm Diff}_0^+(M)\,, \\
 \Bigl(\bigcup_{q\in  \Cal Q_\kappa(M)}  \{q\} \,&\times \bar H_q^{s}(M)\Bigr) /
 \text{\rm Diff}_0^+(M) 
 \end{aligned}
\end{equation}
are well-defined bundles over the stratum $Q_\kappa$ of the Teichm\"uller space of quadratic 
differentials. There is natural action of the mapping class group $\Gamma_g$ on the bundles 
\pref{eq:distbundletwo} induced by the action of $\text{\rm Diff}^+(M)$ on the bundles 
\pref{eq:distbundleone}. The resulting quotient bundles 
\begin{equation}
\label{eq:distbundlethree}
\begin{aligned}
H^{s}_{\kappa}(M)&:= \Bigl(\bigcup_{q\in  \Cal Q_\kappa(M)}  \{q\} \,\times H_q^{s}(M)\Bigr) /
 \text{\rm Diff}^+(M)\,, \\
 \bar H^{s}_{\kappa}(M)&:= \Bigl(\bigcup_{q\in  \Cal Q_\kappa(M)}  \{q\} \,\times H_q^{s}(M)\Bigr) /
 \text{\rm Diff}^+(M)
 \end{aligned}
\end{equation}
are well-defined bundles over the stratum $\Cal M_\kappa$ of the moduli  space.
\smallskip
\noindent We also introduce bundles of $1$-currents over a stratum $\Cal M_\kappa$ as follows:
\begin{equation}
\label{eq:currbundle}
\begin{aligned}
W^{s}_{\kappa}(M):=\Bigl(\bigcup_{q\in  \Cal Q_\kappa(M)}  \{q\} \times W^{s}_q(M)\Bigr)/
\text{\rm Diff}^+(M)\,,\\
\bar W^{s}_{\kappa}(M):=\Bigl(\bigcup_{q\in  \Cal Q_\kappa(M)}  \{q\} \times \bar W^{s}_q(M)\Bigr)/
\text{\rm Diff}^+(M)\,.
\end{aligned}
\end{equation}
Since, for all $q\in \Cal Q_\kappa(M)$ and for all $s\in \R$, the weighted Sobolev space of  $1$-currents ${W}_q^{s}(M)$ is isomorphic to the tensor product ${\R}^2\otimes H^s_q(M)$, 
\begin{equation}
\label{eq:tensor}
\begin{aligned}
W^{s}_{\kappa}(M) \equiv  \R^2 \otimes H^{s}_{\kappa}(M) \,,
\bar W^{s}_{\kappa}(M) \equiv  \R^2 \otimes \bar H^{s}_{\kappa}(M)
\end{aligned}
\end{equation}

\smallskip
\noindent 
Let  $\{G^s_t\}_{t\in\R}$ be the cocycle over the Teichm\"uller flow, defined as the projection onto the bundle $H^{-s}_{\kappa}(M)$ of the trivial skew-product cocycle
\begin{equation}
\label{eq:distcocycle}
 G_t \times\text{\rm id} :  \bigcup_{q\in  \Cal Q_\kappa(M)}  \{q\} \,\times H_q^{-s}(M)  \to 
\bigcup_{q\in  \Cal Q_{\kappa}(M)}  \{q\} \times H_q^{-s}(M) \,\,.
\end{equation}
Let $\{\Phi^s_t\}_{t\in\R}$ be the cocycle over the Teichm\"uller flow, defined as the projection onto the bundle ${W}_{\kappa}^{-s}(M)$ of the trivial skew-product cocycle: 
\begin{equation}
G_t\times \hbox{id} :  \bigcup_{q\in  \Cal Q_\kappa(M)}  \{q\} \,\times W_q^{-s}(M)  \to 
\bigcup_{q\in  \Cal Q_{\kappa}(M)}  \{q\} \times W_q^{-s}(M) \,\,.
\end{equation} 
Such cocycles can be described as the cocycles obtained by parallel transport of distributions and 
$1$-currents with respect to the trivial connection along the orbits of the Teichm\"uller flow. The
following immediate identity allows to express the cocycle $\{\Phi^s_t\}_{t\in\R}$ in terms of the 
cocycle $\{G^s_t\}_{t\in\R}$:
\begin{equation}
\label{eq:currcocycle}
\Phi^s_t:= \hbox{diag}(e^{-t},e^t)\otimes G^s_t \quad \text{ \rm on } \quad {W}_{\kappa}^{-s}(M)=
{\R}^2\otimes H^{-s}_{\kappa}(M)\,.
\end{equation}

\begin{lemma} 
\label{lemma:welldef} For any $s\in \R$, the spaces $H^{s}_\kappa(M)$,  $\bar H^{s}_\kappa(M)$, 
$W^{s}_\kappa(M)$ and $\bar W^{s}_\kappa(M)$ are well-defined Hilbert bundles over a stratum 
$\Cal M^{(1)}_\kappa$ of orientable quadratic differentials. The flows $\{G_t^s\}_{t\in \R}$ and $\{\Phi_t^s\}_{t\in \R}$ are well-defined smooth cocycles on the bundles 
$H^{-s}_\kappa(M)$ and $W^{-s}_\kappa(M)$ respectively.
\end{lemma}
\begin{proof} Let $q_0\in \Cal Q_\kappa(M)$. There exists a neighbourhood $S_0\subset M$ of the
set $\Sigma_{q_0}$ and a neighbourhood $D_0\subset \Cal Q_\kappa(M)$ in the space of quadratic differentials with zeros of multiplicities $\kappa=(k_1, \dots, k_\sigma)$ such that, for all $q\in D_0$, $\Sigma_q\subset S_0$ and the quadratic differential $q$ is isotopic to $q_0$ on $S_0$. Thus there exists a smooth map $f: D_0 \to \text{\rm Diff}^+(M)$  such that $q=f^\ast_q (q_0)$ on $S_0$.  The bundle $H^{s}_\kappa(M)$ and  $\bar H^{s}_\kappa(M)$ are trivialized over $D_0$ 
by the map 
\begin{equation}
\label{eq:triv}
\begin{aligned}
\bar H^{s}_\kappa(M) \vert D_0   \,\, &\to    
\,\,   D_0 \times \bar H^{s}_{q_0}(M); \\
(q, \Cal D)  \,\, &\to \,\,  \left( q, (f^{-1}_q)^\ast (\Cal D)\right),
\end{aligned}
\end{equation}
which maps the subspace $H^{s}_\kappa(M) \vert D_0$ onto $ D_0 \times H^{s}_{q_0}(M)$. 
 It follows that $H^{s}_\kappa(M)$ and $\bar H^{s}_\kappa(M)$ are well-defined Hilbert bundles, hence so are $W^{s}_\kappa(M)$ and  $W^{s}_\kappa(M)$ by formula \pref{eq:tensor}. The dynamical system $\{G^s_t\}_{t\in \R}$ concides  with the product cocycle $\{G_t \times  \text{\rm id }\}$ on  $D_0 \times H^{-s}_{q_0}(M)$ with respect to the trivialization \pref{eq:triv}, hence $\{G_t^s\}_{t\in \R}$ and $\{\Phi_t^s\}_{t\in \R}$ are well-defined smooth cocycles. In fact, the weighted Sobolev spaces
 $H^s_q(M)$ and $W^s_q(M)$ are invariant under the action of the Teichm\"uller flow on the space
 $\Cal Q_\kappa(M)$ (although their Hilbert structure is not). 
\end{proof}
\noindent We point out that for any $s>0$  there is no natural extension of the distributional cocycles $\{G_t^s\}_{t\in \R}$ and $\{\Phi_t^s\}_{t\in \R}$ respectively to the bundles $\bar H^{-s}_\kappa(M)$ and $\bar W^{-s}_\kappa(M)$, since the action of the Teichm\"uller flow on $\Cal Q_\kappa$ does not respect the domain of the Friedrichs Laplacian, hence the trivial cocycle $\{G_t  \times  \text{\rm id }\}$ is not defined on the bundles $D_0 \times \bar H^{-s}_{q_0}(M)$.

\medskip
\noindent Let ${\Cal I}^s_{\kappa,+}(M) [ {\Cal I}^s_{\kappa,-}(M)] \subset {H}^{-s}_{\kappa}(M)$ be the finite dimensional sub-bundle of {\it horizontally [vertically]  invariant distributions}. By definition, the fiber of the bundle ${\Cal I}^s_{\kappa,+}(M)$ [${\Cal I}^s_{\kappa,-}(M)$] at any $q\in {\Cal M}^{(1)}_{\kappa}$ coincides with the vector space  ${\Cal I}^s_{+q}(M)$ [${\Cal I}^s_{-q}(M)$] of horizontally [vertically]  invariant distributions.

\begin{lemma}
\label{lemma:invdistsub} For any $s\geq 0$ and for any $G_t$-invariant ergodic probability measure 
on a stratum $\Cal M^{(1)}_\kappa \subset \Cal M^{(1)}_g$ of orientable quadratic differentials,
 \begin{enumerate}
\item the sub-bundles ${\Cal I}^s_{\kappa,\pm}(M)\subset H^{-s}_{\kappa}(M)$ are 
$G^s_t$-invariant, measurable and of finite, almost everywhere constant rank;
\item the cocycle $\{G^s_t \vert {\Cal I}^s_{\kappa,\pm}(M)\}$ satisfies the Oseledec's theorem.
\end{enumerate}
\end{lemma}
\begin{proof}
\smallskip
\noindent $(1)\,$ The horizontal and vertical  vector fields determined by a quadratic
differential $q \in \Cal Q_\kappa(M)$ are rescaled under the Teichm\"uller geodesic flow and 
the space $\Cal S_q(M)$ of  $q$-tempered currents is invariant. Hence the spaces 
${\Cal I}_{\pm q}(M)$ of all horizontally, respectively vertically, invariant distributions are invariant. 
Since the Sobolev spaces $H^{-s}_q(M)$ are invariant, the spaces ${\Cal I}^s_{\pm q}(M)$ are also invariant under the Teichm\"uller geodesic flow.

\smallskip 
\noindent  For any $q\in \Cal Q_\kappa(M)$, the space  $\Cal I^{-s}_{+ q} (M)$ [$\Cal I^{-s}_{- q} (M)$] can be characterized as the perpendicular of the closure in $H^{s}_q(M)$ of the range of the Lie derivative $\Cal L_{S_q}$ [$\Cal L_{T_q}$]  as a linear operator defined on the space $H^{s+1}_q(M)$:
\begin{equation}
\begin{aligned}
\Cal I^{-s}_{q} (M) &=
\overline{  \{ S_q u\,  \vert \, u\in H^{s+1}_q(M) \} }^\perp \,; \\
[\Cal I^{-s}_{-q} (M) &=
\overline{  \{ T_q u\,  \vert \, u\in H^{s+1}_q(M)\} }^\perp] \,.
\end{aligned}
\end{equation}
Let $q_0\in \Cal Q_\kappa(M)$ and let $f: D_0 \to \text{\rm Diff}^+(M)$ be a map defined on a neighbourhood $D_0\subset \Cal Q_\kappa(M)$ of $q_0$ which trivializes the bundles 
$H^{-s}_\kappa(M)$ as in formula \pref{eq:triv}. For any fixed $v\in H^{s+1}_{q_0}(M)$ the maps
\begin{equation}
\label{eq:Liedermaps}
q\to  (f_q^{-1})^\ast   S_q f_q^\ast (v) \quad \text{ \rm and } \quad  
q\to  (f_q^{-1})^\ast  T_q f_q^\ast (v) 
\end{equation}
are well-defined and continuous on $D_0$ with values in $H^{s}_{q_0}(M)$. Since the Hilbert 
spaces $H^s_q(M)$ are separable, the sub-spaces 
$$(f_q^{-1})^\ast  \Cal I^{-s}_{\pm q} (M)\subset H^{-s}_{q_0}(M)$$ 
 are measurable functions of the quadratic differentials $q\in \Cal Q_\kappa(M)$. In fact, the orthogonal
 projections $I_{\pm q}$ on $(f_q^{-1})^\ast  \Cal I^{-s}_{\pm q} (M)$ can be constructed as follows.
 Let $B_0^{s+1}$ be an orthonormal basis for $H^{s+1}_{q_0}(M)$. If the horizontal [vertical] foliation of  $q\in D_0$ is minimal, the subset of $H^s_q(M)$ defined as
 \begin{equation}
 B_q^{s}:=\{ S_q f_q^\ast (v) \vert v\in B_0^{s+1}\} 
  \quad [ \,  B_{-q}^{s}:=\{ T_q f_q^\ast (v) \vert v\in B_0^{s+1}\} \, ] 
  \end{equation}
  is linerarly independent. Since the set of $q\in D_0$ with non-minimal horizontal or vertical foliation has measure zero, the Gram-Schmidt orthonormalization algorithm applied to the system  $(f_q^{-1})^\ast B^{s}_q$ [ $(f_q^{-1})^\ast B^{s}_{-q}$] yields an orthonormal basis $\{ u^+_k(q) \}_{k\in \N}$ [$\{ u^-_k(q) \}_{k\in \N}$ ]  in $H^{s}_{q_0}(M)$ of the  subspace 
   \begin{equation}
(f_q^{-1})^\ast \overline{  \{ S_q v\,  \vert \, v\in H^{s+1}_q(M) \} } \quad
 [  \,(f_q^{-1})^\ast \overline{  \{ T_q v\,  \vert \, v\in H^{s+1}_q(M) \} }\, ]
  \end{equation}
such that, for all $k\in \N$, the functions $u^\pm _k: D_0 \to H^s_{q_0}(M)$ are defined $\mu$-almost everywhere and are continuous on their domain of definition by the continuity of the maps \pref{eq:Liedermaps}.  
Let  $\{ \Cal D^\pm_k(q) \}_{k\in \N} \subset H^{-s}_{q_0}(M)$ be the (orthonormal) system dual to the orthonormal system $\{ u^\pm_k(q) \}_{k\in \N} \subset H^s_{q_0}(M)$. Then we can write
 \begin{equation}
 \label{eq:projections}
 I_{\pm q} (\Cal D) =  \Cal D \, - \, \sum_{k\in \N} \<\Cal D, \Cal D^\pm_k(q)\>_{\!-s} \,\Cal D^\pm_k(q)\,,
 \,\, \text{ \rm for all } \, \Cal D \in H^{-s}_{q_0}(M).
\end{equation}
It is immediate to verify that formula \pref{eq:projections} yields  the orthogonal projections onto the
subspaces $(f_q^{-1})^\ast  \Cal I^{-s}_{\pm q} (M)$ for any $q\in D_0$. Such projections
 can therefore be obtained as a limit of operators which depend continuously on $q\in \Cal Q_\kappa(M)$ in the complement of the subsets of quadratic differential with non-minimal horizontal
 [vertical] foliation. Since such a set has measure zero with respect to any $G_t$-invariant ergodic
 probability measure, the measurability of the sub-bundles of invariant distributions is proved.
 
 \smallskip
 \noindent The rank of the bundles ${\Cal I}^s_{\kappa,\pm}(M)$ is finite by Lemma \ref{lemma:DtoC} and Theorem \ref{thm:bcstruct} and it is almost everywhere constant with respect to any {\it ergodic }$G_t$-invariant probability measure by definition of ergodicity.
 
 \smallskip
 \noindent $(2)\,$ It follows from the Definition \ref{def:snorm} of the weighted Sobolev norms and from the fundamental theorem of interpolation (see \cite{LiMa68}, Chap. 1, \S 5.1) applied to the interpolation family of Hilbert spaces $\{H^s_q(M)\vert s\in [0,1]\}$ that, for any $s\in \R$ and for any $q\in \Cal Q_\kappa(M)$, the following estimates hold:
 \begin{equation}
 \label{eq:cocyclenormsone}
 \Vert G^s_t :   H^{-s}_q(M) \to H^{-s}_{G_t(q)}(M) \Vert  \,\, \leq \,\, e^{\vert s\vert \,\vert t \vert}  \,.
  \end{equation}
Since the bundles ${\Cal I}^s_{\kappa,\pm}(M)$ are finite dimensional and measurable,
the Oseledec's theorem applies to the cocycle $\{G^s_t\vert {\Cal I}^s_{\kappa,\pm}(M)\}$
with respect to any ergodic $G_t$-invariant probability measure.

\end{proof}
\noindent The measurable dependence of the spaces of invariant distributions on the quadratic
differentials can also holds for distributional solutions of the cohomological equation with
arbitrary data. In order to formulate the result, which will be relevant later on, we introduce the
following:
\begin{definition}  For any $s \in \R$ and any $r>0$, the range (in $ \bar H^{-s}_q(M) $) of the horizontal [vertical] Lie derivative operator (on $ \bar H^{-r}_q(M)$) is the subspace
\begin{equation}
\begin {aligned}
R^{r,s}_q(M)&:= \{ S_q U \in \bar H^{-s}_q(M) \, \vert U\in \bar H^{-r}_q(M) \\
 [ R^{r,s}_{-q}(M)&:= \{ T_q U \in \bar H^{-s}_q(M) \,\vert U\in \bar H^{-r}_q(M)\}]
\end{aligned}
\end{equation} 
The \emph{Green operators }$\Cal U^{r,s}_{\pm q} : R^{r,s}_{\pm q}(M) \to H^{-r}_q(M)$ are defined
 as follows: for any distribution $F\in R^{r,s}_q(M)$ [ $F\in R^{r,s}_{-q}(M)$] :
\begin{equation}
\begin {aligned}
\Cal U^{r,s}_{q} (F) &:=U \,\,,  \quad \text{ if } \,\, S_qU =F  \text{ and } \, U\in \Cal I^r_q(M)^\perp \subset
\bar H^{-r}_q(M) \,; \\
[\Cal U^{r,s}_{-q} (F) &:=U \,\,,  \quad \text{ if } \,\, T_qU =F  \text{ and }  \, U\in \Cal I^r_{-q}(M)^\perp \subset
\bar H^{-r}_q(M)]\,.
\end{aligned}
\end{equation}
\end{definition}
\noindent By Corollary \ref{cor:CEdistribution}, for any $s\in \R$ there exists $r>0$ 
such that, for almost all quadratic differentials in every circle orbit in a stratum 
$\Cal M^{(1)}_\kappa$, the Green operators are {\it bounded }operators defined on the codimension
$1$ subspace $\bar {\Cal H}_q^{-s}(M) \subset \bar H^{-s}_q(M)$ of distributions vanishing  on constant
functions. 

\smallskip
\noindent For any $s\in \R$, let $\bar {\Cal H}_\kappa^s(M) \subset \bar H^s_\kappa(M)$ be the continuous sub-bundle of distributions vanishing on constant functions:
$$
\bar {\Cal H}_\kappa^s(M):= \{ (q, \Cal D) \in \bar H^s_\kappa(M)\, \vert \,  \<\Cal D,1\>_s =0  \} \,.
$$

\begin{lemma} 
\label{lemma:Greenmeas}
Let $\mu$ be any  $SO(2, \R)$-absolutely continous, $G_t$-invariant ergodic probability measure on a stratum $\Cal M^{(1)}_\kappa$ of orientable quadratic differentials. For any $s\in \R$, there exists
$r>0$ such that the  the Green operator $\,\Cal U^{r,s}_{\pm q}\,$ yield a measurable bundle map $\Cal U^{r,s}_{\kappa, \pm}: \bar {\Cal H}^{-s}_\kappa(M) \to H^{-r}_\kappa(M)$. In particular, the operator norm
 $\Vert \Cal U^{r,s}_{\pm q} \Vert$ yields a well-defined measurable real-valued function of the quadratic differential $q\in \Cal M^{(1)}_\kappa$.
\end{lemma}
\begin{proof} 
Let $q_0\in \Cal Q_\kappa(M)$ and let $f: D_0 \to \text{\rm Diff}^+(M)$ be a map defined on a neighbourhood $D_0\subset \Cal Q_\kappa(M)$ of $q_0$ which trivializes the bundles 
$\bar H^{-s}_\kappa(M)$ as in formula \pref{eq:triv}. The argument is similar to the one given in
the proof of the measurability of the sub-bundles of invariant distributions in Lemma \ref{lemma:invdistsub}. Let $B_0^{r+1} \subset H^{\infty}_{q_0}(M)$ be a basis for the Hilbert space $H^{r+1}_{q_0}(M)$. The subset of $H^r_q(M)$ defined as
 \begin{equation}
B_q^{r}:=\{ S_q f_q^\ast (v) \vert v\in B_0^{r+1}\} 
  \quad [ \,  B_{-q}^{r}:=\{ T_q f_q^\ast (v) \vert v\in B_0^{r+1}\} \, ] 
  \end{equation}
  is linerarly independent if the horizontal [vertical] foliation of  $q\in D_0$ is minimal, Hence, by the continuity of the maps \pref{eq:Liedermaps}, the Gram-Schmidt orthonormalization applied to the system  $ B^{r}_q$ [ $B^{r}_{-q}$] yields an orthonormal basis $\{ u^+_k(q) \}_{k\in \N}$ [$\{ u^-_k(q) \}_{k\in \N}$ ]  in $H^{r}_{q}(M)$ of the subspace 
   \begin{equation}
 \overline{  \{ S_q v\,  \vert \, v\in H^{r+1}_q(M) \} } \quad
 [ \overline{  \{ T_q v\,  \vert \, v\in H^{r+1}_q(M) \} }\, ]
  \end{equation}
   such that, for all $k\in \N$, the functions $(f^{-1})^\ast \circ u^\pm _k: D_0 \to H^r_{q_0}(M)$ are defined $\mu$-almost everywhere and continuous on their domain of definition.  In fact, there exists bases $\{v^\pm_k\} _{k\in \N} \subset H^{\infty}_{q_0}(M)$ of the Hilbert space $H^{r+1}_{q_0}(M)$ such that the following holds. For each $k\in\N$ there exists a function $v^\pm_k: D_0 \to \text{\rm span} (v^\pm_1 ,\dots, v^\pm_k)$, defined $\mu$-almost everywhere and continuous on its domain of definition, such that for $\mu$-almost all $q\in D_0$ the set $\{v^\pm_k(q)\} _{k\in \N}$ is a basis of the Hilbert space $H^{r+1}_{q_0}(M)$ and
$$
u^+_k(q) =    S_q (f^\ast \circ v^+_k)(q)   \quad
[ u^-_k(q) =   T_q (f^\ast \circ  v^-_k)(q) ]\,.
$$
Let  $\{ \Cal D^\pm_k(q) \}_{k\in \N} \subset H^{-r}_{q_0}(M)$ be the system dual to the linearly independent system $\{ (f_q^{-1})^\ast  u^\pm_k(q) \}_{k\in \N} \subset H^r_{q_0}(M)$. Such systems 
are in general not orthonormal. Then for any $F \in \bar {\Cal H}^{-s}_{q_0}(M)$, the following formula holds: 
\begin{equation}
\label{eq:greenseries}
(f_q^{-1})^\ast \circ \Cal U^{r,s}_{\pm q}\circ f_q^\ast  (F)  \,:= \,
 \sum_{k\in \N}   \<F, v^\pm_k(q)\> \, \Cal D^\pm_k(q) \,.
\end{equation}
In fact, let $U^\pm_k(q)$ be the series on the right hand side of formula \pref{eq:greenseries}. Since $\<f_q^\ast \Cal D^\pm_k(q), u^\pm_k(q)  \> = \delta_{kh}$ for any $k$, $h\in \N$, it follows that
\begin{equation}
\begin{aligned}
\<f_q^\ast U_k^+(q),  S_q f^\ast_q v^+_k(q) \> = \<F, v^+_k(q)\> 
=  \<f_q^\ast (F), f^\ast_q v^+_k(q)\>  \\
[\<f_q^\ast U_k^-(q),  T_q f^\ast_q v^-_k(q) \> = \<F, v^-_k(q)\> 
=  \<f_q^\ast (F), f^\ast_qv^-_k(q)\>  ]
\end{aligned}
\end{equation}
which implies that  $f_q^\ast U_k^+ (q)$ [$f_q^\ast U_k^- (q)$] is a distributional solution of the 
equation $S_q u=f_q^\ast (F)$ [$T_q u =f_q^\ast (F)$]. Finally, since all the distributions $f_q^\ast \Cal D^+_k(q)$ [$f_q^\ast \Cal D^-_k(q)$] vanish on the orthogonal complement of the space  $\{ S_q v\,  \vert \, v\in H^{r+1}_q(M) \}$ [$\{ T_q v\,  \vert \, v\in H^{r+1}_q(M) \}$] in $H^r_q(M)$, it follows that $f_q^\ast U_k^\pm (q)$ is orthogonal to the space $\Cal I^r_{\pm q}(M)$ of invariant distributions,
hence $\Cal U^{r,s}_{\pm q}\circ f_q^\ast  (F)=f_q^\ast U_k^\pm (q)$. It is also immediate from formula
\pref{eq:greenseries} that,  for any $F\in \bar {\Cal H}^{-s}_{q_0}(M)$, the distribution $(f_q^{-1})^\ast \circ \Cal U^{r,s}_{\pm q}\circ f_q^\ast  (F)$ is a $\mu$-measurable function (defined almost everywhere) of the quadratic differential $q\in D_0$ with values in the Hilbert space $H^{-r}_{q_0}(M)$. In fact, by construction the functions $\Cal D^\pm_k :D_0 \to H^{-r}_{q_0}(M)$ are defined $\mu$-almost everywhere and continuous on their domain of definition.

\end{proof}

\noindent Let ${\Cal Z}^s_{\kappa}(M)\subset {W}^{-s}_{\kappa}(M)$ be the infinite dimensional 
sub-bundle of {\it closed currents} over ${\Cal M}^{(1)}_{\kappa}$. By definition, the fiber of the bundle ${\Cal Z}^s_{\kappa}(M)$ at any $q\in {\Cal M}^{(1)}_{\kappa}$ coincides with the vector space of 
{\it closed currents}:
 \begin{equation}
 \label{eq:Hsclosedcurr}
{\Cal Z}^s_q(M):={\Cal Z}_q(M) \cap  {W}^{-s}_q(M)\,.
\end{equation}
The bundle ${\Cal Z}^s_{\kappa}(M)$ and the sub-bundle $\Cal E^s_{\kappa}(M) \subset 
{\Cal Z}^s_{\kappa}(M)$ of {\it exact }currents are smooth, $\Phi^s_t$-invariant sub-bundles of
the bundle $W^{-s}_q(M)$. The quotient cocycle, defined on the $H^{-s}$ de Rham cohomology 
bundle, is isomorphic to the Kontsevich-Zorich cocycle. The latter isomorphism is the essential motivation for the formulas \pref{eq:distcocycle} and \pref{eq:currcocycle} which define, respectively, 
the cocycles $\{G^s_t\}_{t\in\R}$ and $\{\Phi^s_t\}_{t\in\R}$. In fact, let
\begin{equation}
\label{eq:cohobundlemap}
j_{\kappa}:{\Cal Z}^s_{\kappa}(M)\to {\Cal H}^1_{\kappa}(M,{\R})
\end{equation}
be the natural de Rham cohomology map onto the cohomology bundle ${\Cal H}^1_{\kappa}(M,{\R})$, defined as the restriction to the stratum $\Cal M^{(1)}_\kappa$ of the cohomology bundle ${\Cal H}^1_g(M,{\R})$ introduced in formula \pref{eq:CB}. 

\smallskip
\noindent Let ${\Cal B}^s_{\kappa,\pm}(M)\subset {\Cal Z}^s_{\kappa}(M)$ be the sub-bundles with 
fiber at $q\in {\Cal M}^{(1)}_{\kappa}$ given by the vector spaces ${\Cal B}^s_{\pm q}(M)$ of ${\Cal F}_{\pm q}$-basic currents (defined in \pref{eq:Hsbcurr}). 

\begin{lemma}
\label{lemma:bcurrsub} For any $s\geq 0$ and for any $G_t$-invariant ergodic probability measure 
on a stratum $\Cal M^{(1)}_\kappa \subset \Cal M^{(1)}_g$ of orientable quadratic differentials,
 \begin{enumerate}
\item the identity $j_{\kappa}\circ \Phi^s_t=\Phi_t \circ j_{\kappa}$ holds everywhere
on ${\Cal Z}^s_{\kappa}(M)$;
\item the sub-bundles ${\Cal B}^s_{\kappa,\pm}(M)\subset {\Cal Z}^s_{\kappa}(M)$ are 
$\Phi^s_t$-invariant, measurable and of finite, almost everywhere constant rank;
\item the cocycle $\{\Phi^s_t \vert {\Cal B}^s_{\kappa,\pm}(M)\}$ satisfies the Oseledec's theorem.
\end{enumerate}
\end{lemma}
\begin{proof}
 $(1)\,$ It is an immediate consequence of the definition of the cocycle $\{\Phi^s_t\}_{t\in \R}$
on the bundle ${\Cal Z}^s_{\kappa}(M)$ of closed currents.

\smallskip
\noindent $(2)\,$ The horizontal and vertical measured foliations of a quadratic differential $q \in \Cal Q_\kappa(M)$ are projectively invariant  under the Teichm\"uller geodesic flow and the space $\Cal S_q(M)$ of  $q$-tempered currents is invariant. As a consequence, the spaces ${\Cal B}_{+ q}(M)$ 
and ${\Cal B}_{- q}(M)$  of horizontally, respectively vertically, basic currents are invariant. 
The Sobolev spaces $W^{-s}_q(M)$ are also invariant (although the Hilbert structure is not). Hence
the spaces ${\Cal B}^s_{\pm q}(M)$ are invariant under the Teichm\"uller geodesic flow.

\smallskip 
\noindent By Lemma \ref{lemma:DtoC}, the measurability of the bundles ${\Cal B}^s_{\kappa,\pm }(M)\subset {W}^{-s}_{\kappa}(M)$ of basic currents is equivalent to the measurability of the bundles
${\Cal I}^s_{\kappa,\pm }(M)\subset {H}^{-s}_q(M)$, proved in Lemma \ref{lemma:invdistsub}.

  \smallskip
 \noindent The rank of the bundles ${\Cal B}^s_{\kappa,\pm}(M)$ is finite by Theorem \ref{thm:bcstruct}    and it is almost everywhere constant with respect to any {\it ergodic }$G_t$-invariant probability measure
 by definition of ergodicity.
 
 \smallskip
 \noindent $(3)\,$ It follows immediately from the identities \pref{eq:currcocycle} and from the bound \pref{eq:cocyclenormsone}  that, for any $s\in \R$ and for any $q\in \Cal Q_\kappa(M)$, the following estimates hold:
 \begin{equation}
 \label{eq:cocyclenormstwo}
 \Vert \Phi^s_t :   W^{-s}_q(M) \to W^{-s}_{G_t(q)}(M) \Vert \,\, \leq \,\, e^{(\vert s\vert +1)\vert t \vert}\,.
  \end{equation}
Since the bundles ${\Cal B}^s_{\kappa,\pm}(M)$ are finite dimensional and measurable,
the Oseledec's theorem applies to the cocycle $\{\Phi^s_t\vert {\Cal B}^s_{\kappa,\pm}(M)\}$
with respect to any ergodic $G_t$-invariant probability measure.
\end{proof}

\smallskip
\noindent Lemma \ref{lemma:invdistsub} can be generalized to the bundle of quasi-invariant distributions. Let $\Cal I^s_{\kappa,\pm} (M\setminus\Sigma_\kappa) \subset H^{-s}_{\kappa}(M)$ be the bundle over $\Cal M^{(1)}_\kappa$ defined as follows: its fiber at each $q\in \Cal M^{(1)}_\kappa$ is the vector space $\Cal I^s_{\pm q}(M\setminus\Sigma_q)\subset H^{-s}_q(M)$ of quasi-invariant distributions. 
An argument analogous to the proof of Lemma \ref{lemma:invdistsub} proves the following:
\begin{lemma}
\label{lemma:qinvdistsub} For any $s\geq 0$ and for any $G_t$-invariant ergodic probability measure on a stratum $\Cal M^{(1)}_\kappa \subset \Cal M^{(1)}_g$ of orientable quadratic differentials,
 \begin{enumerate}
\item the sub-bundles ${\Cal I}^s_{\kappa,\pm}(M\setminus\Sigma_\kappa)\subset H^{-s}_{\kappa}(M)$ are $G^s_t$-invariant, measurable and of finite, almost everywhere constant rank;
\item the cocycle $\{G^s_t \vert {\Cal I}^s_{\kappa,\pm}(M\setminus\Sigma_\kappa)\}$ satisfies the Oseledec's theorem.
\end{enumerate}
\end{lemma}

\noindent Lemma \ref{lemma:bcurrsub} can be generalized to the bundle of currents closed on the complement of the singular set and to the sub-bundle of quasi-basic currents. Let $\Cal Z^s_\kappa (M\setminus\Sigma_\kappa) \subset W^{-s}_{\kappa}(M)$ be the bundle over $\Cal M^{(1)}_\kappa$ defined as follows: its fiber at each $q\in \Cal M^{(1)}_\kappa$ is the vector space $\Cal Z^s(M\setminus\Sigma_q)\subset W^{-s}_q(M)$ of closed currents on the open manifold $M\setminus\Sigma_q$. The bundle ${\Cal Z}^s_{\kappa} (M\setminus\Sigma_\kappa)$ and the sub-bundle $\Cal E^s_{\kappa} (M\setminus\Sigma_\kappa) \subset {\Cal Z}^s_{\kappa}(M\setminus\Sigma_\kappa)$ of {\it exact }currents are smooth, $\Phi^s_t$-invariant sub-bundles of the bundle $W^{-s}_q(M)$. The quotient cocycle, defined on the $H^{-s}$ de Rham punctured cohomology bundle, is isomorphic to the punctured Kontsevich-Zorich cocycle. Let
\begin{equation}
\label{eq:pcohobundlemap}
j_{\kappa}:{\Cal Z}^s_{\kappa}(M\setminus\Sigma_\kappa)\to {\Cal H}^1_{\kappa}(M\setminus\Sigma_\kappa,{\R})
\end{equation}
be the natural de Rham cohomology map onto the punctured cohomology bundle ${\Cal H}^1_{\kappa}(M\setminus\Sigma_\kappa,{\R})$ introduced in formula \pref{eq:PCB}. 

\smallskip
\noindent Let ${\Cal B}^s_{\kappa,\pm}(M\setminus\Sigma_\kappa)\subset {\Cal Z}^s_{\kappa}(M\setminus\Sigma_\kappa)$ be the sub-bundles with 
fiber at $q\in {\Cal M}^{(1)}_{\kappa}$ given by the vector spaces ${\Cal B}^s_{\pm q}(M\setminus\Sigma_q)$ of quasi-basic currents for the measured foliation 
$\Cal F_{\pm q}$ (defined in \pref{eq:Hsbcurr}).  An argument analogous to the proof
of Lemma \ref{lemma:bcurrsub} proves the following:

\begin{lemma}
\label{lemma:qbcurrsub} For any $s\geq 0$ and for any $G_t$-invariant ergodic probability measure 
on a stratum $\Cal M^{(1)}_\kappa \subset \Cal M^{(1)}_g$ of orientable quadratic differentials,
 \begin{enumerate}
\item the identity $j_{\kappa}\circ \Phi^s_t=\Psi_t \circ j_{\kappa}$ holds everywhere
on ${\Cal Z}^s_{\kappa}(M\setminus\Sigma_\kappa)$;
\item the sub-bundles ${\Cal B}^s_{\kappa,\pm}(M\setminus\Sigma_\kappa)\subset {\Cal Z}^s_{\kappa}(M\setminus\Sigma_\kappa)$ are $\Phi^s_t$-invariant, measurable and of finite, almost everywhere constant rank;
\item the cocycle $\{\Phi^s_t \vert {\Cal B}^s_{\kappa,\pm}(M\setminus\Sigma_\kappa)\}_{t\in\R}$ 
satisfies the Oseledec's theorem.
\end{enumerate}
\end{lemma}

\noindent It follows immediately from the definitions and from Lemma \ref{lemma:DtoC} that the following cocycle isomorphisms hold:
\begin{equation}
\label{eq:cocycleiso}
\begin{aligned}
\Phi^s_t \vert {\Cal B}^s_{\kappa,\pm}(M)  \,\, &\equiv \,\,     e^{\pm t} \, G^s_t  \vert {\Cal I}^s_{\kappa,\pm}(M) \,;\\
\Phi^s_t \vert {\Cal B}^s_{\kappa,\pm}(M\setminus\Sigma_\kappa)  \,\,&\equiv \,\,    e^{\pm t} \, G^s_t  \vert {\Cal I}^s_{\kappa,\pm}(M\setminus\Sigma_\kappa)\,.
\end{aligned}
\end{equation}
As a consequence of such an isomorphism parts $(2)$ and $(3)$ in Lemma \ref{lemma:bcurrsub}
and Lemma \ref{lemma:qbcurrsub} can be immediately derived from Lemma \ref{lemma:invdistsub}
and Lemma \ref{lemma:qinvdistsub} respectively. In addition, the Lyapunov spectra
and the Oseledec's decomposition of the cocycles $\{\Phi^s_t \vert {\Cal B}^s_{\kappa,\pm}(M)\}$ 
[$\{\Phi^s_t \vert {\Cal B}^s_{\kappa,\pm}(M\setminus\Sigma_\kappa)\}$] and $\{G^s_t  \vert {\Cal I}^s_{\kappa,\pm}(M)\}$  [$\{G^s_t  \vert {\Cal I}^s_{\kappa,\pm}(M\setminus\Sigma_\kappa)\}$] can be immediately derived from one another. By part $(1)$ in Lemma \ref{lemma:bcurrsub} and Lemma
\ref{lemma:qbcurrsub}, by Corollary \ref{cor:basiccohom} and  by the structure theorem for basic currents (Theorem \ref{thm:bcstruct}), informations on the Lyapunov spectrum of the cocycles $\{\Phi^s_t \vert {\Cal B}^s_{\kappa,\pm}(M\setminus\Sigma_\kappa)\}$ and $\{\Phi^s_t \vert {\Cal B}^s_{\kappa,\pm}(M)\}$ can be derived from that of the (punctured) Kontsevich-Zorich cocycle.

\medskip
\noindent For $s=1$, by  \cite{F02}, Lemma 8.1, and by the representation Lemma \ref{lemma:reglowbound}, in the non-uniformly hyperbolic case it is possible to give a quite complete description 
of the distributional cocycle on the bundle of basic currents.

\begin{lemma} 
\label{lemma:B1}
Let $\mu$ be any $SO(2,\R)$-absolutely continuous, KZ-hyperbolic measure on a stratum ${\Cal M}^{(1)}_{\kappa}$ of orientable quadratic differentials.
\begin{enumerate}
\item The cocycle $\{\Phi^1_t\}_{t\in\R}$ has strictly positive [strictly negative]  Lyapunov spectrum,
with respect to the measure $\mu$ on  $\Cal M^{(1)}_\kappa$,
on the invariant sub-bundle ${\Cal B}^1_{\kappa,+}(M)$ [${\Cal B}^1_{\kappa,-}(M)$]. 
\item The sum ${\Cal B}^1_{\kappa}(M):={\Cal B}^1_{\kappa,+}(M) + {\Cal B}^1_{\kappa,-}(M)$ 
of bundles over $\Cal M^{(1)}_\kappa$ is direct and the restriction of the map $j_{\kappa}$ to the 
sub-bundle ${\Cal B}^1_{\kappa}(M)$  is $\mu$-almost everywhere injective.  
\item The coycle $\{\Phi^1_t \vert {\Cal B}^1_{\kappa}(M)\}_{t\in \R}$ is isomorphic to the Kontsevich-Zorich cocycle $\{\Phi_t\}_{t\in \R}$ on the real cohomology bundle $\Cal H^1_\kappa(M,\R)$, hence 
it has the same Lyapunov spectrum.
 \end{enumerate}
\end{lemma}

\subsection{Lyapunov exponents}
\label{ss:Lyapexp}

\noindent We prove below that the the spaces of all horizontally and vertically (quasi)-basic currents
and invariant distributions have well-defined Oseledec decompositions at almost all point of any stratum of the moduli space (with respect to any $G_t$-invariant ergodic measure). We establish a fundamental relation between the Lyapunov exponents and the Sobolev order of basic currents or distributions in
each Oseledec subspace. We conclude with a crucial `spectral gap' result for the distributional cocycle on the bundle of exact currents which is the basis for sharp estimates on the growth of ergodic averages,
hence for the construction of square-integrable solutions of the cohomological equation in \S \ref{ss:ergint}.

\smallskip
\noindent We introduce the following:
\begin{definition}
\label{eq:simple}
A non-zero current $C^\pm \in \Cal B^s_{\kappa,\pm}(M\setminus\Sigma_\kappa)$  will be 
called \emph{ (Oseledec) simple }if it belongs to an Oseledec subspace of the cocycle $\{\Phi_t^s \vert
B^s_{\kappa,\pm}(M\setminus\Sigma_\kappa)\}$. A non-zero invariant distribution$\Cal D^\pm \in \Cal I^s_\kappa(M\setminus\Sigma_\kappa)$ will be called \emph{ (Oseledec) simple }if it belongs to an Oseledec subspace of the cocycle $\{G_t^s \vert \Cal I^s_\kappa(M\setminus\Sigma_\kappa)\}$. 
\end{definition}

\noindent The above definition is well-posed since for any $r\leq s$ the cocycles $\{\Phi^r_t\}$
[$\{G^r_t\}$] are the restrictions of the cocycles $\{\Phi^s_t\}$ [$\{G^s_t\}$]. It is also immediate to prove that the image of any Oseledec simple basic current under the isomorphism $\Cal D_{\pm q} : \Cal B^s_{\pm q} (M\setminus \Sigma_q) \to \Cal I^s_{\pm q}(M\setminus \Sigma_q)$ introduced in formula \pref{eq:DtoC} is an Oseledec simple invariant distribution.

\begin{lemma} 
\label{lemma:currLE}
Let  $\mu$ be any $G_t$-invariant probability measure  on $\Cal M_\kappa^{(1)}$ and let 
$\Cal R_\mu \subset \Cal M_\kappa^{(1)}$ be the set of  all holomorphic orientable quadratic differentials which, for all $s\geq 0$, are Oseledec regular point for the cocycles $\{\Phi^s_t \vert \Cal B^s_{\kappa, +}(M\setminus\Sigma_\kappa)\}$ and $\{\Phi^s_t \vert \Cal B^s_{\kappa, -}(M\setminus\Sigma_\kappa)\}$ over the Teichm\"uller flow $(\{G_t\}, \mu)$. The set $\Cal R_\mu $ has full measure and there exist measurable functions $L^\pm_\mu : \Cal B_{\kappa, \pm}(M\setminus\Sigma_\kappa)  \to \R$ (with almost everywhere constant range) such that, 
for any $q\in \Cal R_\mu$, the number $L^\pm_\mu(C^\pm)$ is equal to the Lyapunov exponent of the Oseledec simple current $C^\pm\in \Cal B_{\pm q}(M\setminus\Sigma_q)$ of Sobolev order ${\Cal O}_q^W(C^\pm)\geq 0$ (see Definition \ref{def:currsobord}),  with respect to the cocycle $\{\Phi^s_t \vert  \Cal B^s_{\kappa, \pm} (M\setminus\Sigma_\kappa)\}$ over the flow $(\{G_t\}, \mu)$ for any $s>{\Cal O}_q^W(C^\pm)$.
\end{lemma}
\begin{proof}
For any $r \leq s$, the embeddings $ \Cal B^{r}_{\kappa, \pm}(M\setminus\Sigma_\kappa) 
\subset \Cal B^{s}_{\kappa, \pm}(M\setminus\Sigma_\kappa)$ hold and the cocycle
$\{\Phi^r_t\}_{t\in \R}$ on the bundle $W^{-r}_\kappa(M)$ coincides with restriction of
the cocycle $\{\Phi^s_t\}_{t\in \R}$, defined on the bundle $W^{-s}_\kappa(M)$. The Oseledec's theorem 
holds for the cocycles $\{\Phi^s_t \vert  \Cal B^{s}_{\kappa, \pm}(M\setminus\Sigma_\kappa)\}_{t\in \R}$
on the measurable, finite dimensional sub-bundles $\Cal B^{s}_{\kappa, \pm}(M\setminus\Sigma_\kappa)\subset W^{-s}_\kappa(M)$ for all $s\geq 0$. Let $\Cal R_\mu
\subset \Cal M^{(1)}_\kappa$ be the set of points which are Oseledec regular for the cocycles $\{\Phi^k_t \vert  \Cal B^{k}_{\kappa, \pm}(M\setminus\Sigma_\kappa)\}_{t\in \R}$ for all $k\in \N$. By the Oseledec theorem the set $\Cal R_\mu$ has full measure and, for any $q\in \Cal R_\mu$ and for all $s\geq r$, the Lyapunov exponents of any Oseledec simple current  $C^\pm\in \Cal B^{r}_{\kappa, \pm}(M\setminus\Sigma_\kappa)$ with respect to the cocycles $\{\Phi^r_t\}_{t\in \R}$ and $\{\Phi^s_t\}_{t\in \R}$ are well-defined and coincide. The common value $L^{\pm}_\mu(C^\pm)$ of all Lyapunov exponents with respect to the cocycles $\{\Phi^s_t\}_{t\in \R}$ for  $s>{\Cal O}_q^W(C^\pm)$ is therefore well-defined.
\end{proof}

\noindent By the isomorphisms \pref{eq:cocycleiso} an analogous statement holds for the cocycles on the bundles of invariant distributions:
\begin{lemma} 
\label{lemma:distLE}
Let  $\mu$ be a $G_t$-invariant probability measure  on $\Cal M_\kappa^{(1)}$. There exist measurable functions $l^{\pm}_\mu: \Cal I_{\kappa, \pm}(M\setminus\Sigma_\kappa)  \to \R$ (with almost everywhere constant range) such that, for any quadratic differential $q\in \Cal R_\mu$, the number $l^{\pm}_\mu(\Cal D^\pm)$ is equal to the Lyapunov exponent of the Oseledec simple invariant distribution $\Cal D^\pm\in \Cal I_{\pm q}(M\setminus\Sigma_q)$ of Sobolev order ${\Cal O}_q^H(\Cal D^\pm)\geq 0$ (see Definition \ref{def:distsobord}),  with respect to the cocycle $\{G^s_t \vert  \Cal I^s_{\kappa, \pm} (M\setminus\Sigma_\kappa)\}$ over the flow $(\{G_t\}, \mu)$ for any $s>{\Cal O}_q^H(\Cal D^\pm)$.  Let $\Cal D_{\pm q} : \Cal B_{\pm q} (M\setminus \Sigma_q) \to \Cal I_{\pm q}(M\setminus \Sigma_q)$ be the isomorphism introduced in  formula \pref{eq:DtoC},  the following identities hold:
\begin{equation}
\label{lemma:currdistLEid}
l^{\pm}_\mu \circ \Cal D_{\pm q} \, =\,  L^\pm_\mu   \, \mp \, 1 \,, \quad \text {\rm on } \,\, \Cal B_{\pm q} (M\setminus \Sigma_q)\,.
\end{equation}
\end{lemma}

\noindent Lyapunov exponents of basic currents impose restrictions on their Sobolev regularity.
In fact, we have:

\begin{lemma} 
\label{lemma:regupperbound}
The following inequalities hold for any quadratic differential $q\in \Cal R_\mu\subset \Cal M^{(1)}_\kappa$. For any Oseledec simple basic current $C^\pm \in \Cal B_{\pm q}(M\setminus\Sigma_q)$ and any Oseledec simple invariant distribution $\Cal D^\pm \in \Cal I_{\pm q}(M\setminus\Sigma_q)$ the Sobolev order functions satisfy the following lower bounds:
\begin{equation}
{\Cal O}_q^W(C^\pm) \geq  \vert L^\pm_\mu (C^\pm) \mp 1\vert \quad \text{ \rm and } \quad
{\Cal O}_q^H(\Cal D^\pm) \geq  \vert l^\pm_\mu (\Cal D^\pm)\vert\,.
\end{equation}
\end{lemma}

\begin{proof}
Since $C^\pm \in \Cal B_{\pm q}(M\setminus\Sigma_q)$, by the identities \pref{eq:currcocycle} and by the bound \pref{eq:cocyclenormsone} the following inequalities hold for
any $s\geq r$:
\begin{equation}
\label{eq:regnormbound}
\vert \Phi_t^s(C^\pm) \vert_{-s}\,  \leq\,\, \vert \Phi_t^r(C^\pm) \vert_{-r} \leq  e^{r\vert t \vert \pm t } 
\vert C^\pm \vert_{-r} \,, \,\, \text{ \rm for all} \,t\in\R  \,.
\end{equation}
On the other hand, by the Oseledec's theorem, for any $q\in \Cal R_\mu$ and any $\epsilon>0$
there exists a constant $K_\epsilon(q)>0$ such that
\begin{equation}
\label{eq:expnormbound}
\begin{aligned}
\vert \Phi_t^s(C^\pm) \vert_{-s}  \,\geq    &\, \,K_\epsilon(q) \, e ^{(L_\mu^\pm(C^\pm) -\epsilon)\, t } \, \vert 
C^\pm \vert_{-s} \,\,,
 \quad \text{ \rm for } \,t\geq 0\,; \\
\vert \Phi_t^s(C^\pm) \vert_{-s} \, \geq     &\, \,K_\epsilon(q) \, e ^{(L_\mu^\pm(C^\pm) +\epsilon)\, t}\, \vert 
C^\pm \vert_{-s} \,\,, \quad \text{ \rm for } \,t< 0\,;
\end{aligned}
\end{equation}
A comparison of the inequalities \pref{eq:regnormbound} and \pref{eq:expnormbound} proves
the statement for the case of basic currents. The statement for invariant distributions follows immediately
by Lemma \ref{lemma:distLE} since the following identity holds:
$$
{\Cal O}_q^H \circ \Cal D_{\pm q} = {\Cal O}_q^W  \quad \text {\rm on } \,\, \Cal B_{\pm q} (M\setminus \Sigma_q)\,.
$$
\end{proof}

\noindent A similar argument, based on the spectral gap for almost all quadratic differentials
in every circle orbit (see Theorem \ref{thm:spectralgaptwo}), proves the following:

\begin{theorem} 
\label{thm:BCgap}
For any quadratic differential $q\in \Cal M^{(1)}_\kappa$ there exists a real 
number $s(q)>0$ such that, for almost all $\theta\in S^1$ and for all $s<s(q)$
\begin{equation}
\label{eq:BCgap}
\begin{aligned}
\text{\rm dim} \,\Cal B^s_{q_\theta}(M\setminus\Sigma_q) &= \text{\rm dim}\, \Cal B^s_{-q_\theta}(M\setminus\Sigma_q) =1 \,;\\
\text{\rm dim} \,\Cal I^s_{q_\theta}(M\setminus\Sigma_q) &= \text{\rm dim}\, \Cal I^s_{-q_\theta}(M\setminus\Sigma_q) =1 \,.
\end{aligned}
\end{equation}
\end{theorem}
\begin{proof}
If there exists $s\in (0,1)$ such that \pref{eq:BCgap} does not hold on a positive measure subset of the 
circle, it follows that there exists a positive measure set $E_s\subset S^1$ such that $\text{\rm dim} \,\Cal B^s_{-q_\theta}(M\setminus\Sigma_q)>1$. Hence, for any $\theta \in E_s$, there exists a {\it non-vanishing }vertically (quasi)-basic current $C_\theta\in \Cal B^s_{-q_\theta}(M\setminus\Sigma_q)$ such that $C_\theta \wedge q^{1/2}_\theta=0$. We claim that for almost all $\theta\in E_s$, the current $C_\theta$ has non-zero cohomology class $c_\theta \in H^1(M\setminus\Sigma_q, \R)$. In fact, by Theorem \ref{thm:bcstruct}, if $C_\theta$ has vanishing cohomology class, then there exists a basic current $\hat C_\theta \in \Cal B^{s-1}_{-q_\theta} (M\setminus\Sigma_q)$ such that $\delta^1_{-q_\theta}(\hat C_\theta) =C_\theta$. Since $s<1$, if the vertical foliation $\Cal F_{q_\theta}$ is ergodic, then the current $\hat C_\theta \in \R\cdot \Re(q_\theta)$ which implies $C_\theta=0$, a contradiction. Since by the Keane conjecture (see \cite{Ma82} or \cite{Ve82}) the vertical foliation $\Cal F_{q_\theta}$ is uniquely ergodic for almost all $\theta\in S^1$, the claim is proved. As in formula \pref{eq:regnormbound}, we have
\begin{equation}
\label{eq:regnormboundbis}
\vert \Phi_t^s(C_\theta) \vert_{-s}  \leq   e^{s \vert t \vert - t } 
\vert C_\theta \vert_{-s} \,, \,\, \text{ \rm for all } \,t\in\R  \,.
\end{equation}
Let $\{\Psi_t\}_{t\in \R}$ be the punctured Kontsevich-Zorich cocycle, introduced in \S \ref{KZcocycle}.
By Lemma \ref{lemma:bcurrsub} the following identity holds: 
$$
[\Phi_t^s(C_\theta)] = \Psi_t(c_\theta) \in H^1(M\setminus\Sigma_q,\R) \,, \,\, \text{ \rm for all} \,t\in\R  \,.
$$
The de Rham punctured cohomology bundle $\Cal H^1_\kappa(M\setminus\Sigma_\kappa,\R)$ can be endowed, for any $s\geq 0$, with the quotient norm $\vert \cdot \vert^s_\kappa $ induced by the Sobolev norm on the bundle $W^{-s}_\kappa(M)$.  It follows by the above discussion and by the estimate \pref{eq:regnormboundbis} that, for each $\theta\in E_s$, there exists a non-vanishing cohomology class $c_\theta \in \Cal H^1_{q_\theta}(M\setminus\Sigma_q,\R)$ such that
\begin{equation}
\label{eq:regul2bound}
\begin{aligned}
&c_\theta \wedge [\Re(q^{1/2}_\theta)]= c_\theta \wedge [\Im(q^{1/2}_\theta)]=0 \,,\\
&\limsup_{t\to +\infty} \frac{1}{\vert t \vert }  \log \vert \Psi_t (c_\theta)\vert^s_\kappa  \,\, \leq \,\, s-1\,.
\end{aligned}
\end{equation}
By Lemma \ref{lemma:cocycleiso} and by the Poincar\'e-Lefschetz duality \pref{eq:PLD} 
between relative and punctured cohomology, the cohomology bundle $\Cal H^1_\kappa(M,\R)$ 
admits a complement in the punctured cohomology bundle $\Cal H^1_\kappa(M\setminus \Sigma_\kappa,\R)$ on which the punctured Kontsevich-Zorich cocycle is isometric (with respect
to a continuous norm). By the estimate in \pref{eq:regul2bound} on the upper Lyapunov exponent, it follows that, since $s<1$, the cohomology class $c_\theta  \in \Cal H^1_{q_\theta}(M,\R)$ for all $\theta \in E_s$. The restriction of the punctured Kontsevich-Zorich cocycle $\{\Psi_t\}_{t\in\R}$ to the bundle $\Cal H^1_\kappa(M,\R)$  concides with the Kontsevich-Zorich cocycle $\{\Phi_t\}_{t\in \R}$. Since the latter is symplectic, it follows from \pref{eq:regul2bound} that its upper Lyapunov exponent on the symplectic subspace $I_{q_\theta}^\perp (M,\R) \subset \Cal H^1_{q_\theta}(M,\R)$ (see formulas \pref{eq:Iq} and \pref{eq:Iqperp}) satisfies the estimate:
\begin{equation}
\lambda^+_2(q_{\theta})\, \geq\, 1-s\,, \quad \text{ \rm for all }\, \theta \in E_s\,.
\end{equation}
By Theorem \ref{thm:spectralgaptwo}  the inequality $\lambda^+_2(q_{\theta}) \leq L_\kappa(q)<1$  holds for all $q\in \Cal M^{(1)}_\kappa$ and for almost all $\theta\in S^1$. Since the set $E_s$ has positive measure, it follows that $s>1-L_\kappa(q)>0$.

\end{proof}

\noindent By Lemma \ref{lemma:reglowbound}, Lemma \ref{lemma:B1} and Lemma \ref{lemma:regupperbound} we derive a fundamental relation between the weighted Sobolev order and the Lyapunov exponents of Oseledec simple basic currents and invariant distributions.

\begin{definition} \label{def:Oseledecbasis}
An \emph{Oseledec basis }for the space $\Cal I^s_{\pm q}(M)$ [ $\Cal B^s_{\pm q}(M)$] of invariant distributions [of basic currents] is a basis contained in the union of all Oseledec subspaces  for the cocycle $\{G_t^s \vert {\Cal I}^s_{\kappa,\pm}(M) \}$ [ $\{\Phi_t^s \vert {\Cal B}^s_{\kappa,\pm}(M) \}$] 
at any Oseledec regular point  $q\in \Cal M^{(1)}_\kappa$. Since the cocycles $\{\Phi_t^s \vert {\Cal I}^s_{\kappa,\pm}(M) \}$ and $\{\Phi_t^s \vert {\Cal B}^s_{\kappa,\pm}(M) \}$ are isomorphic via the maps $\Cal D_{\pm q}: {\Cal B}^s_{\pm q}(M) \to {\Cal I}^s_{\pm q}(M)$ introduced in \pref{eq:DtoC}, any basis $\beta^s_{\pm q}\subset {\Cal I}^s_{\pm q}(M)$  is Oseledec if and only if the basis $\Cal D_{\pm q} (\beta^s_{\pm q}) \subset \Cal B^s_{\pm q}(M)$ is Oseledec as well. 
\end{definition}

\begin{theorem} 
\label{thm:bcreg}
For any $SO(2,\R)$-absolutely continuous, KZ-hyperbolic  measure $\mu$ on a stratum ${\Cal M}^{(1)}_{\kappa}$ of orientable quadratic differentials, for $\mu$-almost all $q\in \Cal M^{(1)}_\kappa$ and for any Oseledec basis $\{C^\pm_1, \dots, C^\pm_g\} \subset \Cal B^1_{\pm q}(M)$:
\begin{equation}
\label{eq:bcreg}
{\Cal O}_q^W(C^\pm_i)= 1-\vert L^\pm_\mu(C^\pm_i)\vert = 1-\lambda^\mu_i  \,,
 \quad \text{ \rm for all }\, i\in \{1, \dots, g\}\,.
\end{equation}
Consequently, for any Oseledec basis $\{\Cal D^\pm_1, \dots, \Cal D^\pm_g\} \subset \Cal I^1_{\pm q}(M)$,
\begin{equation}
\label{eq:idreg}
{\Cal O}_q^H(\Cal D^\pm_i) = \vert l^\pm_\mu(\Cal D^\pm_i)\vert = 1-\lambda^\mu_i 
 \,, \quad \text{ \rm for all }\, i\in \{1, \dots, g\}\,.
\end{equation}
\end{theorem}

\noindent The restriction of the cocycle $\{\Phi^1_t\}_{t\in \R}$ to the invariant sub-bundle $\Cal Z_\kappa^1(M)$ is described by the following Oseledec-type theorem, which is a straightforward generalization of  \cite{F02}, Theorem 8.7, based on Lemma \ref{lemma:B1}.

\begin{theorem} 
\label{thm:Otype}
For any $SO(2,\R)$-absolutely continuous, KZ-hyperbolic  measure on a stratum ${\Cal M}^{(1)}_{\kappa}$ of orientable quadratic differentials, there exists a measurable $\Phi^1_t$-invariant splitting:
\begin{equation}
\label{eq:Z1split}
\Cal Z_\kappa^1(M) =  \Cal B^1_{\kappa,+}(M) \oplus  \Cal B^1_{\kappa,-}(M) \oplus 
\Cal E^1_\kappa(M) \,.
\end{equation}
\begin{enumerate}
\item The restriction of the cocycle $\{\Phi^1_t\}_{t\in \R}$ to the bundle $\Cal B^1_{\kappa}(M):= \Cal B^1_{\kappa,+}(M) \oplus  \Cal B^1_{\kappa,-}(M)$ is (measurably) isomorphic to the Kontsevich-Zorich cocycle, hence it has Lyapunov spectrum \pref{eq:KZexps}. The sub-bundle $\Cal B^1_{\kappa,+}(M)$ 
corresponds to the strictly positive exponents and the sub-bundle $\Cal B^1_{\kappa,-}(M)$ 
to the strictly negative exponents. 
\item The Lyapunov spectrum of restriction of the cocycle $\{\Phi^1_t\}_{t\in \R}$ to the infinite dimensional bundle $\Cal E^1_\kappa(M)$ of exact currents is reduced to the
single Lyapunov exponents $0$ (in fact, the cocycle is isometric with respect to a suitable continuous norm on $\Cal E^1_\kappa(M)$).
\end{enumerate}
\end{theorem}

\noindent Informations on the Lyapunov structure of the restrictions of the cocycles $\{\Phi_t^s\}_{t\in\R}$ to the sub-bundles $\Cal B^s_{\kappa, \pm} (M) \subset \Cal B^s_{\kappa, \pm} (M\setminus\Sigma_\kappa)\subset W^{-s}_\kappa(M)$, for any $s\geq 1$, can be derived from
Theorem \ref{thm:Otype} and from the structure theorem for basic currents (see Theorem 
\ref{thm:bcstruct}) combined with the following result:

\begin{lemma} 
\label{lemma:deltasid}
Let $\mu$ be any $G_t$-invariant ergodic probability measure on a stratum $\Cal M^{(1)}_\kappa$
of orientable quadratic differentials. For any $q\in \Cal R_\mu$, the image under the map $\delta_{\pm q}: \Cal B_{\pm q}(M\setminus \Sigma_q) \to \Cal B_{\pm q} (M\setminus \Sigma_q)$ defined by formulas \pref{eq:deltas} of any simple current $C^\pm\in \Cal B^s_{\pm q}(M\setminus \Sigma_q)$ is a simple current $\delta_{\pm q}(C^\pm)\in \Cal B^{s+1}_{\pm q}(M\setminus \Sigma_q)$.
The following identities hold for the Sobolev order map  and the Lyapunov exponent map on $\Cal B_{\pm q} (M\setminus \Sigma_q)$:
\begin{equation}
\label{eq:deltasid}
\begin{aligned}
{\Cal O}_q^W \circ \delta_{\pm q}  \,&=\, {\Cal O}_q^W  - 1\,, \quad \text{  \rm for all }\, 
q\in \Cal M^{(1)}_\kappa  \,; \\
L^\pm_\mu  \circ \delta_{\pm q}  \,&=\,  L^\pm_\mu \mp 1\,, \quad \text{  \rm for all }\, q\in \Cal R_\mu
\subset \Cal M^{(1)}_\kappa \,.
\end{aligned}
\end{equation}
\end{lemma}
\begin{proof}
For any $C^\pm \in \Cal B_{\pm q}(M\setminus \Sigma_q)$, it follows immediately by the
definitions of the Sobolev spaces and of the maps $\delta_{\pm q}$ on $\Cal B_{\pm q}(M\setminus \Sigma_q)$ that $\delta_{ q}(C^+) \in  W_q^{-s-1}(M)$ if and only if $C^+ \wedge 
\Re(q^{1/2}) \in H^{-s}_q(M)$, hence if and only if $C^+ \in  W_q^{-s}(M)$. Similarly, 
$\delta_{ -q}(C^-) \in W_q^{-s-1}(M)$ if and only if $C^- \wedge \Im(q^{1/2}) 
\in H^{-s}_q(M)$, hence if and only if $C^- \in  W_q^{-s}(M)$. The first identity in \pref{eq:deltasid} 
is therefore proved. In fact, we have proved that the maps $\delta^s_{\pm  q}$ send $\Cal B^s_{\pm q}(M\setminus \Sigma_q)$ onto the subspace of cohomologically trivial currents in $\Cal B^{s+1}_{\pm q}(M\setminus \Sigma_q)$ while the space $\Cal B^s_{\pm q}(M)$ is mapped onto the subspace of cohomologically trivial currents in $\Cal B^{s+1}_{\pm q}(M)$. 

\smallskip
\noindent The following identity follows immediately from the definitions: for any quadratic differential $q\in \Cal M^{(1)}_\kappa$, for each $s\geq 0$ and all $t\in \R$, we have
\begin{equation}
\label{eq:deltasidbis}
\delta^{s}_{ \pm G_t(q)} \circ \Phi^s_t    = e^{\pm t}  \, (\Phi^{s+1}_t \circ \delta^{s}_{\pm q})   \qquad 
\text{ \rm on } \,\, \Cal B^s_{\pm q}(M\setminus\Sigma_q)\,.
\end{equation}
Since the maps $\delta^{s}_{\pm q}:\Cal B^s_{\pm q}(M\setminus \Sigma_q) \to \Cal B^{s+1}_{\pm q}(M\setminus \Sigma_q)$ are embeddings, it follows from the identity \pref{eq:deltasidbis} and from the Oseledec's theorem that the current $\delta^{s}_{\pm q}(C^\pm)\in \Cal B^{s+1}_{\pm q}(M
\setminus\Sigma_q)$ is  simple if $C^\pm \in \Cal B^s_{\pm
q}(M\setminus\Sigma_q)$ is. It also follows that the second indentity in \pref{eq:deltasid} holds.
\end{proof}
\noindent An analogous statement for invariant distributions can be derived. For any $q\in \Cal M^{(1)}_\kappa$,  let $\Cal L_{\pm q}: \Cal S_q(M) \to \Cal S_q(M)$ denote the Lie derivative operators on the space $\Cal S_q(M)$ of all $q$-tempered currents: for any $C\in \Cal S_q(M)$,
\begin{equation}
\label{eq:Lieder}
{\Cal L}_q (C) :=  \Cal L_{S_q} (C) \quad \text {\rm and } \quad {\Cal L}_{-q} (C):=  \Cal L_{T_q} (C)\,.
\end{equation}

\begin{lemma}
\label{lemma:idderid}
For any $q\in \Cal M^{(1)}_\kappa$, the operators $\Cal L_{\pm q} : \Cal I_{\pm q}(M) \to \Cal I_{\pm q}(M)$ are well-defined and injective. Let $\mu$ be any $G_t$-invariant ergodic probability measure on $\Cal M^{(1)}_\kappa$. The following identities hold on the spaces $\Cal I_{\pm q}(M)$:
\begin{equation}
\label{eq:idderid}
\begin{aligned}
{\Cal O}_q^H \circ  \Cal L_{\pm q}  \,&=\, {\Cal O}_q^H  - 1\,, \quad \text{  \rm for all }\, 
q\in \Cal M^{(1)}_\kappa  \,; \\
l^\pm_\mu  \circ \Cal L_{\pm q}  \,&=\,  l^\pm_\mu \mp 1\,, \quad \text{  \rm for all }\, q\in \Cal R_\mu
\subset \Cal M^{(1)}_\kappa \,.
\end{aligned}
\end{equation}
\end{lemma}

\begin{corollary} 
\label{cor:Lspectrum} 
Let $\mu$ be any $G_t$-invariant ergodic probability measure on a stratum $\Cal M^{(1)}_\kappa$
of orientable holomorphic quadratic differentials. For any $s\geq 0$ the (finite) Lyapunov spectrum of the cocycle $\{\Phi^s_t \vert \Cal B^s_{\kappa,\pm}(M)\}$  is a finite subset of the countable set 
\begin{equation}
\label{eq:Lspectrumone}
\{ \pm \lambda^\mu_1\} \,\cup \, \{ \pm \lambda^\mu_i \mp  j\,\vert  1< i <2g, \, j\in \N \}\,;
\end{equation}
hence the Lyapunov spectrum of the cocycle $\{G^s_t \vert \Cal I^s_{\kappa,\pm}(M)\}$ is a finite
subset of the countable set  
 \begin{equation}
 \label{eq:Lspectrumtwo}
\{0\} \,\cup \,  \{\pm \lambda^\mu_i \mp  (j+1)\,\vert  1< i <2g, \, j \in \N\}
\end{equation}
 (each element of the sets \pref{eq:Lspectrumone} and \pref{eq:Lspectrumtwo} is taken with
multiplicity one). 
\end{corollary}
\begin{proof} Let $\Pi^1_{\kappa, +}(M,\R)$ [$\Pi^1_{\kappa, -}(M,\R)$] be the continuous sub-bundles 
of the cohomology bundle $\Cal H^1_\kappa(M,\R)$ over $\Cal M^{(1)}_\kappa$ which fibers at any quadratic differential $q\in {\Cal M}^{(1)}_\kappa$ are given by the spaces $\Pi^1_{+q}(M,\R)$ 
[$\Pi^1_{-q}(M,\R)$]  defined in formula \pref{eq:Piabs}. The sub-bundles $\Pi^1_{\kappa, \pm}(M,\R)$ are invariant under the Kontsevich-Zorich cocycle $\{\Phi_t\}_{t\in\R}$ and the Lyapunov spectrum of
the restriction $\{\Phi_t\vert \Pi^1_{\kappa, \pm}(M,\R)\}$ consists of the set $\{\pm \lambda^\mu_i \vert
1\leq i \leq 2g-1\}$. By formulas \pref{eq:wedgecon},  the image of the bundle $\Cal B^s_{\kappa,\pm}(M)$ under the cohomology map  $j_\kappa: \Cal Z^1_\kappa(M) \to \Cal H^1_\kappa(M,\R)$ is a sub-bundle of $\Pi^1_{\kappa, \pm}(M,\R)$ which is invariant under the Kontsevich-Zorich cocycle. By Lemma \ref{lemma:bcurrsub} the cocycle $\{\Phi^s_t \vert \Cal B^s_{\kappa,\pm}(M)\}$ is mapped under the cohomology map onto a restriction of the 
Kontsevich-Zorich cocycle. Let $\delta^s_{\kappa,\pm}: \Cal B^{s-1}_{\kappa,\pm}(M) \to \Cal B^s_{\kappa,\pm}(M)$ the measurable bundle maps defined fiber-wise for $\mu$-almost all
$q\in \Cal M^{(1)}_\kappa$ as the maps $\delta^s_{\pm q}:  \Cal B^{s-1}_{\pm q}(M) \to \Cal B^s_{\pm q}(M)$, defined in formula \pref{eq:deltas}. The kernel of the cohomology map $j_\kappa$ on $\Cal B^s_{\kappa,\pm}(M)$ is a $\Phi^s_t$-invariant sub-bundle which coincides with the range of the map $\delta^s_{\kappa,\pm}$. By Lemma \ref{lemma:deltasid}, the Lyapunov spectrum \pref{eq:Lspectrumone} of $\Phi^s_t \vert \Cal B^s_{\kappa,\pm}(M)$ can therefore be derived by induction on $[s]\in \N$. The Lyapunov spectrum of  \pref{eq:Lspectrumtwo} of $G^s_t \vert \Cal I^s_{\kappa,\pm}(M)$ can then be derived by the isomorphism \pref{eq:cocycleiso}.
\end{proof}

\noindent By Corollary \ref{cor:basiccohom} for any $s>3$ the cohomology map is surjective for 
almost all quadratic differentials in any circle orbit. The above result can be refined as follows:

\begin{corollary}
Let $\mu$ be any  $SO(2, \R)$-absolutely continous, $G_t$-invariant ergodic probability measure on a stratum $\Cal M^{(1)}_\kappa \subset \Cal M^{(1)}_g$ of orientable holomorphic quadratic differentials.  For any $s>3$, there exists  an integer vector $h^s:=(h^s_2, \dots, h^s_{2g-1}) \in \N^{2g-2}$ such that the Lyapunov spectrum of the cocycle $\{\Phi^s_t \vert \Cal B^s_{\kappa,\pm}(M)\}$  
is the (finite) set
\begin{equation}
\label{eq:Lspectrumonebis}
\{ \pm \lambda^\mu_1\} \,\cup \, \{ \pm \lambda^\mu_i \mp  j\,\vert  1<i <2g, \,0\leq j \leq h^s_i\}\,;
\end{equation}
hence the Lyapunov spectrum of the cocycle $\{G^s_t \vert \Cal I^s_{\kappa,\pm}(M)\}$ is
the (finite) set  
 \begin{equation}
 \label{eq:Lspectrumtwobis}
\{ 0 \} \,\cup \, \{\pm \lambda^\mu_i \mp  (j+1)\,\vert  1< i <2g, \, 0\leq j \leq h^s_i\}
\end{equation}
 (each element of the sets \pref{eq:Lspectrumonebis} and \pref{eq:Lspectrumtwobis} is taken with
multiplicity one).  
\end{corollary}
\noindent The integer vector $h^s$ depends on the Sobolev regularity of basic currents in the Oseledec's spaces related to the Lyapunov exponents $\{\lambda^\mu_2, \dots, \lambda^\mu_{2g-1}\}$ which come from the Kontsevich-Zorich cocycle. Hence, in particular,  in case $\mu$ is the unique absolutely continuous, $SL(2,\R)$-invariant ergodic measure on a connected component of a stratum, by Theorem \ref{thm:bcreg} the following estimate for the numbers $(h^s_2, \dots, h_g)$ holds: for all 
$i\in \{2,\dots,g\}$, 
$$
h^s_i = \max \{ h \vert 1-\lambda^\mu_i+ h <s  \}\,, \quad \text{ \rm if } \, s \not \in \N -\lambda^\mu_i \,.
$$

\noindent The Oseledec-Pesin theory of the cocycles $\{G_t^s \vert {\Cal I}^s_{\kappa,\pm}(M) \}$ has 
crucial implications for the theory of the cohomological equation. In particular, it implies that (if
$\mu$ is a KZ-hyperbolic measure) for $\mu$-almost all $q\in \Cal M^{(1)}_\kappa$, Oseledec bases of the spaces $\Cal I_{\pm q}^s(M)$ of invariant distributions have strong, quantitative linear independence properties.

\begin{theorem} \label{thm:Oseledecbasis} 
Let $\mu$ be any $G_t$-invariant ergodic probability measure on a stratum $\Cal M^{(1)}_\kappa$ of orientable holomorphic quadratic differentials. For any $s>0$ and for $\mu$-almost all $q\in \Cal M^{(1)}_\kappa$, let $b^s_{\pm q} =\{ \Cal D^\pm_1, \dots, 
\Cal D^\pm_{J_\pm(s)}\}$ be a Oseledec basis of the space $\Cal I^s_{\pm q}(M)$ of  invariant distributions and let 
\begin{equation}
(l^\pm_1, \dots , l^\pm_{J_\pm(s)}) := \left(l^\pm_\mu(\Cal D^\pm_1), \dots , l^\pm_\mu(\Cal D^\pm_{J_\pm(s)})\right)
\end{equation}
denote the Lyapunov spectrum of the cocycle $\{G_t^s \vert {\Cal I}^s_{\kappa,\pm}(M) \}$ over 
$(\{G_t\}, \mu)$. For any $\epsilon>0$, there exists a measurable function $K^s_\epsilon: \Cal M^{(1)}_\kappa \to \R^+$ such that the following holds. For every $\tau \in (0,1]$, there exist
linearly independent systems of smooth functions $\{u^\pm_1(\tau), \dots, u^\pm_{J_\pm(s)}(\tau) \} 
\subset H^\infty_q(M)$ such that,  for all $i$, $j\in \{1, \dots, J_\pm(s)\}\,$ and for all $0\leq r \leq s$,
\begin{equation}
\label{eq:Oseledecbasis}
\Cal D^\pm_i\left(u^\pm_j(\tau)\right)\,= \, \delta_{ij} \quad \text{ \rm and } \quad
\vert u^\pm_j(\tau) \vert_{r} \,\leq\, K_\epsilon(q) \, \tau^{\vert l^\pm_j\vert -r -\epsilon}\,.
\end{equation}
\end{theorem}
\begin{proof} By Corollary  \ref{cor:Lspectrum}, for any $G_t$-invariant ergodic probability measure $\mu$ on $\Cal M^{(1)}_\kappa$ and any $s>0$, the Lyapunov exponents of the distributional cocycle $\{G_t^s \vert {\Cal I}^s_{\kappa,\pm}(M) \}$ are all of the same sign, namely
\begin{equation}
\{ l^\pm_1, \dots , l^\pm_{J_\pm(s)}\} \subset  \mp \R^+ \cup \{0\} \,.
\end{equation}
It follows from the Oseledec's theorem that, for any $\epsilon>0$, there exists a strictly positive measurable function $C^{(1)}_\epsilon: \Cal M^{(1)}_\kappa \to \R^+$ such that, for every 
$i\in \{1, \dots, J_\pm(s)\}$, for $\mu$-almost all $q \in \Cal M^{(1)}_\kappa$ and for all $t\geq 0$, 
\begin{equation}
\label{eq:basisgrowth}
\vert G^s_{\mp t}  (\Cal D^\pm_i) \vert_{-s} \, \geq \,  C^{(1)}_\epsilon(q) \, 
e^{ (\vert l^\pm_i\vert-\epsilon)  t }\,.
\end{equation}
It also follows from the Oseledec's theorem  that for any $s>0$, for $\mu$-almost all $q \in \Cal M^{(1)}_\kappa$ and any Oseledec basis $b^s_{\pm q}\subset \Cal I^s_{\pm q}(M)$, the distorsion of the (Oseledec) basis $G^s_t (b^s_{\pm q}) \subset H^s_{G_t(q)}(M)$ grows subexponentially in time. The \emph{distorsion }of a basis $b^s_{\pm q}:=\{\Cal D^\pm_1, \dots, \Cal D^\pm_J\} \subset  \Cal I^s_{\pm q}(M)$ is the number
\begin{equation}
d_q(b^s_{\pm q})   :=   \sup\, \{  \frac{\sum_{i=1}^{J_\pm(s)} \vert c_i \vert 
\,\vert \Cal D^\pm_i\vert_{H^{-s}_q(M)}} { \vert \sum_{i=1}^{J_\pm(s)}
c_i\,\Cal D^\pm_i \vert_{H^{-s}_q(M)}} \,\, \vert \,\,  c \in \C^{J_\pm(s)} \} \,.
\end{equation}
The Oseledec's theorem implies that, if $b^s_{\pm q}\subset \Cal I^s_{\pm q}(M)$ is an Oseledec basis, then for any $\epsilon >0$, there exists a measurable function $C^{(2)}_\epsilon: \Cal M^{(1)}_\kappa \to \R^+$ such that for $\mu$-almost all $q \in \Cal M^{(1)}_\kappa$ and for 
all $t\in \R$,
\begin{equation}
\label{eq:distorsiongrowth}
d_{G_t(q)} \left( G^s_t  (b^s_{\pm q} ) \right) \, \leq \,  C^{(2)}_\epsilon(q) \, 
e^{ \epsilon \vert t \vert }\,.
\end{equation}
By the estimates \pref{eq:basisgrowth} and \pref{eq:distorsiongrowth}, it follows that, for any $\epsilon>0$, there exists a measurable function $C^{(3)}_\epsilon: \Cal M^{(1)}_\kappa \to \R^+$ 
such that the following holds. For every $t\geq 0$, since $H^\infty_q(M)$ is dense in $H^s_q(M)$
for any $s\in \R$, there exists a system $\{v^\pm_1(t), \dots, v^\pm_{J_\pm(s)}(t) \} \subset H^s_q(M)$ such that, for $\mu$-almost all $q\in  \Cal M^{(1)}_\kappa$ and for all $i$, $j \in \{1, \dots, J_\pm(s)\}$,  
\begin{equation}
\Cal D^\pm_i\left(v^\pm_j(t)\right)\,= \, \delta_{ij} \,\, \text{ \rm and } \,\,
\vert v^\pm_j(t) \vert_{H^s_{G_{\mp t}(q)}(M)} \,\leq\, C^{(3)}_\epsilon(q) \, 
e^{-\left(\vert l^\pm_j\vert   -\epsilon\right) t}\,.
\end{equation}
By the bound  \pref{eq:cocyclenormsone} on the norm $\Vert G^s_t: H^s_q(M) \to H^s_{G_t(q)}(M)\Vert$, it follows that, for any $0 \leq r  \leq s$ and for all $j\in \{1, \dots, J_\pm(s)\}$,
\begin{equation}
\vert v^\pm_j(t) \vert_{H^r_{q}(M)} \,\leq\, C^{(3)}_\epsilon(q) \, 
e^{(r-\vert l^\pm_j\vert  + \epsilon) t} \,, \quad \text{ \rm for any } \,t\geq 0\,.
\end{equation}
It follows that the system $\{u^\pm_1(\tau), \dots, u^\pm_{J_\pm(s)}(\tau) \}
\subset H^\infty_q(M)$, defined for any $\tau \in (0,1]$ by the identities
\begin{equation}
u^\pm_j(\tau) := v^\pm_j(-\log \tau) \,, \quad  j\in \{1, \dots, J_\pm(s)\}\,,
\end{equation}
satisisfies the required properties  \pref{eq:Oseledecbasis}. 
\end{proof}

\medskip
\noindent The following theorem, which can be interpreted as spectral gap result for the cocycles
$\{G^0_t\}_{t\in \R}$ on the bundle $H^0_{\kappa}(M)$ of square-integrable functions, is the main
technical result of the paper. 

\smallskip
\noindent For any $(\sigma,\lyap) \in \R^2$, we introduce the {\it upper (forward) Lyapunov norm }of a distribution $U \in H^{-\sigma}_q(M)$ at a quadratic differential $q\in \Cal M^{(1)}_\kappa$ as the non-negative extended real number
\begin{equation}
\label{eq:distLnorm}
\Cal N^{\sigma,\lyap}_q(U) :=  \sup_{t\to +\infty}   e^{-\lyap t} \,\vert G^\sigma_t (U) \vert_{-\sigma} 
\end{equation}
\begin{theorem} 
\label{thm:G0gap}
Let $\mu$ be any  $SO(2, \R)$-absolutely continous, $G_t$-invariant ergodic probability measure on a stratum $\Cal M^{(1)}_\kappa$ of orientable quadratic differentials.  For any any $\sigma>0$ and any $\lyap<1$, there exist a real number $\epsilon:=\epsilon(\sigma,\lyap)>0$  and a measurable function $C_{\sigma,\lyap}: \Cal M^{(1)}_\kappa \to \R^+$ such that, for $\mu$-almost all $q\in \Cal M^{(1)}_\kappa$ and all functions $U \in L^2_q(M)$ of zero average,
\begin{equation}
\label{eq:G0gap}
{\Cal N}^{\sigma,-\epsilon}_q(U) <  C_{\sigma,\lyap} (q) \,  \{ \vert U\vert_0  \,+ \, {\Cal N}^{1,\lyap}_q(S_qU)\}\,.
\end{equation}
\end{theorem}
\begin{proof}
The outline of the argument is as follows.  From the results on the cohomological 
equation and on the Lyapunov spectrum of the coycles $\{G^r_t\}_{t\in\R}$ on the bundles of invariant 
distributions we will derive (exponential) estimates on the norms $\vert G^r_t(U)\vert_{-r}$ for a sufficiently large $r>0$, along a sequence of suitable visiting times. Since $\vert G^0_t(U)\vert_0$
is constant, as a consequence of the invariance of the $L^2_q(M)$ norm under the Teichm\"uller
flow, the interpolation inequality for dual weighted Sobolev norms (Corollary \ref{cor:dualintineq}) implies the required (exponential) estimates on $\vert G^\sigma_t(U)\vert_{-\sigma}$ for any $\sigma>0$.

\smallskip
\noindent  Let $q\in \Cal M^{(1)}_\kappa$ and $U\in L^2_q(M)$ be a function of zero average. 
Let us assume that there exist $\lyap_0<1$ and $\Cal N\in \R^+$ such that the Lyapunov norm $\,\Cal N^{1,\lyap_0}_q(S_qU) \leq \Cal N\,$. By definition \pref{eq:distLnorm} of Lyapunov norms, for any $t\geq 0$, 
\begin{equation}
\label{eq:Sderbound}
\vert G^1_t (S_q U)  \vert _{-1}   \,\leq \,  \Cal N\, e^{\lyap_0 t} \,.
\end{equation}
By the spectral gap Theorem \ref{thm:spectralgapone} and by Corollary \ref{cor:Lspectrum}, 
there exist $C_1>0$, $\lyap_1>0$ and a positive measure set $\Cal P_1\subset 
{\Cal M}^{(1)}_\kappa$ such that, for all $q\in \Cal P_1$ and all $S_q$-invariant distributions 
$\Cal D \in {\Cal H}^{-r}_q(M) \cap \Cal I^r_q(M)$, the following holds:
\begin{equation}
\label{eq:idbound}
\vert G^r_t (\Cal D)  \vert _{-r}   \,\leq \,   C_1 \, e^{-\lyap_1 t}\,  \vert \Cal D\vert_{-r} \,.
\end{equation}
By Lemma \ref{lemma:Greenmeas} on the measurability of the Green operators for the cohomological equation, there exists $r>0$ such that the following holds. There
exists a constant $C_2>0$ and a set $\Cal P_2 \subset \Cal P_1$ of positive
measure such that
\begin{equation}
\label{eq:Greenbound}
\Vert  \Cal U_q^{r,1}: {\Cal H}_q^{-1}(M) \to  H^{-r}_q(M) \Vert 
\,\leq C_2 \,, \quad \text{ for all }\, q\in \Cal P_2 \,.
\end{equation}
For $\mu$-almost all $q\in \Cal M^{(1)}_\kappa$, there exists a sequence $(t_n)_{n\in \N}$ of visiting times of the forward orbit $\{G_t(q) \,\vert\, t\geq 0\}$ to the positive measure set $\Cal P_2$ such that,
for each $n\in \N$, the function $t_n: \Cal M^{(1)}_\kappa \to \R^+$ is measurable  and
\begin{equation}
\label{eq:tn}
\lim_{n\to +\infty}  \, \frac{ t_n(q) } {n}   := p \, > \, \frac{\log C_1}{ \lyap_1} \,.
\end{equation}
For each $n\in \N$ and any $r\in\R$, let us introduce the notations:
\begin{equation*}
\begin{array}{ccc}
G_n:= G_{t_n}   \,,  & q_n:=G_n(q)\,, & \{S_n, T_n\}:= \{S_{q_n}, T_{q_n}\}    \,,\\ 
G_n^r:= G^r_{t_n}\, & U^r_n:= G^r_n(U) \,, & F_n:= G^1_n (S_qU)\,, \\ 
H^r_n(M):= H^r_{q_n}(M)\,, &  {\Cal H}^r_n(M):= {\Cal H}^r_{q_n}(M)\,, & \Cal I ^r_n(M):= \Cal I^r_{q_n}(M)     \,.
\end{array}
\end{equation*}
By the definition of the Teichm\"uller flow $\{G_t\}$, the orthonormal frame 
$\{S_n, T_n\} =\{ e^{-t_n} S_q, e^{t_n} T_q\}$ for all $n\in \N$. Hence the following cohomological equation holds:\begin{equation}
S_n U^r_n= e^{-t_n} F_n   \in H^{-1}_{n}(M) \,.
\end{equation}
For each $n\in\N$ and any $r\in \R$ let $\vert \cdot \vert_{r,n}$ denote the weighted Sobolev norm on the
space $H_n^r(M)$. It follows from equation \pref{eq:Sderbound} that
\begin{equation}
\vert S_n U^r_n \vert_{-1,n}  \, \leq \, \Cal N \, e^{-(1-\lyap_0) t_n}      \,.
\end{equation}
By the definition of the Green operators and by equation \pref{eq:Greenbound}, there exists a
solution $G_n \in H^{-r}_{n}(M)$, orthogonal to $\Cal I^r_{n}(M)$, of the cohomological equation
$S_n G_n = S_n U^r_n \in H^{-1}_{n}(M)$ such that
\begin{equation}
\label{eq:Gnbound}
\vert G_n \vert_{-r,n}\,  \leq \,  C_2\, \Cal N\, e^{-(1-\lyap_0) t_n}\,.
 \end{equation}
It follows that there exists $\Cal D_n \in \Cal I^r_{n}(M) \cap {\Cal H}^{-r}_{n}(M)$ such that  $U_n\in H^0_n(M)$ has the following orthogonal decomposition in the Hilbert space $H^{-r}_{n}(M)$:
\begin{equation}
\label{eq:Unod}
U^r_n =  G_n \, + \, \Cal D_n \,.
\end{equation}
For each $n\in \N$, let $\pi_n: H^{-r}_{n}(M) \to  \Cal I^r_{n}(M)$ be the orthogonal projection onto the subspace of $S_n$-invariant distributions and let $\tau_n:= t_{n+1} -t_n$. Since by the definitions $U^r_{n+1} = G^r_{\tau_n}(U^r_n)$, for all $n\in \N$,  and 
$$
G^r_{\tau_n}\left(\Cal I^r_{n}(M) \cap {\Cal H}^{-r}_{n}(M) \right) \,\subset \Cal I^r_{n+1}(M)
 \cap {\Cal H}^{-r}_{n+1}(M)
 $$ 
 the following recursive identity holds:
\begin{equation}
\Cal D_{n+1} =   G^r_{\tau_n} (\Cal D_n)   \,+\,  \pi_{n+1}\circ  G^r_{\tau_n} (G_n) \,.
\end{equation}
Since $\vert G^r_{\tau_n} (G_n)\vert_{-r,n+1}  \leq  e^{r\tau_n} \vert G_n\vert_{-r,n}$ by the bound 
\pref{eq:cocyclenormsone}  on the norm of $G^r_t : H^{-r}_t(M) \to H^{-r}_{G_t(q)}(M)$, it follows from  \pref{eq:idbound} and  \pref{eq:Gnbound} that 
\begin{equation}
\label{eq:Dnrec}
\vert \Cal D_{n+1} \vert_{-r, n+1}   \,\leq \,  C_1 \,e^{ -\lyap_1 \tau_n} \vert \Cal D_{n} \vert_{-r, n} 
+  C_2 \,\Cal N \,e^{r\tau_n -(1-\lyap_0) t_n} \,,
 \end{equation}
 which implies by induction a bound of the form
 \begin{equation}
 \label{eq:Dnbound}
 \vert \Cal D_{n+1} \vert_{-r, n+1} \, \leq \,  C_1^n \, e^{-\lyap_1 (t_{n+1} -t_0)}   \,
 \vert \Cal D_0 \vert_{-r, 0}\,+\, C_2\,\Cal N\, \sum_{j=0}^{n-1} C_1^j \, e^{s_{n,j}} 
 \end{equation}
 with the sequence $\{ s_{n,j} \vert n\in \N, \,0\leq j \leq n\}  $ given by the identity
 \begin{equation}
 \label{eq:snj}
 s_{n,j} := -\lyap_1 (t_{n+1} -t_{n-j+1}) +
 r\tau_{n-j} -(1-\lyap_0) t_{n-j}\,.
 \end{equation}
 Since $U\in L^2_q(M)$ and $G^r_n(U) =G^0_n(U) \in H^0_n(M)=L^2_q(M) \subset H^{-r}_n(M)$,
the following bound holds:
\begin{equation}
\label{eq:GnU0bound}
\vert G_n^r (U) \vert _{0, n}  \,\leq \, \vert U \vert_{0} \,, \quad \text{ for all } \, n\in \N\,.
\end{equation}
It follows in particular from the decomposition \pref{eq:Unod} that
\begin{equation}
\label{eq:D0bound}
\vert \Cal D_0 \vert_{-r, 0} \,\leq \,\vert U_0^r \vert_{-r, 0}\,\leq  \,\vert G_0^r (U) \vert _{0, 0} 
 \,\leq \, \vert U \vert_{0} \,.
\end{equation}
The main step in the argument is the proof of the following claim. There exist a positive measurable
function $C_3: \Cal P_2 \to \R^+$ and a real number $\lyap_3:=\lyap_3(\lyap_0,\lyap_1)>0$ such that
 \begin{equation}
 \label{eq:sumclaim}
 \sum_{j=0}^{n-1} C_1^j \, e^{s_{n,j}} \leq C_3(q)  \, e^{ -\lyap_3 n} \,, \quad \text{ for all } \, n\in \N\,.
 \end{equation}
 Let $\omega_1<\omega<\omega_2$ be positive real numbers such that
  \begin{equation}
   \begin{aligned}
  &(a) \,\, (1-\lyap_0)\omega_1-(\lyap_1 +r) (\omega_2-\omega_1) >0 \,, \\
   &(b)\,\, \lyap_1 \omega_2 - \lyap_1 (\omega_2-\omega_1)> \log C_2 \,.
    \end{aligned}
  \end{equation}
 By condition \pref{eq:tn} on the sequence $(t_n)_{n\in \N}$ there exists a measurable
 function $n_0: \Cal M^{(1)}_\kappa \to \N$ such that
  \begin{equation}
 \label{eq:om12}
 \omega_1 \,n \leq t_n(q) \leq \omega_2\, n \,, \quad \text{ for all }\, n\geq n_0(q)\,. 
  \end{equation}
  It follows that, for any $q\in \Cal M^{(1)}_\kappa$, all $n\geq n_0(q)$ and all $j \leq n$, 
   \begin{equation}
   \begin{aligned}
 s_{n,j} &\leq  [(\lyap_1 +r) (\omega_2-\omega_1) -(1-\lyap_0)\omega_1] \,n \\&+ [(1-\lyap_0)\omega_1 - r(\omega_2-\omega_1) -\lyap_1\omega_2]\, j
 + \lyap_1 (\omega_2-\omega_1) +r\omega_2 \,.
 \end{aligned}
  \end{equation}
 Let $K:=  C_2\, e^{(1-\lyap_0)\omega_1 - r(\omega_2-\omega_1) -\lyap_1 \omega_2} $. There are two cases to consider: 
  \begin{equation}
 (a) \, K\leq 1 \qquad (b) \, K>1 \,.   
 \end{equation}
 In case $(a)$ we immediately obtain that, for all $n\geq n_0(q)$, 
 \begin{equation}
 \sum_{j=n_0}^{n-1} C_1^j \, e^{s_{n,j}}\, \leq \, n \, e^{\lyap_1 (\omega_2-\omega_1) +r\omega_2}  \, 
 e^{ - [(1-\lyap_0)\omega_1-(\lyap_1 +r) (\omega_2-\omega_1)]\,n}  \,,
 \end{equation}
   in case $(b)$ we obtain instead that
   \begin{equation}
   \sum_{j=n_0}^{n-1} C_1^j \, e^{s_{n,j}}\, \leq \,  
   \frac{ e^{\lyap_1 (\omega_2-\omega_1) +r\omega_2} }{ K-1} \,\, C_1^n\, e^{ -[\lyap_1 \omega_2 - \lyap_1 (\omega_2-\omega_1)]\,n}\,. 
    \end{equation}
  It follows that there exist constants $A>0$ and $\alpha:= \alpha(\lyap_0, \lyap_1)>0$ such that
 \begin{equation}
 \label{eq:sumclaimone}
\sum_{j=n_0}^{n-1} C_1^j \, e^{s_{n,j}} \leq A \, e^{ -\alpha\, n} \,, 
\quad \text{ for all } \, n\geq n_0(q)\,.
 \end{equation}
By condition \pref{eq:tn} on the sequence $(t_n)_{n\in \N}$ and by definition \pref{eq:snj}  \begin{equation}
\frac{s_{n,j}(q)}{n}  \,\to -(1-\lyap_2)\omega \,, \quad \text{ as }  \, n\to +\infty\,,
\end{equation}
for $\mu$-almost all $q\in \Cal M^{(1)}_\kappa$, uniformly with respect to $j\in \{0, \dots, n_0(q)-1\}$.
Hence there exists a measurable function $n_1: \Cal M^{(1)}_\kappa \to \N$ with $n_1\geq n_0$ 
and positive constants $B>0$ and $\beta:= \beta(\lyap_0,\lyap_1)>0$ such that
 \begin{equation}
  \label{eq:sumclaimtwo}
\sum_{j=0}^{n_0-1} C_1^j \, e^{s_{n,j}} \leq B \, e^{ -\beta\, n} \,, 
\quad \text{ for all } \, n\geq n_1(q)\,.
 \end{equation}
The claim \pref{eq:sumclaim} then follows from the estimates \pref{eq:sumclaimone} and \pref{eq:sumclaimtwo}.

\smallskip
\noindent It then follows from the orthogonal decomposition \pref{eq:Unod}, from the claim \pref{eq:sumclaim}, proved above, together with the upper bounds \pref{eq:Gnbound}, \pref{eq:Dnbound},
\pref{eq:D0bound} and the lower bound for visiting times in \pref{eq:om12}, that  there exist a measurable function $C_4:\Cal M^{(1)}_\kappa \to \R^+$ and a constant  $\lyap_4:= 
\lyap_4(\lyap_0, \lyap_1)>0$ such that 
\begin{equation}
\label{eq:GnUrbound}
\vert G_n^r (U) \vert _{-r, n}  \leq C_4(q)  \{\vert U\vert_0 \,+\,\Cal N\} \, e^{-\lyap_4\, n} \,, \quad\text{for all } n\geq n_1(q)\,.
\end{equation}
By the {\it interpolation inequality }for the scale of dual weighted Sobolev norms (see Corollary \ref{cor:dualintineq}), from the upper bounds \pref{eq:GnU0bound} and \pref{eq:GnUrbound} it follows that for any $\sigma>0$ there exist a measurable function $C_\sigma:\Cal M^{(1)}_\kappa \to \R^+$  and a constant $\lyap_\sigma:= \lyap(\sigma,\lyap_0, \lyap_1)>0$ such that 
\begin{equation}
\label{eq:GnUsbound}
\vert G_n^\sigma(U) \vert _{-\sigma, n}  \leq  C_\sigma(q)  \{\vert U\vert_0 \,+\,\Cal N\} \, 
e^{-\lyap_\sigma\, n} \,, \quad \text{for all } n\geq n_1(q)\,.
\end{equation}
It remains to prove that the latter bound implies the statement of the theorem. For $\mu$-almost 
all $q\in \Cal M^{(1)}_\kappa$ and for all $t\geq t_0(q)$, there exists a unique $n(t,q) \in \N$ such
that
\begin{equation}
\label{eq:ntq}
t_{n(t,q)}(q) \, \leq t  <  \, t_{n(t,q)+1}(q) \,.
\end{equation}
Let $\sigma>0$ be fixed and let $\omega_1^{(\sigma)} <\omega< \omega_2^{(\sigma)}$ be positive real numbers such that $\lyap_\sigma - \sigma (\omega_2^{(\sigma)}-\omega_1^{(\sigma)})>0$.  There exists a measurable function $n^{(\sigma)}_2: \Cal M^{(1)}_\kappa \to \R$ such that  $n^{(s)}_2 \geq n_1$  and, for $\mu$-almost 
$q\in \Cal M^{(1)}_\kappa$,
\begin{equation}
\label{eq:om12s}
\omega_1^{(\sigma)} \, n \,  \leq\,  t_n(q) \, \leq   \omega_2^{(\sigma)} \, n \,, \quad  \text{ for all } \, n\geq n^{(\sigma)}_2(q)\,.
\end{equation}
Let $t^{(\sigma)}_1: \Cal M^{(1)}_\kappa \to \R^+$ be a measurable function such that $t^{(\sigma)}_1 \geq t_0$
and $n(t,q)\geq n^{(\sigma)}_2(q)$ if $t\geq t^{(\sigma)}_1(q)$, for $\mu$-almost all $q\in \Cal M^{(1)}_\kappa$. 
It follows from \pref{eq:ntq} and \pref{eq:om12s} that, for any $t\geq t^{(\sigma)}_1(q)$,
\begin{equation}
\label{eq:ntqlbound}
n(t,q)+1\, \geq\, \frac{ t_{n(t,q)+1}(q)}{\omega_2^{(\sigma)}}  \, >  \, \frac{ t }{\omega_2^{(\sigma)}}\,.
\end{equation}
It follows by \pref{eq:om12s} and \pref{eq:ntqlbound}, for $\mu$-almost $q\in \Cal M^{(1)}_\kappa$
and all $t\geq t^{(\sigma)}_1(q)$,
\begin{equation}
\label{eq:taunqt}
\tau_{n(q,t)} \, \leq \, (\omega^{(\sigma)}_2 -\omega^{(\sigma)}_1)\left(n(q,t) +1\right) \,\,+ \,\, \omega^{(\sigma)}_1 \,\,.
\end{equation}
Since, for any $n \in \N$ and any $t_n \leq t <t_{n+1}$, 
\begin{equation}
\begin{aligned}
 \vert G_t^\sigma(U) \vert _{-\sigma}  &= \vert G^\sigma_{t-t_n} \circ G_{t_1}(U)\vert_{-\sigma}  \\ 
 &\leq e^{\sigma(t-t_n)} \, 
 \vert G_n^\sigma(U) \vert _{-\sigma}  \leq  e^{\sigma \tau_n} \, \vert G_n^\sigma(U) \vert _{-\sigma}\,,
 \end{aligned}
\end{equation}
 by the upper bounds \pref{eq:GnUsbound}, \pref{eq:ntqlbound} and \pref{eq:taunqt} the following estimate holds. Let $\epsilon:= \epsilon (\sigma, \lyap_0, \lyap_1)$ be the real number defined as follows:
\begin{equation}
\label{eq:epsilonsnu}
\epsilon: = \frac{ \lyap_\sigma - \sigma (\omega_2^{(\sigma)}-\omega_1^{(\sigma)}) } { \omega_2^{(\sigma)} } \,\, > \,\, 0\,.
\end{equation}
For $\mu$-almost all $q\in \Cal M^{(1)}_\kappa$ and all $t\geq t^{(\sigma)}_1(q)$
\begin{equation}
\label{eq:GtsUboundone}
\vert G_t^\sigma(U) \vert _{-\sigma}  \leq C_\sigma(q) \, e^{\sigma \omega^{(\sigma)}_1 +\lyap_\sigma}  
\{\vert U\vert_0 \,+\, \Cal N\}\, e^{ -\epsilon t }\,.
 \end{equation}
 Finally, for $\mu$-almost all $q\in \Cal M^{(1)}_\kappa$ and all $0\leq t\leq t^{(\sigma)}_1(q)$,
 \begin{equation}
\label{eq:GtsUboundtwo}
\vert G_t^\sigma(U) \vert _{-\sigma}  \leq  \vert U\vert_0   \leq  e^{\epsilon t^{(\sigma)}_1(q)} \vert U\vert_0\, e^{ -\epsilon t }\,.
 \end{equation}
The desired estimate \pref{eq:G0gap} immediately follows from  the upper bounds 
\pref{eq:GtsUboundone} and \pref{eq:GtsUboundtwo}. 
 
\end{proof}

\section{Smooth solutions}
\noindent In this final section we prove our main theorems on Sobolev regularity of smooth 
solutions of the cohomological equation for translation flows. The general result which we are able
to prove for any translation surface (and almost all directions) is a direct consequence of the Fourier 
analysis construction of distributional solutions in \S\S \ref{DS} and of the uniqueness result that 
follows from Theorem \ref{thm:BCgap}. The sharper result which we will prove for almost
all translation surfaces (and almost all directions) requires a deeper analysis of ergodic 
averages of translation flows. In fact, the construction of square integrable (bounded) solutions 
is based (as in \cite{MMY05}) on the Gottschalk-Hedlund theorem.  The required bounds on the 
growth of ergodic integrals are derived from the Oseledec-type result (Theorem \ref{thm:Otype}) 
and the spectral gap result (Theorem \ref{thm:G0gap}) for distributional cocycles, by the methods developed in \cite{F02} in the study of the deviation of ergodic averages. The Oseledec-type theorem was in fact already proved in \cite{F02}, while the spectral gap theorem for distributional cocycles
is new.

\subsection{The general case} \label{ss:general} 
In this section we prove a result on existence of smooth solutions of the cohomological equation for translation flows which holds for {\it any }orientable quadratic differential in almost all directions.  Such a result answers a question of Marmi, Moussa and Yoccoz who asked what is the best, that is, the smallest {\it regularity loss }within reach of the Fourier analysis methods of \cite{F97}. As indicated in \cite{MMY05}, the answer is essentially that the solution loses no more than $3+\epsilon$ derivatives (for any $\epsilon>0$) with respect to the scale of weighted Sobolev spaces introduced in \S \ref{WSS}. We recall that the regularity loss obtained in \cite{MMY05} (which only holds for {\it almost all}  quadratic differentials) is essentially $1+\text{\rm BV}$.

\medskip
\begin{theorem} 
\label{thm:GCEsmooth}
Let $q\in \Cal M^{(1)}_\kappa$ be any quadratic differential. Let $k\in \N$ be
any integer such that $k\geq 3$ and let $s>k$ and $r<k-3$. For almost all $\theta\in S^1$ (with respect  to the Lebesgue measure), there exists a constant $C_{r,s}(\theta)>0$ such that the following holds. 
If $f\in H^{s}_q(M)$ is such that  $\Cal D(f)=0$ for all $\Cal D \in \Cal I^s_{q_\theta}(M)$, the cohomological equation $S_{\theta}u=f$ has a  solution $u\in H^{r}_q(M)$ satisfying the following estimate:
\begin{equation}
\label{eq:GCEest}
\vert u\vert_{r}\leq C_{r,s}(\theta)\, \vert f \vert_{s}\,\,.
\end{equation}
\end{theorem} 
\begin{proof}
As a first step, we prove that, under finitely many distributional conditions on $f\in H^s_q(M)$ 
with $s>3$, the cohomological equation $S_\theta u=f$ has, for almost all $\theta\in S^1$, a solution 
$u\in H^{-r}_q(M)$ for any $r>0$, which satisfies an estimate such as \pref{eq:GCEest}.  Since 
$\Delta^F_q f \in H^{s-2}_q(M)$ with $s-2>1$, by Theorem \ref{thm:CEdistribution}, for any $r>0$ 
and for almost all $\theta\in S^1$, there exists a solution $U \in \bar H^{-2-r}_q(M)$ of the 
cohomological equation $S_\theta U= \Delta^F_q f$, which vanishes on constant functions 
and satisfies the bound
\begin{equation}
\label{eq:GCEUbound}
\Vert U  \Vert_{-2-r}\leq C^{(1)}_{r,s}(\theta)\, \Vert \Delta^F_q f \Vert_{s-2}
\leq C^{(1)}_{r,s}(\theta)\, \Vert f \Vert_{s}\,\,.
\end{equation}
Let $u \in  \bar H^{-r}_q(M)$ be the unique distribution vanishing on constant functions such that
$\Delta_q^F u =U$. Since the commutation relation $S_\theta  \Delta_q v= \Delta_q S_\theta v$
holds for any $v\in H_q^{3+r}(M)$,  the following distributional equation holds:
\begin{equation}
\label{eq:GCELaplkersol}
\Delta_q \left ( S_\theta u -f \right)=0 \,\, \text { in }\,\, H^{-3-r}_q(M)\,.
\end{equation}
In addition, from the estimate \pref{eq:GCEUbound} it follows immediately that
\begin{equation}
\label{eq:GCEubound}
\vert u  \vert_{-r}\, \leq \,  C^{(1)}_{r,s}(\theta)\, \vert f \vert_{s}\,\,.
\end{equation}
Let $N_\theta^r\subset \bar H^{-r}_q(M)$ be the subspace defined as follows: 
$$
N_\theta^r:= \{ u\in H^{-r}_q(M) \,\vert  \Delta_q S_\theta u =0 \in H_q^{-3-r}(M) \} \,.
$$
Since the kernel $\Cal K^r(\Delta_q) \subset H_q^{-1-r}(M)$ (which is equal to the subspace perpendicular to the range of the operator $\Delta_q: H_q^{3+r}(M) \to H_q^{1+r}(M)$)  and the space $\Cal I^{r}_{q_\theta}(M)$ of $S_\theta$-invariant distributions are both finite dimensional, it follows that the subspace $N_\theta^r$ is finite dimensional (for almost all $\theta \in S^1$). As a consequence, the space 
\begin{equation}
S_\theta (N^r_\theta):=\{ S_\theta u \,\vert \, u\in N^r_\theta\}  \subset \Cal K^r(\Delta_q)
\end{equation}
is finite dimensional, hence closed in $H^{-1-r}_q(M)$. We claim that there exists a constant 
$C^{(2)}_{r,s}(\theta)>0$ such that the following holds. For almost all $\theta \in S^1$, there
exists a unique distribution $\Cal U_\theta(f) \in H^{-r}_q(M)$ such that  

\begin{list}{}{}
\item[(a)] $\Cal U_\theta(f)$ is orthogonal to the subspace $\Cal I^r_{q_\theta}(M)$; 
\item[(b)] $S_\theta \Cal U_\theta(f) -f  \in \Cal K^r(\Delta_q)$ is orthogonal to the subspace $S_\theta (N_\theta^r)$;
\item[(c)]$\Cal U_\theta(f)$ satisfies the following bound:
\begin{equation}
\label{eq:GCEUthetaone}
\vert \Cal U_\theta(f)  \vert_{-r}\leq C^{(2)}_{r,s}(\theta)\, \vert  f \vert_{s}\,\,.
\end{equation}
\end{list}
In fact, let $u\in H^r_q(M)$ be any solution of equation \pref{eq:GCELaplkersol} which satisfies the
bound \pref{eq:GCEubound}. Let $u_1 \in  H^{-r}_q(M)$ be the component of $u$ orthogonal in 
$H^{-r}_q(M)$ to the subspace $\Cal I^r_{q_\theta}(M)$ of $S_\theta$-invariant distributions. Since $u-u_1\in  \Cal I^r_{q_\theta}(M)$, the distribution $u_1$ still satisfies the equation \pref{eq:GCELaplkersol}.
In addition, the norm $\vert u_1\vert_{-r} \leq \vert u\vert_{-r}$, hence the bound  \pref{eq:GCEubound}
still holds. Since $N_\theta^r$ is finite dimensional, there exists $u_2 \in H^{-r}_q(M)$, a solution of \pref{eq:GCELaplkersol},  such that $S_\theta u_2-f$ is equal to the component of $S_\theta u_1-f$ orthogonal to $S_\theta (N^r_\theta)$ in $H^{-1-r}_q(M)$. Let $\Cal U_\theta(f)$ be the component of $u_2$ orthogonal to $\Cal I^r_{q_\theta} (M)$ in $H^{-r}_q(M)$. By construction $\Cal U_\theta(f)$ 
satisfies the conditions $(a)$ and $(b)$ above and it is uniquely determined by them. Let us prove condition $(c)$. The linear operator $S_\theta:N^r_\theta \to S_\theta (N^r_\theta)$ is surjective and its domain is a finite dimensional subspace, hence it has a bounded right inverse and there exists $K_r:=K_r(\theta)>0$ such that, for any $u\in N^r_\theta$ orthogonal to $\Cal I^r_{q_\theta}(M)$,
\begin{equation}
\label{eq:finiteinverse}
\vert u \vert_{-r} \,\leq \, K_r \,\vert S_\theta u \vert_{-1-r} \,.
\end{equation}
Since by definition $\Cal U_\theta:=\Cal U_\theta(f)$ and  $u_1\in H^{-r}_q(M)$ are both
orthogonal to $\Cal I^r_{q_\theta} (M)$ and $\Cal U_\theta-u_1 \in N^r_\theta$, by the 
bound \pref{eq:finiteinverse},
\begin{equation}
\label{eq:GCEUthetatwo}
\vert \Cal U_\theta \vert_{-r}  \leq   \vert \Cal U_\theta -u_1 \vert_{-r}  + \vert u_1\vert_{-r}
\leq K_r \vert S_\theta\left (\Cal U_\theta-u_1\right) \vert_{-1-r} + \vert u_1\vert_{-r}\,.
\end{equation}
Since again by definition  $S_\theta \Cal U_\theta=S_\theta u_2$ and $S_\theta u_2-f$ is equal to 
the component of $S_\theta u_1-f$ orthogonal to $S_\theta (N^r_\theta)$,
\begin{equation}
\label{eq:GCEUthetathree}
\begin{aligned}
\vert S_\theta\left (\Cal U_\theta-u_1\right) \vert_{-1-r} &\leq \vert S_\theta u_1-f \vert_{-1-r} 
+ \vert S_\theta u_2-f \vert_{-1-r} \\
 & \leq 2 \vert S_\theta u_1-f \vert_{-1-r}  \leq  2 \vert u_1 \vert_{-r} + 2 \vert f\vert_{-1-r} \,.
\end{aligned}
\end{equation}
It follows from the bounds \pref{eq:GCEUthetatwo}, \pref{eq:GCEUthetathree} and from the bound \pref{eq:GCEubound} for the distribution $u_1\in H^{-r}_q(M)$ that the required bound \pref{eq:GCEUthetaone} holds with 
$$
C^{(2)}_{r,s}(\theta):= [(2K_r(\theta) +1)C^{(1)}_{r,s}(\theta)] +1\,.
$$
Since $S_\theta(N^r_\theta) \subset \Cal K^r(\Delta_q)$ is finite dimensional, there exist a finite (linearly independent) set $\{\Phi_1,\dots, \Phi_K\}$ of bounded linear (real-valued) functionals on the Hilbert space $H_q^{-1-r}(M)$ such that 
$$
\Cal K^r(\Delta_q)  \cap S_\theta(N^r_\theta)^\perp \cap  \text{Ker}(\Phi_1) \cap \dots \cap \text{ Ker }(\Phi_K)  \, =\, \{0\}\,.
$$
Let $\{\Cal D_1, \dots, \Cal D_K\} \subset H^{-s}_q(M)$ be the system of distributions defined as
follows: for each $j\in \{1, \dots, K\}$, 
\begin{equation}
\label{eq:GCEDeltaj}
\Cal D_j (f) :=  \Phi_j \left( S_\theta \Cal U_\theta(f) \,-\, f \right) \,,\quad \text { for all } \, f\in H^s_q(M)\,.
\end{equation}
The system $\{\Cal D_1, \dots, \Cal D_K\}$ has by construction the following property: if
$\Cal D_j(f)=0$ for all $j\in \{1, \dots, K\}$, then $\Cal U_\theta(f)$ is the required solution of
the cohomological equation $S_\theta u=f$. In fact, under such conditions $\Cal U_\theta(f)$ is by construction a solution orthogonal to constant functions. By Theorem \ref{thm:BCgap}, there exists $s(q)>0$ such that for all $0<s<s(q)$ the space $\Cal I^s_{q_\theta}(M)$ is $1$-dimensional 
for almost all $\theta\in S^1$. Thus, if $r <s(q)$, for almost all $\theta \in S^1$ the solution $u\in 
H^{-r}_q(M)$ of the cohomological equation $S_\theta u=f$ is unique (if it exists). It follows that, for 
any $s(q)>r>0$, the distribution $\Cal U_\theta(f)$ is the unique solution in $H^{-r}_q(M)$ of the 
cohomological equation, which implies that $\Cal U_\theta(f) \in H^{-r}_q(M)$ for any $r>0$.

\smallskip
\noindent We claim that $\{\Cal D_1, \dots, \Cal D_K\} \subset \Cal I^s_{q_\theta}(M)$. In fact, by its definition in formula \pref{eq:GCEDeltaj}, the distribution $\Cal D_j \in H^{-s}_q(M)$ for all $j\in \{1,\dots,K\}$. In addition, if there exists $v\in H^{s+1}_q(M)$ such that $f=S_\theta v$ , by conditions $(a)$ and $(b)$ on the distribution $\Cal U_\theta(f)$ we have
$$
S_\theta \left( \Cal U_\theta(f) -v\right)=S_\theta  \Cal U_\theta(f) -f  \in \Cal K^r(\Delta_q) \cap 
S_\theta(N^r_\theta)^\perp\,.
$$ 
However, by definition $S_\theta \left( \Cal U_\theta(f) -v\right) \in S_\theta(N^r_\theta)$, hence $S_\theta  \Cal U_\theta(f) -f  =0$ which implies that $\Cal D_1(f)=\dots = \Cal D_K(f) =0$. Hence by definition all the distributions of the system $\{\Cal D_1, \dots, \Cal D_K\}$ are $S_\theta$-invariant.

\smallskip
\noindent Finally we prove the statement of the theorem by induction on $k\in \N$. For $k=3$ the
statement holds by the previous argument. Let us assume that the statement holds for $k\geq 3$.
By the induction hypothesis, for any $s>k+1$ and any $r<k-2$ there exists $C_{r,s}(\theta)>0$ 
such that,  for almost all $\theta\in S^1$ and for any $f\in H^s_q(M)$ with $\Cal D(f) =0$ 
for all $\Cal D\in \Cal I^s_{q_\theta}(M)$, there exist solution $u$, $u_S$ and $u_T \in H^{r-1}_q(M)$ 
of the cohomological equations $S_\theta u=f$, $S_\theta u_S = Sf$ and $S_\theta u_T =Tf$ respectively such that $u$, $u_S$ and $u_T$ are orthogonal to constant functions and the following bound holds:
\begin{equation}
\label{eq:GCEinductbound}
\vert u \vert^2_{r-1} +  \vert u_S \vert^2_{r-1} +  \vert u_T \vert^2_{r-1}
\leq   C_{r,s}(\theta)\left(\vert f\vert^2_{s-1} +  \vert T f\vert^2_{s-1} + 
\vert S f\vert^2_{s-1}\right)\,.
\end{equation}
In fact, since $f\in H^s_q(M)$ is such that $\Cal D(f) =0$ for all $\Cal D\in \Cal I^s_{q_\theta}(M)$
it follows immediately that $\Cal D(Sf)=\Cal D(Tf)=0$ for all $\Cal D\in \Cal I^{s-1}_{q_\theta}(M)$.
Since $S_\theta$ commutes with $S$, $T$ in the sense of distributions, it follows that the 
distributions $Su-u_S$ and $Tu -u_T \in H^{r-1}_q(M)$ are $S_\theta$-invariant. Let $s(q)>0$ be 
such that for all $0<s<s(q)$ and for almost all $\theta\in S^1$ the space $\Cal I^s_{q_\theta}(M)$ is 
$1$-dimensional (see Theorem \ref{thm:BCgap}).  If $s(q)>1-r$, since $Su-u_S$ and $Tu -u_T
\in H^{r-1}_q(M)$ are $S_\theta$-invariant and orthogonal to constant functions,  it follows that 
$Su=u_S$ and $Tu=u_T$. The latter identities imply that $u \in H^{r}_q(M)$ and by \pref{eq:GCEinductbound} the required bound \pref{eq:GCEest} is satisfied. The proof of the induction step is therefore completed.
\end{proof}

\subsection{Ergodic integrals} 
\label{ss:ergint}
Optimal results on the loss of regularity for almost all orientable quadratic differentials (and almost all directions) will be derived from bounds on ergodic integrals by the Gottschalk-Hedlund theorem. The required bounds will be proved along the lines of  \S 9  in \cite{F02} with the key improvement given by the 'spectral gap'  Theorem \ref{thm:G0gap}.

\smallskip
\noindent The key idea of the argument given in \cite{F02} consists in studying the dynamics of the 
distributional cocycle  $\{\Phi^1_t\}_{t\in\R}$ on the infinite dimensional (non-closed) sub-bundles $\Gamma^{\pm}_{\kappa}\subset W^{-1}_{\kappa}(M)$ of $1$-dimensional currents generated by segments of leaves of the horizontal and vertical foliations. 

\noindent Let $\Cal T>0$. We will denote by $\gamma_{\pm q}^{\Cal T}$ a positively oriented segment of length $\Cal T>0$ of a leaf of the measured foliation ${\Cal F}_{\pm q}$ respectively. By the trace theorems for standard Sobolev spaces and by the comparison Lemma \ref{lemma:comparison}, the vector spaces $\Gamma^{\pm}_q$ generated by the sets $\{ \gamma_{\pm q}^{\Cal T}(p)\vert (p,\Cal T)
\in M\times \R\}$ are subspaces of the weighted Sobolev space $W^{-s}_q(M)$ of $1$-dimensional currents for any $s\geq 1$ and the corresponding sub-bundles $\Gamma^{\pm}_{\kappa}$ are invariant under the action of the cocycle $\{\Phi^s_t\}_{t\in\R}$. Let $F_{q}:M\times \R \to M$ [$F_{-q}:M\times \R \to M$] be the (almost everywhere defined) flow of the horizontal [vertical] vector fields $S_q$ [$T_q$] on $M$. The horizontal and vertical foliations ${\Cal F}_{\pm q}$ coincide almost everywhere with the orbit foliations of the flows $F_{\pm q}$. Let $\gamma_{\pm q}^{\Cal T} \in \Gamma^{\pm}_q$ be a positively oriented segment with initial point $p^{\pm}\in M$ and let $\alpha:=f^+\eta_T +f^-\eta_S \in W^s_q(M)$. The following
identity holds: 
\begin{equation}
\gamma_{\pm q}^{\Cal T}(\alpha)=\int_0^{\Cal T}f^{\pm}\circ F_{\pm q}(p^{\pm},\tau) \,d\tau\,\,,
\end{equation}
The ergodic averages of the functions $f^{\pm}\in H^s_q(M)$ (for $s\geq 1$) can therefore be understood by studying the dynamics of the `renormalization' cocycles  $\{\Phi^s_t\}_{t\in\R}$ on 
the sub-bundle $\Gamma^{\pm}_{\kappa}$.
\smallskip
\noindent In \cite{F02} we have proved the following basic estimate on the weighted Sobolev norm 
of the currents $\gamma_{\pm q}^{\Cal T} \in W^{-1}_q(M)$. Let $R_q$ be the flat metric with conical singularities associated to the quadratic differential $q\in {\Cal M}^{(1)}_{\kappa}$ and let $|\!|q|\!|$ denote the $R_q$-length of the shortest saddle connection of the flat surface $(M,R_q)$. (We recall
that a {\it saddle connection} is a segment joining conical points).

\begin{lemma} \label{lemma:shortgammabound}(\cite{F02}, Lemma 9.2)
There exists a constant $K>0$ such that, for all quadratic differentials $q\in {\Cal M}^{(1)}_{\kappa}$,
\begin{equation}
|\gamma_{\pm q}^{\Cal T}|_{W^{-1}_q(M)} \,\leq\, K(1+ \frac{\Cal T}{|\!|q|\!|})\,\,. 
\end{equation}
\end{lemma}
\noindent The above estimate is a `trivial' bound, linear with respect to the length $\Cal T>0$ of the 
orbit segment, and can be quite easily derived from the Sobolev trace theorem for rectangles in the euclidean plane $\R^2$. The number $|\!|q|\!|>0$, which measures the pinching of the flat surface, gives 
the order of magnitude of the edges of the largest flat rectangle which can be embedded in the flat surface $(M,R_q)$ around an arbitrary regular point. 

\smallskip
\noindent We recall below the terminology and the notations, introduced in  \S 9 of \cite{F02}, 
concerning return trajectories of translation flows. 
\begin{definition}
A point $p\in M$ is \emph{regular }with respect to a measured foliation $\Cal F$ if it belongs to a regular (non-singular) leaf of $\Cal F$. For any quadratic differential $q\in {\Cal Q}(M)$. A point $p\in M$ will 
be said to be $q$-\emph{regular }if it is regular with respect to the horizontal and vertical foliations
${\Cal F}_{\pm q}$. 
\end{definition}
\noindent The set of $q$-regular points is of full measure and it is equivariant under the
actions of the group $\text{Diff}^+(M)$ and of the Teichm\"uller flow on ${\Cal Q}(M)$. 
\begin{definition}
Let $p\in M$ be a $q$-regular point and let $I_{\pm q}(p)$ be the vertical [horizontal] 
segment of length $|\!|q|\!|/2$ centered at $p$. A \emph{forward horizontal [vertical] return time }of 
$p\in M$ is a real number $\Cal T:={\Cal T}_{\pm q}(p)>0$ such that $F_{\pm q}(p,{\Cal T}) \in I_{\mp q}(p)$. If ${\Cal T}>0$ is a horizontal [vertical] return time for a $q$-regular point $p\in M$, 
the horizontal [vertical] forward orbit segment $\gamma_{\pm q}^{\Cal T}(p)$ 
will be called a \emph{forward horizontal [vertical] return trajectory }at $p$. 
\end{definition}
\noindent There is a natural map from the set of horizontal [vertical] return trajectories into the set 
of homotopically non-trivial closed curves.
\begin{definition}
\label{def:closing}
The \emph{closing }of any horizontal [vertical] trajectory segment $\gamma_{\pm q} ^{\Cal T}(p)$ is the
piece-wise smooth homotopically non-trivial closed curve 
\begin{equation}
\label{eq:closing}
\widehat{\gamma}_{\pm q}^{\Cal T}(p):=\gamma_{\pm q}^{\Cal T}(p)\cup\gamma
\bigl(p,F_{\pm q}(p,{\Cal T})\bigr)\,\,, 
\end{equation}
obtained as the union of the trajectory segment $\gamma_{\pm q} ^{\Cal T}(p)$ with the shortest geodesic segment $\gamma\bigl(p,F_{\pm q}(p,{\Cal T})\bigr)$ joining its endpoints $p$ and
$F_{\pm q}(p,{\Cal T})$.   
\end{definition}
\smallskip
\noindent Let ${\Cal T}_{\pm q}^{(1)}(p)>0$ be the \emph{forward horizontal [vertical] first return 
time }of a $q$-regular point $p\in M$, defined to be the real number
\begin{equation}
{\Cal T}_{\pm q}^{(1)}(p):=\min\{{\Cal T}>0\,|\, F_{\pm q}(p,{\Cal T})\in I_{\mp q}(p)\} \,\,. 
\end{equation}
The corresponding forward horizontal [vertical] trajectory $\gamma_{\pm q}^{(1)}(p)$ with initial 
point $p$ will be called the {\it forward horizontal [vertical] first return trajectory }at $p$. The following
bounds hold for first return times:
\begin{lemma} 
\label{lemma:firstreturn}
(\cite{F02}, Lemma 9.2') There exists a measurable function $K_r:{\Cal M}_{\kappa}\to {\R}^+$ such 
that, if ${\Cal T}_{\pm q}^{(1)}(p)$ is the forward horizontal [vertical] first return time of a $q$-regular 
point $p\in M$, then
\begin{equation}
\label{eq:firstreturn}
|\!|q|\!|/2\,\leq \, {\Cal T}_{\pm q}^{(1)}(p)\, \leq\, K_r(q)\,\,. 
\end{equation}
\end{lemma}
\noindent The lower bound in \pref{eq:firstreturn} is an immediate consequence of the definition
of $|\!|q|\!|$ as the length of the shortest saddle connection, while the upper bound depends essentially
on the fact that the first return map of a translation flow to a transverse interval is an interval exchange
transformation, hence the return-time function is piece-wise constant (and bounded). 

\smallskip
\noindent In \S 9.3 of \cite{F02} we have constructed special sequences of horizontal [vertical]
return times for almost all quadratic differentials generated. Such special return times are related
to visiting times of the Teichm\"uller flow, which `renormalizes' translation flows, to appropriate compact subsets of positive measure. Let $q\in {\Cal M}^{(1)}_{\kappa}$ be a Birkhoff generic point of the Teichm\"uller flow $\{G_t\}_{t\in\R}$ and let ${\Cal S}_{\kappa}(q)\subset {\Cal M}^{(1)}_{\kappa}$ be a smooth compact hypersurface of codimension $1$, such that  $q \in {\Cal S}_{\kappa}(q)$ and ${\Cal S}_{\kappa}(q)$ is transverse to the Teichm\"uller flow. Let $(t_k)$ be the sequence of visiting times 
of the orbit $\{G_t(q)\,|\,t\in {\R}\}$ to ${\Cal S}_{\kappa}(q)$. Since, by definition,  for any $t\in \R$,
\begin{equation}
({\Cal F}_{G_t(q)},{\Cal F}_{-G_t(q)})=(e^{-t}{\Cal F}_q, e^t {\Cal F}_{-q}) \,\,, 
\end{equation}
if $t:=t_k<0$ is a {\it backward }visiting time of the orbit $\{G_t(q)\}_{t\in\R}$, any forward first return trajectory of  the horizontal foliation ${\Cal F}_{G_t(q)}$ is a forward return trajectory of the horizontal foliation ${\Cal F}_q$, provided $|t_k|$ is sufficiently large. In a similar way, if $t:=t_k>0$ is a {\it 
forward }visiting time of $\{G_t(q)\}_{t\in\R}$, any forward first return trajectory of the vertical foliation 
${\Cal F}_{-G_t(q)}$ is a forward return trajectory of the vertical foliation ${\Cal F}_{-q}$. 

\smallskip
\noindent By the {\it closing }of the return trajectories of the horizontal [vertical] foliation, as in
\pref{eq:closing}, we obtain {\it closed }currents of Sobolev order $1$. The evolution of such currents under the action of the Teichm\"uller flow is therefore described by the cocycles $\{\Phi^s_t\}_{t\in \R}$ 
on the bundle ${\Cal Z}_{\kappa}^s(M)$ for any $s\geq 1$.  In the case $s=1$ the dynamics of
the cocycle $\{\Phi^1_t \vert {\Cal Z}_{\kappa}^1(M)\}$ is completely described by Theorem 
\ref{thm:Otype}. In the following we will analyze Lyapunov exponents of closed currents 
given by the closing of horizontal  [vertical] trajectories of translation flows under the cocycles 
$\{\Phi^s_t\}_{t\in \R}$ for $s>1$. A complete description of the Lyapunov structure of the cocycles
$\{\Phi^s_t \vert {\Cal Z}_{\kappa}^s(M)\}$ for $s>1$ is not relevant for our purposes and will not be
attempted.

\smallskip
\noindent Let $\mu$ be a $SO(2,\R)$-absolutely continuous, KZ-hyperbolic measure on a stratum $\Cal M^{(1)}_\kappa$ of orientable quadratic differentials (in the sense of Definitions \ref{def:circleabs} and \ref{def:nuhsimple}). We recall that according to \cite{F02} all the canonical absolutely continuous, $SL(2,\R)$-invariant measures on any connected component of any stratum of orientable quadratic differentials are KZ-hyperbolic. A different and stronger proof  which establishes that all the afore-mentioned measures  are KZ-simple has been given in \cite{AV} (see Theorem \ref{thm:simplicity}). The latter theorem is not necessary for any of the results of this paper to hold.

\smallskip
\noindent By Theorem \ref{thm:Otype}, for any $s\geq 1$, there exists a measurable splitting 
\begin{equation}
\label{eq:Zssplit}
\Cal Z^s_\kappa(M) =  \Cal B^1_{\kappa, +}(M)   \oplus \Cal B^1_{\kappa, -}(M) \oplus  \Cal E^s_{\kappa}(M) \,.
\end{equation}
The Lyapunov exponents of the restriction of the distributional cocycle $\{\Phi^s_t\}_{t\in\R}$ to the finite dimensional sub-bundles  $\Cal B^1_{\kappa, +}(M)$ [$\Cal B^1_{\kappa, -}(M)$] are equal to the
non-negative [non-positive] exponents of the Kontsevich-Zorich cocycle. It follows from the non-uniform hyperbolicity hypothesis for the latter that all Lyapunov exponents of the cocycle $\{\Phi^s_t\vert \Cal B^1_{\kappa, +}(M)\}$ are strictly positive while all Lyapunov exponents of  the cocycle $\{\Phi^s_t\vert \Cal B^1_{\kappa, -}(M)\}$ are strictly negative. In particular, by Oseledec's theorem there exists a measurable function $\Lambda^s_\kappa:\Cal M^{(1)}_\kappa \to \R^+$ and a real number $\alpha_\mu >0$ such that, for $\mu$-almost all $q\in \Cal M^{(1)}_\kappa$, for all $C\in \Cal B^1_{-q}(M)$ and all $t\geq 0$, 
\begin{equation}
\label{eq:negLexp}
\vert \Phi^s_t(C)\vert _{-s}   \,\leq  \, \Lambda^s_\kappa(q)  \vert C \vert_{-s} \, e^{-\alpha_\mu t} \,.
\end{equation}

\smallskip
\noindent In \cite{F02}, we have proved that for $s=1$ the restriction $\{\Phi^s_t\vert   \Cal E^s_{\kappa}(M)\}$ has $0$ has the unique Lyapunov exponent. In fact, this cocycle is isometric
with respect to a continuous Lyapunov norm (see Theorem \ref{thm:Otype}). For $s>1$, estimates on Lyapunov exponents and Lyapunov norms for the restriction $\Phi^s_t\vert   \Cal E^s_{\kappa}(M)$ will be derived from Theorem \ref{thm:G0gap}.  Let $\sigma:= s-1 >0$. For any $q\in \Cal M^{(1)}_\kappa$,
any $\lyap<1$ and any $\epsilon >0$, let 
\begin{equation}
\label{eq:Lboundfunct}
{\Cal N}_\kappa^{\sigma,\epsilon,\lyap}(q)\,:= \, \sup_{U\in {\Cal H}^0_q(M)} \, \frac{{\Cal N}^{\sigma,-\epsilon}_q (U)}
{\vert U\vert_0 + {\Cal N}^{1,\lyap}_q (S_qU)}\,.
\end{equation}
By Theorem \ref{thm:G0gap} for any $\lyap <1$  there exists $\epsilon:=\epsilon(\sigma,\lyap)>0$ such that  formula \pref{eq:Lboundfunct} defines a  (measurable) function 
${\Cal N}_\kappa^{\sigma,\epsilon,\lyap} :\Cal M^{(1)}_\kappa \to \R^+$ .

\smallskip
\noindent Let $\Pi_{\pm q} : \Cal Z^1_q(M) \to  \Cal B^1_{\pm q}(M)$ and $\Cal E_q: \Cal Z^1_q(M) \to  \Cal E^1_q(M)$ be the projections determined by the splitting \pref{eq:Zssplit}. For any $s\geq 1$, the
restrictions of  the projections $\Pi_{\pm q}$ and $\Cal E_q$ to the subspace $ \Cal Z^s_q(M) \subset  \Cal Z^1_q(M)$ can be extended to projections $\Pi^s_{\pm q}: W^{-s}_q(M) \to \Cal B^1_{\pm q}(M)$ and $\Cal E^s_q: W^{-s}_q(M) \to \Cal E^s_q(M)$, defined on the Sobolev space $W^{-s}_q(M)$
of $1$-dimensional currents, by composition with the orthogonal projection $W^{-s}_q(M)\to
\Cal Z^s_q(M)$ onto the closed subspace of closed currents. Let $d^s_\kappa: \Cal M^{(1)}_\kappa
\to \R^+$ be the {\it distorsion }of the splitting \pref{eq:Zssplit}, that is, the function defined for $\mu$-almost all $q\in \Cal M^{(1)}_\kappa$ as
\begin{equation}
\label{eq:distorsion}
d^s_\kappa (q)\, := \, \sup_{C\in  \Cal Z^s_q(M) }   \, 
\frac{ \vert \Pi_{q}(C)\vert _{-s} + \vert \Pi_{-q}(C)\vert_{-s} +  \vert \Cal E_q(C) \vert_{-s}} 
{ \vert C \vert_{-s}}\,.
\end{equation}
\smallskip
\noindent Let $s>1$ and $\lyap<1$.  Let ${\Cal P}^{(1)}_{\kappa}\subset{\Cal M}^{(1)}_{\kappa}$ be a compact set satisfying the following conditions:
\begin{enumerate}
\item All $q\in {\Cal P}^{(1)}_{\kappa}$ are Birkhoff generic points for the Teichm\"uller flow $\{G_t\}$ and Oseledec regular points for the cocycle $\{\Phi^s_t \vert \Cal B^1_{\kappa}(M)\}$;
\item the set ${\Cal P}^{(1)}_{\kappa}$ is transverse to the Teichm\"uller flow and has positive 
transverse measure;
\item there exists a constant $K^{(1)}_{r}>0$ such that, for all $q\in {\Cal P}^{(1)}_{\kappa}$ and all $q$-regular $p\in M$, the first return times $\Cal T^{(1)}_{\pm q}(p) \leq K^{(1)}_{r}$;
\item the real-valued functions $\Lambda^s_\kappa$,  ${\Cal N}_\kappa^{\sigma,\epsilon,\lyap}$ and  $\delta^1_{\kappa}$, introduced respectively in formulas \pref{eq:negLexp}, \pref{eq:Lboundfunct} and  \pref{eq:distorsion}, are bounded on ${\Cal P}^{(1)}_{\kappa}$.
\end{enumerate}
\smallskip
\noindent By the ergodicity of the system $(\{G_t\}, \mu)$, by the Oseledec's theorem for the cocycle 
$\Phi^s_t \vert \Cal B^1_{\kappa}(M)$ and by Lemma 4.6, it follows that the union of all sets ${\Cal P}^{(1)}_{\kappa}$ with the properties $(1)-(4)$ is a full measure subset of ${\Cal M}^{(1)}_{\kappa}$. 

\begin{definition}
Let $q\in {\Cal M}^{(1)}_{\kappa}$ and $(t_k)_{k\in {\N}}$ be the sequence of visiting 
times of the {\it backward }orbit $\{G_t(q)\,|\,t\leq t_1=0\}$ to a compact positive measure set ${\Cal P}^{(1)}_{\kappa}\subset {\Cal M}^{(1)}_{\kappa}$.
A \emph{principal sequence } $( {\Cal T}^{(k)}_q(p))_{k\in\N}$ of forward return times for the horizontal foliation ${\Cal F}_q$ at a $q$-regular point $p\in M$ is the sequence
\begin{equation}
\label{eq:principalrt}
{\Cal T}^{(k)}_q(p):= {\Cal T}^{(1)}_{G_t(q)}(p) \exp{|t|}\,\,,\,\,\,\,t=t_k \,\,. 
\end{equation}  
For each $k\in \N$, a \emph{horizontal principal (forward) return trajectory }$\gamma_q^{(k)}(p)$ at a $q$-regular point $p\in M$ is the horizontal forward return trajectory at the point $p$ corresponding to a principal return time ${\Cal T}^{(k)}_q(p)>0$. 
\end{definition}
\noindent  We remark that a horizontal principal return trajectory $\gamma_q^{(k)}(p)$ 
coincides with the horizontal {\it first }return trajectory at $p$ of the quadratic differential 
$G_t(q)$, $t=t_k<0$. A similar definition can be given for the vertical foliation $\Cal F_{-q}$
by considering forward visiting times of the Teichm\"uller flow.

\smallskip
\noindent The following standard splitting lemma allows to reduce the analysis
of arbitrary regular trajectories to that of principal return trajectories.
 
\begin{lemma} \label{lemma:trajsplit} (\cite{F02}, Lemma 9.4) Under conditions $(1)-(3)$ there exists a constant $K_{\Cal P}>0$ such that the following holds. Let $q\in {\Cal P}^{(1)}_{\kappa}$ and, for any ${\Cal T}>0$, let $\gamma^{\Cal T}_q(p)$ be a forward trajectory at a $q$-regular point $p\in M$. There exists a finite 
set of points $\{ p^{(k)}_j \,\vert \,1\leq k\leq n,\,,1\leq j\leq m_k \} \subset\gamma^{\Cal T}_q(p)$, 
such that the principal return trajectories $\gamma_q^{(k)}(p^{(k)}_j)\subset \gamma^{\Cal T}_q(p)$ 
do not overlap and 
\begin{equation}
\label{eq:trajsplit}
\begin{aligned}
 & \gamma^{\Cal T}_q(p)=\sum_{k=1}^n\sum_{j=1}^{m_k}\gamma_q^{(k)}(p^{(k)}_j)\,\,
              +\,\,b^{\Cal T}_q(p)\,\,,\\ 
            & \text{ with }m_k\leq K_{\Cal P}\, e^{|t_{k+1}|-|t_k|}\,\,, \text{ for }1\leq k\leq n\,,  \\
            & \text{and the $R_q$-length } L_q\left( b^{\Cal T}_q(p)\right) \leq K_{\Cal P}\,\,.
\end{aligned}
\end{equation}
\end{lemma}
\begin{proof} We recall the proof given in \cite{F02}, Lemma 9.4. The argument is based on the following estimate on principal return times. By \pref{eq:principalrt},  Lemma \ref{lemma:firstreturn} and condition $(3)$, there exists a constant $K_{pr}>0$ such that, for all $q\in {\Cal P}^{(1)}_{\kappa}$, all $q$-regular points $p\in M$ and all $k\in \N$,
\begin{equation}
\label{eq:prtunifbound}
K_{pr}^{-1}\exp{|t_k|}\leq {\Cal T}^{(k)}_q(p) \leq K_{pr}\exp{|t_k|}\,\,.
\end{equation}
Let $n=\max\{k\in {\N}\,|\, {\Cal T}^{(k)}_q(p)\leq {\Cal T}\,\}$. The maximum exists (finite) by
\pref{eq:prtunifbound}. Let $p^{(n)}_1:=p$. The sequence $(p^{(k)}_j)$ with the properties stated in 
\pref{eq:trajsplit} can be constructed by a finite iteration of the following procedure. Let $p^{(k)}_j$ be 
the last point already determined in the sequence and let 
\begin{equation}
p^{(k)}_{j+}:=F_q\bigl(p^{(k)}_j,{\Cal T}_q^{(k)}(p^{(k)}_j)\bigr)\in \gamma^{\Cal T}_q(p)\,\,.
\end{equation}
Let then $k'\in \{1,\dots ,k\}$ be the largest integer such that 
\begin{equation}
F_q\bigl(p^{(k)}_{j+},{\Cal T}_q^{(k')}(p^{(k)}_{j+})\bigr)\in \gamma^{\Cal T}_q(p)\,\,. 
\end{equation}
If $k'=k$, let $p^{(k)}_{j+1}:=p^{(k)}_{j+}$. If $k'<k$, let $m_k:=j$, $m_h=0$ (no points) for all 
$k'<h<k$ and $p^{(k')}_1:=p^{(k)}_{j+}$. The iteration step is concluded. By \pref{eq:prtunifbound} it
follows that
\begin{equation}
\label{eq:prtrecbound}
K_{pr}^{-1}\,e^{|t_k|}\,m_k \leq {\Cal T}^{(k+1)}_q(p^{(k)}_1)\leq  K_{pr}\,e^{|t_{k+1}|} \,\,. 
\end{equation}
The length of the remainder $b^{\Cal T}_q(p)$ has to be less than any upper bound for the length of first return times. Hence, by \pref{eq:prtrecbound}, Lemma \ref{lemma:firstreturn} and condition $(3)$, the estimates in \pref{eq:trajsplit} are proved and the argument is concluded.
\end{proof}

\noindent The result below gives a fundamental uniform estimate on ergodic integrals.

\begin{theorem} 
\label{thm:GHbound}
Let $\mu$ be an $SO(2,\R)$-absolutely continuous, KZ-hyperbolic measure on a stratum $\Cal M^{(1)}_\kappa$ of orientable quadratic differentials. For any $s>1$ there exists a measurable function $\Gamma^s_\kappa: \Cal M^{(1)}_\kappa \to \R^+$  such that, for $\mu$-almost $q\in \Cal M^{(1)}_\kappa$, for all $q$-regular $p\in M$ and for all $\,\Cal T >0\,$,

\begin{equation}
\label{eq:GHbound}
\vert \Pi^s_{-q}\left( \gamma_q^{\Cal T}(p) \right)\vert_{-s} \, +\, 
\vert {\Cal E}^s_{q}\left( \gamma_q^{\Cal T}(p) \right)\vert_{-s}\,\leq \,   \Gamma^s_\kappa(q)\,.
\end{equation}
\end{theorem}
\begin{proof} Let $\Cal P:={\Cal P}^{(1)}_{\kappa}$ be a compact set satisfying conditions $(1)-(4)$ 
listed above. For $\mu$-almost all $q \in {\Cal P}$, there exists a sequence $(t_k)_{k\in\N}$ 
of backward return times of the Teichm\"uller orbit $G_t(q)$ to the set ${\Cal P}$. By Lemma
 \ref{lemma:trajsplit} the uniform estimate \pref{eq:GHbound} can be reduced 
(exponential) estimates for principal return trajectories $\gamma^{(k)}_q(p)$.  We claim that for
every $s>1$ there exist constants $K^s_\Cal P>0$ and $\alpha_s >0$ such that for 
$\mu$-almost all $q\in \Cal P$, for all $q$-regular $p\in M$ and all $k\in \N$, 
\begin{equation}
\label{eq:prtexpbounds}
\vert \Pi^s_{-q}\left( \gamma^{(k)}_q(p) \right)\vert_{-s} \, + \vert {\Cal E}^s_{-q}\left( \gamma^{(k)}_q(p)\right)\vert_{-s}  \, \leq \, K^s_\Cal P \, e^{-\alpha_s \vert t_k\vert} \,.
\end{equation}
 By Lemma \ref{lemma:shortgammabound}, by conditions $(3)$ and $(4)$ on the set $\Cal P$, in particular the bound on the distorsion, there exists a constant $K^s_1(\Cal P)>0$ such that the following 
holds. For $\mu$-almost all $q\in {\Cal P}$ and all $q$-regular $p\in M$, the closing ${\widehat\gamma}^{(1)}_q(p)$ of the  first return horizontal trajectory $\gamma^{(1)}_q(p)$ 
(see Definition \ref{def:closing}) satisfies the following bound:
\begin{equation}
\label{eq:initialbound}
 |\Pi_{-q}\bigl({\widehat\gamma}^{(1)}_q(p)\bigr)|_{-1} \, +\, |{\Cal E}_q\bigl({\widehat\gamma}^{(1)}_q(p)\bigr)|_{-1} \leq K_1(\Cal P) 
 \end{equation}
 Since by definition of the principal return trajectories
 \begin{equation}
 \label{eq:dprt}
{\widehat\gamma}^{(k)}_q(p)=\Phi^s_{\vert t \vert}\bigl({\widehat\gamma}^{(1)}_{G_t(q)}(p)\bigr)\,\,,\,\,\, \text{ \rm for any } \,t=t_k\leq 0\,\,, 
\end{equation}
it follows from the initial bound \pref{eq:initialbound}, from the invariance of the sub-bundle $\Cal B^1_{\kappa,-}(M)$ under the cocycles $\{\Phi^s_t\}$, from the Lyapunov bound \pref{eq:negLexp} on the cocycle $\{\Phi^s_t\vert  \Cal B^1_{\kappa,-}(M)\}$ and from condition $(4)$ on the set $\Cal P\subset \Cal M^{(1)}_\kappa$ that there exists a constant $K_2^s(\Cal P) >0$ such that the following estimate holds. For $\mu$-almost all $q\in \Cal P$ and for all $q$-regular $p\in M$,
\begin{equation}
\label{eq:Piest}
\vert \Pi_{-q}\left( {\widehat\gamma}^{(k)}_q(p) \right)\vert_{-s} \, \leq \, K_2^s(\Cal P)\, 
e^{-\alpha_\mu \vert t_k\vert} \,, \quad \text{ \rm for all } k\in \N\,.
\end{equation}
\noindent A similar estimate holds for the projections of principal return trajectories on the
sub-bundle $\Cal E^1_\kappa(M)$ of exact currents. In fact, there exist $K_3^s(\Cal P) >0$ and $\epsilon_s>0$ such that, for $\mu$-almost all $q\in \Cal P$ and for all $q$-regular $p\in M$,
\begin{equation}
\label{eq:Eest}
\vert \Cal E_{q}\left( {\widehat\gamma}^{(k)}_q(p) \right)\vert_{-s} \, \leq \, K_2^s(\Cal P)\, 
e^{-\epsilon_s \vert t_k\vert} \,, \quad \text{ \rm for all } k\in \N\,.
\end{equation}
The estimate \pref{eq:Eest} is proved as follows. By definition of the bundle $\Cal E^1_\kappa(M)$
of exact currents, for all $q\in \Cal M^{(1)}_\kappa$ and all $q$-regular $p\in M$, there exists a unique function $U^{(1)}_q(p) \in {\Cal H}^0_q(M) \subset L^2_q(M)$ such that
\begin{equation}
\label{eq:Eid}
\Cal E_q \left({\widehat\gamma}^{(1)}_q(p) \right) =  d U^{(1)}_q(p) \,.
\end{equation}
By the definitions of the distributional cocycles $\{\Phi^s_t\}$ and $\{G^s_t\}$ the following 
identity of cocycles holds for any $s\geq 1$:
\begin{equation}
\label{eq:Ecocycleid}
\{ \Phi^s_t \circ d \} = \{ d \circ G^{s-1}_t\}     \quad \text { \rm on } \, H^{-s+1}_q(M) \,.
\end{equation}
The exponential estimate \pref{eq:Eest}  will therefore follow from Theorem \ref{thm:G0gap} if we can prove that there exist a positive number $\lyap<1$ and a constant $K_3(\Cal P)>0$ such that, for $\mu$-almost all $q\in \Cal P$ and all $q$-regular $p\in M$, the following holds:
\begin{equation}
\label{eq:initialboundbis}
\vert U^{(1)}_q(p) \vert_0 \, + \, {\Cal N}^{1,\lyap}_q\left[S_qU^{(1)}_q(p)) \right] \leq  K_3^s(\Cal P)\,.
\end{equation}
Let $\vert \cdot \vert_\kappa$ be the norm on the bundle $\Cal E^1_\kappa(M)$ of exact currents defined as follows: for any $q\in \Cal M^{(1)}_\kappa$ and $C\in \Cal E^1_q(M)$,
\begin{equation}
\label{eq:Enorm}
\vert C \vert_\Cal E :=   \vert U \vert_{L^2_q(M)} \,, \quad  \text{ \rm if }\, \, C=dU\,\,  \text{ \rm with }\,\, U \in {\Cal H}^0_q(M)\,.
\end{equation}
By the definitions of weighted Sobolev norms, it is immediate to prove that the map $d: {\Cal H}^0_q(M) \to 
\Cal E^1_q(M)$ is a continous bijective map of Hilbert spaces. Hence, by definition \pref{eq:Enorm}
and by the open mapping theorem, for each $q\in \Cal M^{(1)}_\kappa$ there exists $K(q)>0$ such that
\begin{equation}
\label{eq:Enorms}
K(q)  \vert \cdot \vert_\Cal E \, \leq \, \vert \cdot \vert_{-1} \,\leq\,\vert \cdot \vert_\Cal E 
\quad \text{ \rm on }\,\, \Cal E^1_q(M)\,.
\end{equation}
Since the sub-bundle ${\Cal H}^0_\kappa(M)  \subset H_\kappa^0(M)$ and the norms $\vert \cdot \vert_\Cal E$, $\vert \cdot \vert_{-1}$ on the bundle $\Cal E^1_\kappa(M)$ are all continuous, the function 
$K:\Cal M^{(1)}_\kappa \to \R^+$ can be chosen continuous (see Lemma 9.3 in \cite{F02}
for a detailed argument).  Hence by the estimates \pref{eq:initialbound} and by identity \pref{eq:Eid}, there exists a constant $C_1(\Cal P)>0$ such that, for $\mu$-almost all $q\in \Cal P$ and all $q$-regular $p\in M$,
\begin{equation}
\label{eq:L2initialbound}
 \vert U^{(1)}_q(p) \vert_{L^2_q(M)} \, \leq \, K(q)^{-1} \vert \Cal E_q \left({\widehat\gamma}^{(1)}_q(p) \right)  \vert_{-1} \leq C_1(\Cal P)\,.
 \end{equation}
 The proof of the estimate \pref{eq:initialboundbis} is completed as follows. Since $\gamma^{(1)}_q(p)$
 is a horizontal trajectory, by the splitting \pref{eq:Zssplit} and by the identity \pref{eq:Eid}, the following formula holds:
 \begin{equation}
 \label{eq:SUid}
 S_q U^{(1)}_q(p) =  \imath_{S_q} [{\widehat\gamma}^{(1)}_q(p) - \gamma^{(1)}_q(p)] \,-\,
 \imath_{S_q} \Pi_{-q}[{\widehat\gamma}^{(1)}_q(p)]\,.
 \end{equation}
Since the distribution $\imath_{S_q}[{\widehat\gamma}^{(1)}_q(p) - \gamma^{(1)}_q(p)]$ is given by integration along a vertical arc of unit length, by the Sobolev embedding Lemma 
\ref{lemma:shortgammabound} and by the logarithmic law of geodesics for the Teichm\"uller geodesic flow on moduli spaces (see \cite{Ma93}, Prop. 1.2, or  \S \ref{KZcocycle},  formula \pref{eq:Loglaw}) it follows that for any $\lyap>0$ there exists a constant $C_\lyap>0$ such that, for $\mu$-almost all $q\in \Cal M^{(1)}_\kappa$ and all $q$-regular $p\in M$, the Lyapunov norm
\begin{equation}
\label{eq:Lyapnormone}
{\Cal N}_q^{1,\lyap} \left( \imath_{S_q} [{\widehat\gamma}^{(1)}_q(p) - \gamma^{(1)}_q(p)] \right)\,\leq \, C_\lyap\,.
\end{equation}
By the cocycle isomorphism \pref{eq:cocycleiso}, by the Lyapunov estimate \pref{eq:negLexp} and
 by condition $(4)$ on the set $\Cal P$, since the distribution $ \imath_{S_q} \Pi_{-q} [{\widehat\gamma}^{(1)}_q(p)] \in \Cal I^1_{-q}(M)$, there exists a constant $C_2(\Cal P)>0$ such that, for any $1>\lyap \geq 1-\alpha_\mu$, for $\mu$-almost all $q\in \Cal M^{(1)}_\kappa$ and all $q$-regular $p\in M$, the Lyapunov norm
 \begin{equation}
 \label{eq:Lyapnormtwo}
{\Cal N}_q^{1,\lyap} \left( \imath_{S_q} \Pi_{-q}[{\widehat\gamma}^{(1)}_q(p)]  \right)\,\leq \, 
C_2(\Cal P)\,.
\end{equation}
The bound \pref{eq:initialboundbis} follows immediately from the bound \pref{eq:L2initialbound},
from the identity \pref{eq:SUid} and from the bounds \pref{eq:Lyapnormone} and \pref{eq:Lyapnormtwo}.

\smallskip
\noindent
As hinted above, the estimate \pref{eq:initialboundbis} implies the required exponential estimate 
\pref{eq:Eest} on the projections of principal return trajectories on the space of exact currents. In fact, by Theorem \ref{thm:G0gap} and by condition $(4)$ on the set $\Cal P$, for any $1>\lyap\geq 1-\alpha_\mu$ and any $\sigma>0$, there exist $K_4^\sigma(\Cal P)>0$ and $\epsilon_s:=\epsilon(\lyap, \sigma)>0$ such that, for $\mu$-almost all $q\in\Cal P$ and all $q$-regular $p\in M$,
\begin{equation}
\label{eq:Uexpbound}
{\Cal N}_q^{\sigma,\epsilon} \left (U^{(1)}_q(p)\right) \, \leq \, K_4^s(\Cal P)\,.
\end{equation}
Let us introduce the following notation: for any $k\in \N$, let
\begin{equation}
\label{eq:qkUknot}
q_k := G_{t_k}(q) \in \Cal P \quad \text { \rm and } \quad U^{(k)}_q(p):= G_{\vert t_k \vert}^0 
\left ( U^{(1)}_{q_k}(p) \right) \in L^2_q(M)\,.
\end{equation}
 Since the splittings \pref{eq:Zssplit} are $\{\Phi^s_t\}$-invariant (in particular for $s=1$), by the ddefinition of principal return trajectories \pref{eq:dprt}, by the identity \pref{eq:Eid} and the cocycle identity \pref{eq:Ecocycleid}, the following holds:
\begin{equation}
\begin{aligned}
\Cal E_q \left(  {\widehat\gamma}^{(k)}_q(p) \right) &= \Cal E_q\circ \Phi^1_{\vert t_k\vert} 
\left(  {\widehat\gamma}^{(1)}_{q_k} (p) \right)= \\ &=  \Phi^1_{\vert t_k\vert} \left( dU^{(1)}_{q_k}(p) \right)
= d U^{(k)}_q(p)
\end{aligned}
\end{equation}
Let $s>1$ and $\sigma:=s-1$. It follows from \pref{eq:Uexpbound} (by the definition \pref{eq:distLnorm} of Lyapunov norms) that, for $\mu$-almost all $q\in \Cal P$ and all $q$-regular $p\in M$,
\begin{equation}
\label{eq:Uexpboundbis}
\vert \Cal E_q \left(  {\widehat\gamma}^{(k)}_q(p) \right) \vert_{-s} 
\leq \vert G_{\vert t_k \vert}^0 \left ( U^{(1)}_{q_k}(p) \right)\vert_{-\sigma}\, \leq \, K_4^s(\Cal P) \,e^{-\epsilon_s \vert t_k\vert }\,.
\end{equation}
The crucial exponential estimate \pref{eq:prtexpbounds} on projections of principal return trajectories, claimed above, can now be derived from estimates \pref{eq:Piest} and \pref{eq:Eest} together with
the remark that, for any $\lyap>0$, there exists a constant  $C_\lyap>0$ such that, for $\mu$-almost all $q\in \Cal P$ and all $q$-regular $p\in M$,
\begin{equation}
\vert \gamma^{(k)}_q(p) - {\widehat\gamma}^{(k)}_q(p) \vert_{-1} \, \leq \, C_\lyap \,e^{-\lyap \vert t_k \vert}\,.
 \end{equation}
By the definition of principal return trajectories and closing, the latter estimate follows immediately from \pref{eq:Lyapnormone}, since $\gamma^{(k)}_q(p) - {\widehat\gamma}^{(k)}_q(p)$ are currents of integration along a vertical $q$-regular  arc.

\smallskip
\noindent The required estimate \pref{eq:GHbound} finally follows from estimate \pref{eq:prtexpbounds} by the trajectory splitting Lemma \ref{lemma:trajsplit}. In fact, for any $s>1$ there exist positive real numbers $\omega^{(1)}_s < \omega^{(2)}_s$ and a measurable map $k_s: \Cal P \to \N$ such that $\alpha_s\omega^{(1)}_s > \omega^{(2)}_s - \omega^{(1)}_s$ and, for $\mu$-almost all $q \in \Cal P$,
\begin{equation}
\label{eq:alphaom12}
\omega^{(1)}_s \, k  \, \leq \,  \vert t_k(q)  \vert   \, \leq \,  \omega^{(2)}_s \, k \,, \quad
 \text{ \rm for all }\, k\geq k_s(q) \,.
\end{equation}
Let $\{ p^{(k)}_j \vert 1\leq k\leq n\,, 1\leq j \leq m_k\} \subset \gamma^{\Cal T}(p)$ be the sequence of $q$-regular points constructed in Lemma \ref{lemma:trajsplit}. Since $\alpha_s\omega^{(1)}_s -
(\omega^{(2)}_s - \omega^{(1)}_s)>0$, by estimates \pref{eq:prtexpbounds} and \pref{eq:alphaom12} , there exists a constant $K^s_5(\Cal P)>0$ such that, for $\mu$-almost all $q \in \Cal P$, all $q$-regular $p\in M$
and all $\Cal T>0$, 
\begin{equation}
\label{eq:PiEfinest1}
\begin{aligned}
\sum_{k=k_s(q)}^n e^{\vert t_{k+1}\vert - \vert t_k\vert }\, \vert \Pi^s_{-q}\left ( \gamma_q^{(k)}(p^{(k)}_j) \right)\vert_{-s}   \,\, &\leq \,\, K^s_5(\Cal P)\,; \\
\sum_{k=k_s(q)}^ne^{\vert t_{k+1}\vert - \vert t_k\vert }\, \vert {\Cal E}^s_{q} \left ( \gamma_q^{(k)}(p^{(k)}_j) \right)\vert_{-s}   \,\, &\leq \,\, K^s_5(\Cal P)\,.
\end{aligned}
\end{equation}
In addition, by Lemma \ref{lemma:trajsplit} there exists a constant $K_6(\Cal P)>0$ such that 
\begin{equation}
\label{eq:PiEfinest2}
L_q \bigl[  \sum_{k=1}^{k_s(q)}  \sum_{j=1}^{m_k} \gamma_q^{(k)}(p_k^{(j)}) \bigr] \,\leq \, K_6(\Cal P)\,
\sum_{k=1}^{k_s(q)}  e^{ \vert t_{k+1}(q) \vert } \,,
\end{equation}
hence the estimate \pref{eq:GHbound} follow from estimates \pref{eq:PiEfinest1}, from
the trajectory splitting Lemma  \ref{lemma:trajsplit}, from estimate \pref{eq:PiEfinest2} and finally 
from Lemma \ref{lemma:shortgammabound}, which yields a bound on weighted Sobolev norms of (horizontal and vertical) trajectories in terms of their $R_q$-length.
\end{proof}
\subsection{The generic case} \label{ss:gencase}
The above Theorem \ref{thm:GHbound} implies, by a standard Gottschalk-Hedlund argument, the following sharp result on the existence of a {\it Green operator }for the (horizontal) cohomological equation in the case of generic orientable quadratic differentials.
\begin{theorem} 
\label{thm:CEsharp1}
For any $SO(2,\R)$-absolutely continuous, KZ-hyperbolic measure $\mu$ on a stratum $\Cal M^{(1)}_\kappa$ of orientable quadratic differentials and for $\mu$-almost all $q\in \Cal M^{(1)}_\kappa$, there exists a densely defined, anti-symmetric Green operator ${\Cal U}_q:{\Cal H}^0_q(M) \to {\Cal H}^0_q(M)$ for the horizontal cohomological equation $S_qu=f$. For any $s>1$, the maximal domain of the Green operator $\,{\Cal U}_q\,$ on the Hilbert space ${\Cal H}^0_q(M)\subset L^2_q(M)$ contains the dense subspace 
\begin{equation}
\label{eq:CEsharp1Greend}
[\Cal I^s_q(M)]^\perp:= \{ f\in H^s_q(M) \vert \Cal D(f) =0\, \text{ for all }\Cal D \in \Cal I^s_q(M)\}\,,
\end{equation}
and there exists a measurable function $C^s_\kappa:\Cal M^{(1)}_\kappa \to \R^+$  such that the following holds. For any $f\in [\Cal I^s_q(M)]^\perp$, the Green solution  ${\Cal U}_q(f) \in L^{\infty}(M)$ 
has zero-average and satisfies the estimate:
\begin{equation}
\label{eq:CEsharp1}
\vert {\Cal U}_q(f) \vert _{L^{\infty}(M)}   \, \leq  \,  C^s_\kappa(q)\, \vert f \vert_s \,.
\end{equation}
\end{theorem}
\begin{proof} Let $\{u_{\Cal T}\}_{\Cal T \in \R}$ be the $1$-parameter family of measurable 
functions defined almost everywhere on $M$ as follows: for all $q$-regular $p\in M$,
\begin{equation}
\label{eq:uT}
u_{\Cal T}(p) :=  \frac{1}{\Cal T} \, \int_0^{\Cal T} \int_0^{\tau}   f\circ F_q(p,s) \, ds \,d\tau \,.
\end{equation}
By Theorem  \ref{thm:GHbound}  (and Lemma \ref{lemma:shortgammabound}) for any $s>1$ there
exists a measurable function $C^s_\kappa:\Cal M^{(1)}_\kappa \to \R^+$ such that, for $\mu$-almost all $q\in \Cal M^{(1)}_\kappa$, for all $q$-regular $p\in M$ and all $\Cal T >0$,
\begin{equation}
\label{eq:CEsharpest1}
\vert   \gamma_q^{\Cal T}(p) -  \Pi^s_q\left (\gamma_q^{\Cal T}(p)\right) \vert_{-s} \,\leq \, C^s_\kappa(q)\,.
\end{equation}
We recall that the projection $\Pi^s_q: W^{-s}_q(M) \to \Cal B^1_q(M)$ is defined as the composition of the orthogonal projection $W^{-s}_q(M)  \to \Cal Z^s_q(M)$ with the projection $\Cal Z^s_q(M) \to \Cal B^1_q(M)$ determined by the splitting \pref{eq:Zssplit}.

\smallskip
\noindent 
If $f \in [\Cal I^s(M)]^\perp$, the $1$-form $\alpha_f:= f\eta_T \in W^s_q(M)$ is such that $C(\alpha)=0$ for any horizontally basic current $C\in \Cal B^1_q(M)$. It follows from the estimate \pref{eq:CEsharpest1} that, for $\mu$-almost all $q\in \Cal M^{(1)}_\kappa$, for all $q$-regular $p\in M$ and all $\Cal T >0$,
\begin{equation}
\label{eq:Linftybound}
\vert u_{\Cal T}(p) \vert \leq  \vert \gamma_q^{\Cal T}(\alpha) \vert  =
\vert \{\gamma_q^{\Cal T}(p) -  \Pi^s_q\left (\gamma_q^{\Cal T}(p)\right)\} (\alpha)\vert
 \leq  C^s_\kappa (q) \, \vert f \vert_s \,,
\end{equation}
hence the family  $\{u_{\Cal T}\}_{\Cal T \in \R}$ is uniformly bounded in the Hilbert space $L^2_q(M)$ for $\mu$-almost all $q\in \Cal M^{(1)}_\kappa$. In addition, a computation shows that if the horizontal flow $F_q$ is ergodic, as $\Cal T \to +\infty$, 
\begin{equation}
S_q u_{\Cal T} = f  \,-\,  \frac{1}{\Cal T} \, \int_0^{\Cal T}   f\circ F_q(\cdot,\tau) \, d\tau \,\,
\to \,\, f \quad \text{ in } \,\, L^2_q(M) \,,
\end{equation}
hence any weak limit $u\in L^2_q(M)$ of the family $\{u_{\Cal T}\}_{\Cal T \in \R}$ is a 
solution of the cohomological equation $S_qu=f$.  Since any function $f\in [\Cal I^s(M)]^\perp$ has zero
average, it follows by the definition \pref{eq:uT} that $u_\Cal T$ has zero average for all $\Cal T>0$,
hence any weak limit of the family $\{u_{\Cal T}\}_{\Cal T \in \R}$ has zero average. By the ergodicity 
of the horizontal flow $F_q$, the cohomological equation $S_qu=f$ has a unique zero average solution in $\Cal U_q(f)\in L^2_q(M)$. It follows that the operator $\Cal U_q : [\Cal I^s_q(M)]^\perp \to {\Cal H}_q^0(M)\subset L^2_q(M)$ is well-defined and linear, for $\mu$-almost all $q\in \Cal M^{(1)}_\kappa$. In addition, by the uniform bound \pref{eq:Linftybound} the function $\Cal U_q(f)\in L^{\infty}(M)$ and satisfies the required bound \pref{eq:CEsharp1}. 

\smallskip
\noindent It remains to be proven that the linear operator $\Cal U_q: {\Cal H}_q^0(M) \to {\Cal H}_q^0(M)$ is
anti-symmetric. For any $f\in [\Cal I^s_q(M)]^\perp$ and $v\in {\Cal H}^{s+1}_q(M)\subset H^{s+1}_q(M)$, since $\Cal U_q(f)\in L^2_q(M)$ is a weak-solution of the equation $S_q u=f$, 
\begin{equation}
\label{eq:Greenantisymm1}
\<\Cal U_q(f), S_q v \>_q = -\<f, v\>_q = - \<f, \Cal U_q(S_q v )\>_q \,.
\end{equation}
Since the subspace $\{ S_q v \in H^s_q(M) \vert v \in H^{s+1}_q(M) \}$ is dense in the subspace 
$[\Cal I^s_q(M)]^\perp \subset H^s_q(M)$ and $\Cal U_q: [\Cal I^s_q(M)]^\perp \to L^2_q(M)$ is continuous, it follows from identity \pref{eq:Greenantisymm1} by a density argument that
\begin{equation}
\label{eq:Greenantisymm2}
\<\Cal U_q(f), g \>_q = - \<f, \Cal U_q(g )\>_q \,, \quad \text{ \rm for all } f, \,g \in [\Cal I^s_q(M)]^\perp\,.
\end{equation}
\end{proof}

\medskip
\noindent Optimal results on the regularity of solutions of the cohomological equation for intermediate `fractional' regularity require smoothing and interpolation techniques in the presence of distributional obstructions. The following definition appears to capture the relevant condition on the obstructions which allow for effective smoothing and interpolation results.

\begin{definition}
\label{def:regularsystem} 
Let  $\{H^s\vert s\geq 0\}$ be a $1$-parameter family of normed spaces such that the following embeddings hold:
\begin{equation}
\label{eq:normedspaces}
H^\infty := \bigcap_{s>0} H^s \subset   H^s  \subset H^r \subset H^0  \,, \quad \text{ \rm for all } \, s\geq r \,.
\end{equation}
The \emph{order }${\Cal O}^H(\Cal D)$ (with respect to $\{H^s\vert s\geq 0\}$)  of any linear functional 
$\Cal D \in (H^\infty)^\ast$ is the non-negative real number
\begin{equation}
\label{eq:order}
{\Cal O}^H(\Cal D) := \inf \{ s\geq 0 \,\vert \,   \Cal D \in (H^s)^\ast \}\,.
\end{equation}
A finite system of linear functionals $\{\Cal D_1, \dots, \Cal D_J\} \subset (H^\sigma)^\ast$ will be called 
\emph{$\sigma$-regular }(with respect to the family $\{H^s\vert s\geq 0\}$) if, for any $\tau\in (0,1]$ there exists a \emph{dual }system $\{u_1(\tau), \dots, u_J(\tau)\} \subset H^\sigma$ such that the following estimates hold. For all $0\leq r  \leq \sigma$  and all $\epsilon>0$, there exists a constant $C^\sigma_r(\epsilon)>0$ such that, for all $i$, $j\in \{1,\dots, J\}$,
\begin{equation}
\label{eq:regularsystem}
\vert u_j(\tau) \vert_r \leq C^\sigma_r(\epsilon)\, \tau^{{\Cal O}^H(\Cal D_j)-r-\epsilon}\,.
\end{equation}
A finite system  $\{\Cal D_1, \dots, \Cal D_J\} \subset (H^s)^\ast$ will be called \emph{regular} if it is $\sigma$-regular  for any $\sigma\geq s$. A finite dimensional subspace $\Cal I \subset (H^\infty)^\ast$ will be called \emph{$\sigma$-regular  [regular]} if it admits a $\sigma$-regular [regular]  basis. 
\end{definition}
\medskip
\noindent It was proved in \S 1.4, in particular in Corollary \ref{cor:deltareg}, that for any $s\geq 0$ the closure $\overline{ H^s_q(M)} \subset \bar H^s_q(M)$ of the weighted Sobolev space as a subspace of the {\it Friedrichs }weighted Sobolev space, is equal to the perpendicular of a subspace $\Cal D^s_q \subset \bar H^{-s}_q(M)$, introduced in \pref{eq:Dsq}, of distributions supported on the singular set $\Sigma_q \subset M$. It follows by Theorem \ref{thm:SO} and by Corollary \ref{cor:deltareg} that all the spaces $\Cal D^s_q$ are regular with respect to the family $\{\bar H^s_q(M)\vert s\geq 0\}$ of Friedrichs weighted Sobolev spaces. This regularity result is the key to the existence of smoothing operators for the family $\{H^s_q(M)\vert s\geq 0 \}$ of weighted Sobolev spaces. In fact, their construction in Theorem \ref{thm:smoothingop} is based on the regularity of the subspaces $\Cal D^s_q \subset \bar H^{-s}_q(M)$ and on the existence of smoothings for the family of Friedrichs weighted Sobolev spaces, which can be defined by appropriate truncations of the Fourier expansion.

\smallskip 
\noindent Our goal in this section is to implement a similar strategy  in order to construct smoothing families with values in the perpendicular of the subspaces of invariant distributions. We prove below, under the hypothesis that the Kontsevich-Zorich cocycle is non-uniformly hyperbolic, a preliminary result on the regularity of the spaces of horizontal [vertical] invariant distributions  with respect to the family of weighted Sobolev spaces. The basic criterion for the regularity of spaces in invariant distributions
is based on Theorem \ref{thm:Oseledecbasis} and the following notion:

\begin{definition}
\label{def:coherent}
Let $\mu$ be any $G_t$-invariant ergodic probability measure on any stratum ${\Cal M}^{(1)}_{\kappa}\subset {\Cal M}^{(1)}_g$ of orientable holomorphic quadratic differentials.
For any $q\in \Cal R_\mu$, a simple invariant distribution $\Cal D^\pm \in \Cal I_{\pm q}(M\setminus\Sigma_q)$ will be called \emph{coherent }(with respect to the family $\{H^s_q(M)\vert s\geq 0\}$) if its weighted Sobolev order and  its Lyapunov exponent are related:
$$
{\Cal O}_q^H(\Cal D^\pm) = \vert l^\pm_\mu (\Cal D^\pm)\vert \,.
$$
A finite dimensional space of  invariant distributions will be called coherent (with respect to the family $\{H^s_q(M)\vert s\geq 0\}$) if it has a basis of simple coherent elements.
\end{definition}

\noindent By the definitions of coherence and regularity, the following result follows immediately from Theorem \ref{thm:Oseledecbasis}:
\begin{lemma}
\label{lemma:cohreg} 
Let $\mu$ be any $G_t$-invariant ergodic probability measure on any stratum ${\Cal M}^{(1)}_{\kappa}\subset {\Cal M}^{(1)}_g$ of orientable holomorphic quadratic differentials.
For any $q\in \Cal R_\mu$, any coherent finite dimensional space of  invariant distributions 
is regular (with respect to the family $\{H^s_q(M)\vert s\geq 0\}$).
\end{lemma}

\smallskip
\noindent For any  $q\in \Cal M^{(1)}_\kappa$ and for any $k\in \N\setminus\{0\}$,  let $\Cal J^k_q(M) \subset \Cal I^k_q(M)$ be the space of horizontally invariant distributions defined as
follows:
\begin{equation}
\label{eq:regularsubspace}
\Cal J^k_q(M) := \oplus_{j=0}^{k-1} \,  {\Cal L}^j_{T_q} \left[\Cal I^1_q(M)\right]\,.
\end{equation}

\begin{corollary}
\label{cor:Jkregular} For any $SO(2,\R)$-absolutely continuous, KZ-hyperbolic measure $\mu$ on any stratum ${\Cal M}^{(1)}_{\kappa} \subset {\Cal M}^{(1)}_g$ of orientable quadratic differentials, for $\mu$-almost all orientable quadratic differential $q\in \Cal M^{(1)}_\kappa$ and for any $k\in \N$,  the space $\Cal J^k_q(M) \subset \Cal I^k_q(M)$  is coherent, hence regular, with respect to the family $\{H^s_q(M)\vert s\geq 0\}$ of weighted Sobolev spaces. 
\end{corollary}
\begin{proof} It follows immediately from Theorem \ref{thm:bcreg} for the case $k=0$ and from Theorem \ref{thm:bcreg} and Lemma \ref{lemma:idderid} in  the general case.
\end{proof}

\noindent It was proved in Theorem \ref{thm:smoothingop} that the family $\{H^s_q(M)\vert s\geq 0\}$
of weighted Sobolev spaces admits a family of smoothing operators. The existence of a smoothing family together with the regularity of the distributional obstructions are the key elements of the interpolation theory for solutions of the cohomological equation. We formalize below the basic construction of smoothing projectors onto the perpendicular of a regular subspace.

\begin{definition} 
\label{def:smoothproj}
Let $\{H^s \vert s\geq 0\}$ be a $1$-parameter family of normed spaces such that the embeddings \pref{eq:normedspaces} hold. A \emph{smoothing projection} of degree $\sigma\in \overline\R^+$ relative to a subspace $\Cal I^\sigma \subset (H^\sigma)^\ast$ is a family $\{P^\sigma(\tau) \vert \tau\in (0,1]\}$ of linear operators such that the operator $P^\sigma(\tau):H^0 \to [\Cal I^\sigma]^\perp\subset H^\sigma$ is bounded for all $\tau\in (0,1]$ and the following estimates hold. For any $r$, $s\in [0, \sigma]$ and for any $\epsilon>0$, there exists a constant $C^\sigma_{r,s}(\epsilon)>0$ such that, for all $u\in [\Cal I^\sigma \cap (H^s)^\ast ]^\perp \subset H^s$ and for all $\tau \in (0,1]$,
\begin{equation}
\label{eq:smoothproj}
\begin{aligned}
\vert P^\sigma(\tau)(u) -u \vert_r \,&\leq\,  C^\sigma_{r,s}(\epsilon) \,\vert u \vert_s \, \tau^{s-r-\epsilon} \,, 
\quad \text{ \rm if } \, s>r\,; \\
\vert P^\sigma(\tau)(u) \vert_r \,&\leq\,  C^\sigma_{r,s}(\epsilon) \,\vert u \vert_s \, \tau^{s-r-\epsilon} \,, 
\quad \text{ \rm if } \, s\leq r\,. 
\end{aligned}
\end{equation}
A smoothing projection relative to the trivial subspace $\Cal I^\sigma=\{0\}\subset H^\sigma$ will be called a \emph{smoothing} of degree $\sigma\in \overline\R^+$ for the family $\{H^s \vert s\geq 0\}$.
\end{definition} 
\noindent The following result is a straightforward generalization of Theorem \ref{thm:smoothingop},
which implies the existence of smoothings of any finite degree for the family $\{ H^s_q(M)\vert s\geq 0\}$ of weighted Sobolev spaces.

\begin{theorem} 
\label{thm:smoothrel}
Let $\{H^s \vert s\geq 0\}$ be a $1$-parameter family of normed spaces such that
the embeddings \pref{eq:normedspaces} hold. If the family $\{H^s \vert s\geq 0\}$  has a smoothing of degree $\sigma\in \R^+$ and the subspace $\Cal I^\sigma\subset  (H^{\sigma})^\ast$ is $\sigma$-regular, there exists a smoothing projection of degree $\sigma\in\R^+$ relative to $\,\Cal I^\sigma\,$.
\end{theorem}
\begin{proof}  Let $\{{\Cal S}^\sigma(\tau) \vert \tau\in (0,1]\}$ be a smoothing of degree $\sigma>0$
for the family $\{H^s \vert s\geq 0\}$ and let $\{\Cal D_1, \dots \Cal D_{J}\}$ be a $\sigma$-regular basis for the subspace $\Cal I^\sigma\subset (H^\sigma)^\ast$. By the Definition \ref{def:regularsystem} of regularity, for any $\tau \in (0,1]$, there exists a dual basis $\{u_1(\tau), \dots, u_{J}(\tau)\}\subset H^\sigma$ such that the estimates \pref{eq:regularsystem} hold for all $0\leq r \leq\sigma$ and all $\epsilon >0$. For any $u \in H^0$, we define 
\begin{equation}
\label{eq:Psigmadef}
P^\sigma(\tau)(u) :=   {\Cal S}^\sigma(\tau)(u) - \sum_{j=1}^{J} 
\Cal D_j \left( {\Cal S}^\sigma(\tau)(u) \right) u_j(\tau)  \,.
\end{equation}
We claim that the family $\{P^\sigma(\tau) \vert \tau\in (0,1]\}$ is a smoothing projection (of degree $\sigma>0$) relative to the subspace $\Cal I^\sigma\subset (H^\sigma)^\ast$. It follows immediately
from the definition  that $P^\sigma(\tau)(u) \in [\Cal I^\sigma]^\perp
\subset H^\sigma$ for any $u\in H^0$, hence the operators $P^\sigma(\tau):H^0 \to [\Cal I^\sigma]^\perp \subset H^\sigma$ are well-defined, linear  and bounded. We claim that, for any $s\in [0,\sigma]$, for any $j\in \{1, \dots, J\}$ and for any $\epsilon>0$, there exists a constant $C^s_j(\epsilon)>0$ such that, for all vectors $u\in [\Cal I^\sigma \cap (H^s)^\ast ]^\perp \subset H^s$, the following estimate holds:
\begin{equation}
\label{eq:Diestfinal}
\vert  \Cal D_j \left( {\Cal S}^\sigma(\tau)(u) \right)\vert  \, \leq \, C^s_j(\epsilon) \, \vert u \vert_s \, 
\tau^{s-{\Cal O}^H(\Cal D_j)-\epsilon} \,, \quad \text{ \rm for all }\tau \in (0,1] \,.
\end{equation}
In fact, if the order ${\Cal O}^H(\Cal D_j) \geq s$, since $\Cal D_j \in (H^{r_j})^\ast$ for any $r_j>{\Cal O}^H(\Cal D_j)$ and the family $\{\Cal S^\sigma(\tau) \vert \tau\in (0,1]\}$ is a smoothing for $\{H^s\vert s\geq 0\}$, for any $\epsilon>0$ there exists a constant $A^s_j(\epsilon)>0$ such that for all $u \in H^s$ and all $\tau\in (0,1]$,
\begin{equation}
\label{eq:Diestone}
\vert  \Cal D_j \left( {\Cal S}^\sigma(\tau)(u) \right)\vert  \, \leq \,\vert \Cal D_j \vert_{r_j}^\ast \,
\vert {\Cal S}^\sigma(\tau)(u) \vert _{r_j} \leq A^s_j(\epsilon) \, \vert u\vert_s \, \tau^{s-r_j-\epsilon/2} \,.
\end{equation}
If ${\Cal O}^H(\Cal D_j) <s$, since $\Cal D_j \in (H^{r_j})^\ast \subset (H^s)^\ast$, for any ${\Cal O}^H(\Cal D_j) <r_j<s$
 and the family $\{\Cal S^\sigma(\tau) \vert \tau\in (0,1]\}$ is a smoothing for $\{H^s\vert s\geq 0\}$,  there exists $B^s_j(\epsilon)>0$ such that, for all $u\in [\Cal I^\sigma \cap (H^s)^\ast ]^\perp$ and all $\tau\in (0,1]$,
 \begin{equation}
\label{eq:Diesttwo}
\begin{aligned}
\vert  \Cal D_j \left( {\Cal S}^\sigma(\tau)(u) \right)\vert  &= \vert  \Cal D_j 
\left( {\Cal S}^\sigma(\tau)(u) -u \right)\vert \\  
&\leq \vert \Cal D_j \vert_{r_j}^\ast \,
\vert {\Cal S}^\sigma(\tau)(u) -u \vert _{r_j} \leq B^s_j(\epsilon) \, \vert u\vert_s \, \tau^{s-r_j-\epsilon/2} \,.
\end{aligned}
\end{equation}
The estimate \pref{eq:Diestfinal} follows immediately from estimates \pref{eq:Diestone} and
\pref{eq:Diesttwo} by taking $r_j \in ({\Cal O}^H(\Cal D_j), {\Cal O}^H(\Cal D_j) +\epsilon/2)$ for all $j\in\{1, \dots, J\}$ 
in both cases.  

\smallskip
\noindent By estimate \pref{eq:Diestfinal} (just proved) and by the estimates on $\{u_1(\tau), \dots, u_J(\tau)\}$ in formula \pref{eq:regularsystem}, for any $r \in [0, \sigma]$ and for any $\epsilon>0$, there exists a constant $C^{(1)}_{r,s}(\epsilon)>0$ such that, for all $u\in [\Cal I^\sigma \cap (H^s)^\ast ]^\perp$ and all $\tau\in (0,1]$,
\begin{equation}
\label{eq:sumDiest}
\sum_{j=1}^J \vert \Cal D_j \left( {\Cal S}^\sigma(\tau)(u) \right)\vert \, \vert u_j(\tau) \vert_r  \,\,
\leq  \, \,C^{(1)}_{r,s}(\epsilon)\,  \vert u\vert_s\,\tau^{s-r-\epsilon}\,.
\end{equation}
Finally, the required estimates \pref{eq:smoothproj} for the family $\{P^\sigma(\tau) \vert \tau\in(0,1]\}$, defined in \pref{eq:Psigmadef}, follow from \pref{eq:sumDiest}, since $\{{\Cal S}^\sigma(\tau)\vert \tau\in(0,1]\}$ is a smoothing (of degree $\sigma>0$) for the family $\{H^s\vert s\geq 0\}$.
 \end{proof}

\noindent The smoothness of the solutions of the cohomological equation is best expressed with respect to the following uniform norms. 
\begin{definition}
\label{def:uniformnorms}
For any $k\in \N$, let $B^k_q(M)$ be the space of all functions $u\in H^{k}_q(M)$ such that
$S_q^i T_q^j u=T_q^i S_q^j u \in L^{\infty}(M)$ for all pairs of integers $(i,j)$ such that
$0\leq i+j \leq k$. The space $B^k_q(M)$ is endowed with the norm defined as follows: for any 
$u \in B^k_q(M)$, 
\begin{equation}
\label{eq:uniformnorms1}
\vert u \vert_{k,\infty} \,:= \,  \left[\sum_{i+j \leq k}  \vert S_q^i T_q^j u\vert^2_{\infty}\right]^{1/2}
\,=\,\left[\sum_{i+j \leq k}  \vert T_q^iS_q^j u\vert^2_{\infty} \right]^{1/2} \,.
\end{equation}
For $s\in [k,k+1)$, let $B^s_q(M):= B^k_q(M) \cap H^s_q(M)$ endowed with the norm defined
as follows: for any  $u \in B^s_q(M)$, 
\begin{equation}
\label{eq:uniformnorms2}
\vert u \vert_{s,\infty} :=  \left(\,\vert u \vert^2_{k,\infty}\,+\,  \vert u \vert^2_s \,\right)^{1/2}\,.
\end{equation}
\end{definition}

\begin{theorem}
\label{thm:CEsharpfundth}
Let $\mu$ be a $SO(2,\R)$-absolutely continuous, KZ-hyperbolic measure on a stratum $\Cal M^{(1)}_\kappa\subset \Cal M^{(1)}_g$  of orientable quadratic differentials. For $\mu$-almost all $q\in \Cal M^{(1)}_\kappa$ and for any $k\in\N$, the space $\Cal I^k_q(M)$ is coherent, hence regular (with respect to the family $\{H^s_q(M)\vert s\geq 0\}$) and, for any $s>k$, there exists a measurable function $C^{k,s}_\kappa: \Cal M^{(1)}_\kappa \to \R^+$ such that the following holds. For any function $f\in  [\Cal I_q^k(M)]^\perp \cap H^s_q(M)$ the Green solution $\,\Cal U_q(f)\,$ of the cohomological equation $S_qu=f$ belongs to the space $B^{k-1}_q(M)$ and satisfies the estimates:
\begin{equation}
\label{eq:CEsharpfundth}
\vert \Cal U_q(f) \vert _{k-1}\,\leq\,  \vert \Cal U_q(f) \vert_{k-1,\infty} \,  \leq \,  
C^{k,s}_\kappa(q)\, \vert f\vert_s 
\end{equation}
\end{theorem}

\begin{proof} For $k=1$ the statement reduces to Theorem \ref{thm:CEsharp1} and Corollary 
\ref{cor:Jkregular}. Since  $f\in [\Cal J^2_q(M)]^\perp$ implies by definition that $f$ and
$T_q f\in [\Cal I^1_q(M)]^\perp$, by Theorem \ref{thm:CEsharp1},  for $\mu$-almost all $q\in \Cal M^{(1)}_\kappa$ and for any $s>k=2$, the  cohomological equations $\,S_qu=f\,$ 
and $\,S_qu_T=T_q f\,$  have Green (zero average) solutions $u(f)$ and $u_T(f)\in L^{\infty}(M) 
\subset L^2_q(M)$ which satisfy the bounds,
\begin{equation}
\label{eq:CEsharpest3}
\begin{aligned}
\vert u(f) \vert_\infty \,&\leq \,  C^{s-1}_\kappa (q) \,  \vert f \vert_{s-1}  \,,\\
\vert u_T(f) \vert_\infty \,&\leq \,  C^{s-1}_\kappa (q) \,  \vert  f \vert_{s}  \,.
\end{aligned}
\end{equation}
It follows that the distribution $u_T(f) - T_q u(f) \in H^{-1}_q(M)$ is horizontally invariant. Let  $\{\Cal D_1, \dots, \Cal D_g\} \subset \Cal I^1_q(M)$ be any regular basis such that 
\begin{equation}
\label{eq:I1soborders}
\Cal O^H_q(\Cal D_j)= 1-\lambda_j^\mu  \,, \quad \text{ \rm for all }\, j\in \{1, \dots, g\}\,.
\end{equation}
It is no restrictive to assume that $\Cal D_1$ is the average. Since Green solutions have zero average
by definition, there exists $F_2(f) , \dots F_g(f) \in \C$ such that
\begin{equation}
\label{eq:CEsharpFjs}
u_T(f) - T_q u(f) =   \sum_{j=2}^g  F_j(f) \, \Cal D_j \,.
\end{equation}
It follows from the bounds \pref{eq:CEsharpest3} that the maps $F_j: [\Cal J^2_q(M)]^\perp \to \C$ are linear bounded functionals, for all $j\in \{2, \dots, g\}$, defined on the closed subspace $[\Cal J^2_q(M)]^\perp \subset H^s_q(M)$. In fact, let $\{u_1(\tau) ,\dots, u_g(\tau)\} \subset H^1_q(M)$ be any dual basis of the regular basis $\{\Cal D_1, \dots, \Cal D_g\}$ as in Definition \ref{def:regularsystem}. By Theorem \ref{thm:CEsharp1}, for any $s>2$ and for any $(\tau_2, \dots, \tau_g) \in (0,1]^{g-1}$, 
\begin{equation}
\label{eq:Fjest1}
\begin{aligned}
\vert F_j(f) \vert   &= \vert \<u_T(f), u_j(\tau_j) \>_q \vert  \,+\, \vert \<u(f), T_q u_j(\tau_j)\>_q \vert \\
&\leq  \, C^{s-1}_\kappa(q) \, \left\{ \vert u_j(\tau_j) \vert_0 \, \vert f \vert_{s}  \, +\,
 \vert u_j(\tau_j) \vert_1 \, \vert f \vert_{s-1} \right\}\,.
\end{aligned}
\end{equation}
We claim that, for each $j\in \{2, \dots, g\}$, the linear functional $F_j$ extends to a horizontally invariant distribution $\Phi_j  \in H^{-s}_q(M)$ such that 
\begin{equation}
\label{eq:I2soborders}
\Cal O^H_q(\Phi_j)\,\leq \, 1+\lambda_j^\mu  \,, \quad \text{ \rm for all }\, j\in \{2, \dots, g\}\,.
\end{equation}
In fact, since by Corollary \ref{cor:Jkregular} the space $\Cal J^2_q(M)$ is regular, by Theorem \ref{thm:smoothrel} there exists a smoothing projection $\{ P_J^\sigma(\tau) \vert \tau\in (0,1]\}$ of degree $\sigma>2$ relative to the subspace $\Cal J^2_q(M) \subset H^{-2}_q(M)$. Hence by definition, for any $r$, $s\in [0,\sigma]$ and any $\epsilon>0$,  there exists a constant $C^\sigma_{r,s}(\epsilon)>0$ such that, for all $f\in [\Cal J^2_q(M)]^\perp \cap H^s_q(M)$ and for all $\tau \in (0,1]$, 
\begin{equation}
\label{eq:Jsmoothproj}
\begin{aligned}
\vert P^\sigma_J(\tau)(f) -f\vert_r \,&\leq\,  C^\sigma_{r,s}(\epsilon) \,\vert u \vert_s \, \tau^{s-r-\epsilon} \,, 
\quad \text{ \rm if } \, s>r\,; \\
\vert P^\sigma_J(\tau)(f) \vert_r \,&\leq\,  C^\sigma_{r,s}(\epsilon) \,\vert f \vert_s \, \tau^{s-r-\epsilon} \,, 
\quad \text{ \rm if } \, s\leq r\,. 
\end{aligned}
\end{equation}
Let $(f_n)_{n\in\N}$ be the sequence of functions defined as follows:
\begin{equation}
\label{eq:Jfn}
f_n := P^\sigma_J(2^{-n})(f) \in [\Cal J^2_q(M)]^\perp \cap  H^\sigma_q(M)\,,  
\quad \text{ \rm for all } \, n\in \N\,.
\end{equation}
It follows from estimates \pref{eq:Jsmoothproj} that, if  $s-1<1+\lambda^\mu_j <s_j <s$, for any
$\epsilon>0$, there exists a constant $C^\sigma_{s_j,s}(\epsilon)>0$ such that
\begin{equation}
\begin{aligned}
\label{eq:Jfnest}
\vert f_{n+1} - f_n \vert_{s-1}  \, &\leq \, C^\sigma_{s_j,s}(\epsilon) \, \vert f \vert_{s_j} \, 
2^{-n(s_j-s+1-\epsilon)}\,;\\
\vert f_{n+1} - f_n \vert_{s}  \, &\leq \, C^\sigma_{s_j,s}(\epsilon)\, \vert f \vert_{s_j} \, 
2^{n(s -s_j+\epsilon)}\,.
\end{aligned}
\end{equation}
Let $\{u_1^{(n)}, \dots, u_g^{(n)}\}\subset H^1_q(M)$ be the sequence of dual basis of the regular basis $\{\Cal D_1, \dots, \Cal D_g\} \subset \Cal I^1_q(M)$ defined as follows: for each $j\in \{1, \dots, g\}$, 
\begin{equation}
u_j^{(n)}  :=  u_j (2^{-n})  \, ,\quad \text{ for all } \, n\in \N\,.
\end{equation}
By the estimate \pref{eq:regularsystem}  in Definition \pref{def:regularsystem} and by the identities
\pref{eq:I1soborders}, for any $\epsilon >0$
there exists a constant $C(\epsilon)>0$ such that, for all $j\in \{1, \dots, g\}$,
\begin{equation}
\label{eq:dualunest}
\begin{aligned}
\vert u_j^{(n)} \vert_0  \, &\leq  \,  C(\epsilon) \,  2^{-n(1-\lambda_j^\mu-\epsilon)} \,; \\
\vert u_j^{(n)} \vert_1  \, &\leq  \, C(\epsilon) \,  2^{n(\lambda_j^\mu + \epsilon)} \,.
\end{aligned}
\end{equation}
For any  $s_j >1+\lambda_j$, there exists $s\in (2, 2+\lambda_j^\mu)$ and $\epsilon >0$ such that
\begin{equation}
\label{eq:poscond}
s_j +1 -s-\lambda_j^\mu  -2\epsilon := \alpha_j >0 \,.
\end{equation}
Hence by the estimate \pref{eq:regularsystem} and  \pref{eq:Fjest1}, for $\mu$-almost $q\in \Cal M^{(1)}_\kappa$, there exists a constant $C_q^j(\epsilon):=C_q(\epsilon,s,s_j, \alpha_j)>0$ such that 
\begin{equation}
\label{eq:Fjest2}
\vert F_j (f_{n+1}-f_n) \vert \, \leq \,C^j_q(\epsilon) 2^{-\alpha_j n}  \, \vert f \vert _{s_j}\,.
\end{equation}
By \pref{eq:Jsmoothproj} and \pref{eq:Jfn} the sequence $f_n \to f$ in $H^s_q(M)$ for any
$s<\sigma$. Since \pref{eq:Fjest1} holds for any $s>2$, it follows from \pref{eq:Fjest2}
that, for $\mu$-almost $q\in \Cal M^{(1)}_\kappa$,
\begin{equation}
\label{eq:Fjest3}
\vert F_j( f ) \vert \, \leq \,\frac {C^j_q(\epsilon)}{1- 2^{-\alpha_j}}  \, \vert f \vert _{s_j}  \,\, +\,\, 
\vert P^\sigma_J(1)(f) \vert_\sigma\,.
\end{equation}
Since the operator $P_J^\sigma(1): L^2_q(M) \to H^\sigma_q(M)$ is bounded, we conclude that  
for any $s_j> 1+ \lambda^\mu_j$ there exists a continuous extension $\Phi^{s_j}_j \in H^{-s_j}_q(M)$
of the linear functional $F_j: [\Cal J^2_q(M)]^\perp \cap H^s_q(M) \to \C$.  By construction, for any $r_j$, $s_j>1+ \lambda^\mu_j$, the extensions $\Phi^{r_j}_j = \Phi^{s_j}_j$ (mod.  $\Cal J ^2_q(M)$). Since $\Cal J ^2_q(M)$ is finite dimensional, for each $j\in\{2, \dots,g\}$ there exists a distribution $\Phi_j \in H^{-2}_q(M)$, which extends the linear functional $F_j$, such that  the Sobolev order $\Cal O^H_q(\Phi_j) \leq 1+\lambda^\mu_j$ as claimed in \pref{eq:I2soborders}. Finally, we prove that, by construction, the distributions $\Phi_2, \dots, \Phi_g \in H^2_q(M)$ are horizontally invariant. In fact, for any $s>2$ and any $v \in \Cal H^{s+1}_q(M)$, the function $f_v :=S_q v \in [\Cal I^2_q(M)]^\perp \cap H^s_q(M)$, hence, by Theorem \ref{thm:CEsharp1}, the cohomological equations $\,S_q u =f_v\,$ and $S_q u_T=T_q f_v$ have unique Green solutions $u(f_v)$ and $u_T(f_v) \in L^{\infty}(M)$ respectively. By the ergodicity of the horizontal foliation, since $v$ and $T_q v \in L^2_q(M)$ are also zero-average solutions, the identities $u(f_v)=v$ and $u_T(f_v)=T_q v$ hold. It follows that $u_T(f_v) - T_q u(f_v)=0$, hence by \pref{eq:CEsharpFjs}, for all $j\in \{2,\dots, g\}$,
\begin{equation}
\Phi_j(S_q v)\, = \,F_j (S_q v) \,=\, 0\,, \quad \text {for all }\, v\in H^{s+1}_q(M)\,,
\end{equation}
thus $\{\Phi_2, \dots\Phi_g\} \subset \Cal I^s_q(M)$ and the claim is completely proved.

\smallskip
\noindent For $\mu$-almost all $q\in \Cal M^{(1)}_\kappa$ and for all integers $k\geq 2$,  let
$\hat {\Cal I}^k_q(M) \subset \Cal I^k_q(M)$ be the subspaces defined as follows:
\begin{equation}
\label{eq:hatIk}
\hat {\Cal I}^k_q(M) :=    \bigcup_{h=0}^{k-2} \Cal L_{T_q} ^h\left[ \Cal J^2_q(M)\, \oplus\, \bigoplus_{j=2}^g \C\cdot \Phi_j\right] \,.
\end{equation}
It follows from the above construction that the following holds: for any integer $k\geq 2$, for
any $s>k$ and for any function $f \in H^s_q(M) \cap [\hat {\Cal I}^k_q(M)]^\perp$, the Green
solution $\Cal U_q(f) \in B^{k-1}_q(M)$ and satisfies the required estimate \pref{eq:CEsharpfundth}.
In fact, for $k=2$ the statement follows by the construction of the system of invariant distributions 
$\{\Phi_2, \dots, \Phi_g\}\subset H^{-2}_q(M)$. For $k\geq 3$, the statement can be proved by induction. In fact, by the induction hypothesis, for any $s>k$ and any $f\in H^s_q(M)\cap  [\hat {\Cal I}^k_q(M)]^\perp$, the cohomological equations $\,S_q u =f\,$ has a unique solution 
$u\in B^{k-2}_q(M)$ such that 
\begin{equation}
\label{eq:CEsharpinductest1}
 \vert u \vert_{k-2,\infty} \,  \leq \,  C^{k-1,s-1}_\kappa(q)\, \vert f\vert_{s-1} \,.
\end{equation}
In addition, the function $u_T:=T_q u \in H^{k-3}_q(M)$ is the unique solution of the cohomological equation $\,S_qu_T= T_qf\,$. Since  $T_q f \in H^{s-1}_q(M) \cap [\hat {\Cal I}^{k-1}_q(M)]^\perp$, by the induction hypothesis, the following estimate holds:
\begin{equation}
\label{eq:CEsharpinductest2}
\vert T_q u \vert_{k-2,\infty}  \leq C^{k-1,s-1}_\kappa(q)\, \vert  T_q f\vert_{s-1}  \leq  
C^{k-1,s-1}_\kappa(q)\, \vert f\vert_{s}\,.
\end{equation}
Finally, by the Sobolev embedding theorem, there exists a continuous function  $\tilde C_\kappa:\Cal M^{(1)}_\kappa \to \R^+$ such that the following estimate holds:
\begin{equation}
\label{eq:CEsharpinductest3}
\left[ \sum_{i=1}^{k-1} \vert S_q^i u \vert_{\infty}^2 \right]^{1/2}  =  \left[\sum_{i=0}^{k-2}
 \vert S_q^i f \vert_{\infty}^2  \right]^{1/2} \leq \tilde C_\kappa(q) \, \vert f \vert_s \,.
\end{equation}
The required estimate  follows from \pref{eq:CEsharpinductest1}, \pref{eq:CEsharpinductest2} and \pref{eq:CEsharpinductest3}.

\smallskip
\noindent It remains to be proven that, for $\mu$-almost $q\in \Cal M^{(1)}_\kappa$ and all integers
$k\geq 2$, the space $\Cal I^k_q(M)$ is coherent.  From the above argument it follows that
\begin{equation}
\label{eq:Iidentity1}
\Cal I_q(M) \,=\,  \bigoplus_{k=2}^{+\infty} \hat {\Cal I}^{k}_q(M) \,.
\end{equation}
In fact, if $f \in H^{\infty}_q(M)$ is such that $f \in [\hat {\Cal I}^{k}_q(M)]^\perp$ for all $k\in \N$, there exists $u\in H^{\infty}_q(M)$ such that $S_qu=f$, hence $f\in \Cal I_q(M)^\perp$. The identity
\pref{eq:Iidentity1} immediately implies that 
\begin{equation}
\label{eq:Iidentity2}
\Cal I^k_q(M) \,=\,  \hat {\Cal I}^{k}_q(M) \,, \quad \text{ \rm for all integers } k\geq 2  \,.
\end{equation}
The identities \pref{eq:Iidentity2} in turn imply that the system $\{\Phi_2, \dots, \Phi_g\}$ 
constructed above is linearly independent over $\Cal J^2_q(M)$, in particular 
\begin{equation}
\label{eq:I2maxdim}
\text{ \rm dim}_\C \, \Cal I^2_q(M) = \text{ \rm dim}_\C \, \hat{\Cal I}^2_q(M) 
= 3g -2\,.
\end{equation}
In fact, by the above construction the following identity holds:
\begin{equation}
\label{eq:identity3}
\Cal I_q(M) =  \Cal I^1_q(M) \, + \,\bigoplus_{j=2}^g  \C\cdot \Phi_j   \, +\,
 \Cal L_{T_q} [\Cal I_q(M)] \,.
\end{equation}
 Let $\Cal D_q : \Cal B_q(M) \to \Cal I_q(M)$ the isomorphism between the space $\Cal B_q(M)$ of horizontally basic currents and the space of horizontally invariant distributions defined in \pref{eq:DtoC}.
 Let $\delta_q: \Cal B_q(M) \to \Cal B_q(M)$ be the differential map on the space of basic currents introduced in \pref{eq:deltas}.  From \pref{eq:identity3} it follows immediately that
 \begin{equation}
\label{eq:identity4}
\Cal B_q(M) =  \Cal B^1_q(M)\, + \, \bigoplus_{j=2}^g  \C\cdot {\Cal D}_q^{-1}( \Phi_j) \,  
+\, \delta_q [\Cal B_q(M)]\,.
\end{equation}
By the structure theorem for (real) basic currents (Theorem \ref{thm:bcstruct}),  the cohomology map $j_q:\Cal B_q(M) \to H^1(M,\R)$ vanishes on the space $\delta_q [\Cal B_q(M)]$. On the other hand, by Corollary \ref{cor:basiccohom}, for $\mu$-almost all $q\in \Cal M^{(1)}_\kappa$, the cohomology map on $\Cal B_q(M)$ has rank of codimension $1$ in the homology $H^1(M,\R)$, hence of dimension $2g-1$. Since $ \Cal B^1_q(M)$ has dimension $g$ and the map $\Cal D_q : \Cal B_q(M) \to \Cal I_q(M)$ is an isomorphism, it follows that the system $\{\Phi_2, \dots, \Phi_g\}$ is linearly independent over $\Cal J^2_q(M)$ and \pref{eq:I2maxdim} holds.

\smallskip
\noindent We claim that, for $\mu$-almost all $q\in \Cal M^{(1)}_\kappa$, the space $\hat {\Cal I}^2_q(M)= \Cal I^{2}_q(M)$ is coherent. For any Lyapunov exponent $l<0$ of the cocycle $\Phi^2_t\vert {\Cal I}^2_q(M)$, let $E^2_q(l) \subset \Cal I^{2}_q(M)$ denote the corresponding Oseledec subspace and let $F^2_q(l) \subset E^2_q(l)$ be the subspace of coherent distributions.  Let $$
l_1 <\dots < l_d
$$ 
be the distinct Lyapunov exponents of the cocycle $\Phi^2_t\vert {\Cal I}^2_q(M)$ on the the Oseledec complement of the subspace $\Cal J^{2}_q(M)$. Since $\Cal J^2_q(M)$ is coherent and 
$$
\Cal I^{2}_q(M)=\hat {\Cal I}^2_q(M)= \Cal J^{2}_q(M) \oplus \bigoplus_{i=1}^d E_q^2(l_i)\,,
$$
it is sufficient to prove the identities:
\begin{equation}
\label{eq:coherence}
F^2_q(l_i) \,=\, E^2_q(l_i) \,, \quad \text{ \rm for all }\, i\in \{1, \dots, d\} \,.
\end{equation}
By Theorem \ref{thm:Otype} and Lemma \ref{lemma:idderid}, the Lyapunov spectrum of $\{\Phi^2_t\vert {\Cal J}^2_q(M)\}$ is the ordered set 
$$
0 > \lambda_2^\mu-1 \geq \dots \geq \lambda^\mu_g -1 \geq \lambda_2^\mu-2 \geq \dots \geq
\lambda^\mu_g -2 \,.
$$ 
Hence, by the description of the Lyapunov spectrum of $\{\Phi^2_t\vert {\Cal I}^2_q(M)\}$ given in Corollary \ref{cor:Lspectrum},  the set
$$
\{l_1, \dots, l_d\} =  \{-1-\lambda_2^\mu, \dots, -1-\lambda_g^\mu \}\,.
$$
For any $s\geq 0$, let $\Cal I_q(s)\subset E^2_q(l_1)\oplus \dots \oplus E^2_q(l_d)$ be the subset of horizontally invariant distributions of Sobolev order less or equal to $s\geq 0$. It follows from Lemma \ref{lemma:regupperbound} that the following inclusions hold: 
\begin{equation}
\label{eq:IEinclusion}
\Cal I_q(\vert l_i\vert )  \subset \bigoplus_{j\leq i} E^2_q(l_j)\,, \quad \text{ \rm for all }\, i\in \{1, \dots, d\} \,.
\end{equation}
By the estimate \pref{eq:I2soborders} on Sobolev orders of the distribuitons in the system $\{\Phi_2, \dots,\Phi_g\}$, the following lower bounds hold:
\begin{equation}
\label{eq:IEdimension}
\text{ \rm dim}_\C I_q(\vert l_i\vert ) \,\geq \,  \sum_{j=1}^i \text{ \rm dim}_\C E^2_q(l_j)\,, 
\quad \text{ \rm for all }\, i\in \{1, \dots, d\} \,.
\end{equation}
It follows from  \pref{eq:IEinclusion} and \pref{eq:IEdimension} that the inclusions in
\pref{eq:IEinclusion} are in fact identities, for all $i\in \{1, \dots, d\}$, and by Lemma \ref{lemma:regupperbound} the claim \pref{eq:coherence} holds. 

\smallskip
\noindent We have thus proved that the space ${ \Cal I}^2_q(M)=\hat{\Cal I}^2_q(M)$ is coherent.
It follows  from definition \pref{eq:hatIk} and Lemma \ref{lemma:idderid} that the space $\hat {\Cal I}^k_q(M)$ is coherent for any integer $k\geq 3$. Since by \pref{eq:Iidentity2} the identity ${ \Cal I}^k_q(M)=\hat { \Cal I}^k_q(M)$ holds for all $k\geq 2$, the space ${ \Cal I}^k_q(M)$ is coherent and the proof is complete.
\end{proof}

\begin{theorem}
\label{thm:CEsharpth}
Let $\mu$ be a $SO(2,\R)$-absolutely continuous, KZ-hyperbolic measure on a stratum $\Cal M^{(1)}_\kappa\subset \Cal M^{(1)}_g$  of orientable quadratic differentials. For $\mu$-almost all $q\in \Cal M^{(1)}_\kappa$ and for any $s\in \R^+$, the space $\Cal I^s_q(M)$ is coherent, hence regular (with respect to the family $\{H^s_q(M)\vert s\geq 0\}$) and, for any $0<r <s-1$ there exists a measurable function $C^{r,s}_\kappa: \Cal M^{(1)}_\kappa \to \R^+$ such that the following holds. For any function $f\in  [\Cal I_q^s(M)]^\perp \cap H^s_q(M)$ the Green solution $\,\Cal U_q(f)\,$ of the cohomological equation $S_qu=f$ belongs to the space $H^r_q(M)$ and satisfies the estimates:
\begin{equation}
\label{eq:CEsharpth}
\vert \Cal U_q(f) \vert _r \,\leq\,  C^{r,s}_\kappa(q)\, \vert f\vert_s 
\end{equation}
\end{theorem}
\begin{proof}
By Theorem \ref{thm:CEsharpfundth}, the subspace $\Cal I^k_q(M)$ is coherent  for $\mu$-almost 
all $q\in \Cal M^{(1)}_\kappa$ and for any $k\in \N$.  Since for any $s<k$ the sub-bundle $\Cal I^s_{\kappa,+}(M) \subset \Cal I^k_{\kappa,+}(M)$ is $\{\Phi^k_t\}$-invariant, it follows that $\Cal I^s_q(M)$ is coherent, hence regular by Lemma \ref{lemma:cohreg}, for $\mu$-almost all $q\in \Cal M^{(1)}_\kappa$.

\smallskip
\noindent By Theorem \ref{thm:smoothrel}, for $\mu$-almost all $q\in \Cal M^{(1)}_\kappa$  there 
exists a smoothing projection $\{P^\sigma(\tau) \vert \tau\in (0,1]\}$ of any given degree $\sigma>0$ relative to the subspace $\Cal I^\sigma_q(M)$. Let $s>1$ and let $\sigma\in \R^+$ and $k\in \N$ be such that $\sigma>k+1\geq s\geq \sigma-1 >k$. Let  $f \in [\Cal I^s_q(M)]^\perp \cap H^s_q(M)$.
By Definition \ref{def:smoothproj} of a smoothing projection, for any  $f \in [\Cal I^s_q(M)]^\perp \cap H^s_q(M)$, the function $P^\sigma(\tau)(f) \in [\Cal I^\sigma_q(M)]^\perp\cap H^\sigma_q(M)$ and satisfies the Sobolev bounds \pref{eq:smoothproj}. By Theorem \ref{thm:CEsharpfundth} the cohomological equation $S_qu = P^\sigma(\tau)(f)$ has a (unique) Green solution $u(\tau)\in 
H^{k-1}_q(M)$ and there exists a measurable function $C_\kappa: \Cal M^{(1)}_\kappa \to \R^+$ 
such that , for all $\tau\in (0,1]$, 
\begin{equation}
\label{eq:CEsharpth1}
\begin{aligned}
\vert u(\tau)-u(\tau/2) \vert_{k} &\leq C_\kappa(q) \, \vert P^\sigma(\tau)(f) 
- P^\sigma(\tau/2)(f) \vert_{\sigma} \,, \\
\vert u(\tau)-u(\tau/2) \vert_{k-1} &\leq C_\kappa(q) \, \vert P^\sigma(\tau)(f) - P^\sigma(\tau/2)(f)\vert_{\sigma-1} \,. \\
\end{aligned}
\end{equation}
By the interpolation inequality proved in Lemma \ref{lemma:intineq}, for any $r\in [k-1,k]$ there
exists $C_{k,r}>0$ such that, for all $\tau\in (0,1]$, 
\begin{equation}
\label{eq:CEsharpth2}
\vert u(\tau)-u(\tau/2) \vert_r \leq C_{k,r}\, \vert u(\tau)-u(\tau/2) \vert^{k-r}_{k-1} \, 
\vert u(\tau)-u(\tau/2) \vert^{r-k+1}_{k} \,.
\end{equation}
By the bounds  \pref{eq:smoothproj}, it follows from \pref{eq:CEsharpth1} and \pref{eq:CEsharpth1} 
that, for any $\epsilon >0$ there exists $C^\sigma_{r,s}(\epsilon)>0$ such that
\begin{equation}
\label{eq:CEsharpth3}
\vert u(\tau)-u(\tau/2) \vert_r  \leq  C^\sigma_{r,s}(\epsilon) \, \vert f \vert_s \, 
\tau^{(k-r)(s-\sigma+1-\epsilon)} \tau^{(r-k+1)(s-\sigma-\epsilon)}\,.
\end{equation}
Since $r<s-1$, it is possible to choose $\sigma\in \R^+$, $\epsilon>0$ and $k\in \N$ so that
$$
\alpha = (k-r)(s-\sigma+1-\epsilon) + (r-k+1)(s-\sigma-\epsilon) >0\,.
$$
It follows then from the bound \pref{eq:CEsharpth3}, that the sequence $\{u(2^{-n})\,\vert\, n\in\N\}$
is Cauchy in $H^r_q(M)$, hence it converges to a function $u\in H^r_q(M)$ of zero average. Since  $P^\sigma(\tau)(f)\to f$ in $H^s_q(M)$ as $\tau \to 0^+$, it follows that the function $u\in H^r_q(M)$ is the unique zero-average (Green) solution of the cohomological equation $S_q u=f$. The required Sobolev bound \pref{eq:CEsharpth} also follows from \pref{eq:CEsharpth3}. In fact, by the interpolation inequality (Lemma \ref{lemma:intineq}), by Theorem \ref{thm:CEsharpfundth} and by the bounds \pref{eq:smoothproj} for the smoothing projection, there exists a measurable function $C'_\kappa:
\Cal M^{(1)}_\kappa \to \R^+$ such that, 
\begin{equation}
\label{eq:CEsharpth4}
\vert u(1) \vert_r  \leq C'_\kappa(q) \, \vert f \vert_s \,.
\end{equation}
The bound \pref{eq:CEsharpth} can then be derived from \pref{eq:CEsharpth3} and \pref{eq:CEsharpth4}.
\end{proof}

\end{document}